\documentclass[oneside,12pt]{article}
\usepackage{graphicx}
\begin{document}
\title{Stabilizations of reducible Heegaard splittings}

\author {Ruifeng Qiu \thanks{Work supported in part by NSFC}  }
\date{}
\maketitle
\begin{abstract}Cameron Gordon([Problem 3.91 in [Ki]) conjectured that a connected sum
of two Heegaard splittings is stabilized if and only if one of the
two factors is stabilized. In this paper, we shall prove this
conjecture.
\end{abstract}

{\bf Keywords}: Abstract tree, band sum, connected sum,
stabilization.

AMS Classification: 57M25

\section{Introduction}

\ \ \ \ \ Let $M$ be a compact 3-manifold such that $\partial M$
has no 2-sphere components. A Heegaard splitting of $M$ is a pair
$(\cal V,\cal W)$, where $\cal V$ and $\cal W$ are compression
bodies such that $\cal V\cup \cal W$$=M$, and $\cal V\cap\cal
W=\partial_{+}\cal V=\partial_{+}\cal W=$$F$. $F$ is called a
Heegaard surface of $M$. The splitting is often denoted as $\cal
V$$\cup_F$$\cal W$ or $\cal V\cup W$. It is known that any compact
3-manifold has a Heegaard splitting.

Let $\cal V\cup\cal W$ be a Heegaard splitting of $M$. $\cal V\cup
\cal W$ is said to be reducible if there exist essential disks
$B_{\cal V}\subset\cal V$ and $B_{\cal W}\subset\cal W$ with
$\partial B_{\cal V}=
\partial B_{\cal W}$. Otherwise, it is said to be irreducible.
W. Haken[H] showed that any Heegaard splitting
of a reducible, compact 3-manifold is reducible, Kneser[Kn] and
Milnor[M] showed that any orientable, compact 3-manifold is a
connected sum of $n$ irreducible 3-manifolds $M_{1},\ldots, M_{n}$
where $M_{1},\ldots,M_{n}$ are unique up to isotopy. Thus any
Heegaard splitting of $M$ is the connected sums of the Heegaard
splittings of $n$ factors $M_{1},\ldots, M_{n}$ defined as
following:

Now let $M=\cal V\cup\cal W$ be a reducible Heegaard splitting.
Then $P=B_{\cal V}\cup B_{\cal W}$ is a 2-sphere. Suppose that $P$
cuts $M$ into $M_{+}^{*}$ and $M_{-}^{*}$. Then $B_{\cal V}$
separates $\cal V$ into $\cal V_{+}$ and $\cal V_{-}$, $B_{\cal
W}$ separates $\cal W$ into $\cal W_{+}^{*}$ and $\cal W_{-}^{*}$.
We may assume that $\cal V_{+}, W_{+}^{*}\subset$$ M_{+}$ and
$\cal V_{-},W_{-}^{*}$$\subset M_{-}$. Let
$M_{+}=M_{+}^{*}\cup_{P} H^{3}_{+}$ and $M_{-}=M_{-}^{*}\cup_{P}
H^{3}_{-}$ where $H^{3}_{+}$ and $H^{3}_{-}$ are two 3-balls. Then
$M$ is the connected sum of $M_{+}$ and $M_{-}$, denoted by
$M=M_{+}\sharp M_{-}$. Let $\cal W_{+}=\cal W_{+}^{*}$$\cup
H^{3}_{+}$, and $\cal W_{-}=\cal W_{-}^{*}$$\cup H^{3}_{-}$. Then
$\cal W_{+}$ and $\cal W_{-}$ are two compression bodies such that
$\partial_{+} \cal V_{+}=\partial_{+} W_{+}$ and $\partial_{+}
\cal V_{-}=\partial_{+} W_{-}$. Hence $M_{+}=\cal V_{+}\cup W_{+}$
is a Heegaard splitting of $M_{+}$ and $M_{-}=\cal V_{-}\cup
W_{-}$ is a Heegaard splitting of $M_{-}$. In this case, $\cal
V\cup W$ is called the connected sum of $\cal V_{+}\cup W_{+}$ and
$\cal V_{-}\cup W_{-}$.

A Heegaard splitting $M=\cal V\cup W$ is said to be stabilized if
there are two properly embedded disks $V\subset \cal V$ and
$W\subset \cal W$ such that $V$ intersects $W$ in one point;
otherwise, it is said to be unstabilized. Some important results
on stabilizations of Heegaard  splittings have been given in [RS],
[S], [ST] and [W].

An interesting problem on stabilizations of Heegaard splittings
offered by C. Gordon is the following:

{\bf \em  Gordon's conjecture.} \ The connected sum of two
Heegaard splittings is stabilized if and only if one of the two
factors is stabilized. (See Problem 3.91 in [Ki].)

In this paper, we shall give a  proof to Gordon's conjecture. The
main result is the following theorem:

{\bf \em Theorem 1.} \ The connected sum of two Heegaard
splittings is stabilized if and only if one of the two factors is
stabilized.

{\bf Comments on Theorem 1.}

(1) \ By Haken's lemma,  the connected sum of the minimal Heegaard
splittings of $M_{+}$ and $M_{-}$ is unstabilized;  but there are
many manifolds which have unstabilized Heegaard splittings of
distinct genera. There are examples, given by A. Casson and C.
Gordon[CG1], independently by T. Kobayashi[Ko], which have
irreducible Heegaard splittings of arbitrarily high genera. Now
let $M_{1}$ be such a manifold, and $M_{2}$ be any compact
3-manifold. Then, by Theorem 1, $M_{1}\sharp M_{2}$ has
unstabilized Heegaard splittings of arbitrarily high genus.

(2) \  David Bachman[B] announced that if $M_{+}$ and $M_{-}$ are
closed and irreducible then Gordon's conjecture is true. In this
case, $\cal V_{+}\cup W_{+}$ and $\cal V_{-}\cup W_{-}$ are
irreducible.

(3) \ Two applications of Theorem 1 have been given in [QM]:

(a) \ A Heegaard splitting $M=$$\cal W$$\cup_{S}$$\cal V$ is said
to be boundary reducible if there is an essential disk $D$ of $M$
which intersects $S$ in an essential simple curve in $S$. A.
Casson and C. Gordon proved that any Heegaard splitting of a
boundary reducible 3-manifold is boundary reducible in [CG2]. In
[QM], we shall prove that any Heegaard splitting $\cal W\cup V$ of
a boundary reducible, irreducible 3-manifold is obtained by doing
boundary connected sums and self-boundary connected sums from
Heegaard splittings of $n$ 3-manifolds $M_{1},\ldots, M_{n}$,
where $M_{i}$ is either boundary irreducible or a solid torus.
Furthermore, $\cal W\cup V$ is unstabilized if and only if  one of
the factors is unstabilized. This result can be taken as the disk
version of Gordon's conjecture.

(b) \ Suppose that $M_{1}$ and $M_{2}$ are two compact 3-manifolds
with boundary. Let $A_{i}$ be an incompressible annulus in
$\partial M_{i}$, and $M=M_{1}\cup_{A_{1}=A_{2}} M_{2}$. Let
$M_{i}=$$\cal W$$_{i}\cup$$\cal V$$_{i}$ be an unstabilized
Heegaard splitting of $M_{i}$.  Then $M$ has a natural Heegaard
splitting $\cal W$$\cup$$\cal V$ induced by $\cal W$$_{1}\cup$$
\cal V$$_{1}$ and $\cal W$$_{2}\cup$$\cal V$$_{2}$ such that
$g($$\cal W$$)=g($$\cal W$$_{1})+g($$\cal W$$_{2})$. Without loss
of generality, we may assume that $A_{1}\subset \partial_{-}$$
\cal V$$_{i}$. We denote by $M_{i}(A_{i})$ the manifold obtained
by attaching a 2-handle to $M_{i}$ along $A_{i}$, $\cal
V$$_{i}(A_{i})$ the manifold obtained by attaching a 2-handle to
$\cal V$$_{i}$ along $A_{i}$. Then $M_{i}(A_{i})=$$\cal
W$$_{i}\cup$$\cal V$$_{i}(A_{i})$ is a Heegaard splitting of
$M_{i}(A_{i})$. In [QM], we shall prove that if
$M_{1}(A_{1})=$$\cal W$$_{1}\cup$$\cal V$$_{1}(A_{1})$ and
$M_{2}(A_{2})=$$\cal W$$_{2}\cup$$\cal V$$_{2}(A_{2})$ are
unstabilized, then $\cal W$$\cup$$\cal V$ is unstabilized. \vskip
0.5mm

We shall use a basic tool  in 3-manifold theory, called band sums
of disks, to prove Theorem 1. The argument in this paper is
self-contained. We shall give an outline of the proof of Theorem 1
in Chapter 2.

\section{The outline of Theorem 1}

\ \ \ \ \ In this chapter, we shall introduce the ideas of the
proof of Theorem 1. Before doing this, we first give some
notations and related basic observations.

\subsection{The surface generated by an abstract tree}

\ \ \ \ \ {\bf \em Definition 2.1.1.} \ (1) \ Let $I=[-1,1]$.

(2) \ Suppose that  $a$ is an arc in a surface $F$, $a\times I$ be
a neighborhood of $a$ in $F$ such that $a=a\times\bigl\{0\bigr\}$.
For a sub-arc $b$ of $a$, we denote $b\times I\subset a\times I$
by $(b\times I)_{a}$. In this case, $b\times I$ and $a\times I$
have the same width. \vskip 0.5mm

Let $m(l)=\bigl\{m^{1},\ldots, m^{l}\bigr\}$ be a subset of
$\bigl\{1,\ldots, n\bigr\}$ possibly not in  a natural order. This
means that it is possible $m^{i}>m^{j}$ when $i<j$.

{\bf \em Definition 2.1.2.} \ Suppose that $F$ is a closed
surface, $P_{0},\ldots, P_{k}$ are pairwise disjoint disks in $F$,
and $e_{m^{1}}\ldots,e_{m^{l}}$ are pairwise disjoint arcs in $F$
satisfying the following conditions:

(1) \ $\partial e_{\gamma}\subset \cup_{f=0}^{k}\partial P_{f}$.

(2) \ If we denote by $e_{\gamma}^{'}$ the arc obtained by pushing
$inte_{\gamma}$ off $\cup_{f} P_{f}$ in $F\times I$, then each
component of $\cup_{f} P_{f}\cup_{\gamma} e_{\gamma}^{'}$ is a
tree when we take $P_{f}$ as a fat vertex and $e_{\gamma}^{'}$ as
an edge. \\
Then we say $\cup_{f} P_{f}\cup_{\gamma} e_{\gamma}$ is an
abstract tree.\vskip 0.5mm

Let $\cup_{f} P_{f}\cup_{\gamma} e_{\gamma}$ be an abstract tree
in a closed surface, and $e_{\gamma}\times I$ be a regular
neighborhood of $e_{\gamma}$ in $F$ satisfying the conditions:

($1^{*}$) \ If $\lambda>\gamma\in m(l)$, then either
$e_{\lambda}\cap e_{\gamma}\times I=\emptyset$ or each component
of $e_{\lambda}\cap e_{\gamma}\times I$ is an arc $c\subset
inte_{\lambda}$ which is a core of $e_{\gamma}\times (0,1)$, and
each component of $e_{\lambda}\times I\cap e_{\gamma}\times I$ is
$(c\times I)_{\lambda}\subset e_{\gamma}\times (0,1)$.

($2^{*}$) \ Each component of $inte_{\lambda}\times I\cap(\cup_{f}
P_{f}\cup_{\gamma<\lambda} e_{\gamma}\times I)$ is $(b\times
I)_{\lambda}$ where $b\subset inte_{\gamma}$ is a properly
embedded arc in $\cup_{f} P_{f}\cup_{\gamma<\lambda}
e_{\gamma}\times I$.

($3^{*}$) \ For $\gamma\in m(l)$, $(\partial e_{\gamma})\times
I\subset \cup_{f} \partial P_{f}$, for $\lambda\neq\gamma\in
m(l)$, $(\partial e_{\gamma})\times I\cap (\partial
e_{\lambda})\times I=\emptyset$.\vskip 0.5mm

{\bf \em Lemma 2.1.3.} \  If Conditions ($1^{*}$), ($2^{*}$) and
($3^{*}$) are satisfied, then $S=\cup_{f} P_{f}\cup_{\gamma\in
m(l)} e_{\gamma}\times I$ is a compact surface.

{\bf \em Proof.} \ Let $\gamma=Minm(l)$. Then,  by Condition
($2^{*}$), each component of $inte_{\gamma}\times I\cap (\cup_{f}
P_{f})$ is $(b\times I)_{\gamma}$ where $b\subset int e_{\gamma}$
is  a properly embedded arc in $\cup_{f} P_{f}$. Hence $b$ is
compact and $\partial b\cap
\partial e_{\gamma}=\emptyset$. We denote by $b_{1},\ldots, b_{\alpha}$ the
components of $e_{\gamma}-\cup_{f} intP_{f}$. Then $b_{i}$ is an
arc with $\partial b_{i}\subset \cup_{f} \partial P_{f}$. Hence
$\cup_{f} P_{f}\cup e_{\gamma}\times I=\cup_{f} P_{f}\cup
\cup_{i=1}^{\alpha} (b_{j}\times I)_{\gamma}$  is a surface. By
Conditions ($1^{*}$),($2^{*}$) and ($3^{*}$) and induction on
$m(l)=\bigl\{m^{1},\ldots, m^{l}\bigr\}$, we can prove Lemma
2.1.3.\qquad Q.E.D.

{\bf \em Definition 2.1.4.} \ Let $\cup_{f} P_{f}\cup_{\gamma}
e_{\gamma}$ be an abstract tree in a closed surface. If Conditions
($1^{*}$), ($2^{*}$) and ($3^{*}$) are satisfied, then $S=\cup_{f}
P_{f}\cup_{\gamma\in m(l)} e_{\gamma}\times I$ is called  a
surface generated by $\cup_{f} P_{f}\cup_{\gamma} e_{\gamma}$.

\begin{center}
\includegraphics[totalheight=4cm]{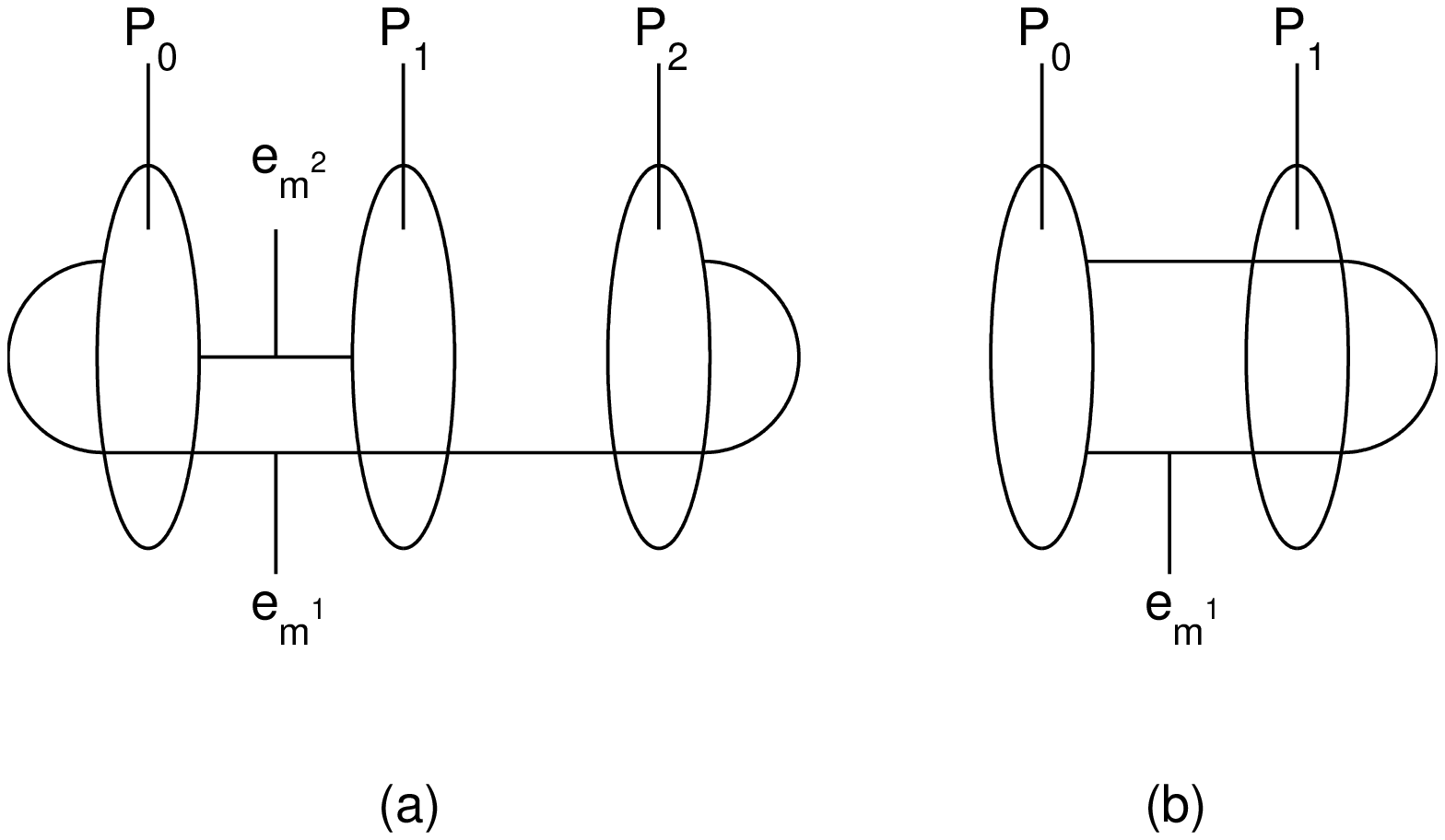}
\begin{center}
Figure 1
\end{center}
\end{center}

\begin{center}
\includegraphics[totalheight=4.5cm]{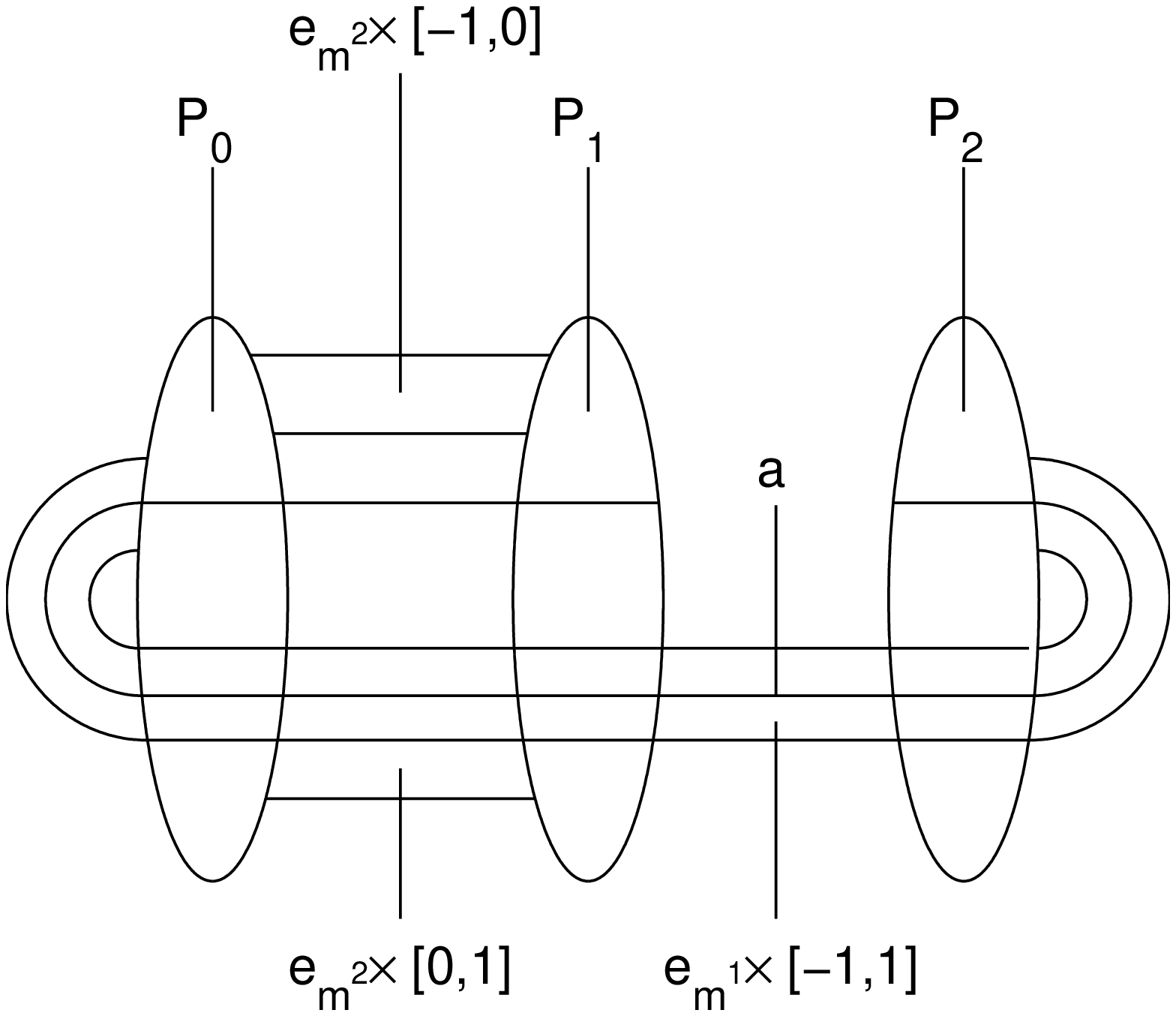}
\begin{center}
Figure 2
\end{center}
\end{center}

{\bf \em Example 1:}

(1). \ Suppose that $P_{0}, P_{1}, P_{2}$ are three disks in a
surface $F$ and $e_{m^{1}}, e_{m^{2}}$ are two arcs in $F$ as in
Figure 1(a). Then $\cup_{f} P_{f}\cup_{\gamma} e_{\gamma}\times I$
is an abstract tree.

(2). \ Suppose that $P_{0}, P_{1}$ and $e_{m^{1}}$ are as in
Figure 1(b). Then $P_{1}\cup P_{2}\cup e_{m^{1}}$ is not an
abstract tree. In this case, if $e^{'}_{m^{1}}$ is an arc  such
that $\partial e_{m^{1}}=\partial e_{m^{1}}^{'}$ and
$inte_{m^{1}}^{'}$ is disjoint from $P_{0}\cup P_{1}$, then
$P_{0}\cup P_{1}\cup e_{m^{1}}^{'}$ is not a tree.\vskip 0.5mm

(3). \ Figure 2 is a surface generated by an abstract tree
$\cup_{f} P_{f}\cup_{\gamma} e_{\gamma}$ as in Figure 1(a).  In
this case, $m^{1}>m^{2}$. Note that the surface generated by a
fixed abstract tree is not unique.

Now we consider a kind of special arcs in a surface generated by
an abstract tree.

{\bf \em Definition 2.1.5.} \ Let $S$  be a surface generated by
an abstract tree $\cup_{f} P_{f}\cup_{\gamma} e_{\gamma}$. If $a$
is a properly embedded arc in $S$ such that

1) \ for $\gamma\in m(l)$, each component of $a\cap
e_{\gamma}\times I$ is a core of $e_{\gamma}\times (0,1)$ which is
contained in $inta$,

2) \ each component of $a\cap (\cup_{\gamma} e_{\gamma}\times I)$
is a core of $e_{\lambda}\times (0,1)$ for some $\lambda\in m(l)$,

3) \ for each $\lambda\in m(l)$, there is at most one component of
$a\cap (\cup_{\gamma} e_{\gamma}\times I)$ which  is a core of
$e_{\lambda}\times (0,1)$. \\
Then $a$ is said to be regular in $S$.\vskip 0.5mm

{\bf \em Example 2:}

The arc $a$ in Figure 2 is regular in $P_{0}\cup P_{1}\cup
P_{2}\cup e_{m^{1}}\times I\cup e_{m^{2}}\times I$. In this case,
$a$ intersects $e_{m^{2}}\times I$ in two cores of
$e_{m^{2}}\times (0,1)$, $e_{m^{1}}\times I$ in one core of
$e_{m^{1}}\times (0,1)$, $e_{m^{1}}\times I\cup e_{m^{2}}\times I$
in two components, one of which is a core of $e_{m^{1}}\times
(0,1)$ and the other is a core of $e_{m^{2}}\times (0,1)$.
$a-(inte_{m^{1}})\times I\cup (inte_{m^{2}}\times I)$ contains
three components $a_{0}, a_{1}, a_{2}$ with $a_{i}\subset
P_{i}$.\vskip 0.5mm

{\bf \em Lemma 2.1.6.} \ Let $a$ be a regular arc in a surface $S$
generated by an abstract tree $\cup_{f} P_{f}\cup_{\gamma}
e_{\gamma}$. Then $a=\cup_{i=1}^{\theta(a)}
a_{f_{i}}\cup_{i=1}^{\theta(a)-1} e_{\gamma_{i},a}$ satisfying the
following conditions:

1) \ $0\leq f_{i}\leq k$, $a_{f_{i}}$ is a properly embedded arc
in $P_{f_{i}}$ which is disjoint from $\cup_{\gamma\in m(l)
}inte_{\gamma}\times I$.

2) \  $\gamma_{i}\in m(l)$, and $e_{\gamma_{i},a}$ is a core of
$e_{\gamma_{i}}\times (0,1)$.

3) \ For $1\leq i\leq \theta(a)-1$, $\partial_{2}
a_{f_{i}}=\partial_{1} e_{\gamma_{i},a}$ and  for $2\leq i\leq
\theta(a)-1$,  $\partial_{1} a_{f_{i+1}}=\partial_{2}
e_{\gamma_{i},a}$, $\partial_{1} a=\partial_{1} a_{f_{1}}$,
$\partial_{2} a=\partial_{2} a_{f_{\theta(a)}}$.

4) \ For $1\leq i\neq r\leq \theta(a)$, $f_{i}\neq f_{r}$, for
$1\leq i\neq r\leq \theta(a)-1$, $\gamma_{i}\neq \gamma_{r}$.

5) \ For each $f$, $P_{f}\cap (a-\cup_{\gamma} inte_{\gamma}\times
I)$ contains at most one component.

{\bf \em Proof.} \  Since $a$ is regular in $S$, By Definition
2.1.5(1) and (2), $\partial a$ is disjoint from $\cup_{\gamma}
e_{\gamma}\times I$, and  each component of $a\cap (\cup_{\gamma}
e_{\gamma}\times I)$ is a core of $e_{\lambda}\times (0,1)$ for
some $\lambda\in m(l)$. By Definition 2.1.4, $(\partial
e_{\gamma})\times I\subset \cup_{f} \partial P_{f}$, each
component of $a-\cup_{\gamma} inte_{\gamma}\times I$ is an
properly embedded arc in $P_{f}$ for some $0\leq f\leq k$. Hence
$a=\cup_{j=1}^{\theta(a)} a_{f_{j}}\cup_{j=1}^{\theta(a)-1}
e_{\gamma_{j},a}$,  and (1), (2), (3) holds.

By Definition 2.1.5(3), $\gamma_{i}\neq\gamma_{r}$ for $1\leq
i\neq r\leq \theta(a)-1$. Now if $f_{i}=f_{r}$ for $1\leq i\neq
r\leq \theta(a)$, then $\cup_{f} P_{f}\cup_{\gamma} e_{\gamma}$ is
not an abstract tree, a contradiction. Hence (4) holds. (5)
follows  from (1) and (4). \qquad Q.E.D.\vskip 0.5mm

{\bf \em Remark 2.1.7.} \ The properties of regular arcs in Lemma
2.1.6 are important in the proof of Theorem 1. Lemmas 3.1.2 and
3.3.1 follow from Lemma 2.1.6 and they take roles in the inductive
proof of Theorem 1. \vskip 0.5mm

{\bf \em Lemma 2.1.8.} \ Let $S=\cup_{f} P_{f}\cup_{\gamma\in
m(l)} e_{\gamma}\times I$ be a surface generated by an abstract
tree. Then $S_{j}=\cup_{f} P_{f}\cup_{\gamma<j} e_{\gamma}\times
I$  is also a surface generated by an abstract tree where  $j$ is
any integer.

{\bf \em Proof.} \ Since $\cup_{f} P_{f}\cup_{\gamma\in
m(l)}e_{\gamma}$ is an abstract tree, $\cup_{f}
P_{f}\cup_{\gamma<j} e_{\gamma}$ is also an abstract tree.
Obviously, Conditions ($1^{*}$), ($2^{*}$) and ($3^{*}$) are
satisfied. Hence $S_{j}=\cup_{f} P_{f}\cup_{\gamma<j}
e_{\gamma}\times I$ is a surface generated by $\cup_{f}
P_{f}\cup_{\gamma<j} e_{\gamma}$.\qquad Q.E.D.

\subsection{The element of the induction}

\ \ \ \ \ {\bf It is easy to see that if one of $\cal V_{+}\cup
W_{+}$ and $\cal V_{-}\cap W_{-}$ is stabilized then $\cal V\cup
W$ is stabilized. So in order to obtain a contradiction, we may
assume that  each of $\cal V_{+}\cup W_{+}$ and $\cal V_{-}\cup
W_{-}$ is unstabilized and $\cal V\cup W$ is stabilized.}\vskip
3mm

{\bf Assumption(*).} (1) Let $(V, W)$ be a pair of stabilized
disks such that $V\subset\cal V$, $W\subset\cal W$.

Now  $V$ intersects $W$ in only one point $x$. Recalling the two
disks $B_{\cal V}$ and $B_{\cal W}$ defined in Chapter 1. We may
assume that each component of $V\cap B_{\cal V}$ is an arc in both
$V$ and $B_{\cal V}$, each component of $W\cap B_{\cal W}$ is an
arc in both $W$ and $B_{\cal W}$. It is easy to see that if $V\cap
B_{\cal V}=\emptyset$ or $W\cap B_{\cal W}=\emptyset$, then one of
$\cal V_{+}\cup W_{+}$ and $\cal V_{-}\cap W_{-}$ is stabilized.
So we may assume that $V\cap B_{\cal V}\neq\emptyset$ and  $W\cap
B_{\cal W}\neq\emptyset$.

{\bf \em Assumption (*).}  (2) \ $x\in\partial_{+}{\cal
V_{+}}-B_{\cal V}$.

(3) \ $|V\cap B_{\cal V}|=m$ and $|W\cap B_{\cal W}|=n$.\vskip
1.5mm

Now consider $V\cap B_{\cal V}$. By assumption, each component $e$
of $V\cap B_{\cal V}$ separates $V$ into two disks $V_{e}^{'}$ and
$V_{e}^{''}$ such that $x\in
\partial V_{e}^{'}$. Now we denote by
$V_{e}$ the disk in $V_{e}^{''}$ which is bounded by $e$, with
some components in $V\cap B_{\cal V}$ and some arcs in $\partial
V$, such that $\textrm{int}V_{e}$ is disjoint from $B_{\cal V}$ as
in Figure 3. Then $V_{e}$ is a properly embedded disk in $\cal
V_{+}$ or $\cal V_{-}$.
\begin{center}
\includegraphics[totalheight=4cm]{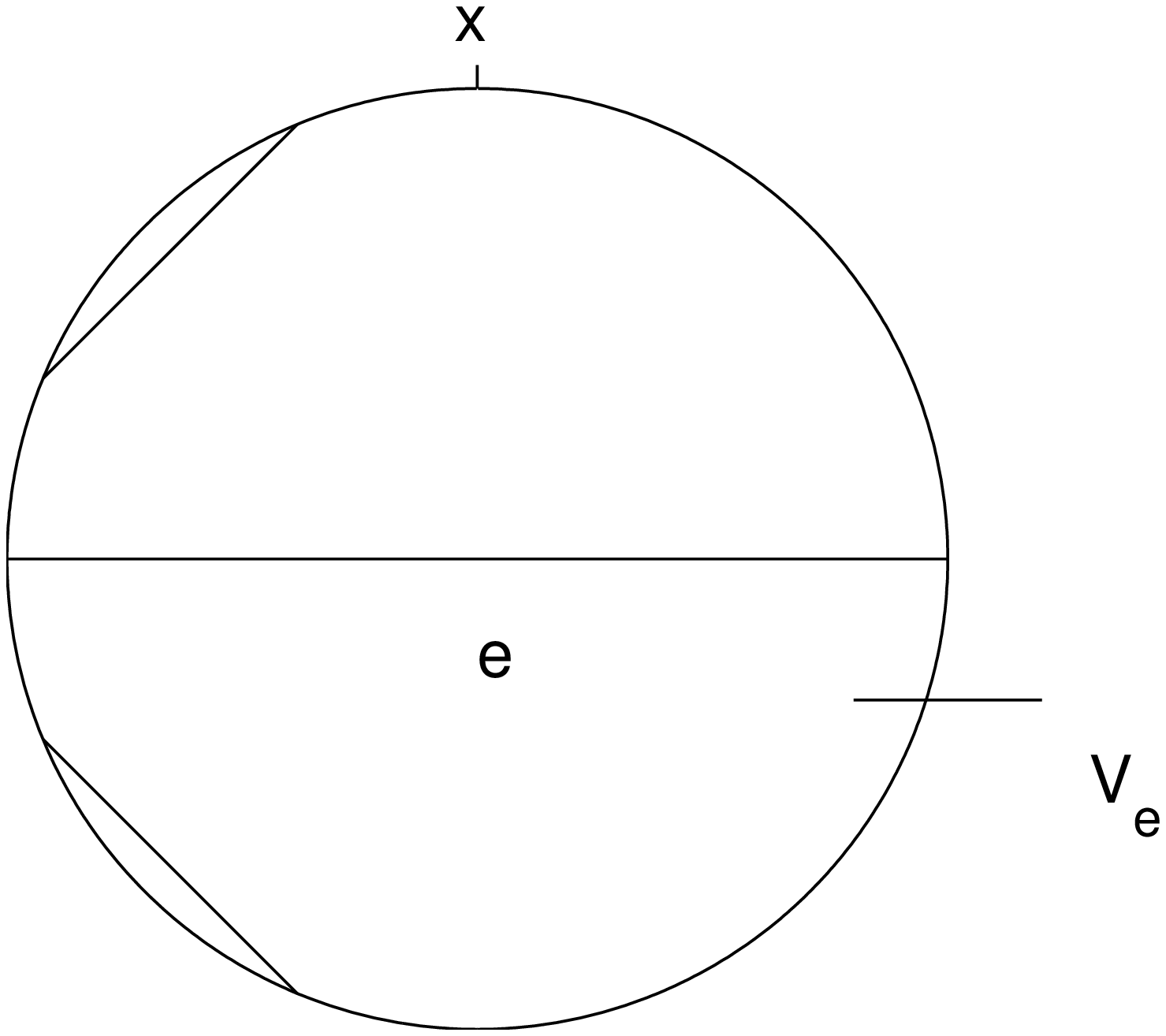}
\begin{center}
Figure 3
\end{center}
\end{center}

{\bf \em Definition 2.2.1.} \ A component $e$ in $V\cap B_{\cal
V}$ is labeled a symbol $s(e)$ where $s(e)=+$  if $V_{e}\subset
\cal V_{+}$ and $s(e)=-$ if $V_{e}\subset\cal V_{-}$.\vskip 0.5mm

Similarly, each component $e$ of $W\cap B_{\cal W}$ separates $W$
into two disks $W_{e}^{'}$ and $W_{e}^{''}$ such that $x\in
\partial W_{e}^{'}$. Now we denote by
$W_{e}$ the disk in $W_{e}^{''}$ which is bounded by $e$ with some
components in $W\cap B_{\cal W}$ and some arcs in $\partial W$
such that $\textrm{int}W_{e}$ is disjoint from $B_{\cal W}$. Then
$W_{e}$ is a properly embedded disk in $\cal W_{+}^{*}$ or $\cal
W_{-}^{*}$ defined in Section 1.\vskip 0.5mm

{\bf \em Definition 2.2.2.} \ A component $e$ in $W\cap B_{\cal
W}$ is labeled a symbol $s(e)$ where $s(e)=+$  if $W_{e}\subset
\cal W_{+}^{*}$ and $s(e)=-$ if $W_{e}\subset\cal
W_{-}^{*}$.\vskip 0.5mm

By Assumption(*), we number the components of $V\cap B_{\cal V}$,
$v_{1},\ldots,v_{m}$, and the components of $W\cap B_{\cal W}$,
$w_{1},\ldots, w_{n}$, so that if $V_{v_{i}}^{''}\subset
V_{v_{k}}^{''}$ and $W_{w_{j}}^{''}\subset W_{w_{l}}^{''}$, then
$i<k$, $j<l$. Now we denote by $V_{i}$ the disk $V_{v_{i}}$,
$W_{j}$ the disk $W_{w_{j}}$ for $1\leq i\leq m$ and $1\leq j\leq
n$.\vskip 0.5mm

{\bf \em Definition 2.2.3.} \ For each $v_{i}$ in $V\cap B_{\cal
V}$ and each $w_{j}$ in $W\cap B_{\cal W}$, let $I(v_{i})=\bigl\{r
\ | \ v_{r}\neq v_{i}\subset
\partial V_{i} \bigr\},$ and
$I(w_{j})=\bigl\{r \ | \ w_{r}\neq w_{j}\subset
\partial W_{j} \bigr\}.$

The following two lemmas are immediately from definitions.\vskip
3mm

{\bf \em Lemma 2.2.4.} \ 1) \ If $r\in I(v_{i})$, then $r<i$,  and

2) \ if  $r\in I(w_{j})$, then $r<j$.

{\bf \em Proof.} \ We need only to prove (1).

If $r\in I(v_{i})$, then $V_{r}^{''}\subset V_{i}^{''}$. By
definitions, $r<i$. \vskip 0.5mm

{\bf \em Lemma 2.2.5.} \ (1) $s(v_{i})=+$ if and only if
$s(v_{r})=-$ for each $r\in I(v_{i})$.

(2) If $r\in I(v_{i}), I(v_{k})$, then $k=i$.

(3) $s(w_{j})=+$ if and only if $s(w_{r})=-$ for each $ r\in
I(w_{j})$.

(4) If $r\in I(w_{j}), I(w_{k})$, then $k=j$.\vskip 0.5mm

{\bf \em Proof.} \ This lemma is immediately from definitions and
simple observations.\qquad Q.E.D.\vskip 0.5mm

Since $x\in \partial_{+} \cal V_{+}$, each component of $V\cap
\cal V_{+}$ is either the disk $V_{i}$ for some $1\leq i\leq m$
with $s(v_{i})=+$ or a disk containing $x$, denoted by $V_{x}$. By
definition of $V_{i}$, $V_{x}$ is disjoint from $V_{i}$ for
$s(v_{i})=+$. Now each component of $V\cap \cal V_{-}$ is the disk
$V_{i}$ for some $1\leq i\leq m$ with $s(v_{i})=-$. Similarly,
each component of $W\cap \cal W_{+}^{*}$ is either the disk
$W_{j}$ for some $1\leq j\leq n$ with $s(w_{j})=+$ or a disk
containing $x$, denoted by $W_{x}$, each component of $W\cap \cal
W_{-}^{*}$ is the disk $W_{j}$ for some $1\leq j\leq n$ with
$s(w_{j})=-$.\vskip 0.5mm

{\bf \em Definition 2.2.6.} \  Let $I(v)=\bigl\{r \ | \ v_{r}\in
V_{x}\cap B_{\cal V}\bigr\}$ and $I({w})=\bigl\{r \ | \ w_{r}\in
W_{x}\cap B_{\cal W} \bigr\}$.\vskip 0.5mm

{\bf \em Lemma 2.2.7.} \  1) If $r\in I(v)$, then $s(v_{r})=-$.

2) \ If $r\in I(w)$, then $s(w_{r})=-$.

3) \ $m\in I(v)$, $n\in I(w)$.

{\bf \em Proof.} \ By assumption, $x\in
\partial_{+}\cal V_{+}$. Hence $V_{x}\subset \cal V_{+}$,
$W_{x}\subset \cal W_{+}$. If $r\in I(v)$, then $V_{r}\subset \cal
V_{-}$. If $r\in I(w)$, then $W_{r}\subset \cal W_{-}$.

Now $v_{m}$ separates $V$ into $V_{v_{m}}^{'}$ and
$V_{v_{m}}^{''}$ such that $x\in V_{v_{m}}^{'}$. Suppose,
otherwise, $m\notin I(v)$. By definition, $v_{m}$ is disjoint from
$ \partial V_{x}$. This means that there is an integer $1\leq i<m$
such that $v_{i}$ separates $v_{m}$ and the point $x$ in
$V_{v_{m}}^{'}$. By definition, $V_{v_{m}}^{''}\subset
V_{v_{i}}^{''}$. Hence $m<i$, a contradiction. \qquad Q.E.D.\vskip
0.5mm

Recalling the definitions of $\cal V_{+}\cup W_{+}$ and $\cal
V_{-}\cup W_{-}$ in section 1. Now $W_{j}$ and $W_{x}$ are
properly embedded in $\cal W_{+}^{*}$ or $\cal W_{-}^{*}$ for
$1\leq j\leq n$. It is easy to see that for each $j$, there is an
arc $w_{j,v}$ in $B_{\cal V}$ such that $w_{j,v}\cup w_{j}$ bounds
a disk in $H_{+}^{3}$, denoted by $W_{j,+}$, and   a disk in
$H_{-}^{3}$, denoted by $W_{j,-}$ as in Figure 4. Note that for
each $j\neq k$, $W_{{j},+}\cap W_{{k},+}=\emptyset$ and
$W_{{j},-}\cap W_{{k},-}=\emptyset$. Thus if $s(w_{j})=+$, then
$W_{j,v}=W_{j}\cup W_{j,+}\cup_{r\in I(w_{j}) } W_{r,+}$ is a
properly embedded disk in $\cal W_{+}$, and if $s(w_{j})=-$, then
$W_{j,v}=W_{j}\cup W_{j,-}\cup_{r\in I(w_{j})} W_{r,-}$ is a
properly embedded disk in $\cal W_{-}$. Specially,
$W_{x,v}=W_{x}\cup_{r\in I(w)} W_{r,+}$ is a properly embedded
disk in $\cal W_{+}$.

\begin{center}
\includegraphics[totalheight=4cm]{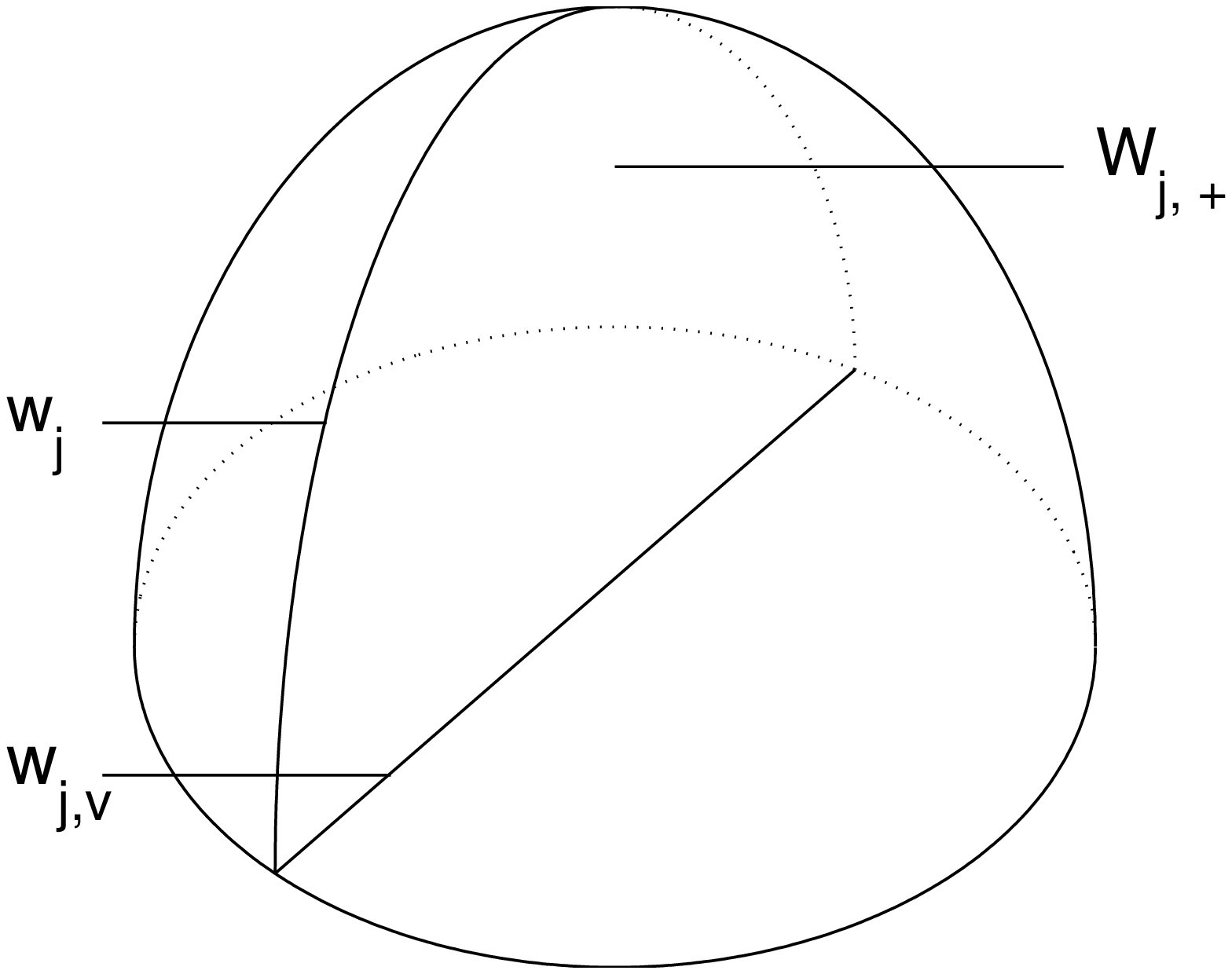}
\begin{center}
Figure 4
\end{center}
\end{center}

\begin{center}
\includegraphics[totalheight=4cm]{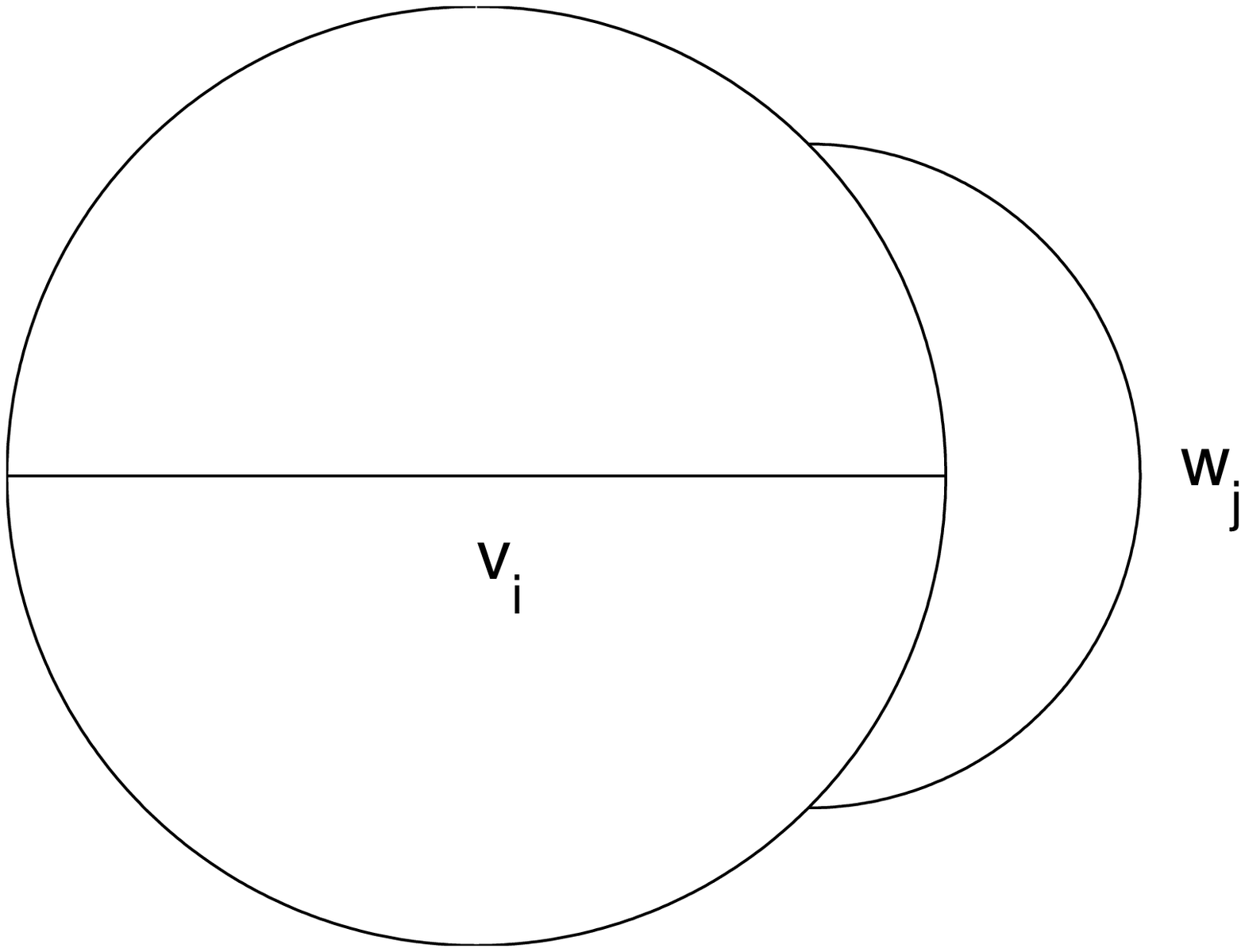}
\begin{center}
Figure 5
\end{center}
\end{center}

{\bf \em Lemma 2.2.8.} \ $\bigl\{v_{1},\ldots, v_{m}\bigr\}$ and
$\bigl\{w_{1,v},\ldots,w_{n,v}\bigr\}$ are two sets of pairwise
disjoint arcs properly embedded in $B_{\cal V}$ such that

(1) \ $v_{i}$ intersects $w_{j,v}$ in at most one point, and

(2) \ $w_{j,v}\cap\partial v_{i}=\emptyset$ for each $1\leq i\leq
m$ and $1\leq j\leq n$.

{\bf \em Proof.} \ Since $v_{i}\subset V\cap B_{\cal V}$,
$w_{j}\subset W\cap B_{\cal W}$ and $W\cap V=\bigl\{x\bigr\}$, By
Assumption(*),  $v_{i}$ is properly embedded in $B_{\cal V}$ and
$w_{j}$ is properly embedded in $B_{\cal W}$ such that $v_{i}$ is
disjoint from $w_{j}$. Note that $\partial B_{\cal V}=\partial
B_{\cal W}$. Hence either the two end points of $w_{j}$ lie in the
same component of $\partial B_{\cal V}-\partial v_{i}$ or the two
end points of $w_{j}$ lie in the distinct components of $\partial
B_{\cal V}-\partial v_{i}$ as in Figure 5. Hence the lemma
holds.\qquad Q.E.D.

{\bf \em  Lemma 2.2.9.} \ There are four sets of pairwise disjoint
properly embedded disks $\bigl\{V_{i} \ | \
s(v_{i})=+\bigr\}\cup\bigl\{V_{x}\bigr\}$ in $\cal V_{+}$,
$\bigl\{V_{i} \ | \ s(v_{i})=-\bigr\}$ in $\cal V_{-}$,
$\bigl\{W_{j,v} \ | \ s(w_{j})=+\bigr\}\cup\bigr\{W_{x,v}\bigr\}$
in $\cal W_{+}$, and $\bigl\{W_{j,v} \ | \ s(w_{j})=-\bigr\}$ in
$\cal W_{-}$ satisfying the following conditions:

(1) \ $V_{i}\cap B_{\cal V}= v_{i}\cup_{r\in I(v_{i})} v_{r}$,
$W_{j,v}\cap B_{\cal V}=w_{j,v}\cup_{r\in I(w_{j})} w_{r,v}$,
$V_{x}\cap B_{\cal V }=\cup_{r\in I(v)} v_{r}$, $W_{x,v}\cap
B_{\cal V}=\cup_{r\in I(w)} w_{r,v}$.

(2) \ If $s(v_{i})=+, s(w_{j})=+$, then $V_{i}\cap
W_{j,v}=V_{i}\cap W_{j,v}\cap B_{\cal V}$, $V_{i}\cap
W_{x,v}=V_{i}\cap W_{x,v}\cap B_{\cal V}$, $V_{x}\cap
W_{j}=V_{x}\cap W_{j}\cap B_{\cal V}$, $V_{x}\cap
W_{x,v}=(V_{x}\cap W_{x,v}\cap B_{\cal V})\cup\bigl\{x\bigr\}$.

3) \ If $s(v_{i})=-, s(w_{j})=-$, then $V_{i}\cap
W_{j,v}=V_{i}\cap W_{j,v}\cap B_{\cal V}$.

{\bf \em Proof.} \ (1) follows from the construction of $w_{j,v}$.

Since $W\cap V=W_{x}\cap V_{x}=\bigl\{x\bigr\}$, $W_{j}\cap
V_{i}=\emptyset$ for each $1\leq i\leq m$ and $1\leq j\leq n$. By
the constructions of $W_{j,v}$, $W_{x,v}$, $w_{j,v}$, (2) and (3)
holds.\qquad Q.E.D.

\subsection {Outline of the proof of Theorem 1}

{\bf The idea of the Proof of Theorem 1}

Now let $\cal V$$^{1}_{+}=\partial_{+} \cal V$$_{+}\times
I\cup_{s(v_{i})=+} N(V_{i})\cup
N(V_{x})\cup\bigl\{3-handles\bigr\}$ where $N(V_{i})$ and
$N(V_{x})$ are regular neighborhoods of $V_{i}$ and $V_{x}$ in
$\cal V_{+}$, $\cal V$$^{1}_{-}=\partial_{+} \cal V$$_{-}\times
I\cup_{s(v_{i})=-} N(V_{i})\cup\bigl\{3-handles\bigr\}$ where
$N(V_{i})$ is a regular neighborhood of $V_{i}$ in $\cal V_{-}$.
Then $\cal V$$^{1}_{+}\subset $$V_{+}$ and $\cal
V$$_{-}^{1}\subset $$V_{-}$ are two compression bodies. Let $\cal
W$$^{1}_{+}=\partial_{+}\cal W$$_{+}\times I\cup_{s(w_{j})=+}
N(W_{j,v}) \cup N(W_{x,v})\cup\bigl\{3-handles\bigr\}$, where
$N(W_{j,v})$ and $N(W_{x,v})$ are regular neighborhoods of
$W_{j,v}$ and $W_{x,v}$ in $\cal W_{+}$,  $\cal
W$$^{1}_{-}=\partial_{+}\cal W$$_{-}\times I\cup_{s(w_{j})=-}
N(W_{j,v})\cup\bigl\{3-handles\bigr\}$ where $N(W_{j,v})$ is a
regular neighborhood of $W_{j,v}$ in $\cal W_{-}$.

Now $W$ and $V$ defined in Section 2.2 are a pair of stabilized
disks of the connected sum of $\cal V$$_{+}^{1}\cup$$\cal
W$$^{1}_{+}$ and $\cal V$$_{-}^{1}\cup$$\cal W$$^{1}_{-}$. If
Theorem 1 is true, then one of $\cal V$$_{+}^{1}\cup$$\cal
W$$^{1}_{+}$ and $\cal V$$_{-}^{1}\cup$$\cal W$$^{1}_{-}$, say
$\cal V$$_{+}^{1}\cup$$\cal W$$^{1}_{+}$, is stabilized. This
means that there are two essential disks $W^{'}\subset \cal
W$$^{1}_{+}$ and $V^{'}\subset \cal V$$^{1}_{+}$ such that $W^{'}$
intersects $V^{'}$ in only one point. Hence $W^{'}$ is obtained by
doing band sums from $W_{j,v}$ with $s(w_{j})=+$ and $W_{x,v}$,
and $V^{'}$ is obtained by doing band sums from $V_{i}$ with
$s(v_{i})=+$ and $V_{x}$. We do want to do this.

Recalling $I(v_{i}), I(w_{j}), I(v), I(w)$, $s(v_{i}), s(w_{j})$
which are defined in Section 2.2.

{\bf \em Definition 2.3.1.} \  For $0\leq k\leq m$ and
$m(k)=\bigl\{m^{0},\ldots, m^{k}\bigr\}\subset\bigl\{1,\ldots,
n\bigr\}$. Let $I(v_{i},k)=I(v_{i})-\bigl\{1,\ldots,k\bigr\}$,
$I(w_{j},k)=I(w_{j})-m(k)$, $I(v,k)=I(v)-\bigl\{1,\ldots,
k\bigr\}$ and $I(w,k)=I(w)-m(k)$.\vskip 0.5mm

By definition, if $k>l$, then $I(v_{i},k)\subset I(v_{i},l)$,
$I(v,k)\subset I(v,l)$. Furthermore,
$I(v_{i},m)=I(v,m)=\emptyset$. By definitions, $I(v_{i}),
I(v)\subset\bigl\{1,\ldots,m\bigr\}$. Hence
$I(v_{i},k)=I(v_{i})-\bigl\{1,\ldots,k\bigr\}$,
$I(v,k)=I(v)-\bigl\{1,\ldots, k\bigr\}$. Similarly, $I(w_{j}),
I(w)\subset\bigl\{1,\ldots,n\bigr\}$.  \vskip 0.5mm

{\bf To prove Theorem 1, we only need to prove the following
propositions:}\vskip 0.5mm

{\bf \em Proposition 1.} \  For each $0\leq k\leq m$, there are a
surface $P^{k}$ in $\partial_{+}\cal V_{+}$, a subset of
$\bigl\{1,\ldots, n\bigr\}$, say $m(k)=\bigl\{m^{0}, m^{1},\ldots,
m^{k}\bigr\}$, and two sets of pairwise disjoint arcs
$\bigl\{v_{i}^{k} \ | \ k+1\leq i\leq m \bigr\}$ and
$\bigl\{w_{j}^{k} \ | \ j\in \bigl\{1,\ldots,n\bigr\}-m(k)\bigr\}$
in $P^{k}$ satisfying the following conditions:

(1) \ $P^{k}=\cup_{f=0}^{k_{}} D^{k}_{f}\cup_{\gamma\in m(k)}
b_{\gamma}^{k}\times I$ is a surface generated by an abstract tree
$\cup_{f=0}^{k} D^{k}_{f}\cup_{\gamma\in m(k)} b_{\gamma}^{k}$
such that

(i) \ for  each $\lambda\subset m(k)$ with $s(w_{\lambda})=+$,
$(\textrm{int}b^{k}_{\lambda}\times I)\cap (\cup_{\gamma<\lambda}
b_{\gamma}^{k}\times I\cup_{f} D_{f}^{k})=\cup_{r\in
I(w_{\lambda},k)} (w_{r}^{k}\times I)_{\lambda}$;

(ii) for each $\lambda\subset m(k)$ with $s(w_{\lambda})=-$,
$\textrm{int}b^{k}_{\lambda}\times I$ is disjoint from
$\cup_{f=0}^{k} D^{k}_{f}\cup_{\gamma<\lambda}
b^{k}_{\gamma}\times I$.

(2) \ For $j\notin m(k)$, $w_{j}^{k}$ is regular in
$\cup_{\gamma<j} b_{\gamma}^{k}\times I\cup_{f} D_{f}^{k}$.

(3) \ For $j\notin m(k)$ and $\gamma\in m(k)$, if $j<\gamma$, then
either $j\in I(w_{\gamma},k)$ with $s(w_{\gamma})=+$ or
$w_{j}^{k}$ is disjoint from $b^{k}_{\gamma}\times I$.

(4)  \ For each $i\geq k+1$, $v_{i}^{k}$ is a properly embedded
arc in $P^{k}$  lying in $D^{k}_{f}$ for some $f$. Furthermore,
for each $j\in\bigl\{1,\ldots,n\bigr\}-m(k)$,
$w_{j}^{k}-\cup_{\gamma<j} intb^{k}_{\gamma}\times I$ intersects
$v_{i}^{k}$ in at most one point.

(5) \ For each $i\geq k+1$, $j\in\bigl\{1,\ldots,n\bigr\}-m(k)$
and $\gamma\in m(k)$, $\partial v_{i}^{k}\cap
w_{j}^{k}=\emptyset$, $\partial v_{i}^{k}\cap b_{\gamma}^{k}\times
I=\emptyset$.\vskip 0.5mm

{\bf \em Definition 2.3.2.} \ If $w_{j}^{k}-\cup_{\gamma<j} int
b^{k}_{\gamma}\times I$ intersects $v_{i}^{k}$ in one point, then
we say  $i\in L(w_{j}^{k})$ and $j\in L(v_{i}^{k})$.\vskip 0.5mm

{\bf \em Proposition 2.} \ For each $1\leq k\leq m$ and $j\in
\bigl\{1,\ldots,n\bigr\}- m(k)$.

(1) If $j\notin L(v_{k}^{k-1})$, then $L(w_{j}^{k})=
L(w_{j}^{k-1})$.

(2) \ If $j\in L(v_{k}^{k-1})$, then $L(w_{j}^{k})=
L(w_{j}^{k-1})\cup L(w_{m^{k}}^{k-1})- L(w_{j}^{k-1})\cap
L(w_{m^{k}}^{k-1})$.

(3) \ $m^{0}=\emptyset$ and $m^{k}=MinL(c_{k}^{k-1})$. \\ Where
$m(k)=\bigl\{m^{0},\ldots, m^{k}\bigr\}$ is as in Proposition 1.
\vskip 0.5mm

{\bf \em Proposition 3.} \ For each $0\leq k\leq m$, there are two
sets of pairwise disjoint disks $\bigl\{V_{i}^{k} \ | \ k+1\leq
i\leq m \ with \ s(v_{i})=+ \bigr\}\cup\bigl\{V_{x}^{k}\bigr\}$
properly embedded in $\cal V_{+}$ and $\bigl\{W_{j}^{k} \ | \ j\in
\bigl\{1,\ldots,n\bigr\}-m(k) \ with \ s(w_{j})=+
\bigr\}\cup\bigl\{W_{x}^{k}\bigr\}$ properly embedded in $\cal
W_{+}$ satisfying the following conditions:

(1) \ $V_{i}^{k}\cap P^{k}=v_{i}^{k}\cup_{r\in I(v_{i},k)}
v_{r}^{k}$, $W_{j}^{k}\cap P^{k}=w_{j}^{k}\cup_{r\in I(w_{j},k)}
w_{r}^{k}$, $V_{x}^{k}\cap P^{k}=\cup_{r\in I(v,k)} v_{r}^{k}$,
$W_{x}^{k}\cap P^{k}=\cup_{r\in I(w,k)} w_{r}^{k}$.

(2) \ $V_{i}^{k}\cap W_{j}^{k}=V_{i}^{k}\cap W_{j}^{k}\cap P^{k}$,
$V_{i}^{k}\cap W_{x}^{k}=V_{i}^{k}\cap W_{x}^{k}\cap P^{k}$,
$V_{x}^{k}\cap W_{j}^{k}=V_{x}^{k}\cap W_{j}^{k}\cap P^{k}$,
$V_{x}^{k}\cap W_{x}^{k}=(V_{x}^{k}\cap W_{x}^{k}\cap
P^{k})\cup\bigl\{x\bigr\}$ where $x\in W_{x,v}\cup V_{x}$ in Lemma
2.2.9. \vskip 0.5mm

Recalling the set $\bigl\{m^{0}, m^{1},\ldots, m^{k}\bigr\}$ in
Proposition 1.\vskip 0.5mm

{\bf \em Proposition 4.} \  For each $1\leq k\leq m$, there are a
surface $F^{k}$ in $\partial_{+}\cal V_{-}$,  and two sets of
pairwise disjoint arcs $\bigl\{c_{i}^{k} \ | \ k+1\leq i\leq m
\bigr\}$ and $\bigl\{d_{j}^{k} \ | \ j\in
\bigl\{1,\ldots,n\bigr\}-m(k)\bigr\}$ in $F^{k}$ satisfying the
following conditions:

(1) \ $F^{k}=\cup_{f=0}^{k} E^{k}_{f}\cup_{\gamma\in m(k)}
e_{\gamma}^{k}\times I$ is a surface generated by an abstract tree
$\cup_{f=0}^{k} E^{k}_{f}\cup_{\gamma\in m(k)} e_{\gamma}^{k}$
such that

(i) \ for each $\lambda\subset m(k)$,  if $s(w_{\lambda})=-$, then
$(\textrm{int}e^{k}_{\lambda}\times I)\cap (\cup_{\gamma<\lambda}
e_{\gamma}^{k}\times I\cup_{f} E_{f}^{k})=\cup_{r\in
I(w_{\lambda},k)} (d_{r}^{k}\times I)_{\lambda}$;

(ii) \ if $s(w_{\lambda})=+$, $\textrm{int}e^{k}_{\lambda}\times
I$ is disjoint from $\cup_{f=0}^{k} E^{k}_{f}\cup_{\gamma<\lambda}
e^{k}_{\gamma}\times I$.

(2) \ For $j\notin m(k)$,  $d_{j}^{k}$ is regular in
$\cup_{\gamma<j} e_{\gamma}^{k}\times I\cup_{f} E_{f}^{k}$.

(3) \ For $j\notin m(k)$ and $\gamma\in m(k)$, if $j<\gamma$, then
either $j\in I(w_{\gamma},k)$ with $s(w_{\gamma})=-$ or
$d_{j}^{k}$ is disjoint from $e^{k}_{\gamma}\times I$.

(4)  \ For each $i\geq k+1$, $c_{i}^{k}$ is a properly embedded
arc in $F^{k}$ lying in $E^{k}_{f}$ for some $f$. Furthermore, for
each $j\in\bigl\{1,\ldots,n\bigr\}-m(k)$,
$d_{j}^{k}-\cup_{\gamma<j} \textrm{int}e_{\gamma}^{k}\times I$
intersects $c_{i}^{k}$ in at most one point.

(5) \ For each $i\geq k+1$, $j\in\bigl\{1,\ldots,n\bigr\}-m(k)$
and $\gamma\in m(k)$, $\partial c_{i}^{k}\cap
d_{j}^{k}=\emptyset$, $\partial c_{i}^{k}\cap e_{\gamma}^{k}\times
I=\emptyset$. \vskip 0.5mm

{\bf \em Definition 2.3.3.} \ If $d_{j}^{k}-\cup_{\gamma<j}
\textrm{int}e_{\gamma}^{k}\times I$ intersects $c_{i}^{k}$ in one
point, then we say $i\in L(d_{j}^{k})$ and $j\in L(c_{i}^{k})$.
\vskip 0.5mm

{\bf \em Proposition 5.} \ $L(c^{0}_{i})=L(v_{i}^{0})$,
$L(d_{j}^{0})=L(w_{j}^{0})$. For each $1\leq k\leq m$ and $j\notin
m(k)$.

(1) If $j\notin L(c_{k}^{k-1})$, then $L(d_{j}^{k})=
L(d_{j}^{k-1})$.

(2) \ If $j\in L(c_{k}^{k-1})$, then $L(d_{j}^{k})=
L(d_{j}^{k-1})\cup L(d_{m^{k}}^{k-1})- L(d_{j}^{k-1})\cap
L(d_{m^{k}}^{k-1})$.\vskip 0.5mm

{\bf \em Proposition 6.} \ For each $0\leq k\leq m$, there are two
sets of pairwise disjoint disks $\bigl\{V_{i}^{k} \ | \ k+1\leq
i\leq m \ with \ s(v_{i})=- \bigr\}$ properly embedded in $\cal
V_{-}$ and $\bigl\{W_{j}^{k} \ | \ j\in
\bigl\{1,\ldots,n\bigr\}-m(k) \ with \ s(w_{j})=- \bigr\}$
properly embedded in $\cal W_{-}$ such that

(1) \ $V_{i}^{k}\cap F^{k}=c_{i}^{k}\cup_{r\in I(v_{i},k)}
c_{r}^{k}$, $W_{j}^{k}\cap F^{k}=d_{j}^{k}\cup_{r\in I(w_{j},k)}
d_{r}^{k}$;

(2) \ $V_{i}^{k}\cap W_{j}^{k}=V_{i}^{k}\cap W_{j}^{k}\cap F^{k}$.
\vskip 0.5mm

{\bf Now we prove Theorem 1 under the assumption that Propositions
1-6 are true.}

{\bf \em The proof of Theorem 1.} Suppose that $k=m$. By the
definition of $I(v,k)$,
$I(v,m)=I(v)-\bigl\{1,\ldots,m\bigr\}=\emptyset$. By Proposition
3(2), $V_{x}^{m}\cap W_{x}^{m}=(V_{x}^{m}\cap W_{x}^{m}\cap
P^{m})\cup\bigl\{x\bigr\}$. By Proposition 3(1), $V_{x}^{m}\cap
P^{m}=\cup_{r\in I(v,m)} v_{r}^{m}$, $W_{x}^{m}\cap
P^{m}=\cup_{r\in I(w,m)} w_{r}^{m}$. Hence $V_{x}^{m}\cap
W_{x}^{m}\cap P^{m}=\emptyset$ and $V_{x}^{m}\cap
W_{x}^{m}=\bigl\{x\bigr\}$. This means that $\cal V_{+}\cup W_{+}$
is stabilized.\qquad Q.E.D.\vskip 4mm

{\bf Remark.} \ Though Theorem 1 follows immediately from
Proposition 3. But Propositions 1, 2, 4, 5 and 6 are necessary in
the inductive proof of Proposition 3. See Section 2.5

\subsection{ The organizing of the inductive proofs of
Propositions 1-6}

We organize the proofs of Propositions 1-6 as follows:

In Section 2.5, we shall prove  Propositions 1-6 for $k=0$.
Furthermore, we shall introduce the ideas of the inductive proofs
of Propositions 1-6.

In Chapter 3, we shall study properties of $\bigl\{v_{i}^{k},
w_{j}^{k}, b^{k}_{\gamma}\times I, P^{k}\bigr\}$ and
$\bigl\{c^{k}_{i}, d_{j}^{k}, e^{k}_{\gamma}\times I,
F^{k}\bigr\}$ under the assumptions that Propositions 1-6 are true
for $k\leq l$.

In Chapter 4, we shall prove that Propositions 4-6 are true for
$k=l+1$ under the assumptions: $s(v_{l+1})=+$ and Propositions 1-6
are true  for $k\leq l$.

In Chapter 5, we shall prove that Propositions 4-6 are true for
$k=l+1$ under the assumptions: $s(v_{l+1})=-$ and Propositions 1-6
are true  for $k\leq l$.

In Chapter 6, we shall finish the proofs of Propositions 1-6.

\subsection{ The proofs of Propositions 1-6 for $k=0$}

{\bf \em Lemma 2.5.1.} \   Propositions 1-3 are true for $k=0$.

{\bf \em Proof.} \ By the definition, $B_{\cal V}$ is a disk in
$\partial_{+} \cal V_{+}$. Hence $B_{\cal V}$ is a surface
generated  by an trivial abstract tree $B_{\cal V}\cup \emptyset$.
Now Let $v_{i}^{0}=v_{i}$ $w_{j}^{0}=w_{j,v}$ for $1\leq i\leq m$
and $1\leq j\leq n$. By Lemma 2.2.8, $v_{i}^{0}$ intersects
$w_{j}^{0}$ in at most one point and $\partial v_{i}^{0}\cap
w_{j}^{0}=\emptyset$. Specially, let
$m(0)=\bigl\{m^{0}\bigr\}=\emptyset$. Thus Proposition 1 holds.

We denote by $W_{j}^{0}$ the disk $W_{j,v}$ for $s(w_{j})=+$,
$W_{x}^{0}$ the disk  $W_{x,v}$, $V_{i}^{0}$ the disk $V_{i}$ for
$s(v_{i})=+$, $V_{x}^{0}$ the disk $V_{x}$. Then, by Lemma 2.2.9,
we obtain two sets of pairwise disjoint properly embedded disks
$\bigl\{V_{i}^{0} \ | \ m\geq i\leq 1,
s(v_{i})=+\bigr\}\cup\bigl\{V_{x}^{0}\bigr\}$ in $\cal V_{+}$ and
$\bigl\{W_{j}^{0} \ | \ j\in\bigl\{1,\ldots,n\bigr\}- m(0),
s(w_{j})=+\bigr\}\cup\bigl\{W_{x}^{0}\bigr\}$ in $\cal W_{+}$. By
Definition 2.3.1, $I(v_{i},0)=I(v_{i})$, $I(v,0)=I(v)$,
$I(w_{j},0)=I(w_{j})-\emptyset=I(w_{j})$,
$I(w,0)=I(w)-\emptyset=I(w)$. By Lemma 2.2.9, Proposition 3
holds.\qquad Q.E.D.

{\bf \em Lemma 2.5.2.} Propositions 4-6 are true for $k=0$.

{\bf \em Proof.} \ By the definitions of $\cal W_{-}$ and $\cal
V_{-}$, $B_{\cal V}$ is a disk in $\partial_{+} \cal V_{-}$. Hence
$B_{\cal V}$ is a surface generated by an trivial abstract tree
$B_{\cal V}\cup \emptyset$ in $\partial_{+} V_{-}$. Now Let
$c_{i}^{0}=v_{i}$ $d_{j}^{0}=w_{j,v}$. By Lemma 2.2.8, $c_{i}^{0}$
intersects $d_{j}^{0}$ in at most one point and $\partial
c_{i}^{0}\cap d_{j}^{0}=\emptyset$. It is easy to see that
$L(c_{i}^{0})=L(v_{i}^{0})$ and $L(d_{j}^{0})=L(w_{j}^{0})$. Hence
Propositions 4 and 5 holds.

We denote by $W_{j}^{0}$ the disk $W_{j,v}$ for $s(w_{j})=-$,
$V_{i}^{0}$ the disk $V_{i}$ for $s(v_{i})=-$. Then, by Lemma
2.2.9, we obtain two sets of pairwise disjoint properly embedded
disks $\bigl\{V_{i}^{0} \ | \ m\geq i\leq 1, s(v_{i})=-\bigr\}$ in
$\cal V_{-}$ and $\bigl\{W_{j}^{0} \ | \
j\in\bigl\{1,\ldots,n\bigr\}- m(0), s(w_{j})=-\bigr\}$ in $\cal
W_{-}$. Note that $m(0)=\emptyset$. By Definition 2.3.1,
$I(v_{i},0)=I(v_{i})$, $I(w_{j},0)=I(w_{j})$. By Lemma 2.2.9,
Proposition 6 holds.\qquad Q.E.D.

{\bf \em The ideas of the proofs of Propositions 1-6.}

By Lemmas 2.5.1 and 2.5.2, Propositions 1-6 are true for $k=0$.
Now we may assume that Propositions 1-6 are true for $k\leq l$. We
only need to prove Propositions 1-6 are true for $k=l+1$. The
inductive proofs depend on $s(v_{l+1})$.

Case 1. $s(v_{l+1})=+$.

By Proposition 1, $v_{l+1}^{l}$ is a properly embedded arc in
$P^{l}$. Let $m^{l+1}=MinL(v_{l+1}^{l})$. By Definition 2.3.2 and
Proposition 1, either $m^{l+1}=\emptyset$ or $m^{l+1}\notin m(l)$.
Now let $m(l+1)=m(l)\cup\bigl\{m^{l+1}\bigr\}$. By Proposition 3,
there is a properly embedded disk $V_{l+1}^{l}$ in $\cal V_{+}$
such that $V_{l+1}^{l}\cap P^{l}=c^{l}_{l+1}$. By Proposition 1,
$w_{m^{l+1}}^{l}$  is an arc in $P^{l}$ if $m^{l+1}\neq\emptyset$.
We can prove that $c_{l+1}^{l}$ intersects $w_{m^{l+1}}^{l}$ in
only one point, and $s(w_{m^{l+1}})=-$. Let $N(\partial
V_{l+1}^{l})$ be a regular neighborhood of $\partial V_{l+1}^{l}$
in $\partial \cal V_{+}$, and $N(w_{m^{l+1}}^{l})$ be a regular
neighborhood of $w_{m^{l+1}}^{l}$ in $P^{l}$. In fact,
$P^{l+1}=(P^{l}-N(w_{m^{l+1}}^{l}))\cup N(\partial V_{l+1}^{l})$,
$w_{j}^{l+1}=w_{j}^{l}$ for $j\notin m(l+1)$, and
$v^{l+1}_{i}$($i\geq l+2$) is the band sum of $v^{l}_{i}$ with
some copies of $\partial V_{l+1}^{l}$ along $w_{m^{l+1}}^{l}$.
Specially, $W_{x}^{l+1}=W_{x}^{l}$ and $W_{j}^{l+1}=W_{j}^{l}$ for
$j\notin m(l+1)$ with $s(w_{j})=+$; $V_{x}^{l+1}$  is the
connected sum of $V_{x}$ with some copies of $V_{l+1}^{l}$,  and
$V_{i}^{l+1}$ is the connected sum of $V^{l}_{i}$ with some copies
of $V_{l+1}^{l}$ for $i\geq l+2$ with $s(v_{i})=+$.

Note that Propositions 1-6 are true for $k\leq l$. By Propositions
2, 4 and Lemmas 2.5.1, 2.5.2, we can prove that
$L(v_{i}^{l})=L(c^{l}_{i})$ and $L(d_{j}^{l})=L(w_{j}^{l})$. Since
$s(w_{m^{l+1}})=-$, by Proposition 6, $W_{m^{l+1}}^{l}$ is a
properly embedded disk in $\cal W_{-}$ such that
$W_{m^{l+1}}^{l}\cap F^{l}=d_{m^{l+1}}^{l}\cup_{r\in
I(w_{m^{l+1}},l)} d_{r}^{l}$. We can also prove that
$d_{m^{l+1}}^{l}$ intersects $c_{l+1}^{l}$ in only one point if
$m^{l+1}\neq\emptyset$. Let $N(\partial W_{m^{l+1}}^{l})$ be a
regular neighborhood of $\partial W_{m^{l+1}}^{l}$ in $\partial
\cal V_{-}$, and $N(c_{l+1}^{l})$ be a regular neighborhood of
$c_{l+1}^{l}$ in $F^{l}$. In fact,
$F^{l+1}=(F^{l}-N(c_{l+1}^{l})\cup N(\partial W_{m^{l+1}}^{l})$,
$c_{i}^{l+1}=c^{l}_{i}$ for $i\geq l+2$, and $d_{j}^{l+1}$ is the
connected sum of $d_{j}^{l}$ with some copies of $\partial
W_{m^{l+1}}^{l}$ along $c_{l+1}^{l}$ for $j\notin m(l+1)$.
Specially, $V_{i}^{l+1}=V_{i}^{l}$ for $i\geq l+2$ with
$s(v_{i})=-$, and $W_{j}^{l+1}$ is the connected sum of
$W^{l}_{j}$ with some copies of $W_{m^{l+1}}^{l}$ for $j\notin
m(l+1)$ with $s(w_{j})=-$.

Case 2. $s(v_{l+1})=-$.

Now let $m^{l+1}=MinL(c_{l+1}^{l})$. By Proposition 4, we can
prove that $s(w_{m^{l+1}})=+$. Thus we have an alternating proof
with the one of Case 1. So Propositions 4-6 are necessary  to show
$s(w_{m^{l+1}})=+$ when $s(v_{l+1})=-$.\qquad Q.E.D.

\section{Some properties of  $\bigl\{c^{l}_{i},
d_{j}^{l}, e^{l}_{\gamma}\times I, F^{l}\bigr\}$}

\ \ \ \ \ By the argument in Section 2.5, Propositions 1-6 holds
for $k=0$. In this chapter, we shall give some properties of
$\bigl\{v_{i}^{l}, w_{j}^{l}, b^{l}_{\gamma}\times I,
P^{l}\bigr\}$ and $\bigl\{c^{l}_{i}, d_{j}^{l},
e^{l}_{\gamma}\times I, F^{l}\bigr\}$ under the  assumptions that
Propositions 1-6 hold for $k\leq l$.

We first consider $\bigl\{c^{l}_{i}, d_{j}^{l},
e^{l}_{\gamma}\times I, F^{l}\bigr\}$.

\subsection{The intersection of $c^{l}_{i}$ and $d_{j}^{l}$}

\ \ \ \ \ {\bf \em Lemma 3.1.1.} \ Suppose that  $j\in
\bigl\{1,\ldots,n\bigr\}-m(l)$, and $j\notin I(w_{\gamma},l)$ for
each $\gamma\in m(l)$ with $s(w_{\gamma})=-$. Then $d_{j}^{l}$ is
regular in $F^{l}$.

{\bf \em Proof.}  \ Now by Proposition 4(3), if $j<\gamma$, then
$d_{j}^{l}$ is disjoint from $e^{l}_{\gamma}\times I$. By
Proposition 4(2), the lemma holds.\qquad Q.E.D.

{\bf \em Lemma 3.1.2.} \ For each $j\in
\bigl\{1,\ldots,n\bigr\}-m(k)$,
$d_{j}^{l}=\cup_{i=1}^{\theta(j)}d_{j,f_{i,j}}^{l}
\cup_{i=1}^{\theta(j)-1} e_{\gamma_{i,j}}$ satisfying the
following conditions:

(1) \ $j>\gamma_{i,j}\in m(l)$, and $e_{\gamma_{i,j}}$ is a core
of $e^{l}_{\gamma_{i,j}}\times (0,1)$ for $1\leq
i\leq\theta(j)-1$.

(2) \ $0\leq f_{i,j}\leq l$, $d_{j,f_{i,j}}^{l}$ is a properly
embedded arc in $E^{l}_{f_{i,j}}$ which is disjoint from
$\cup_{\gamma<j} inte^{l}_{\gamma}\times I$ for $1\leq i\leq
\theta(j)$.

(3) \ For $i\neq r$, $\gamma_{i,j}\neq \gamma_{r,j}$, $f_{i,j}\neq
f_{r,j}$.

(4) \  $\partial_{1} e_{\gamma_{i,j}}=\partial_{2} d_{j,
f_{i,j}}^{l}$, $\partial_{2} e_{\gamma_{i,j}}=\partial_{1} d_{j,
f_{i+1,j}}^{l}$, $\partial_{1} d_{j}^{l}=\partial_{1} d_{j,
f_{1,j}}^{l}$, $\partial_{2} d_{j}^{l}=\partial_{2} d_{j,
f_{\theta(j),j}}^{l}$.

(5) \ For each $f$, $d_{j}^{l}-\cup_{\gamma<j}
inte^{l}_{\gamma}\times I$ intersects $E^{l}_{f}$ in at most one
component, denoted by $d_{j,f}^{l}$.  Furthermore, if $f=f_{i,j}$,
then $d_{j,f}^{l}=d_{j,f_{i,j}}^{l}$; if $f\neq f_{i,j}$ for each
$1\leq i\leq \theta(j)$, then $d_{j,f}^{l}=\emptyset$.

{\bf \em Proof.} \ By Proposition 4(2), $d_{j}^{l}$ is regular in
$\cup_{\gamma<j} e^{l}_{\gamma}\times I\cup_{f} E_{f}^{l}$. Hence
the lemma follows from Lemma 2.1.6. \qquad Q.E.D.\vskip 0.5mm

{\bf \em Remark.} \ It is possible that $\theta(j)=1$. In this
case, $d_{j}^{l}=d_{j,f_{1,j}}^{l}$. Hence $d_{j}^{l}$ is disjoint
from $e^{l}_{\gamma}\times I$ for $\gamma<j$.\vskip 0.5mm

Now by Proposition 4(4), $d_{j,f_{i,j}}^{l}$ intersects
$c^{l}_{i}$ in at most one point.

{\bf \em Definition 3.1.3.} \ We say $i\in L(d_{j,f_{i,j}}^{l})$
if $d_{j,f_{i,j}}^{l}$ intersects $c^{l}_{i}$ in  one point.\vskip
3mm

{\bf \em Lemma 3.1.4.} \ (1) If $r\neq i$, then
$L(d_{j,f_{i,j}}^{l})\cap L(d_{j,f_{r,j}}^{l})=\emptyset$.

(2) \ $L(d_{j}^{l})=\cup_{i} L(d_{j,f_{i,j}}^{l})$.

{\bf Proof.} By Proposition 4(4),  if $d_{j,f_{i,j}}^{l}$
intersects $c^{l}_{i}$ in  one point, then $c^{l}_{i}\subset
E^{l}_{f_{i,j}}$. By Lemma 3.1.2(3), (1) holds.

By Lemma 3.1.2(5), each component of
$d_{j}^{l}-\cup_{\gamma<j}inte^{l}_{\gamma}\times I$ is
$d_{j,f_{i,j}}^{l}$ for some $f_{i,j}$. Hence (2) holds. \qquad
Q.E.D. \vskip 0.5mm

{\bf \em Lemma 3.1.5.} \ $L(v_{i}^{l})=L(c^{l}_{i})$ and
$L(d_{j}^{l})=L(w_{j}^{l})$.

{\bf \em Proof.} \   By Lemmas 2.5.1 and 2.5.2,
$L(w_{j}^{0})=L(d_{j}^{0})$. $L(c_{i}^{0})=L(v_{i}^{0})$. Now by
Propositions 2 and 5, the lemma holds. \qquad Q.E.D. \vskip 0.5mm

{\bf \em Lemma 3.1.6.} \  (1) \ Each component of $d_{j}^{l}\cap
E^{l}_{f}$ is either $d_{j,f}^{l}$  or a core of
$d_{r,f}^{l}\times (0,1)$ for some $\gamma\in m(l)$ and $r\in
I(w_{\gamma},l)$. Furthermore, $r<\gamma<j$, $s(w_{\gamma})=-$.

(2) \ If $d_{j}^{l}\cap c^{l}_{i}\neq \emptyset$, then either
$j\in L(c^{l}_{i})$ or $r\in L(c^{l}_{i})$ for some $r<j$.

{\bf \em Proof.} \ (1) \ By Lemma 3.1.2,
$d_{j}^{l}=\cup_{i=1}^{\theta(j)}d_{j,f_{i,j}}^{l}
\cup_{i=1}^{\theta(j)-1} e_{\gamma_{i,j}}$ such that
$\gamma_{i,j}<j$ and $\gamma_{i,j}\in m(l)$. Now if one component
$c$ of $d_{j}^{l}\cap E^{l}_{f}$ is not $d_{j,f}^{l}$, then, by
Lemma 3.1.2(2) and (4), $c\subset inte_{\gamma_{i,j}}\cap
E^{l}_{f}$ for some $i$ where $e_{\gamma_{i,j}}$ is a core of
$e^{l}_{\gamma_{i,j}}\times (0,1)$. Since $
inte_{\gamma_{i,j}}\cap E^{l}_{f}\neq\emptyset$. By Proposition
4(1), $s(w_{\gamma_{i,j}})=-$ and $inte^{l}_{\gamma_{i,j}}\times
I\cap (\cup_{f} E^{l}_{f}\cup_{\gamma<\gamma_{i,j}}
e^{l}_{\gamma}\times I)=\cup_{r\in I(w_{\gamma_{i,j}},l)}
d_{r}^{l}\times I$. Hence $c\subset e_{\gamma_{i,j}}\cap
d_{r}^{l}\times I$ for $r\in I(w_{\gamma_{i,j}},l)$. By
Proposition 4(2), $e_{\gamma_{i,j}}\cap d_{r}^{l}\times I$ is a
core of $d_{r}^{l}\times (0,1)$.  Note that $r<\gamma_{i,j}<j$. By
induction, (1) holds.

(2) \ Let $p$ be a point in $d_{j}^{l}\cap c^{l}_{i}$. Then, by
(1), either $p\in d_{j,f}^{l}\cap c^{l}_{i}$ or $p\in c\cap
c^{l}_{i}$ where $c$ is a core of $d_{r}^{l}\times (0,1)$ for some
$\gamma\in m(l)$ and $r\in I(w_{\gamma},l)$. By Proposition 4(5),
$c^{l}_{i}$ intersects $c$ in one point if and only if it
intersects $d_{r}^{l}=d_{r}^{l}\times\bigl\{0\bigr\}$ in one
point. By induction, (2) holds. \qquad Q.E.D. \vskip 0.5mm

\subsection{The intersection of  $d_{j}^{l}$ and $e^{l}_{\gamma}
\times I$}\vskip 0.5mm

\ \ \ \ \ {\bf For each $j\notin m(l)$ and $\gamma\in m(l)$,
either $j\in I(w_{\gamma},l)$ or $j\notin I(w_{\gamma},l)$. We may
assume that $j\in I(w_{\gamma_{j}},l)$. In this case, if $j\notin
I(w_{\gamma},l)$ for each $\gamma\in m(l)$, then
$\gamma_{j}=\emptyset$.}

{\bf \em Definition 3.2.1.} \ (1) \ Let $(d_{j}^{l}\times
I)_{\gamma_{j}}=d_{j}^{l}\times I\subset e^{l}_{\gamma_{j}}\times
I$ if $j\in I(w_{\gamma_{j}},l)$ and $s(w_{\gamma_{j}})=-$,

(2) \ let $(d_{j}^{l}\times I)_{\gamma_{j}}=d_{j}^{l}$ if $j\notin
I(w_{\gamma},l)$ for each $\gamma\in m(l)$ or $j\in
I(w_{\gamma_{j}},l)$ with $s(w_{\gamma_{j}})=+$.\vskip 4mm

Definition 3.2.1 means that $(d_{j}^{l}\times I)_{\gamma_{j}}$ is
either a disk or an arc.

Recalling the equality:
$d_{j}^{l}=\cup_{i=1}^{\theta(j)}d_{j,f_{i,j}}^{l}
\cup_{i=1}^{\theta(j)-1} e_{\gamma_{i,j}}$ for each $j\notin m(l)$
in Lemma 3.1.2, and $d_{j,f}^{l}$ in Lemma 3.1.2(5).\vskip 0.5mm

{\bf \em Lemma 3.2.2.} \ If $r\neq j$, then $(d_{j,f}^{l}\times
I)_{\gamma_{j}}$ is disjoint from $(d_{r,f}^{l}\times
I)_{\gamma_{r}}$ for each $0\leq f\leq l$.

{\bf \em Proof.} \   Without loss of generality, we may assume
that $j\in I(w_{\gamma_{j}},l)$,  $r\in I(w_{\gamma_{r}},l)$ and
$s(w_{\gamma_{j}})=s(w_{\gamma_{r}})=-$.

Assume that $\gamma_{j}=\gamma_{r}$. Since $j\neq r$,
$(d_{j,f}^{l}\times I)_{\gamma_{j}}$ is disjoint from
$(d_{r,f}^{l}\times I)_{\gamma_{r}}$.

Assume  now that $\gamma_{j}>\gamma_{r}$.  If $(d_{j,f}^{l}\times
I)_{\gamma_{j}}\cap(d_{r,f}^{l}\times I)_{\gamma_{r}}\neq
\emptyset$, then $e_{\gamma_{j}}^{l}\times I\cap
e_{\gamma_{r}}^{l}\times I\neq\emptyset$. By Proposition 4(1),
$F^{l}$ is generated by an abstract tree. By Definition 2.1.4,
each component of $e^{l}_{\gamma_{j}}\times I\cap
e^{l}_{\gamma_{r}}\times I$ is $(c\times I)_{\gamma_{j}}\subset
e^{l}_{\gamma_{r}}\times (0,1)$ where $c\subset
inte_{\gamma_{j}}^{l}$ is a core of $e^{l}_{\gamma_{r}}\times
(0,1)$. Since $d_{j,f}^{l}$ and $d_{r,f}^{l}$ are properly
embedded in $E^{l}_{f}$,  $(d_{j,f}^{l}\times
I)_{\gamma_{j}}\subset (d_{r,f}^{l}\times (0,1))_{\gamma_{r}}$.
Hence $d_{j,f}^{l}=d_{j,f}^{l}\times\bigl\{0\bigr\}\subset
e^{l}_{\gamma_{r}}\times (0,1)$. By Lemma 2.2.4, $j\notin
I(w_{\gamma_{r}},l)$. By Proposition 4(3), $j>\gamma_{r}$. Note
that $d_{j,f}^{l}=d_{j,f}^{l}\times\bigl\{0\bigr\}\subset
(d_{r,f}^{l}\times (0,1))_{\gamma_{r}}\subset
e^{l}_{\gamma_{r}}\times (0,1)$, contradicting Lemma
3.1.2(2).\qquad Q.E.D.\vskip 0.5mm

{\bf \em Lemma 3.2.3.} \  If $d_{j}^{l}\cap e^{l}_{\gamma}\times
I\neq \emptyset$, then either $\gamma\leq
Max\bigl\{\gamma_{1,j},\ldots,\gamma_{\theta(j)-1,j}\bigr\}$ or
$j\in I(w_{\gamma},l)$ with $s(w_{\gamma})=-$.

{\bf \em Proof.} \  Suppose that $\gamma\neq \gamma_{i,j}$ for
$1\leq i\leq \theta(j)-1$ and $j\notin I(w_{\gamma},l)$. If
$\gamma>j$, then, by Proposition 4(3), $d_{j}^{l}$ is disjoint
from $e^{l}_{\gamma}\times I$.

Suppose that $\gamma<j$. Then, by  Lemma 3.1.2,
$d_{j,f_{i,j}}^{l}$ is disjoint from
$\textrm{int}e^{l}_{\gamma}\times I$ for each $1\leq i\leq
\theta(j)$.  By Proposition 4(2), each component of $d_{j}^{l}\cap
e^{l}_{\gamma}\times I$  is a core of $e^{l}_{\gamma}\times
(0,1)$.  Since $F^{l}$ is generated by an abstract tree, $\partial
e^{l}_{\gamma}\times I\cap
\partial e^{l}_{\lambda}\times I=\emptyset$ for each $\gamma\neq\lambda\in m(l)$.
Since $\gamma\neq\gamma_{i,j}$, each component of $d_{j}^{l}\cap
e^{l}_{\gamma}\times I$ lies in $inte_{\gamma_{i,j}}$. (See the
equality of $d_{j}^{l}$ in Lemma 3.1.2.) Thus each component of
$e^{l}_{\gamma_{i,j}}\cap e^{l}_{\gamma}\times I$ is a core of
$e^{l}_{\gamma}\times I$. By Definition 2.1.4,
$\gamma<\gamma_{i,j}$. \qquad Q.E.D. \vskip 0.5mm

{\bf \em Lemma 3.2.4.} \ If $\textrm{int}e^{l}_{\gamma}\times
I\cap E^{l}_{f}\neq\emptyset$ for some $\gamma\in m(l)$ and some
$0\leq f\leq l$, then each component of
$\textrm{int}e^{l}_{\gamma}\times I\cap E^{l}_{f}$ is either
$(d_{r,f}^{l}\times I)_{\gamma}$ for some $r\in I(w_{\gamma},l)$
or $(c\times I)_{\gamma}\subset (d_{j,f}^{l}\times
(0,1))_{\gamma_{j}}$ where $c\subset inte^{l}_{\gamma}$ is a core
of $d_{j,f}^{l}\times (0,1)$ for some $\gamma_{j}\in m(l)$ and
$j\in I(w_{\gamma_{j}},l)$. Furthermore, $j<\gamma_{j}<\gamma$,
$s(w_{\gamma})=s(w_{\gamma_{j}})=-$.

{\bf \em Proof.} \ Now by Proposition 4(1), $s(w_{\gamma})=-$ and
if $C$ is a component of $\textrm{int}e^{l}_{\gamma}\times I\cap
E^{l}_{f}$, then $C\subset (d_{r}^{l}\times I)_{\gamma}$ for some
$r\in I(w_{\gamma},l)$. Hence $r<\gamma$. If $C\neq
(d_{r,f}^{l}\times I)_{\gamma}$. Then, by Lemma 3.1.2, $C\subset
\textrm{int}e^{l}_{\gamma_{i,r}}\times I\cap E^{l}_{f}$ for some
$1\leq i\leq \theta(r)-1$. By induction, the lemma holds.\qquad
Q.E.D.

{\bf \em Corollary 3.2.5.} \ If $L(c_{l+1}^{l})=\emptyset$, then
$d_{j}^{l}$ and $e^{l}_{\gamma}\times I$ are disjoint from
$c_{l+1}^{l}$ for each $j\in\bigl\{1,\ldots,n\bigr\}-m(l)$ and
$\gamma\in m(l)$.

The corollary follows from Lemma 3.1.6 and Lemma 3.2.4 and
Proposition 4(5).

{\bf \em Lemma 3.2.6.} \ (1) \ Suppose that $j<\gamma<\gamma_{j}$
and $j\in I(w_{\gamma_{j}},l)$. Then $(d_{j}^{l}\times
I)_{\gamma_{j}}$ is disjoint from $e^{l}_{\gamma}\times I$.

(2) \ Suppose that $j>\gamma$. Then $(intd_{j,f}^{l}\times
I)_{\gamma_{j}}$ is disjoint from $e^{l}_{\gamma}\times I$.

{\bf \em Proof.} \ (1) \ Suppose that $s(w_{\gamma_{j}})=+$. Then,
by Proposition 4(3) and Definition 3.2.1, $(d_{j}^{l}\times
I)_{\gamma_{j}}=d_{j}^{l}$ is disjoint from $e^{l}_{\gamma}\times
I$.

Now suppose that $s(w_{\gamma_{j}})=-$. By Proposition 4(1) and
Definition 2.1.4, $e^{l}_{\gamma_{j}}\times I\cap
e^{l}_{\gamma}\times I$ is $(c\times I)_{\gamma_{j}}$ where
$c\subset inte^{l}_{\gamma_{j}}$ is a core of
$e^{l}_{\gamma}\times (0,1)$. By Lemma 3.2.3, $(d_{j}^{l}\times
I)_{\gamma_{j}}$ is disjoint from $e^{l}_{\gamma}\times I$.

(2) \ Suppose that $s(w_{\gamma_{j}})=+$. Then, by Lemma 3.1.2 and
Definition 3.2.1, $(intd_{j}^{l}\times
I)_{\gamma_{j}}=intd_{j}^{l}$ is disjoint from
$e^{l}_{\gamma}\times I$.

Now suppose that $s(w_{\gamma_{j}})=-$. By Proposition 4(1) and
Definition 2.1.4, each component of $e^{l}_{\gamma_{j}}\times
I\cap e^{l}_{\gamma}\times I$ is $(c\times I)_{\gamma_{j}}\subset
e^{l}_{\gamma}\times (0,1)$, where $c\subset
inte^{l}_{\gamma_{j}}$ is a core of $e^{l}_{\gamma}\times (0,1)$.
Since $intd_{j}^{l}$ is disjoint from $e^{l}_{\gamma}\times I$,
(2) holds. \qquad  Q.E.D.

\subsection{Properties of $d_{j}^{l}$  for $j\in L(c^{l}_{l+1})$}\vskip 0.5mm

\ \ \ \ \ {\bf \em Definition 3.3.0.} \ Let
$m^{l+1}=MinL(c_{l+1}^{l})$.\vskip 0.5mm

By Proposition 4(4), $c^{l}_{l+1}$ lies in $E^{l}_{f}$ for some
$0\leq f\leq l$.  Without loss of generality, we may assume that
$c_{l+1}^{l}\subset E^{l}_{0}$. Recalling
$d_{j,f}^{l}=(d_{j}^{l}-\cup_{\gamma<j} inte^{l}_{\gamma}\times
I)\cap E^{l}_{f}$ defined in Lemma 3.1.2(5). Now if $j\in
L(c_{l+1}^{l})$, then, by Proposition 4(4) and Lemma 3.1.2,
$d_{j,0}^{l}$ intersects $c_{l+1}^{l}$ in one point and $j\geq
m^{l+1}$.

Now we rearrange all the elements in $L(c^{l}_{l+1})$ as $\ldots,
j_{-1}, j_{0}=m^{l+1}, j_{1},\ldots$ according to the order of
$(\cup_{j\in L(c_{l+1}^{l})} d_{j,0}^{l})\cap c_{l+1}^{l}$ lying
in $c_{l+1}^{l}$ as in Figure 6.
\begin{center}
\includegraphics[totalheight=5cm]{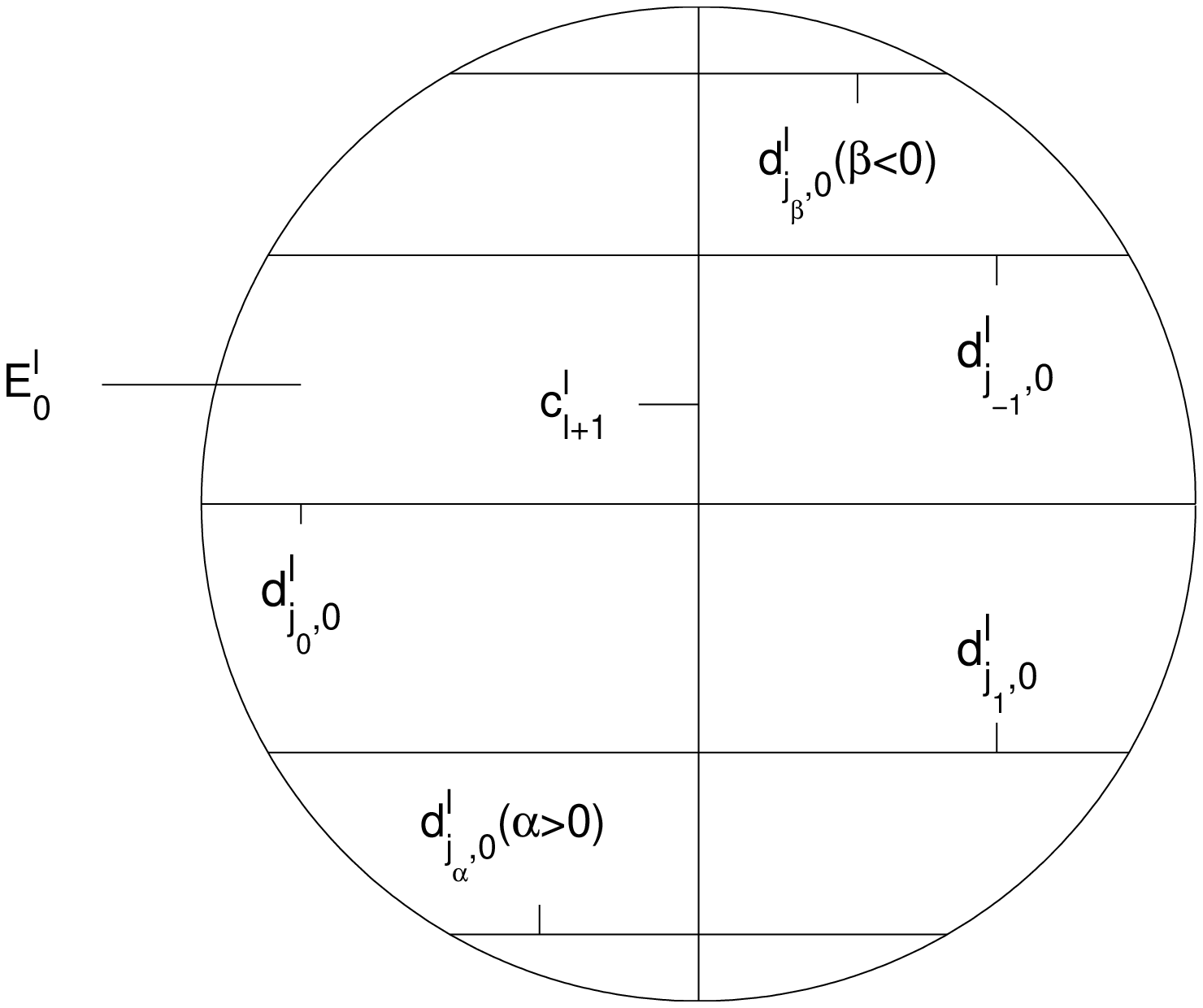}
\begin{center}
Figure 6
\end{center}
\end{center}

{\bf \em Lemma 3.3.1.} \ For each $j_{\alpha}\in L(c_{l+1}^{l})$,
$d_{j_{\alpha}}^{l}=d_{0,\alpha}\cup_{i=\delta(\alpha)}^{\theta(\alpha)}d_{i,\alpha}
\cup_{i=\delta(\alpha)}^{\theta(\alpha)} e_{i,\alpha}$ satisfying
the following conditions:

(1) \ For $\delta(\alpha)\leq i\neq 0\leq \theta(\alpha)$,
$e_{i,\alpha}$ is a core of $e^{l}_{\gamma_{i,\alpha}}\times
(0,1)$ for some $\gamma_{i,\alpha}<j_{\alpha}$.

(2) \ For $\delta(\alpha)\leq i\neq 0\leq \theta(\alpha)$ and
$d_{i,\alpha}$ is a properly embedded arc in
$E^{l}_{f_{i,\alpha}}$ which is disjoint from
$\cup_{\gamma<j_{\alpha}} inte^{l}_{\gamma}\times I$.

(3)  \ $d_{0,\alpha}$ is a properly embedded arc in $E^{l}_{0}$.

(4) \ For $i\neq r\neq 0$, $\gamma_{i,\alpha}\neq
\gamma_{r,\alpha}$, $f_{i,\alpha}\neq f_{r,\alpha}$, and
$f_{i,\alpha}\neq 0$.

(5) \ For $\delta(\alpha)\leq i\neq 0\leq \theta(\alpha)$,
$\partial_{1} e_{i,\alpha}=\partial_{2} d_{i,\alpha}$,
$\partial_{2} e_{i,\alpha}=\partial_{1} d_{i+1,\alpha}$,
$\partial_{1} e_{1,\alpha}=\partial_{2} d_{0, \alpha}$,
$\partial_{2} e_{-1,\alpha}=\partial_{1} d_{ 0,\alpha}$,
$\partial_{1} d_{j_{\alpha}}^{l}=\partial_{1}
d_{\delta(\alpha),\alpha}$, $\partial_{2}
d_{j_{\alpha}}^{l}=\partial_{2} d_{\theta(\alpha),\alpha}$.

{\bf \em Proof.} \ Since $j_{\alpha}\in L(c_{l+1}^{l})$,
$d_{j_{\alpha},0}^{l}\neq\emptyset$, by Lemma 3.1.2,
$d_{j_{\alpha}}^{l}=\cup_{i=\delta(j_{\alpha})}^{-1}
d_{j_{\alpha},f_{i,j_{\alpha}}}^{l}\cup
d_{j_{\alpha},0}^{l}\cup_{i=1}^{\theta(j_{\alpha})}d_{j_{\alpha},f_{i,j_{\alpha}}}^{l}\cup_{i=\delta(j_{\alpha})}^{-1}
e_{\gamma_{i,j_{\alpha}}}\cup_{i=1}^{\theta(j_{\alpha})}
e_{\gamma_{i,j_{\alpha}}}$. Now we denote by $\delta(\alpha)$ the
integer $\delta(j_{\alpha})$,  $\theta(\alpha)$ the integer
$\theta(j_{\alpha})$, $d_{i,\alpha}$ the arc
$d_{j_{\alpha},f_{i,j_{\alpha}}}^{l}$, $e_{i,\alpha}$ the arc
$e_{\gamma_{i,j_{\alpha}}}$ for $i\neq 0$. Let
$f_{i,\alpha}=f_{i,j_{\alpha}}$ and
$\gamma_{i,\alpha}=\gamma_{i,j_{\alpha}}$. In particular, we
denote by $d_{0,\alpha}$ the arc $d_{j_{\alpha},0}^{l}$. By Lemma
3.1.2, the lemma holds.\qquad Q.E.D.\vskip 0.5mm

{\bf \em Remark 3.3.2.} \ 1) \ In order to simplify the
formulation, we shall write as
$d_{j_{\alpha}}^{l}=\cup_{i=\delta(\alpha)}^{\theta(\alpha)}d_{i,\alpha}
\cup_{i=\delta(\alpha)}^{\theta(\alpha)} e_{i,\alpha}$. In the
formulation,
"$\cup_{i=\delta(\alpha)}^{\theta(\alpha)}d_{i,\alpha}$" means
"the union index from  ${\delta(\alpha)}$ to $-1$, then 0, then
$1$, and then to ${\theta(\alpha)}$", but
"$\cup_{i=\delta(\alpha)}^{\theta(\alpha)} e_{i,\alpha}$" means
"the union index from  ${\delta(\alpha)}$ to $-1$, then $1$, and
then to ${\theta(\alpha)}$".

2) \ If $\delta(\alpha)=\theta(\alpha)=0$, then
$d_{j_{\alpha}}^{l}=d_{0,\alpha}$.

{\bf \em Lemma 3.3.3.} \ (1) \ $d_{j}^{l}\cap c^{l}_{l+1}\subset
(\cup_{\alpha} (d_{0,\alpha}\times I)_{\gamma_{\alpha}})\cap
c_{l+1}^{l}$ where $j\in\bigl\{1,\ldots,n\bigr\}-m(l)$,
$j_{\alpha}\in L(c_{l+1}^{l})$, and $(d_{0,\alpha}\times
I)_{\gamma_{\alpha}}$ is defined in Definition 3.2.1.

(2) \ $e^{l}_{\gamma}\times I\cap c^{l}_{l+1}\subset
(\cup_{\alpha} (d_{0,\alpha}\times I)_{\gamma_{\alpha}})\cap
c_{l+1}^{l}$ for $\gamma\in m(l)$.

{\bf \em Proof.} \ (1) Now each component of $d_{j}^{l}\cap
c^{l}_{l+1}$ is a point. By assumption,  $c_{l+1}^{l}\subset
E^{l}_{0}$. If $d_{j}^{l}\cap c^{l}_{l+1}\neq\emptyset$, then
$d_{j}^{l}\cap c^{l}_{l+1}=(d_{j}^{l}\cap E^{l}_{0})\cap
c^{l}_{l+1}$. By Lemma 3.1.6, each component of $d_{j}^{l}\cap
E^{l}_{0}$ is either $d_{j,0}^{l}$  or a core of
$d_{r,0}^{l}\times (0,1)$ for some $\gamma\in m(l)$ and $r\in
I(w_{\gamma},l)$. Suppose that $d$ is a component of
$d_{j}^{l}\cap E^{l}_{0}$ such that $d\cap
c_{l+1}^{l}\neq\emptyset$. If $d=d_{j,0}^{l}$, then $j\in
L(c_{l+1}^{l})$. If $d=d_{r,0}^{l}$, then $r\in L(c_{l+1}^{l})$.

(2) follows from Lemma 3.2.4 and Proposition 4(5).\qquad
Q.E.D.\vskip 0.5mm

\begin{center}
\includegraphics[totalheight=4cm]{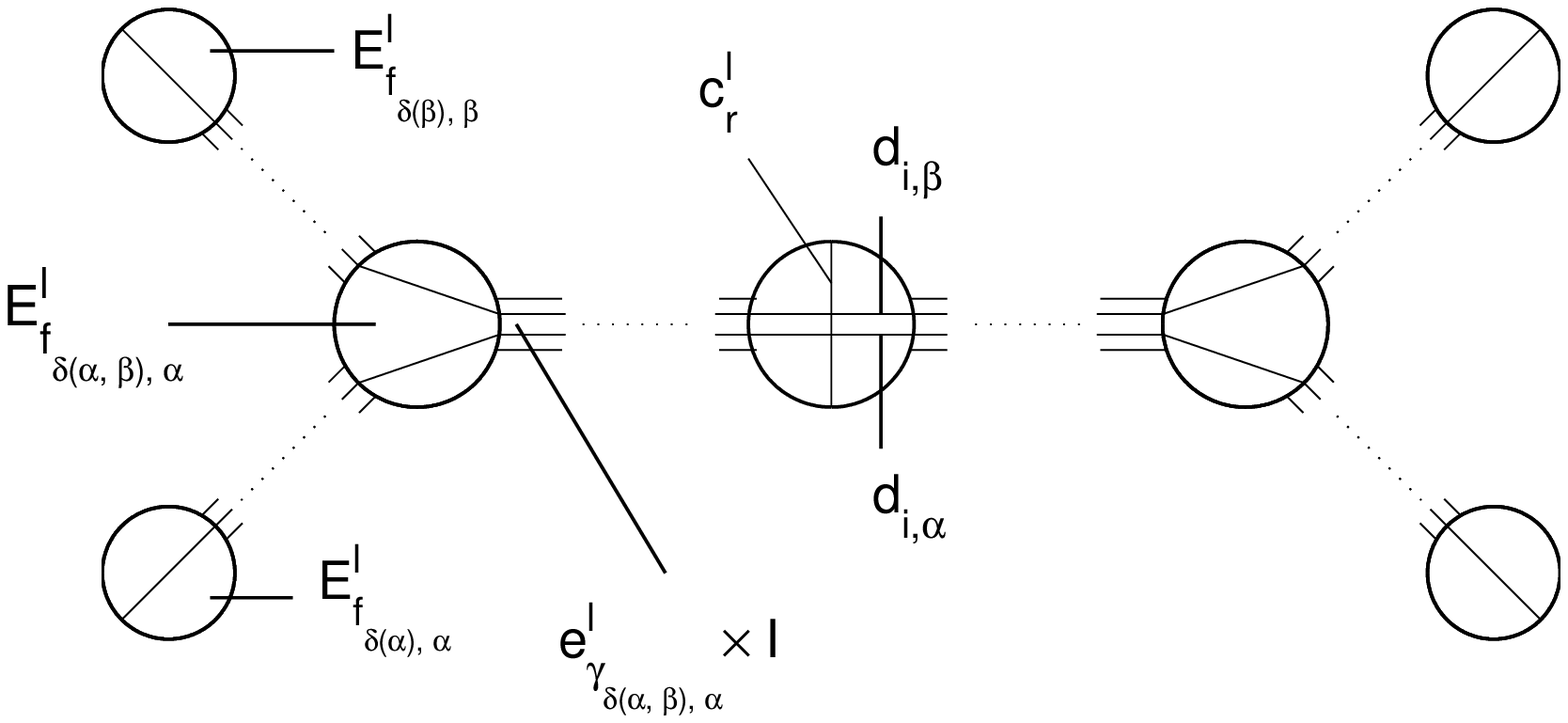}
\begin{center}
Figure 7
\end{center}
\end{center}
\vskip 4mm

{\bf \em Lemma 3.3.4.} \ If $j_{\alpha},j_{\beta} \in
L(c^{l}_{l+1})$, then there are two integers
$\delta(\alpha,\beta)\leq 0$  and $\theta(\alpha,\beta)\geq 0$
satisfying the following conditions:

(1) \ $\delta(\alpha),\delta(\beta)\leq
\delta(\alpha,\beta)\leq\theta(\alpha,\beta)\leq
\theta(\alpha),\theta(\beta)$.

(2) \ If $\delta(\alpha,\beta)\leq i\leq \theta(\alpha,\beta)$,
then $\gamma_{i,\alpha}=\gamma_{i,\beta}$ and
$f_{i,\alpha}=f_{i,\beta}$.

(3) \ $\bigl\{\gamma_{i,\alpha} \ | \ i>\theta(\alpha,\beta) \ or
\ i<\delta(\alpha,\beta)\bigr\}\cap \bigl\{\gamma_{i,\beta} \ | \
i>\theta(\alpha,\beta) \ or \
i<\delta(\alpha,\beta)\bigr\}=\emptyset$.

(4) \ $\bigl\{f_{i,\alpha} \ | \ i>\theta(\alpha,\beta) \ or \
i<\delta(\alpha,\beta)\bigr\}\cap \bigl\{f_{i,\beta} \ | \
i>\theta(\alpha,\beta) \ or \
i<\delta(\alpha,\beta)\bigr\}=\emptyset$.

(5) \ If $\delta(\alpha,\beta)+1\leq i\leq
\theta(\alpha,\beta)-1$, then $L(d_{i,\alpha})=L(d_{i,\beta})$.

{\bf \em Proof.} \  Since $j_{\alpha},j_{\beta}\in
L(c_{l+1}^{l})$, $d_{0,\alpha}, d_{0,\beta}\neq\emptyset$. Suppose
that $\delta(\alpha,\beta)$ and $\theta(\alpha,\beta)$ are two
integers such that $f_{i,\alpha}=f_{i,\beta}$ for
$\delta(\alpha,\beta)\leq i\leq \theta(\alpha,\beta)$. Since
$\cup_{f} E^{l}_{f}\cup e_{\gamma}^{l}$ is an abstract tree,
$\gamma_{i,\alpha}=\gamma_{i,\beta}$ for $\delta(\alpha,\beta)\leq
i\leq \theta(\alpha,\beta)$. Now $\delta(\alpha),\delta(\beta)\leq
\delta(\alpha,\beta)\leq\theta(\alpha,\beta)\leq
\theta(\alpha),\theta(\beta)$.

Assume now that $\gamma_{\delta(\alpha,\beta)-1,\alpha}\neq
\gamma_{\delta(\alpha,\beta)-1,\beta}$ and
$\gamma_{\theta(\alpha,\beta)+1,\alpha}\neq
\gamma_{\theta(\alpha,\beta)+1,\beta}$, and
$\gamma_{i,\alpha}=\gamma_{r,\beta}$ for some $i<
\delta(\alpha,\beta)-1$ or $i>\theta(\alpha,\beta)+1$ and some $r<
\delta(\alpha,\beta)-1$ or $r>\theta(\alpha,\beta)+1$.  Since
$d_{0,\alpha},d_{0,\beta}\neq\emptyset$, $\cup_{\gamma\in
m(l)}e_{\gamma}\cup_{f=0}^{l} E_{f}^{l}$ is not an abstract tree
even if $\delta(\alpha,\beta)=\theta(\alpha,\beta)=0$,
contradicting Proposition 4(1). Hence (1), (2), (3) and (4) hold.

Suppose that $\delta(\alpha,\beta)+1\leq i\leq -1$. Now
$f_{i,\alpha}=f_{i,\beta}$ and
$\gamma_{i-1,\alpha}=\gamma_{i-1,\beta}$,
$\gamma_{i,\alpha}=\gamma_{i,\beta}$. Hence $\partial
d_{i,\alpha}\cup \partial d_{i,\beta}\subset (\partial
e^{l}_{\gamma_{i-1,\alpha}}\cup \partial
e^{l}_{\gamma_{i,\alpha}})\times I$. By Proposition 4(4) and (5),
$\partial c_{r}^{l}$ is disjoint from $(
e^{l}_{\gamma_{i-1,\alpha}}\cup e^{l}_{\gamma_{i,\alpha}})\times
I$. Hence $c_{r}^{l}$ intersects $d_{i,\alpha}$ in one point if
and only if $c_{r}^{l}$ intersects $d_{i,\beta}$ in one point as
in Figure 7. Hence (5) holds. \qquad Q.E.D.

\begin{center}
\includegraphics[totalheight=7cm]{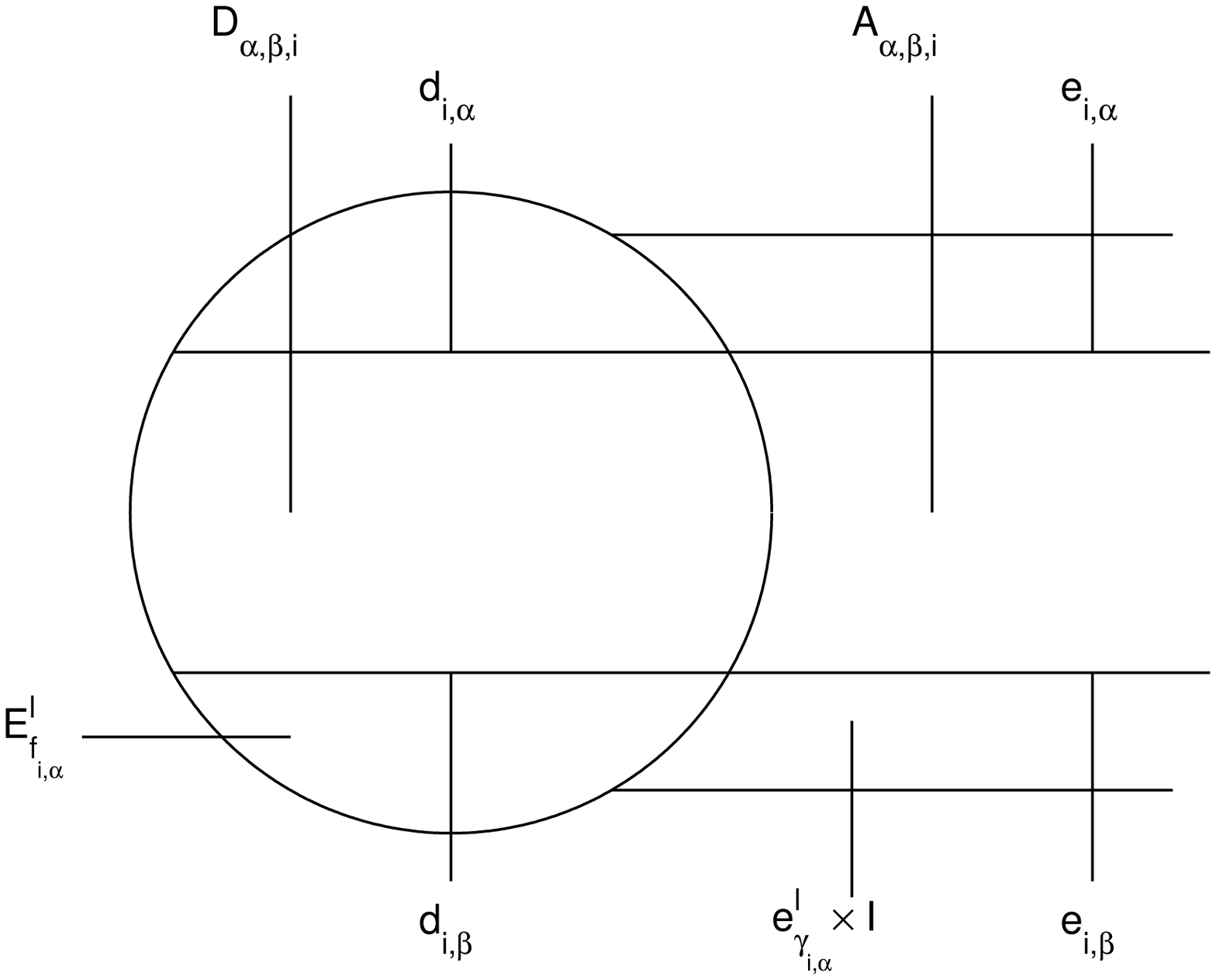}
\begin{center}
Figure 8(a)
\end{center}
\end{center}

\subsection{ $\delta(\alpha,\beta)$ and $\theta(\alpha,\beta)$}\vskip 0.5mm

\ \ \ \ \ Recalling the equality
$d_{j_{\alpha}}^{l}=\cup_{i=\delta(\alpha)}^{\theta(\alpha)}d_{i,\alpha}
\cup_{i=\delta(\alpha)}^{\theta(\alpha)} e_{i,\alpha}$ in Remark
3.3.2, and the two integers $\delta(\alpha,\beta)$ and
$\theta(\alpha,\beta)$ in Lemma 3.3.4. In this case,
$d_{i,\alpha}$ is an arc in $E^{l}_{f_{i,\alpha}}$ and
$e_{i,\alpha}$ is a core of $e^{l}_{\gamma_{i,\alpha}}\times
(0,1)$ for $i\neq 0$. In particular, $d_{0,\alpha}$ is an arc in
$E^{l}_{0}$. By Lemma 3.3.4, $f_{i,\alpha}=f_{i,\beta}$ and
$\gamma_{i,\alpha}=\gamma_{i,\beta}$ for $\delta(\alpha,\beta)\leq
i\neq 0\leq \theta(\alpha,\beta)$.

{\bf \em Definition 3.4.1.} \ For each $\delta(\alpha,\beta)\leq
i\neq 0 \leq \theta(\alpha,\beta)$.

(1) \ Let $D_{\alpha,\beta,i}$ be the disk  in
$E^{l}_{f_{i,\alpha}}$ which is bounded by $d_{i,\alpha}$ and
$d_{i,\beta}$ with two arcs in $\partial E^{l}_{f_{i,\alpha}}$.

(2) \ Let $A_{\alpha,\beta,i}$ be the disk in
$e^{l}_{\gamma_{i,\alpha}}\times I$ which is  bounded by
$e_{i,\alpha}$ and $e_{i,\beta}$  with two arcs in $(\partial
e^{l}_{\gamma_{i,\alpha}})\times I$ as in Figure 8(a).

(3) \ Let $D_{\alpha,\beta,0}$ be the disk in $E^{l}_{0}$ which is
 bounded by $d_{0,\alpha}$ and $d_{0,\beta}$ with two arcs in $\partial E^{l}_{0}$ as in Figure 8(a).

(4) \ Let
$D_{\alpha,\beta}=D_{\alpha,\beta,0}\cup_{i=\delta(\alpha,\beta)+1}^{\theta(\alpha,\beta)-1}D_{\alpha,\beta,i}\cup_{i=\delta(\alpha,\beta)}^{\theta(\alpha,\beta)}
A_{\alpha,\beta,i}$. Now we rewrite $D_{\alpha,\beta}$ as
$\cup_{i=\delta(\alpha,\beta)+1}^{\theta(\alpha,\beta)-1}D_{\alpha,\beta,i}\cup_{i=\delta(\alpha,\beta)}^{\theta(\alpha,\beta)}
A_{\alpha,\beta,i}$.

Now suppose that $\alpha<\lambda<\beta$.  In this section, we
shall prove that $\delta(\alpha,\lambda)\leq
\delta(\alpha,\beta)\leq \theta(\alpha,\beta)\leq
\theta(\alpha,\lambda)$, and $D_{\alpha,\beta}$ is a disk in
$F^{l}$. \vskip 0.5mm

%{\bf \em Lemma 3.4.1.} \ If $j_{\lambda}\in L(c_{k+1}^{l})$ and
%$\lambda\neq \alpha, \alpha+1$, then $d_{0,\lambda}$ is disjoint
%from $D_{\alpha,\alpha+1,0}$.

%{\bf \em Proof.} \ The lemma is immediately from Lemma 3.2.2 and
%the order of the elements in $L(c_{k+1}^{l})$ we rearranged in
%Section 3.3.\qquad Q.E.D. \vskip 0.5mm

{\bf \em Lemma 3.4.2.} \ Let $j\in\bigl\{1,\ldots,n\bigr\}-m(l)$.
If $d_{j}^{l}\cap D_{\alpha,\alpha+1,i}\neq\emptyset$ for some
$\delta(\alpha,\alpha+1)<i<\theta(\alpha,\alpha+1)$, then
$d_{j}^{l}\cap D_{\alpha,\alpha+1,i}\subset (d_{i,\alpha}\times
I)_{\gamma_{\alpha}}\cup (d_{i,\alpha+1}\times
I)_{\gamma_{\alpha+1}}$. Furthermore, if $j\neq j_{\alpha},
j_{\alpha+1}$, then $j>{\gamma_{\alpha}}$ or
${\gamma_{\alpha+1}}$.

{\bf \em Proof.} \   Assume that $d_{j}^{l}\cap
D_{\alpha,\alpha+1,i}\neq \emptyset$ for some
$\delta(\alpha,\alpha+1)<i\leq 0$. By Lemma 3.1.6, each component
of $d_{j}^{l}\cap D_{\alpha,\alpha+1,i}$ is either
$d_{j,f_{i,\alpha}}^{l}$ or a copy of $d_{r,f_{i,\alpha}}^{l}$ for
some $r\in I(w_{\gamma},l)$ and $\gamma\in m(l)$ with $\gamma<j$.
There are two cases:

Case 1. one component $c$ of $d_{j}^{l}\cap D_{\alpha,\alpha+1,i}$
is  $d_{j,f_{i,\alpha}}^{l}$.

Now if $j=j_{\alpha}$ or $j_{\alpha+1}$, then, by Lemma 3.3.2,
either $d_{j,f_{i,\alpha}}^{l}=d_{i,\alpha}$ or
$d_{j,f_{i,\alpha}}^{l}=d_{i,\alpha+1}$.

Suppose that $j\neq j_{\alpha},j_{\alpha+1}$. By Lemma 3.2.2,
$d_{j,f_{i,\alpha}}^{l}=c$ is disjoint from $(d_{i,\alpha}\times
I)_{\gamma_{\alpha}}\cup (d_{i,\alpha+1}\times
I)_{\gamma_{\alpha+1}}$. Since $\delta(\alpha,\alpha+1)<i\leq 0$,
$\gamma_{i-1,\alpha}=\gamma_{i-1,\alpha+1}$ and
$\gamma_{i,\alpha}=\gamma_{i, \alpha+1}$.  Now $\partial
d_{j,f_{i,\alpha}}^{l}\subset (\partial
e^{l}_{\gamma_{i-1,\alpha}}\cup
\partial e^{l}_{\gamma_{i,\alpha}})\times I$.

Now we claim that $\partial_{1} d_{j,f_{i,\alpha}}^{l}\subset
\partial e^{l}_{\gamma_{i-1,\alpha}}\times I$ and $\partial_{2}
d_{j,f_{i,\alpha}}^{l}\subset (\partial
e^{l}_{\gamma_{i,\alpha}})\times I$.

Suppose that $\partial d_{j,f_{i,\alpha}}^{l}\subset (\partial
e^{l}_{\gamma_{i-1,\alpha}})\times I$, then $d_{j}^{l}\cap
(\partial e^{l}_{\gamma_{i-1,\alpha}})\times I\neq \emptyset$. By
Proposition 4(1), $j\notin I(w_{\gamma_{i-1,\alpha}},l)$. By
Proposition 4(3), $j>\gamma_{i-1,\alpha}$. By Proposition 4(2),
$d_{j}^{l}$ is regular in $\cup_{f} E^{l}_{f}\cup_{\gamma<j}
e^{l}_{\gamma}\times I$.  This means that $d_{j}^{l}-\cup_{f}
d_{j,f}^{l}$ intersects $e^{l}_{\gamma_{i-1,\alpha}}\times I$ in
at least two cores of $e^{l}_{\gamma_{i-1,\alpha}}\times I$,
contradicting Lemma 3.1.2(3).

Now $j>\gamma_{i-1,\alpha}, \gamma_{i,\alpha}$, $\partial_{1}
d_{j,f_{i,\alpha}}^{l}\subset
\partial e^{l}_{\gamma_{i-1,\alpha}}\times I$ and $\partial_{2}
d_{j,f_{i,\alpha}}^{l}\subset (\partial
e^{l}_{\gamma_{i,\alpha}})\times I$. Hence $d_{j}^{l}$ intersects
$A_{\alpha,\alpha+1,i}$ in at least a core of
$A_{\alpha,\alpha+1,i}$. By Lemma 3.1.2,
$d_{j,f_{i+1,\alpha}}^{l}\neq\emptyset\subset
D_{\alpha,\alpha+1,i+1}$ as in Figure 8(b). By induction,
$d_{j,0}^{l}\neq \emptyset\subset D_{\alpha,\alpha+1,0}$, and
$j\in L(c_{l+1}^{l})$. Contradicting the order of $j_{\alpha},
j_{\alpha+1}$ we rearrange.

Case 2.  one component $c$ of  $d_{j}^{l}\cap
D_{\alpha,\alpha+1,i}$ is a copy of $d_{r,f_{i,\alpha}}^{l}$.

By the argument in Case 1, $r=j_{\alpha}$ or $j_{\alpha+1}$ and
$c\subset (d_{i,\alpha}\times I)_{\gamma_{\alpha}}\cup
(d_{i,\alpha+1}\times I)_{\gamma_{\alpha+1}}$. Assume now  that
$c\subset(d_{i,\alpha}\times I)_{\gamma_{\alpha}}$. Since
$c\subset d_{j}^{l}$ and $j\neq j_{\alpha}$. Then, by Proposition
4(3), $j>\gamma_{\alpha}$.\qquad Q.E.D.
\begin{center}
\includegraphics[totalheight=4cm]{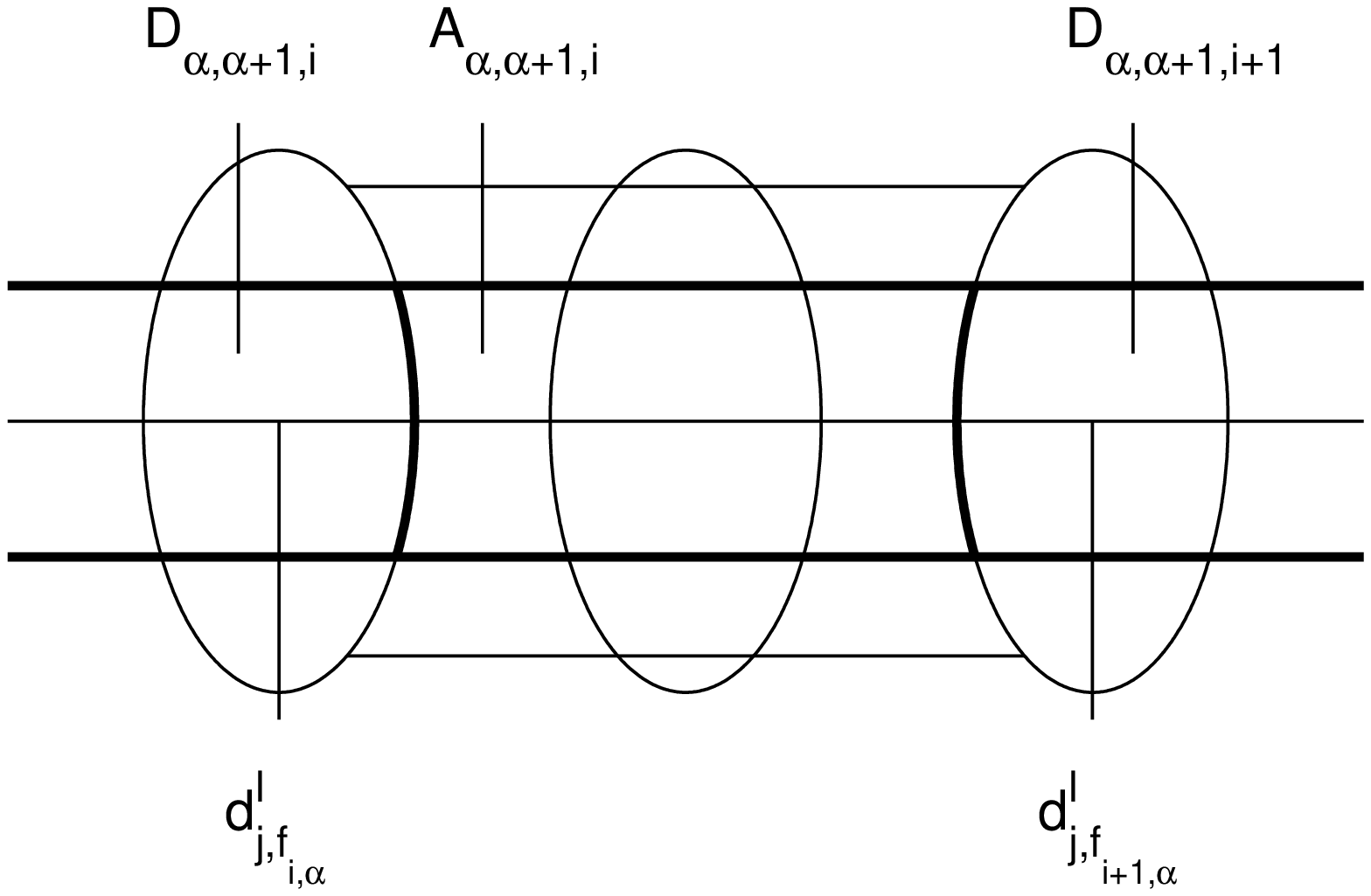}
\begin{center}
Figure 8(b)
\end{center}
\end{center}

{\bf \em Lemma 3.4.3.} \  If $\textrm{int}e^{l}_{\gamma}\times
I\cap D_{\alpha,\alpha+1,i}\neq \emptyset$ for some
$\delta(\alpha,\alpha+1)<i<\theta(\alpha,\alpha+1)$,  then each
component of $\textrm{int}e^{l}_{\gamma}\times I\cap
D_{\alpha,\alpha+1,i}$ lies in one of $(d_{i,\alpha}\times
I)_{\gamma_{\alpha}}$ and $(d_{i,\alpha+1}\times
I)_{\gamma_{\alpha+1}}$. Furthermore, if one component of
$\textrm{int}e^{l}_{\gamma}\times I\cap D_{\alpha,\alpha+1,i}$
lies in $(d_{i,\alpha}\times I)_{\gamma_{\alpha}}$, then
$\gamma\geq\gamma_{\alpha}$.

{\bf \em Proof.} \  Suppose   that
$\textrm{int}e^{l}_{\gamma}\times I\cap D_{\alpha,\alpha+1,i}\neq
\emptyset$. By Lemma 3.2.4, each component of
$\textrm{int}e^{l}_{\gamma}\times I\cap E^{l}_{f_{i,\alpha}}$ is
contained in $(d_{r,f_{i,\alpha}}^{l}\times I)_{\gamma_{r}}$ for
some $\gamma_{r}\in m(l)$ and $r\in I(w_{\gamma_{r}},l)$.
Furthermore, $\gamma\geq\gamma_{r}>r$. If $r\neq
j_{\alpha},j_{\alpha+1}$, then, by Lemma 3.2.2 and Definition
3.2.1, $(d_{i,\alpha}\times I)_{\gamma_{\alpha}}$ and
$(d_{i,\alpha+1}\times I)_{\gamma_{\alpha+1}}$ are disjoint from
$(d_{r,f_{i,\alpha}}^{l}\times I)_{\gamma_{r}}$. In this case, if
$(d_{r,f_{i,\alpha}}^{l}\times I)_{\gamma_{r}}\cap
D_{\alpha,\alpha+1,i}\neq\emptyset$, then
$d_{r,f_{i,\alpha}}^{l}\subset D_{\alpha,\alpha+1,i}$. By the
proof of Lemma 3.4.2, this is impossible.  \qquad Q.E.D.

{\bf \em Lemma 3.4.4.} \ If $\delta(\alpha,\alpha+1)\leq i\neq
r\neq 0\leq \theta(\alpha,\alpha+1)$, then $A_{\alpha,\alpha+1,i}$
is disjoint from $A_{\alpha,\alpha+1,r}$.

{\bf \em Proof.} \ By Lemma 3.3.1,
$\gamma_{i,\alpha}\neq\gamma_{r,\alpha}$. We may assume that
$\gamma_{i,\alpha}<\gamma_{r,\alpha}$. If
$e^{l}_{\gamma_{r,\alpha}}\times I$ is disjoint from
$e^{l}_{\gamma_{i,\alpha}}\times I$, then
$A_{\alpha,\alpha+1,r}\subset e^{l}_{\gamma_{r,\alpha}}\times I$
is disjoint from $A_{\alpha,\alpha+1,i}\subset
e^{l}_{\gamma_{i,\alpha}}\times I$.

Suppose now that $e^{l}_{\gamma_{r,\alpha}}\times I\cap
e^{l}_{\gamma_{i,\alpha}}\times I\neq \emptyset$. Now by
Proposition 4(1) and Definition 2.1.4, each component of
$A_{\alpha,\alpha+1,r}\cap e^{l}_{\gamma_{i,\alpha}}\times I$ is a
disk in $e^{l}_{\gamma_{i,\alpha}}\times I$. If one component of
$A_{\alpha,\alpha+1,r}\cap e_{\gamma_{i,\alpha}}^{l}\times I$, say
$A$, is not disjoint from $A_{\alpha,\alpha+1,i}$, then either one
of the two arcs $A\cap e_{r,\alpha}$ and $A\cap e_{r,\alpha+1}$
lies in $A_{\alpha,\alpha+1,i}$ as in Figure 9(a) or one of
$e_{i,\alpha}$ and $e_{i,\alpha+1}$ lies in $A$ as in Figure 9(b).
Here $e_{r,\alpha}$ is defined in Definition 3.4.1.

We first suppose that   one of $A\cap e_{r,\alpha}$ and $A\cap
e_{r,\alpha+1}$, say $A\cap e_{r,\alpha}$, lies in
$A_{\alpha,\alpha+1,i}$. By Definition 2.1.4, $(\partial
e^{l}_{\gamma_{r,\alpha}})\times I\cap (\partial
e^{l}_{\gamma_{i,\alpha}})\times I=\emptyset$. Hence
$e_{r,\alpha}\cap intD_{\alpha,\alpha+1,i+1}\neq \emptyset$. By
Lemma 3.2.4, there is $j$ such that  a copy of
$d_{j,f_{i+1,\alpha}}^{l}$ is contained in
$D_{\alpha,\alpha+1,i+1}$, and
$j<\gamma_{i,\alpha}<j_{\alpha},j_{\alpha+1}$, contradicting Lemma
3.4.2.

Assume that one of $e_{i,\alpha}$ and $e_{i,\alpha+1}$, say
$e_{i,\alpha}$, lies in $A$ as in Figure 9(b). Then
$d_{i+1,\alpha}\subset A_{\alpha,\alpha+1,r}\subset
e^{l}_{\gamma_{r,\alpha}}\times I$. Note that
$j_{\alpha}>\gamma_{r,\alpha}$. By  Lemma 3.1.2(2), this is
impossible. \qquad Q.E.D.\vskip 0.5mm
\begin{center}
\includegraphics[totalheight=4cm]{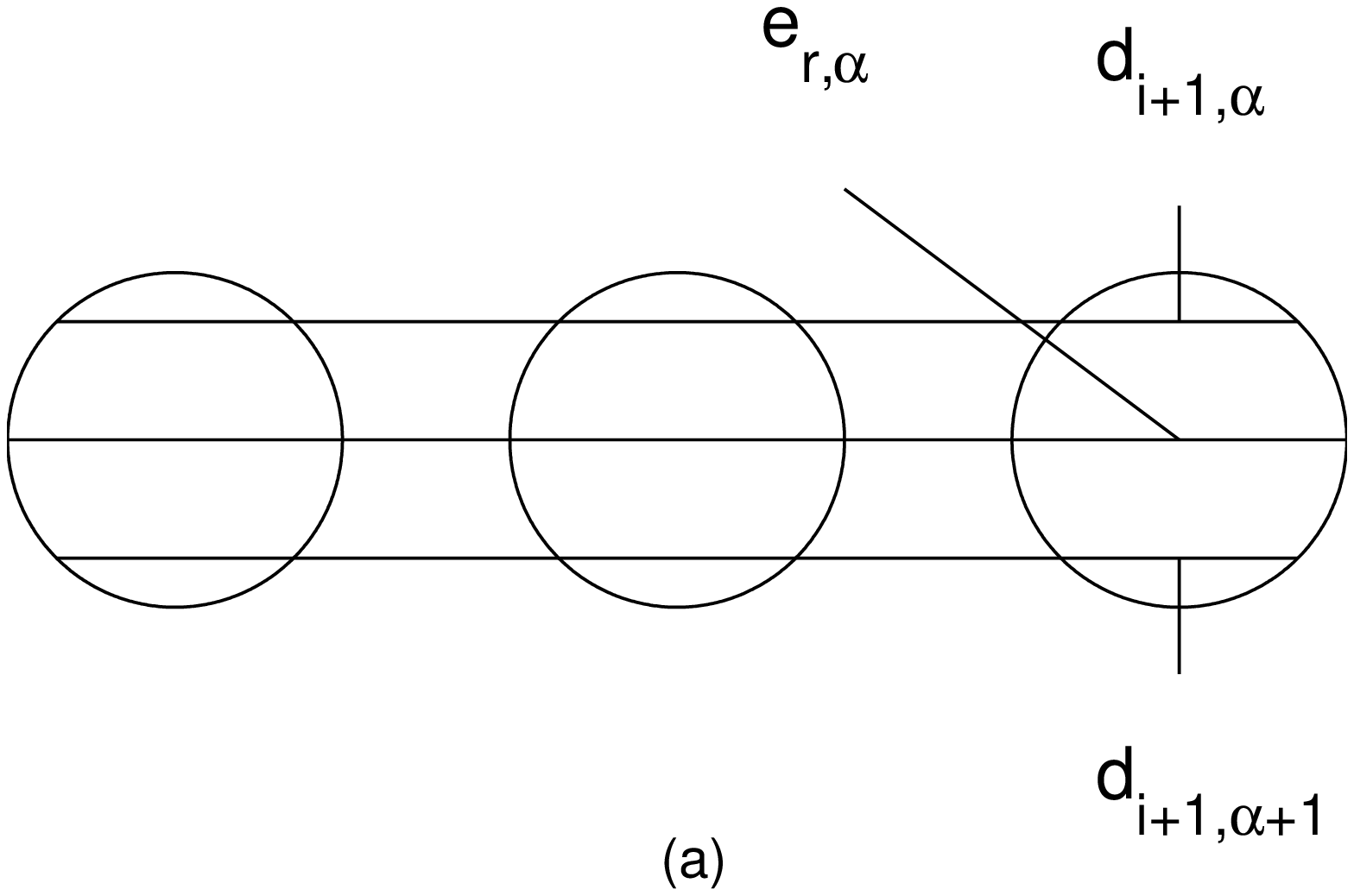}
\end{center}
\begin{center}
\includegraphics[totalheight=4cm]{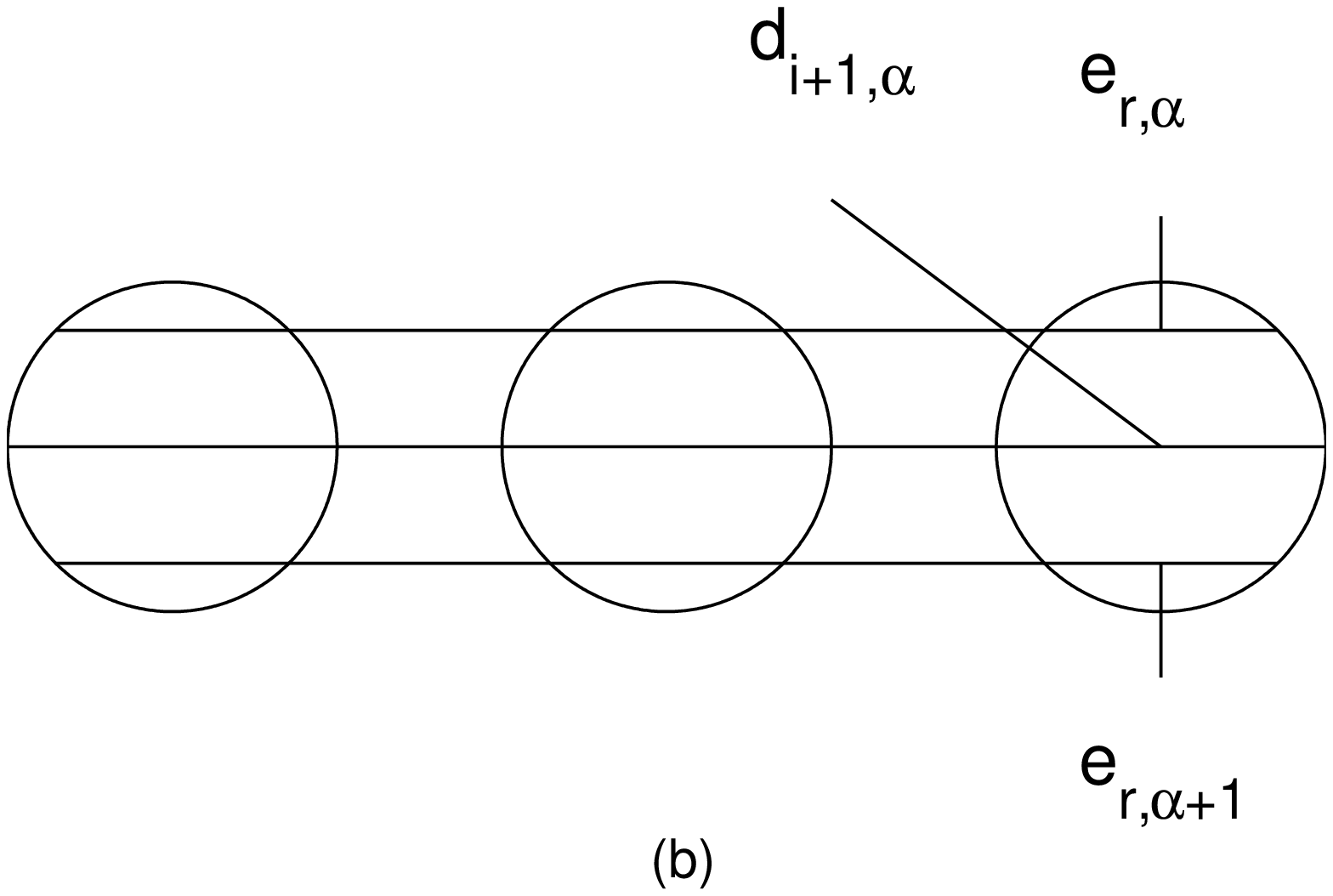}
\begin{center}
Figure 9
\end{center}
\end{center}

{\bf \em Lemma 3.4.5.} For each $\delta(\alpha,\alpha+1)\leq
i\leq\theta(\alpha,\alpha+1)$ and $\delta(\alpha,\alpha+1)+1\leq
r\neq 0\leq\theta(\alpha,\alpha+1)-1$,
$inte^{l}_{\gamma_{i,\alpha}}\times I\cap A_{\alpha,\alpha+1,i}$
is disjoint from $D_{\alpha,\alpha+1,r}$. Furthermore, if $i\neq
r-1,r$, then $A_{\alpha,\alpha+1,i}$ is disjoint from
$D_{\alpha,\alpha+1,r}$.

{\bf \em Proof.} \ By Lemma 3.2.4, each component of
$\textrm{int}e^{l}_{\gamma_{i,\alpha}}\times I\cap
E^{l}_{f_{r,\alpha}}$ is contained in $(d_{j,f_{r,\alpha}}^{l}
\times I)_{\gamma}$, where $\gamma\leq \gamma_{i,\alpha}\in m(l)$
and $j\in I(w_{\gamma},l)$. Hence $j<\gamma_{i,\alpha}$. Since
$\gamma_{i,\alpha}<j_{\alpha}, j_{\alpha+1}$, by Lemma 3.3.1,
$d_{r,\alpha}$ and $d_{r,\alpha+1}$ are disjoint from
$\textrm{int}e^{l}_{\gamma_{i,\alpha}}\times I$. Thus if
$\textrm{int}e^{l}_{\gamma_{i,\alpha}}\times I\cap
D_{\alpha,\alpha+1,r}\neq \emptyset$, then one component of
$\textrm{int}e^{l}_{\gamma_{i,\alpha}}\times I\cap
E^{l}_{f_{r,\alpha}}$ lies in $D_{\alpha,\alpha+1,r}$. Hence
$(d_{j,f_{r,\alpha}}^{l} \times I)_{\gamma}$ lies in
$D_{\alpha,\alpha+1,r}$ for some $j<j_{\alpha}, j_{\alpha+1}$. By
Lemma 3.4.2, this is impossible.

Note that $\partial_{1} d_{r,\alpha}, \partial_{1}
d_{r,\alpha+1}\subset (\partial e^{l}_{\gamma_{r-1,\alpha}})\times
I$ and $\partial_{2} d_{r,\alpha}, \partial_{2}
d_{r,\alpha+1}\subset (\partial e^{l}_{\gamma_{r,\alpha}})\times
I$. If $i\neq r,r-1$, then, by Lemma 2.1.4, $(\partial
e^{l}_{\gamma_{i,\alpha}})\times I\cap ((\partial
e^{l}_{\gamma_{r-1,\alpha}})\times I\cup(\partial
e^{l}_{\gamma_{r,\alpha}})\times I)=\emptyset$. \qquad
Q.E.D.\vskip 0.5mm

{\bf \em Lemma 3.4.6.} \ (1) \ If $d_{j}^{l}\cap
A_{\alpha,\alpha+1,i}\neq \emptyset$, then $d_{j}^{l}\cap
A_{\alpha,\alpha+1,i}\subset (e_{i,\alpha}\times
I)_{\gamma_{\alpha}}\cup (e_{i,\alpha+1}\times
I)_{\gamma_{\alpha+1}}$. Furthermore, if $j\neq
j_{\alpha},j_{\alpha+1}$, then $j>\gamma_{\alpha}$ or
$\gamma_{\alpha+1}$.

(2) \ If $\gamma> \gamma_{i,\alpha}$ for
$\delta(\alpha,\alpha+1)\leq i\leq \theta(\alpha,\alpha+1)$ and
$e^{l}_{\gamma}\times I\cap A_{\alpha,\alpha+1,i}\neq \emptyset$,
then $e^{l}_{\gamma}\times I\cap A_{\alpha,\alpha+1,i}\subset
(e_{i,\alpha}\times I)_{\gamma_{\alpha}}\cup (e_{i,\alpha+1}\times
I)_{\gamma_{\alpha+1}}$. Furthermore, $\gamma\geq \gamma_{\alpha}$
or $\gamma_{\alpha+1}$.

{\bf \em Proof}. \ (1) \ By Lemma 3.3.1 and Definition 3.4.1,
$e_{i,\alpha},e_{i,\alpha+1}, A_{\alpha,\alpha+1,i}\subset
e^{l}_{\gamma_{i,\alpha}}\times (0,1)$. If $d_{j}^{l}\cap
A_{\alpha,\alpha+1,i}\neq \emptyset$, then $j\notin
I(w_{\gamma_{i,\alpha}},l)$; otherwise, $d_{j}^{l}\subset
e^{l}_{\gamma_{i,\alpha}}\times\bigl\{0\bigr\}$. By Proposition
4(3), $j>\gamma_{i,\alpha}$. By Proposition 4(2) and Definition
2.1.5, each component of $d_{j}^{l}\cap
e^{l}_{\gamma_{i,\alpha}}\times I$ is a core of
$e^{l}_{\gamma_{i,\alpha}}\times I$ which lies in $intd_{j}^{l}$.
Thus $d_{j}^{l}\cap D_{\alpha,\alpha+1,i+1}\neq \emptyset$. By
Lemma 3.4.2, (1) holds.

(2) \ Now suppose that $\gamma>\gamma_{i,\alpha}$ and
$e^{l}_{\gamma}\times I\cap A_{\alpha,\alpha+1,i}\neq \emptyset$.
Then each component of $e^{l}_{\gamma}\times I\cap
e^{l}_{\gamma_{i,\alpha}}\times I$ is $(c\times I)_{\gamma}\subset
e^{l}_{\gamma_{i,\alpha}}\times I$ where $c$ is a core of
$e^{l}_{\gamma_{i,\alpha}}\times I$. By Definition 2.1.4,
$(\partial e^{l}_{\gamma})\times I\cap(\partial
e^{l}_{\gamma_{i,\alpha}})\times I=\emptyset$. Hence
$\textrm{int}e^{l}_{\gamma}\times I\cap
D_{\alpha,\alpha+1,i+1}\neq \emptyset$. By Lemma 3.4.3,
$\textrm{int}e^{l}_{\gamma}\times I\cap
D_{\alpha,\alpha+1,i+1}\subset(d_{i,\alpha}\times
I)_{\gamma_{\alpha}}\cup(d_{i,\alpha+1}\times
I)_{\gamma_{\alpha+1}}$. Furthermore, If $e^{l}_{\gamma}\times
I\cap (d_{i,\alpha}\times I)_{\gamma_{\alpha}}\neq\emptyset$, then
$\gamma\geq \gamma_{\alpha}$. In this case,   each component of
$e^{l}_{\gamma}\times I\cap e^{l}_{\gamma_{\alpha}}\times I$ is
$(c\times I)_{\gamma}\subset e^{l}_{\gamma_{\alpha}}\times (0,1)$
when $\gamma>\gamma_{\alpha}\in m(k)$. Hence (2) holds. \qquad
Q.E.D.\vskip 0.5mm

{\bf \em Lemma 3.4.7.} \ (1) Suppose that $\alpha<\lambda<\beta$.
Then $\delta(\alpha,\beta)\geq
\delta(\alpha,\lambda),\delta(\lambda,\beta)$ and
$\theta(\alpha,\beta)\leq
\theta(\alpha,\lambda),\theta(\lambda,\beta)$.

(2) \ $D_{\alpha,\beta}$ is a disk in $F^{l}$.

(3) \ For $\alpha\leq \lambda\leq\beta$, let
$h^{\lambda}_{\alpha,\beta}=\cup_{i=\delta_{0}+1}^{\theta_{0}-1}d_{i,\lambda}\cup_{i=\delta_{0}}^{\theta_{0}}e_{i,\lambda}$
where $\theta_{0}=\theta(\alpha,\beta)$ and
$\delta_{0}=\delta(\alpha,\beta)$.

(i) if $d_{j}^{l}\cap D(\alpha,\beta)\neq \emptyset$, then
$d_{j}^{l}\cap D(\alpha,\beta)\subset
\cup_{\lambda}(h^{\lambda}_{\alpha,\beta}\times
I)_{\gamma_{\lambda}}$. Furthermore, if $j\neq j_{\lambda}$ for
each $\alpha\leq\lambda\leq\beta$, then $j>\gamma_{\lambda}$ for
some $\alpha\leq\lambda\leq\beta$,

(ii) \ If $\gamma> \gamma_{i,\alpha}$ for
$\delta(\alpha,\beta)\leq i\leq \theta(\alpha,\beta)$ and
$e^{l}_{\gamma}\times I\cap D(\alpha,\beta)\neq \emptyset$, then
$e^{l}_{\gamma}\times I\cap
D(\alpha,\beta)\subset\cup_{\lambda}(h^{\lambda}_{\alpha,\beta}\times
I )_{\gamma_{\lambda}}$. Furthermore, $\gamma\geq
\gamma_{\lambda}$ for some $\alpha\leq\lambda\leq\beta$.

{\bf \em Proof.} \  By Lemma 3.3.1, if $i\neq r$, then
$f_{i,\lambda}\neq f_{r,\lambda}$. Hence
$D_{\lambda,\lambda+1,i}\cap D_{\lambda,\lambda+1,r}=\emptyset$ if
$i\neq r$. By Lemmas 3.4.4 and 3.4.5, $D_{\lambda,\lambda+1}$ is a
disk in $F^{l}$. See in Definition 3.4.1.

We first prove that $\delta(\alpha,\alpha+2)\geq
\delta(\alpha,\alpha+1),\delta(\alpha+1,\alpha+2)$ and
$\theta(\alpha,\alpha+2)\leq
\theta(\alpha,\alpha+1),\theta(\alpha+1,\alpha+2)$.

Suppose, otherwise, that $\delta(\alpha,\alpha+2)<
\delta=\delta(\alpha,\alpha+1)$. Then
$f_{i,\alpha}=f_{i,\alpha+1}=f_{i,\alpha+2}$ and
$\gamma_{i,\alpha}=\gamma_{i,\alpha+1}=\gamma_{i,\alpha+2}$ for
$\delta\leq i\leq -1$. Furthermore,
$\gamma_{i,\alpha}<j_{\alpha},j_{\alpha+1}, j_{\alpha+2}$ for
$\delta\leq i\leq -1$.  Now let $D_{1}=\cup_{i=\delta+1}^{0}
D_{\alpha,\alpha+1,i}\cup_{i=\delta}^{-1}A_{\alpha,\alpha+1,i}$,
$D_{2}=\cup_{i=\delta+1}^{0}
D_{\alpha+1,\alpha+2,i}\cup_{i=\delta}^{-1}A_{\alpha+1,\alpha+2,i}$.
Then $D_{1}\subset D_{\alpha,\alpha+1}$ and $D_{2}\subset
D_{\alpha+1,\alpha+2}$ are two disks.

Now by Lemma 3.2.2 and Lemma 3.4.2, $d_{i,\alpha+2}$ is disjoint
from $D_{\alpha,\alpha+1,i}$, otherwise, $d_{0,\alpha+2}\subset
D_{\alpha,\alpha+1,0}$. Hence $D_{\alpha,\alpha+1,i}\cap
D_{\alpha+1,\alpha+2,i}=d_{i,\alpha+1}$. Since
$\gamma_{i,\alpha}<j_{\alpha},j_{\alpha+1}, j_{\alpha+2}$, by
Lemma 3.4.3, $A_{\alpha+1,\alpha+2,i}$ is disjoint from
$intD_{\alpha,\alpha+1,r}$. Similarly, $A_{\alpha,\alpha+1,i}$ is
disjoint from $intD_{\alpha+1,\alpha+2,r}$. By the proof of Lemma
3.4.4, $A_{\alpha,\alpha+1,i}$ is disjoint from
$A_{\alpha+1,\alpha+2,r}$ for $i\neq r$. Note that
$A_{\alpha,\alpha+1,i}\cap
A_{\alpha+1,\alpha+2,i}=e_{i,\alpha+1}$; Otherwise,
$e_{i,\alpha+2}$ separates $e_{i,\alpha},e_{i,\alpha+1}$ in
$e^{l}_{\gamma_{i,\alpha}}\times I$, and $d_{i+1,\alpha+2}$
separates $d_{i+1, \alpha}$ and $d_{i+1,\alpha+1}$ in
$E^{l}_{f_{i+1,\alpha}}$, but by the proof of Lemma 3.4.2, this is
impossible. Now $D_{1}\cup D_{2}$ is a disk in $F^{l}$. Hence
$d_{\delta,\alpha+1}$ separates $d_{\delta,\alpha}$ and
$d_{\delta,\alpha+2}$ in $E^{l}_{f_{\delta,\alpha}}$. Since
$F^{l}$ is generated by $\cup_{f}E^{l}_{f}\cup_{\gamma}
e^{l}_{\gamma}$. So if $(\partial_{1} e^{l}_{\gamma})\times
I\subset E^{l}_{f}$ for some $f$, then $(\partial_{2}
e^{l}_{\gamma})\times I$ is disjoint from $E^{l}_{f}$. Now if
$\gamma_{\delta-1,\alpha}=\gamma_{\delta-1,\alpha+2}$, then
$\gamma_{\delta-1,\alpha}=\gamma_{\delta-1,\alpha+1}$. This means
that $\delta(\alpha,\alpha+2)\geq \delta=\delta(\alpha,\alpha+1)$,
$\delta(\alpha+1,\alpha+2)$, a contradiction.

Similarly, $\theta_{0}=\theta(\alpha,\alpha+2)\leq
\theta(\alpha,\alpha+1),\theta(\alpha+1,\alpha+2)$.

Let $\delta_{0}=\delta(\alpha,\alpha+2)$ and
$\theta_{0}=\theta(\alpha,\alpha+2)$. Let
$D_{*}=\cup_{i=\delta_{0}+1}^{\theta_{0}-1}
D_{\alpha,\alpha+1,i}\cup_{i=\delta_{0}}^{\theta_{0}}A_{\alpha,\alpha+1,i}$,
$D_{**}=\cup_{i=\delta_{0}+1}^{\theta_{0}-1}
D_{\alpha+1,\alpha+2,i}\cup_{i=\delta_{0}}^{\theta_{0}}A_{\alpha+1,\alpha+2,i}$.
Then $D_{*}\subset D_{\alpha,\alpha+1}$ and $D_{**}\subset
D_{\alpha+1,\alpha+2}$.

By the above argument, $D(\alpha,\alpha+2)=D_{*}\cup D_{**}$ is a
disk in $F^{l}$. By Lemma 3.4.2, Lemma 3.4.3 and Lemma 3.4.6, (3)
and (4) hold.

By induction, $\alpha<\lambda<\beta$. Then
$\delta(\alpha,\beta)\geq
\delta(\alpha,\lambda),\delta(\lambda,\beta)$ and
$\theta(\alpha,\beta)\leq
\theta(\alpha,\lambda),\theta(\lambda,\beta)$. Furthermore,
$D_{\alpha,\beta}=\cup_{\lambda=\alpha}^{\beta-1}(\cup_{i=\delta(\alpha,\beta)+1}^{\theta(\alpha,\beta)-1}
D_{\lambda,\lambda+1,i}\cup_{i=\delta(\alpha,\beta)}^{\theta(\alpha,\beta)}A_{\lambda,\lambda+1,i})$
is a disk as in Figure 10(a). (3) is immediately from Lemmas
3.4.2, 3.4.3 and 3.4.6. \qquad Q.E.D.
\begin{center}
\includegraphics[totalheight=3cm]{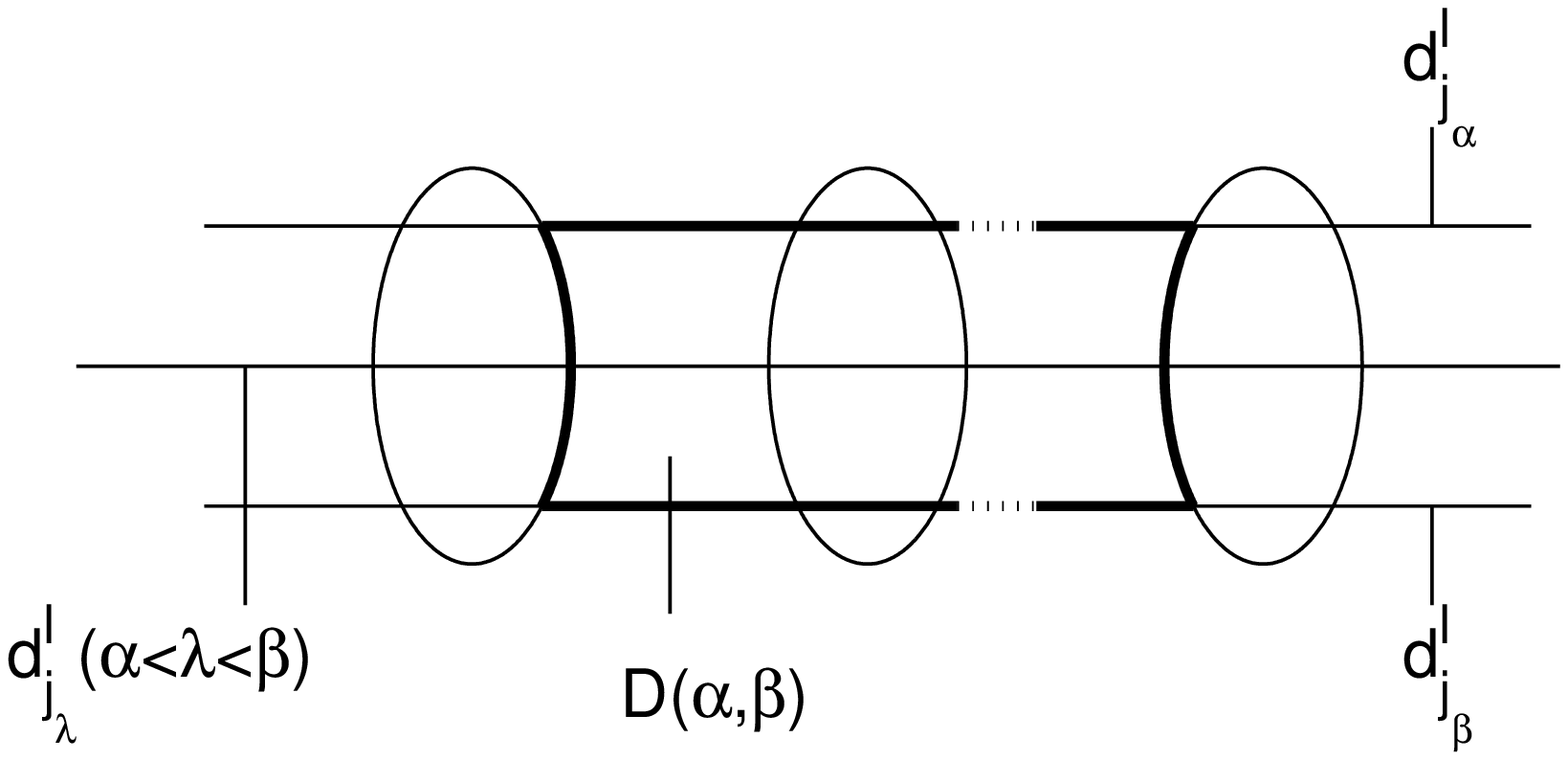}
\begin{center}
Figure 10(a)
\end{center}
\end{center}

\begin{center}
\includegraphics[totalheight=4.5cm]{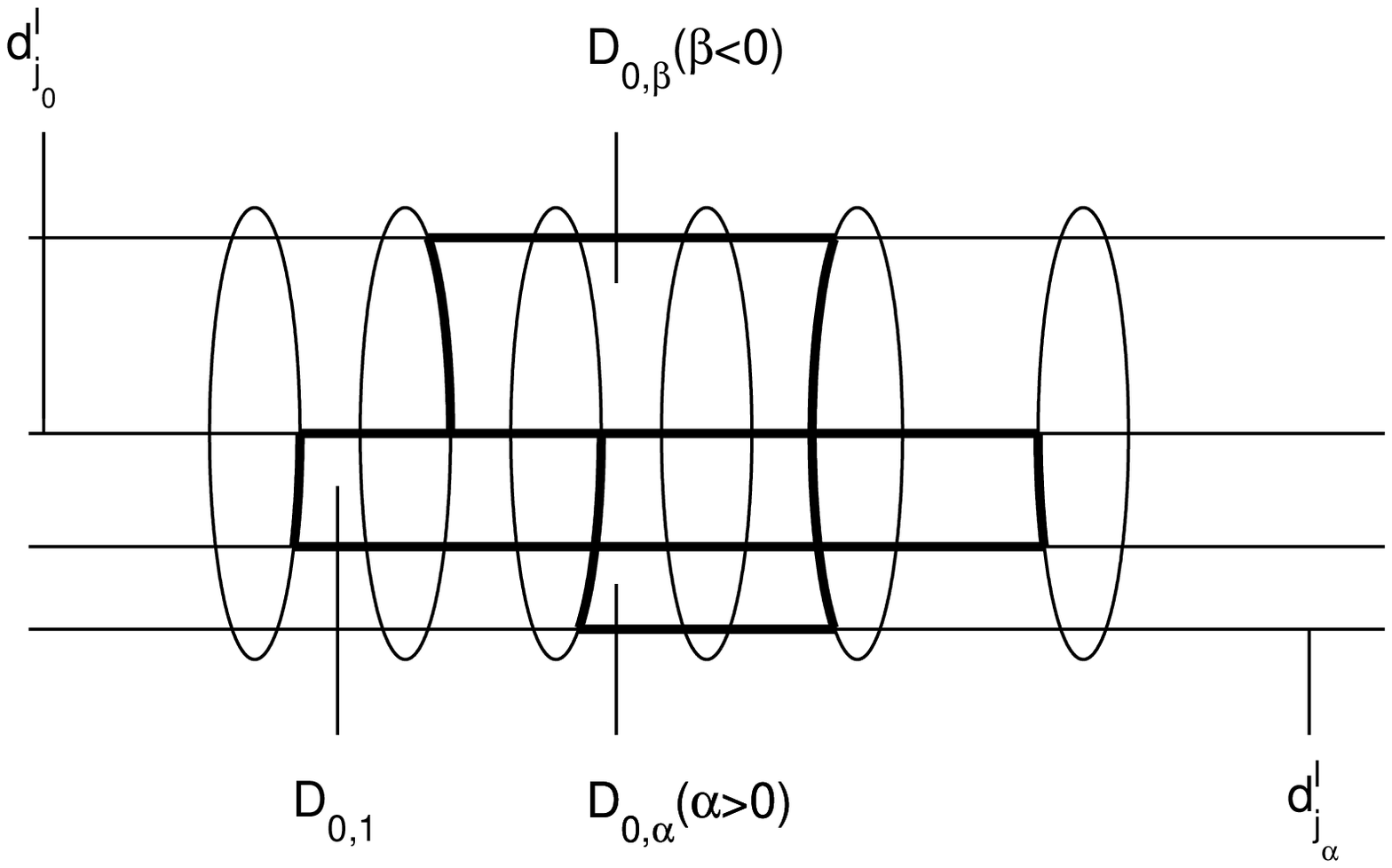}
\begin{center}
Figure 10(b)
\end{center}
\end{center}

{\bf \em Lemma 3.4.8.} \ $\cup_{\alpha\neq 0} D_{0,\alpha}$ is a
disk.

{\bf \em Proof.} \ By Lemma 3.4.7,
$\delta(0,\alpha)\leq\delta(0,\beta)\leq\theta(0,\beta)\leq\theta(0,\alpha)$
if $0<\alpha<\beta$ or $\beta<\alpha<0$.  Now suppose that
$L(c_{l+1}^{l})=\bigl\{j_{\beta_{1}},\ldots,j_{-1},j_{0},j_{1},\ldots,j_{\beta_{2}}\bigr\}$.
Then

\begin{center}$\cup_{\alpha\neq 0}
D_{0,\alpha}=\cup_{\alpha=\beta_{1}}^{\beta_{2}-1}(\cup_{i=\delta(0,\alpha)+1}^{\theta(0,\alpha)-1}D_{\alpha,\alpha+1,i}
\cup_{i=\delta(0,\alpha)}^{\theta(0,\alpha)}A_{\alpha,\alpha+1,i})$.\end{center}

Since $j_{0}<j_{\alpha}$, $\gamma_{i,0}<j_{0}<j_{\alpha}$ for
$\delta(0)\leq i\leq \theta(0)$ and $\alpha\neq 0$.

Now if $\beta\neq \alpha,\alpha+1$, then $intD_{\beta,\beta+1,i}$
is disjoint from $D_{\alpha,\alpha+1,i}$. Otherwise, one of
$d_{i,\beta}$ and $d_{i,\beta+1}$, say $d_{i,\beta}$, separates
$d_{i,\alpha}$ and $d_{i,\alpha+1}$ in $E^{l}_{f_{i,\alpha}}$, but
by the proof of Lemma 3.4.2, $d_{0,\beta}$ separates
$d_{0,\alpha}$ and $d_{0,\alpha+1}$ in $E^{l}_{0}$, a
contradiction. Similarly, if $\bigl\{\beta,\beta+1\bigr\}\cap
\bigl\{\alpha,\alpha+1\bigr\}=\emptyset$, then
$A_{\beta,\beta+1,i}$ is disjoint from $A_{\alpha,\alpha+1,i}$. If
$\bigl\{\beta,\beta+1\bigr\}\cap
\bigl\{\alpha,\alpha+1\bigr\}=\bigl\{\beta\bigr\}$, then
$A_{\beta,\beta+1,i}\cap A_{\alpha,\alpha+1,i}=e_{i,\beta}$.

Since $\gamma_{i,\alpha}=\gamma_{i,0}<j_{0}$ for each $\alpha$ and
$\delta(0,\alpha)\leq i\leq \theta(0,\alpha)$,  by the proof of
Lemma 3.4.4, $int A_{\alpha,\alpha+1,i}$ is disjoint from
$A_{\beta,\beta+1,r}$. By Lemma 3.4.3, $int A_{\alpha,\alpha+1,i}$
is disjoint from $D_{\beta,\beta+1,r}$. Now
 $\cup_{\alpha\neq 0} D_{0,\alpha}$ is a
disk as in Figure 10(b).\qquad Q.E.D.

\subsection{$a_{\alpha}$ and $a_{\alpha}^{0}$}

\ \ \ \ \  For each $j_{\alpha}\in L(c_{l+1}^{l})$, recalling the
equality:
$d_{j_{\alpha}}^{l}=d_{0,\alpha}\cup_{i=\delta(\alpha)}^{\theta(\alpha)}d_{i,\alpha}
\cup_{i=\delta(\alpha)}^{\theta(\alpha)} e_{i,\alpha}$.

{\bf \em Definition 3.5.1.} \ For $\alpha\neq 0$,

(1) \ let
$a_{\alpha}=\cup_{i=\delta(0,\alpha)}^{\theta(0,\alpha)}d_{i,\alpha}
\cup_{i=\delta(0,\alpha)}^{\theta(0,\alpha)}e_{i,\alpha}$,

(2) \ let
$a_{\alpha}^{0}=\cup_{i=\delta(0,\alpha)}^{\theta(0,\alpha)}d_{i,0}
\cup_{i=\delta(0,\alpha)}^{\theta(0,\alpha)}e_{i,0}$.\vskip 0.5mm

{\bf \em Lemma 3.5.2.} \ (1) \ $a_{\alpha}\subset
d_{j_{\alpha}}^{l}$ and $a_{\alpha}^{0}\subset d_{j_{0}}^{l}$.

(2) \ If $0<\alpha<\beta$, then $a_{\beta}^{0}\subset
a_{\alpha}^{0}$.

(3) \ If $\beta<\alpha<0$, then $a_{\beta}^{0}\subset
a_{\alpha}^{0}$.

{\bf \em Proof.} \ The lemma follows from Definition 3.5.1 and
Lemma 3.4.7.\qquad Q.E.D.\vskip 0.5mm

{\bf \em Lemma 3.5.3.} \ (1) \ Suppose that
$s(w_{\gamma_{\alpha}})=+$ or $j_{\alpha}\notin I(w_{\gamma},l)$
for each $\gamma\in m(l)$. If $e^{l}_{\gamma}\times
I\cap(inta_{\alpha}\times I)_{\gamma_{\alpha}}\neq \emptyset$,
then  $\gamma\leq Max\bigl\{\gamma_{i,0} \ | \ \delta(0)\leq i\leq
\theta(0)\bigr\}$.

(2) \ Suppose that $s(w_{\gamma_{\alpha}})=-$. If
$e^{l}_{\gamma}\times I\cap(inta_{\alpha}\times
I)_{\gamma_{\alpha}}\neq \emptyset$, then either $\gamma\leq
Max\bigl\{\gamma_{i,0} \ | \ \delta(0)\leq i\leq \theta(0)\bigr\}$
or $\gamma\geq\gamma_{\alpha}$.

{\bf \em Proof.} \ By Definition 3.5.1,
$a_{\alpha}=\cup_{i=\delta(0,\alpha)}^{\theta(0,\alpha)}d_{i,\alpha}
\cup_{i=\delta(0,\alpha)}^{\theta(0,\alpha)}e_{i,\alpha}$. By
Lemma 3.3.4, $\gamma_{i,\alpha}=\gamma_{i,0}$,
$f_{i,\alpha}=f_{i,0}$ for $\delta(0,\alpha)\leq i\leq
\theta(0,\alpha)$.

(1) Suppose that $s(w_{\gamma_{\alpha}})=+$ or $j_{\alpha}\notin
I(w_{\gamma},l)$ for each $\gamma\in m(l)$. In this case, by
Definition 3.2.1, $(inta_{\alpha}\times
I)_{\gamma_{\alpha}}=inta_{\alpha}$. If $\gamma>j_{\alpha}$, then,
by Proposition 4(3), $a_{\alpha}$ is disjoint from
$e^{l}_{\gamma}\times I$. So $j_{\alpha}>\gamma$. By Lemma 3.3.1,
$intd_{i,\alpha}$ is disjoint from $e^{l}_{\gamma}\times I$. Note
that $\gamma_{i,\alpha}=\gamma_{i,0}$ for $\delta(0,\alpha)\leq
i\leq \theta(0,\alpha)$. Now by the proof of Lemma 3.2.3,
$inta_{\alpha}$ is disjoint from $e^{l}_{\gamma}\times I$ for
$\gamma>Max\bigl\{\gamma_{i,0} \ | \ \delta(0)\leq i\leq
\theta(0)\bigr\}$.

(2) \ Suppose that $s(w_{\gamma_{\alpha}})=-$,
$\gamma_{\alpha}>\gamma> Max\bigl\{\gamma_{i,0} \ | \
\delta(0)\leq i\leq \theta(0)\bigr\}$ and $e^{l}_{\gamma}\times
I\cap(inta_{\alpha}\times I)_{\gamma_{\alpha}}\neq\emptyset$.

If $j_{\alpha}>\gamma$, then, by (1), $inta_{\alpha}$ is disjoint
from $e^{l}_{\gamma}\times I$. Now suppose that
$j_{\alpha}<\gamma$. Since $j_{\alpha}\in
I(w_{\gamma_{\alpha}},l)$,  by Lemma 2.2.5, $j_{\alpha}\notin
I(w_{\gamma},l)$. By Proposition 4(3), $a_{\alpha}\subset
d_{j_{\alpha}}^{l}$ is disjoint from $e^{l}_{\gamma}\times I$.

Since $\gamma_{\alpha}>\gamma$, each component of
$e^{l}_{\gamma}\times I\cap e^{l}_{\gamma_{\alpha}}\times I$ is
$(c\times I)_{\gamma_{\alpha}}$ where $c\subset
inte^{l}_{\gamma_{\alpha}}$. Now if $e^{l}_{\gamma}\times
I\cap(inta_{\alpha}\times I)_{\gamma_{\alpha}}\neq\emptyset$, then
$inta_{\alpha}\cap e^{l}_{\gamma}\times I\neq \emptyset$, a
contradiction.\qquad Q.E.D.\vskip 0.5mm

{\bf \em Lemma 3.5.4.} \ If $\alpha,\beta\neq 0$ and
$\alpha\neq\beta$, then $(a_{\alpha}\times I)_{\gamma_{\alpha}}$
is disjoint from $(a_{\beta}\times I)_{\gamma_{\beta}}$.

{\bf \em Proof.} \ There are three cases:

Case 1.  \ $j_{\alpha}, j_{\beta}\notin I(w_{\gamma},l)$ for each
$\gamma\in m(l)$ with $s(w_{\gamma})=-$.

Now by Definition 3.2.1, $(a_{\alpha}\times
I)_{\gamma_{\alpha}}=a_{\alpha}\subset d_{j_{\alpha}}^{l}$ and
$(a_{\beta}\times I)_{\gamma_{\beta}}=a_{\beta}\subset
d_{j_{\beta}}^{l}$. Since $j_{\alpha}\neq j_{\beta}$,  the lemma
holds.

Case 2. \ $j_{\alpha}\notin I(w_{\gamma},l)$ for each $\gamma\in
m(k)$ with $s(w_{\gamma})=-$, and $j_{\beta}\in
I(w_{\gamma_{\beta}},l)$ for some $\gamma_{\beta}\in m(l)$ with
$s(w_{\gamma_{\beta}})=-$.

Since $\gamma_{\beta}>j_{\beta}$, $\gamma_{\beta}>j_{0}=m^{l+1}$.
Since $\gamma_{i,\alpha}=\gamma_{i,0}$ for $\delta(0,\alpha)\leq
i\leq \theta(0,\alpha)$.  Hence $\gamma_{\beta}>\gamma_{i,0}$ for
$\delta(0,\alpha)\leq i\leq \theta(0,\alpha)$, and
$\gamma_{\beta}>Max\bigl\{\gamma_{i,0} \ | \ \delta(0)\leq
i\leq\theta(0)\bigr\}$.  By the definition,
$a_{\alpha}=\cup_{i=\delta(0,\alpha)}^{\theta(0,\alpha)}d_{i,\alpha}
\cup_{i=\delta(0,\alpha)}^{\theta(0,\alpha)}e_{i,\alpha}$ and
$a_{\beta}=\cup_{i=\delta(0,\beta)}^{\theta(0,\beta)}d_{i,\beta}
\cup_{i=\delta(0,\beta)}^{\theta(0,\beta)}e_{i,\beta}$.  Now if
$j_{\alpha}<\gamma_{\beta}$, then, by Proposition 4(3),
$a_{\alpha}$ is disjoint from $(a_{\beta}\times
I)_{\gamma_{\beta}}$.

Suppose that $j_{\alpha}>\gamma_{\beta}$. By Proposition 4(2),
either $d_{j_{\alpha}}^{l}\cap e^{l}_{\gamma_{\beta}}\times I
=\emptyset$ or each component  of $d_{j_{\alpha}}^{l}\cap
e^{l}_{\gamma_{\beta}}\times I$ is a core of
$e^{l}_{\gamma_{\beta}}\times I$. Hence either
$d_{j_{\alpha}}^{l}\cap (a_{\beta}\times I)_{\gamma_{\beta}}
=\emptyset$ or each component  of $d_{j_{\alpha}}^{l}\cap
(a_{\beta}\times I)_{\gamma_{\beta}}$ is a core of
$(a_{\beta}\times I)_{\gamma_{\beta}}$. By Lemma 3.2.2,
$d_{\delta(0,\alpha),\alpha}\cup d_{\theta(0,\alpha),\alpha}$ is
disjoint from $(d_{\delta(0,\beta),\beta}\times
I)_{\gamma_{\beta}}\cup (d_{\theta(0,\beta),\beta})\times
I)_{\gamma_{\beta}}$. Now if $a_{\alpha}\cap (a_{\beta}\times
I)_{\gamma_{\beta}}\neq\emptyset$, then $inta_{\alpha}\cap
(a_{\beta}\times I)_{\gamma_{\beta}}\neq\emptyset$. By  Lemma
3.5.3(1), this is impossible.

Case 3. \ $s(w_{\gamma_{\alpha}})=-$, $s(w_{\gamma_{\beta}})=-$.

Now if $\gamma_{\alpha}=\gamma_{\beta}$, then the lemma holds.

Suppose that $\gamma_{\alpha}<\gamma_{\beta}$.  By Proposition
4(1), each component of $e^{l}_{\gamma_{\beta}}\times I\cap
e^{l}_{\gamma_{\alpha}}\times I$ is $(c\times I)_{\gamma_{\beta}}$
where $c\subset e^{l}_{\gamma_{\beta}}$ is a core of
$e^{l}_{\gamma_{\alpha}}\times (0,1)$. By (2), $a_{\beta}$ is
disjoint from $(a_{\alpha}\times I)_{\gamma_{\alpha}}$. Hence the
lemma holds. \qquad Q.E.D.\vskip 0.5mm

{\bf \em Lemma 3.5.5.}  \ For each $\alpha\neq 0$, there are an
arc $a_{0,\alpha}$ in $E_{f_{\delta(0,\alpha),\alpha}}^{l}$ and an
arc $b_{0,\alpha}$ in $E_{f_{\theta(0,\alpha),\alpha}}^{l}$ such
that

(1) \ $\partial_{1} a_{0,\alpha}=\partial_{1}a_{\alpha}$ and
$\partial_{2} a_{0,\alpha}=\partial_{1} a_{\alpha}^{0}$,
$\partial_{1} b_{0,\alpha}=\partial_{2} a_{\alpha}$ and
$\partial_{2} b_{0,\alpha}=\partial_{2} a_{\alpha}^{0}$,

(2) \ $a_{0,\alpha}\cup a_{\alpha}\cup b_{0,\alpha}\cup
a_{\alpha}^{0}$ bounds a disk $D^{*}_{0,\alpha}$ in $F^{l}$.

(3) \ $d_{j}^{l}\cap a_{0,\alpha} \subset
\cup_{\lambda=0}^{\alpha}(d_{\delta(0,\alpha),\lambda}\times
I)_{\gamma_{\lambda}}$,  $d_{j}^{l}\cap b_{0,\alpha} \subset
\cup_{\lambda=0}^{\alpha}(d_{\theta(0,\alpha),\lambda}\times
I)_{\gamma_{\lambda}}$.

(4) \ $e^{l}_{\gamma}\times I\cap a_{0,\alpha} \subset
\cup_{\lambda=0}^{\alpha}(d_{\delta(0,\alpha),\lambda}\times
I)_{\gamma_{\lambda}}$ and $e^{l}_{\gamma}\times I\cap
b_{0,\alpha} \subset
\cup_{\lambda=0}^{\alpha}(d_{\theta(0,\alpha),\lambda}\times
I)_{\gamma_{\lambda}}$. Furthermore, if  $e^{l}_{\gamma}\times
I\cap (inta_{0,\alpha}\cup intb_{0,\alpha})\neq\emptyset$, then
$\gamma\geq\gamma_{\lambda}$ for some
$\alpha\leq\lambda\leq\beta$.

{\bf \em Proof.} \ We may assume that $\alpha>0$. By Lemma 3.4.7,
 for $0\leq \lambda \leq \alpha$,
$\delta(0,\alpha)\leq\delta(0,\lambda)\leq\theta(0,\lambda)\leq\theta(0,\alpha)$.
By Lemma 3.4.7(3) and Definition 3.5.1, $h_{0,\alpha}^{\alpha}\cup
d_{\delta(0,\alpha),\alpha}\cup
d_{\theta(0,\alpha),\alpha}=a_{\alpha}$. By Lemma 3.4.7 and Lemma
3.5.4, $d_{\delta(0,\alpha),\lambda}\cap D_{0,\alpha}=\partial_{2}
d_{\delta(0,\alpha),\lambda}$ for $0\leq \lambda\leq \alpha$.

Now let $a_{0,\alpha}$  be in $E_{f_{\delta(0,\alpha)}}^{l}$ such
that $\partial_{1} a_{0,\alpha}=\partial_{1}
d_{\delta(0,\alpha),\alpha}=\partial_{1} a_{\alpha}$ and
$\partial_{2} a_{0,\alpha}=\partial_{1}
d_{\delta(0,\alpha),0}=\partial_{1} a_{\alpha}^{0}$ as in Figure
11(a). In fact, $a_{0,\alpha}$ can be obtained by pushing
$d_{\delta(0,\alpha),0}\cup a^{*}\cup d_{\delta(0,\alpha),\alpha}$
slightly where $a^{*}\subset \partial D_{0,\alpha}\cap
E^{l}_{f_{\delta(0,\alpha),0}}$ as in Figure 11(a). Then
$a_{0,\alpha}$ is disjoint from $D_{0,\alpha}$. Similarly, let
$b_{0,\alpha}$  be in $E_{f_{\theta(0,\alpha),0}}^{l}$ such that
$\partial_{1} b_{0,\alpha}=\partial_{2}
d_{\theta(0,\alpha),\alpha}=\partial_{2} a_{\alpha}$ and
$\partial_{2} b_{0,\alpha}=\partial_{2}
d_{\theta(0,\alpha),0}=\partial_{2} a_{\alpha}^{0}$ as in Figure
11(a). Hence $a_{0,\alpha}\cup a_{\alpha}\cup b_{0,\alpha}\cup
b_{\alpha}$ bounds a disk $D_{0,\alpha}^{*}$ in $F^{l}$.
\begin{center}
\includegraphics[totalheight=4cm]{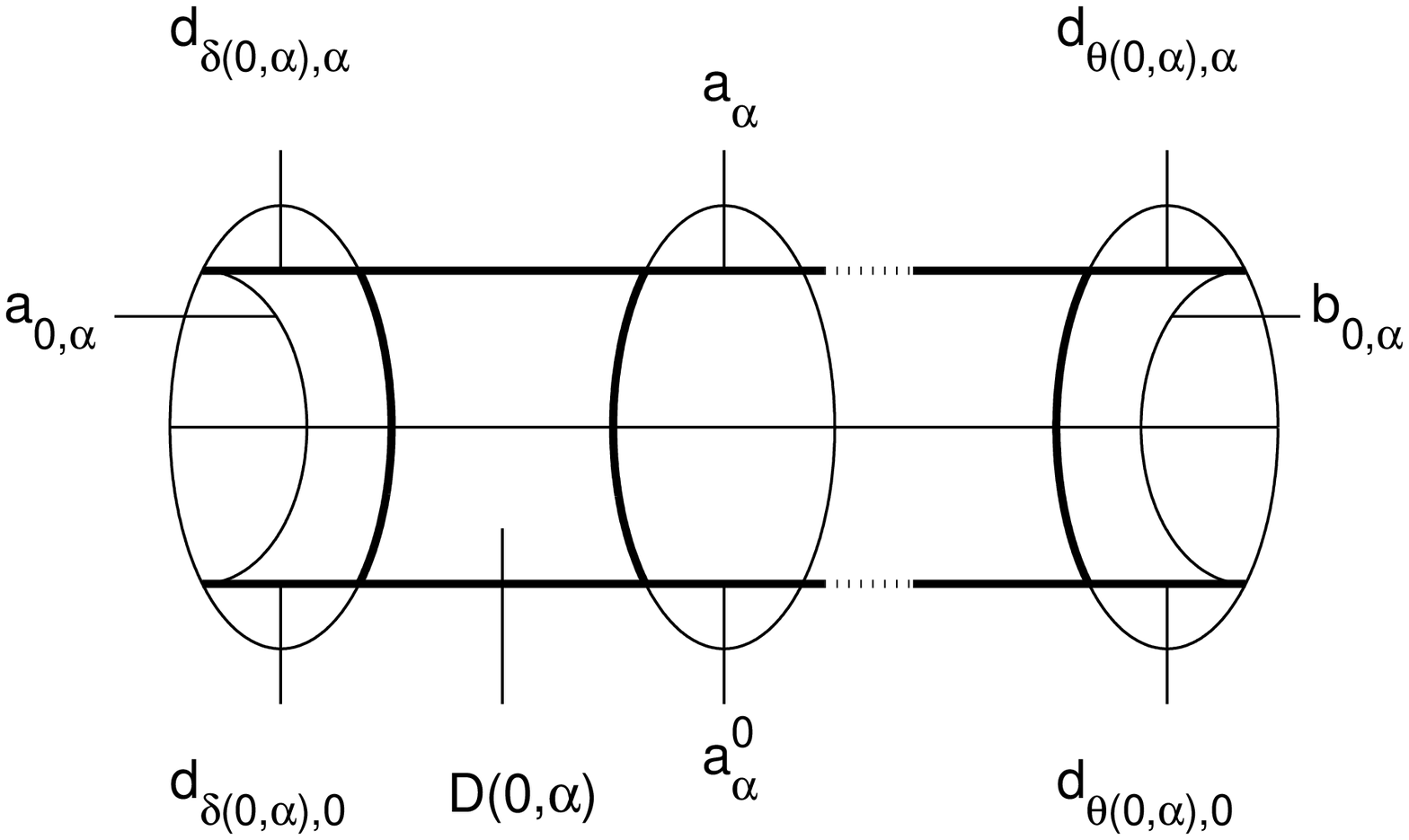}
\begin{center}
Figure 11(a)
\end{center}
\end{center}

Suppose that $d_{j}^{l}\cap a_{0,\alpha}\neq \emptyset$. Then, by
Lemma 3.1.6 and Lemma 3.2.2,  one component of $d_{j}^{l}\cap
E^{l}_{f_{\delta(0,\alpha),\alpha}}$, say $c$, lies in
$D_{0,\alpha,\delta(0,\alpha)}$. If $\partial c$ is disjoint from
$A_{0,\alpha,\delta(0,\alpha)}$, then $a_{0,\alpha}$ can be
isotoped to be disjoint from $c$. If one end point of $c$ lies in
$A_{0,\alpha,\delta(0,\alpha)}$, then, by the proof of Lemma
3.4.1, (2) holds.  In this case, either $j=j_{\lambda}$ or
$j>\gamma_{\lambda}$ for some $\alpha\leq\lambda\leq \beta$.

(3) follows from Lemma 3.2.4 and the proof of Lemma 3.4.3. \qquad
Q.E.D.

\begin{center}
\includegraphics[totalheight=4.5cm]{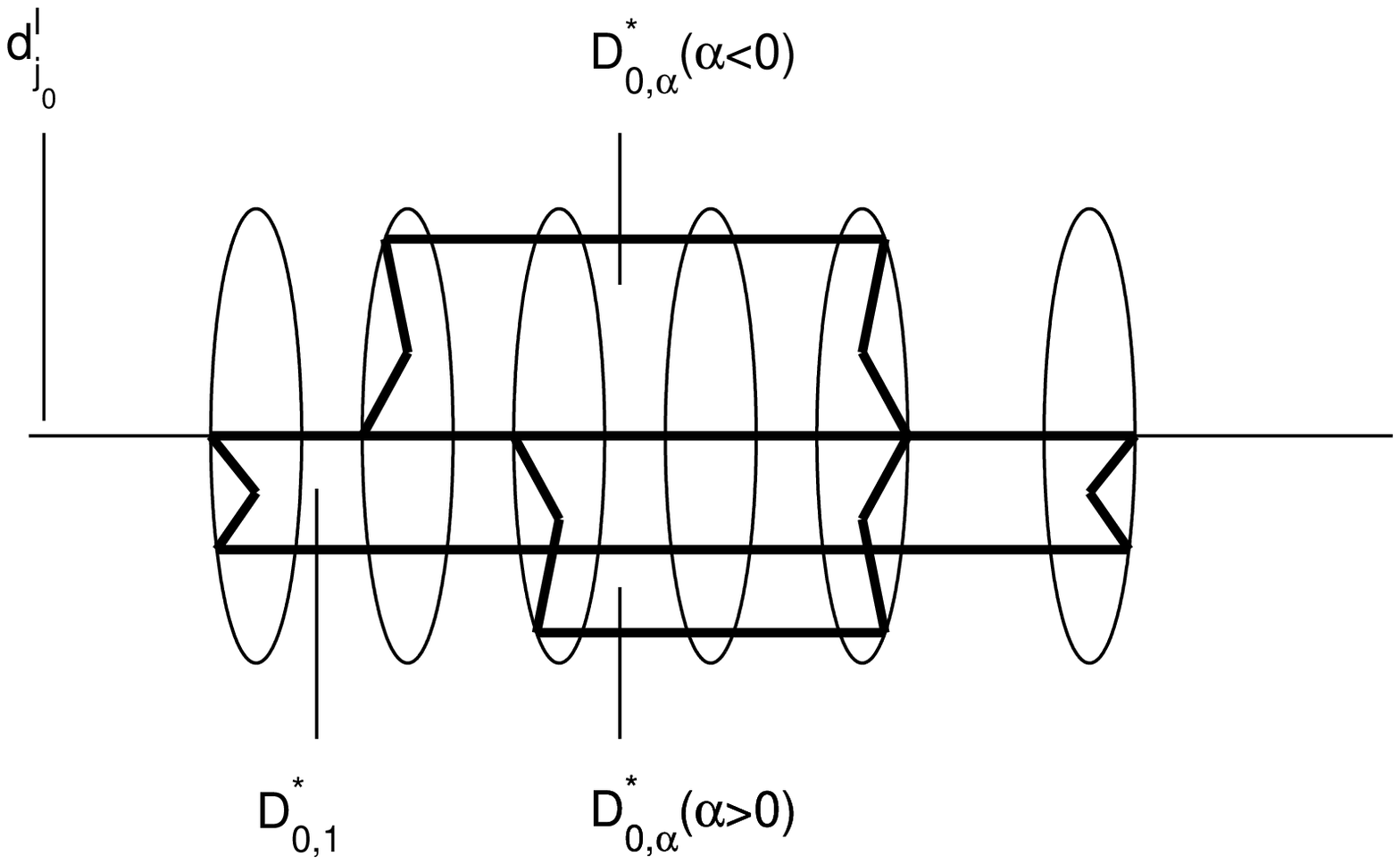}
\begin{center}
Figure 11(b)
\end{center}
\end{center}

{\bf \em Lemma 3.5.6.} \ $\cup_{\alpha\neq 0} D_{0,\alpha}^{*}$ is
a disk.

{\bf \em Proof.} \ By Lemma 3.4.8, $\cup_{\alpha\neq 0}
D_{0,\alpha}$ is a disk. By Lemma 3.5.4, Lemma 3.4.7,
$d_{\delta(0,\beta),\beta}\cap (\cup_{\alpha\neq 0}
D_{0,\alpha})=\partial_{2} d_{\delta(0,\beta),\beta}$ for
$\beta\neq 0$ as in Figure 10(b).  We may assume that
$\delta(0,-1)\geq \delta(0,1)$. Then, by Lemma 3.4.7,
$\delta(0,\alpha)\geq\delta(0,1)$ for each $\alpha$. Then
$d_{\delta(0,1),0}\cap(\cup_{\alpha\neq 0}
D_{0,\alpha})=\partial_{2} d_{\delta(0,1),0}$. Furthermore,
$d_{\delta(0,\alpha),0}\subset a_{1}^{0}$ for each $\alpha$.   By
the proof of Lemma 3.5.5, $\cup_{\alpha\neq 0} D_{0,\alpha}^{*}$
is a disk as in Figure 11(b).\qquad Q.E.D.

{\bf \em Lemma 3.5.7.} \ (1) For each $\alpha\neq 0$,
$d_{j_{0}}^{l}\cap D_{0,\alpha}^{*}=a^{0}_{\alpha}$.

(2) \ If $r<j_{0}$, then $d_{r}^{l}$ is disjoint from
$D_{0,\alpha}^{*}$.

{\bf \em Proof.} \ By Lemma 3.5.4, $a^{0}_{\alpha}\subset
\partial D_{0,\alpha}^{*}$. Since $j_{0}<j_{\alpha}$ for each
$\alpha\neq 0$. Hence, by Proposition 4(3), $d_{j_{0}}^{l}$ is
disjoint from $(a_{\alpha}\times I)_{\gamma_{\alpha}}$ for
$\alpha\neq 0$. By Lemma 3.5.5(3) and Lemma 3.4.7,
$d_{j_{0}}^{l}-a_{\alpha}^{0}$ is disjoint from
$D_{0,\alpha}^{*}$.

Suppose that $r<j_{0}$. Then $r<j_{\alpha}$ for each $\alpha$.
Hence $d_{r}^{l}$ is disjoint from $(a_{\alpha}\times
I)_{\gamma_{\alpha}}$. By Lemma 3.5.5(3) and Lemma 3.4.7, (2)
holds.\qquad Q.E.D.

\begin{center}
\includegraphics[totalheight=5cm]{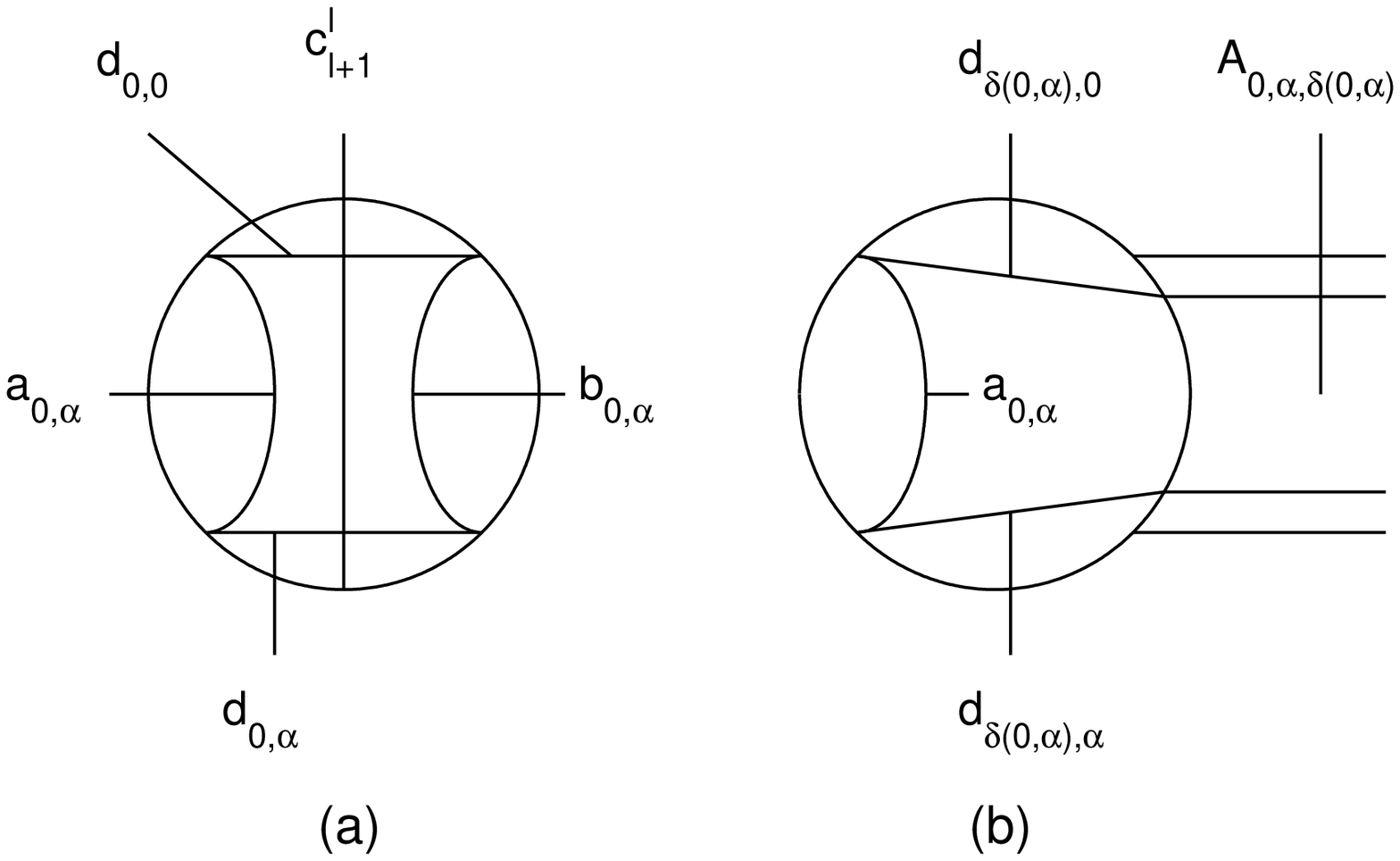}
\begin{center}
Figure 12
\end{center}
\end{center}

{\bf \em Lemma 3.5.8.} \ (1) If
$\delta(0,\alpha)=\theta(0,\alpha)=0$, $c^{l}_{i}$ intersects
$a_{0,\alpha}\cup b_{0,\alpha}$ in one point if and only if
$c^{l}_{i}$ intersects $d_{0,\alpha}\cup d_{0,0}$ in one point.

(2) \ If $\delta(0,\alpha)\neq \theta(0,\alpha)$, then $c^{l}_{i}$
intersects $a_{0,\alpha}$ in a one point if and only if
$c^{l}_{i}$ intersects $d_{\delta(0,\alpha),\alpha}\cup
d_{\delta(0,\alpha),0}$ in one point,   $c^{l}_{i}$ intersects
$b_{0,\alpha}$ in one point if and only if $c^{l}_{i}$ intersects
$d_{\theta(0,\alpha),\alpha}\cup d_{\theta(0,\alpha),0}$ in one
point.

{\bf \em Proof.} \ (1) \ Since
$\delta(0,\alpha)=\theta(0,\alpha)=0$, $a_{\alpha}=d_{0,\alpha}$
and $a_{\alpha}^{0}=d_{0,0}$. Hence $a_{0,\alpha}\subset
E^{l}_{0}$ and $b_{0,\alpha}\subset E^{l}_{0}$ as in Figure 12(a).
By Lemma 3.1.2, $\partial d_{0,0}\cup \partial d_{0,\alpha}\subset
\partial d_{j_{0}}^{l}\cup \partial d_{j_{\alpha}}^{l}\cup (\partial e^{l}_{\gamma})\times
I$. By Proposition 4(4),(5), $d_{0,0}\cup d_{0,\alpha}$ is
disjoint from $\partial c^{l}_{i}$,  and $c^{l}_{i}$ is properly
embedded in $E^{l}_{f}$. Since $j_{0},j_{\alpha}\in
L(c_{l+1}^{l})$,  $c_{l+1}^{l}$ intersects each of $d_{0,0}$ and
$d_{0,\alpha}$ in one point, we can moved $c_{l+1}^{l}$ so that
$a_{0,\alpha}\cup b_{0,\alpha}$ is disjoint from $c_{l+1}^{l}$.
Similarly, if $c^{l}_{i}$ intersects each of $d_{0,0}$ and
$d_{0,\alpha}$ in one point, then $c^{l}_{i}$ is disjoint from
$a_{0,\alpha}\cup b_{0,\alpha}$.  Hence (1) holds.

(2) By Lemma 3.3.1, $d_{\delta(0,\alpha),0}$ and
$d_{\delta(0,\alpha),\alpha}$ are properly embedded arcs in
$E^{l}_{f_{\delta(0,\alpha),\alpha}}$ and $d_{\theta(0,\alpha),0}$
and $d_{\theta(0,\alpha),\alpha}$ are properly embedded arcs in
$E^{l}_{f_{\theta(0,\alpha),\alpha}}$. Since $\delta(0,\alpha)\neq
\theta(0,\alpha)$, by Lemma 3.3.1,
$f_{\delta(0,\alpha),\alpha}\neq f_{\theta(0,\alpha),\alpha}$. By
Proposition 4(5), $\partial c^{l}_{i}$ is disjoint from
$e^{l}_{\gamma}\times I$ and $d_{j}^{l}$. Hence $c^{l}_{i}$ is
disjoint from $A_{0,\alpha,\delta(0,\alpha)}\cup
A_{0,\alpha,\theta(0,\alpha)}$. Note that $\partial_{2}
d_{\delta(0,\alpha),0}\cup
\partial_{2} d_{\delta(0,\alpha),\alpha}\subset
A_{0,\alpha,\delta(0,\alpha)}$ and  $\partial_{1}
d_{\theta(0,\alpha),0}\cup
\partial_{1} d_{\theta(0,\alpha),\alpha}\subset
A_{0,\alpha,\theta(0,\alpha)}$ as in Figure 12(b). By the argument
in (1), (2) holds. \qquad Q.E.D.

{\bf \em Remark 3.5.9.} \ Since Propositions 1-6 hold for $k\leq
l$. By Lemma 3.1.5, $m^{k}=MinL(v_{k+1}^{k})$ for $k\leq l$.  By
Propositions 1 and 2, all the arguments in Sections 3.1-3.5 are
true when we take place of $c^{l}_{i}, d^{l}_{j}, e^{l}_{\gamma},
F^{l}$ with $v^{l}_{i}, w^{l}_{j}, b^{l}_{\gamma}, P^{l}$,
$s(w_{\gamma})=-$ with $s(w_{\gamma})=+$ in Lemma 3.1.6 and Lemma
3.2.3, $s(w_{\gamma})=s(w_{\gamma_{j}})=-$ with
$s(w_{\gamma})=s(w_{\gamma_{j}})=+$ in Lemma 3.2.4.

\section{The Proofs of Propositions 4-6 for the case: $k=l+1$ and $s(v_{l+1})=+$}

Assume now that Propositions 1-6 hold for each $0\leq k\leq l.$
Now we only need to prove that Propositions 1-6 hold for $k=l+1$.

There are two cases:

1. \ $s(v_{l+1})=+$.

2. \ $s(v_{l+1})=-$.

In this chapter, we shall prove that Propositions 4-6 hold for
Case 1. Hence, in this chapter, we assume that $s(v_{l+1})=+$. We
first construct $c^{l+1}_{i}, d_{j}^{l+1}, e^{l+1}_{\gamma}\times
I, F^{l+1}$ from $c^{l}_{i}, d_{j}^{l}, e^{l}_{\gamma}\times I,
F^{l}$, then we shall prove  Propositions 4-6 for this case.

\subsection{The element of the construction}

\ \ \ \ \ Since Propositions 1-6 hold for $k\leq l$. So
$L(c^{l}_{i})$, $L(d_{j}^{l})$, and $m(k)$ are well defined for
$k\leq l$. Recall that $m^{l+1}$ defined in Definition 3.3.0.

{\bf \em Definition 4.1.1.} \ Let
$m(l+1)=m(l)\cup\bigl\{m^{l+1}\bigr\}$. In particular, if
$m^{l+1}=\emptyset$, then let $m(l+1)=m(l)$.\vskip 2mm

Without loss of generality, we assume that $c_{l+1}^{l}\subset
E^{l}_{0}$. Recalling the argument in Section 3.3, we can number
all the elements in $L(c^{l}_{l+1})$ as $\ldots, j_{-1},
j_{0}=m^{l+1}, j_{1},\ldots$ according to the order of
$(\cup_{j\in L(c_{l+1}^{l})} d_{j,0}^{l})\cap c_{l+1}^{l}$ lying
in $c_{l+1}^{l}$. In some time, we shall use $j_{0}$ to take place
of $m^{l+1}$. By Lemma 3.3.1 and Remark 3.3.2, for each
$j_{\alpha}$,
$d_{j_{\alpha}}^{l}=\cup_{i=\delta(\alpha)}^{\theta(\alpha)}d_{i,\alpha}\cup_{i=\delta(\alpha)}^{\theta(\alpha)}
e_{i,\alpha}$.\vskip 0.5mm

{\bf \em Lemma 4.1.2.} \ Suppose that
$j_{0}=m^{l+1}\neq\emptyset$.

(1) \ In $F^{l}$, $c_{l+1}^{l}$ intersects $d_{j_{0}}^{l}$ in one
point. Furthermore, if $r<j_{0}$, then $d_{r}^{l}$ is disjoint
from $c_{l+1}^{l}$; if $\gamma<j_{0}$, then $e^{l}_{\gamma}\times
I$ is disjoint from $c_{l+1}^{l}$.

(2) \ In $P^{l}$, $v_{l+1}^{l}$ intersects $w_{j_{0}}^{l}$ in one
point. Furthermore, if $r<j_{0}$, then $v_{l+1}^{l}$ is disjoint
from $w_{r}^{l}$; if $\gamma<j_{0}$, then $b^{l}_{\gamma}\times I$
is disjoint from $v_{l+1}^{l}$.

{\bf \em Proof.} \ By assumptions, Propositions 1-6 hold for
$k=l$.

(1) \ By Proposition 4(4), $c_{l+1}^{l}$ is a properly embedded
arc in $F^{l}$ lying in $E^{l}_{0}$. Since $j_{0}=m^{l+1}\in
L(c_{l+1}^{l})$, $d_{0,0}$ intersects $c_{l+1}^{l}$ in one point.
By Proposition 4(5) and Lemma 3.1.2, $d_{0,0}\cap
c_{l+1}^{l}=intc_{l+1}^{l}\cap intd_{0,0}$. Suppose, otherwise,
$d_{j_{0}}^{l}$ intersects $c_{l+1}^{l}$ in at least two points.
Then, by Lemma 3.1.6, $r\in L(c_{l+1}^{l})$ for some $r<j_{0}$,
contradicting the minimality of $j_{0}$. Similarly,  if $r<j_{0}$,
then $c_{l+1}^{l}$ is disjoint from $d_{r}^{l}$. By Lemma 3.2.4,
if $\gamma<j_{0}$, then $e^{l}_{\gamma}\times I$ is disjoint from
$c_{l+1}^{l}$.

(2) follows from Remark 3.5.9 and Propositions 1 and 2. \qquad
Q.E.D.\vskip 0.5mm

{\bf \em Lemma 4.1.3.} \ (1) \ If $j_{0}=m^{l+1}\neq\emptyset$,
then $s(w_{j_{0}})=-$.

(2) \ $W_{j_{0}}^{l}$ is a properly embedded disk in $\cal W_{-}$
such that $\partial W_{j_{0}}^{l}$ intersects $c_{l+1}^{l}$ in one
point which lies in $intd_{0,0}$.

{\bf \em Proof.} \ (1) \ Now consider $P^{l}$. By assumptions,
Propositions 1-3 hold for $k=l$. Suppose, otherwise,
$s(w_{j_{0}})=+$. By assumption, $s(v_{l+1})=+$. By Proposition 3,
$V_{l+1}^{l}$ is a properly embedded disk in $\cal V_{+}$ and
$W_{j_{0}}^{l}$ is a properly embedded disk in $\cal W_{+}$. By
Definition 2.3.1 and Lemma 2.2.4,
$I(v_{l+1},l)=I(v_{l+1})-\bigl\{1,\ldots, l\bigr\}= \emptyset$. By
Lemma 3.1.5, $m^{l+1}=MinL(v_{l+1}^{l})$. Note that
$j_{0}=m^{l+1}$. By Proposition 3, $V_{l+1}^{l}\cap
W_{j_{0}}^{l}=(v_{l+1}^{l}\cup_{r\in
I(v_{l+1},l)}v_{r}^{l})\cap(w_{j_{0}}^{l}\cup_{r\in
I(w_{j_{0}},l)} w_{r}^{l})=v_{l+1}^{l}\cap(w_{j_{0}}^{l}\cup_{r\in
I(w_{j_{0}},l)} w_{r}^{l})$. If $r\in I(w_{j_{0}},l)$, then, by
Lemma 2.2.4, $r<j_{0}$. By Lemma 4.1.2(2), $w_{r}^{l}$ is disjoint
from $v_{l+1}^{l}$. By also Lemma 4.1.2(2), $V_{l+1}^{l}$
intersects $W_{j_{0}}^{l}$ in only one point. Hence $\cal
V_{+}\cup W_{+}$ is stabilized, a contradiction.

(2) \ By (1) and Proposition 6, $W_{j_{0}}^{l}$ is a properly
embedded disk in $\cal W_{-}$.  By Proposition 6,
$W_{j_{0}}^{l}\cap c_{l+1}^{l}=(d_{j_{0}}^{l}\cup_{r\in
I(w_{j_{0}},l)} d_{r}^{l})\cap c_{l+1}^{l}$. Since $r\in
I(w_{j_{0}},l)$, $r<j_{0}$.  By Lemma 4.1.2, $\partial
W_{j_{0}}^{l}$ intersects $c_{l+1}^{l}$ in one point lying in
$intd_{0,0}$. \qquad Q.E.D.\vskip 0.5mm

{\bf \em Lemma 4.1.4.} \ Suppose that
$j_{0}=m^{l+1}\neq\emptyset$.

(1)  \ $j_{0}\notin I(w_{\gamma},l)$ for each $\gamma\in
m(l)=\bigl\{m^{1},\ldots,m^{l}\bigr\}$ with $s(w_{\gamma})=-$.
Furthermore, if $r\in I(w_{j_{0}},l)$, then $r\notin
I(w_{\gamma},l)$ for each $\gamma\in m(l)$.

(2) \ $d_{j_{0}}^{l}$ is properly embedded in $F^{l}$,
$d_{r}^{l}$ is  properly embedded in $F^{l}$ for $r\in
I(w_{j_{0}},l)$.

{\bf \em Proof.} \  By Lemma 4.1.3, $s(w_{j_{0}})=-$. Now if
$j_{0}\in I(w_{\gamma},l)$ for some $\gamma\in m(l)$, then, by
Lemma 2.2.5, $s(w_{\gamma})=+$. By Lemma 3.1.1, $d_{j_{0}}^{l}$ is
regular in $F^{l}$. By Definition 2.1.5, $d_{j_{0}}^{l}$ is
properly embedded in $F^{l}$.

Suppose that $r\in I(w_{j_{0}},l)$. By Definition 2.3.3, $j_{0}\in
\bigl\{1,\ldots,n\bigr\}-m(l)$. Hence $j_{0}\notin m(l)$.  By
Lemma 2.2.5, $r\notin I(w_{\gamma},l)$ for each $\gamma\in m(l)$.
By Lemma 3.1.1 and Definition 2.1.5, $d_{r}^{l}$ is properly
embedded in $F^{l}$.\qquad Q.E.D. \vskip 0.5mm

{\bf \em Lemma 4.1.5.} \ Suppose that
$j_{0}=m^{l+1}\neq\emptyset$. Then there is a regular neighborhood
of $\partial W_{j_{0}}^{l}$ in $\partial_{+}\cal V_{-}$, say
$\partial W_{j_{0}}^{l}\times I$, satisfying the following
conditions:

(1) \ $\partial W_{j_{0}}^{l}\times I\cap
F^{l}=(d_{j_{0}}^{l}\times I)_{j_{0}}\cup_{r\in I(w_{j_{0}},l)}
(d_{r}^{l}\times I)_{j_{0}}$, $(\partial d_{j_{0}}^{l})\times
I,(\partial d_{r}^{l})\times I\subset \cup_{f} \partial
E^{l}_{f}-\cup_{\gamma\in m(l)} e^{l}_{\gamma}\times I$.

(2) \ If $j_{0}<\gamma$,  then $(d_{j_{0}}^{l}\times I)_{j_{0}}$
is disjoint from $e^{l}_{\gamma}\times I$, if $r\in
I(w_{j_{0}},l)$ and $r<\gamma$, then $(d_{r}^{l}\times I)_{j_{0}}$
is disjoint from $e^{l}_{\gamma}\times I$.

(3) \ If $j_{0}>\gamma$,  then each component of
$(d_{j_{0}}^{l}\times I)_{j_{0}}\cap e^{l}_{\gamma}\times I$ is
$(c\times I)_{j_{0}}\subset e^{l}_{\gamma}\times (0,1)$ where
$c\subset intd_{j_{0}}^{l}$ is a core of $e^{l}_{\gamma}\times I$,
if $r\in I(w_{j_{0}},l)$ and $r>\gamma$, then each component of
$(d_{r}^{l}\times I)_{j_{0}}\cap (e^{l}_{\gamma}\times I)$ is
$(c\times I)_{j_{0}}\subset e^{l}_{\gamma}\times (0,1)$ where
$c\subset intd_{r}^{l}$ is a core of $e^{l}_{\gamma}\times I$.

(4) \ If $j\neq j_{0}$ and $j\notin I(w_{j_{0}},l)$, then
$d_{j}^{l}$ is disjoint from $\partial W_{j_{0}}^{l}\times I$.

(5) \ For each $i\geq l+1$, $\partial c_{i}^{l}$ is disjoint from
$\partial W_{j_{0}}^{l}\times I$.

(6) \ If $j_{0}>\gamma$, then $(intd_{i,0}\times I)_{j_{0}}$ is
disjoint from $e^{l}_{\gamma}\times I$ for $\delta(0)\leq i\leq
\theta(0)$.

(7) \ For each $i\geq l+2$ with $s(v_{i})=-$, $\partial
W_{j_{0}}^{l}\times I$ is disjoint from $\partial
V_{i}^{l}-c_{i}^{l}\cup_{r\in I(v_{i},l)} v_{r}^{l}$, for each
$j\notin m(l+1)$ with $s(w_{j})=-$, $\partial W_{j_{0}}^{l}\times
I$ is disjoint from $\partial W_{j}^{l}$.

(8) \ $\partial W_{j_{0}}^{l}\times I\cap D_{0,\alpha}^{*}\subset
a_{\alpha}^{0}\times I$.

{\bf \em Proof.} \ Now $F^{l}=\cup_{f} E^{l}_{f}\cup_{\gamma\in
m(l)} e^{l}_{\gamma}\times I$. By Lemma 3.1.1 and Lemma 4.1.4(1),
$d_{j}^{l}, d_{r}^{l}$ are regular in $F^{l}$. By Definition 2.1.5
and Lemma 2.1.6, $\partial d_{j}^{l},
\partial d_{r}^{l}$ are disjoint from $e^{l}_{\gamma}\times I$ for each $\gamma\in m(l)$.
Hence (1) follows from Proposition 6(1). (2) follows from
Proposition 4(3) and Lemma 4.1.4(1). (3) follows from Proposition
4(2) and Definition 2.1.5. By Proposition 4, $\bigl\{d_{j}^{l} \ |
\ j\notin\bigl\{1,\ldots,n\bigr\}-m(l)\bigr\}$ is a set of
pairwise disjoint arcs in $F^{l}$. Since $j\neq j_{0}$ and
$j\notin I(w_{j_{0}},l)$, $d_{j}^{l}\cap (d_{j_{0}}^{l}\cup_{r\in
I(w_{j_{0}},l)} d_{r}^{l})=\emptyset$. By (1), (4) holds. Since
$c^{l}_{i}\subset F^{l}$, by Proposition 4(5), $\partial
c_{i}^{l}$ is disjoint from $\partial W_{j_{0}}^{l}$ for each
$i\geq l+1$. Hence (5) holds. (6) follows from Lemma 3.1.2 and
Lemma 3.3.1. (7) follows  from Proposition 6. (8) follows from
Lemma 3.5.7. \qquad Q.E.D.\vskip 0.5mm

\subsection{ The proofs of Propositions 4-6 for  one special case}

\ \ \ \ \ In this section, we shall prove Propositions 4-6 for the
speical case: $k=l+1$ and $L(c_{l+1}^{l})=\emptyset$.

{\bf \em The proof of Proposition 4.} \  Suppose that
$L(c_{l+1}^{l})=\emptyset$. Now by Corollary 3.2.5, for each
$j\notin m(l)$ and $\gamma\in m(l)$, $d_{j}^{l}$ and
$e^{l}_{\gamma}\times I$ are disjoint from $c_{l+1}^{l}\times I$
where $c_{l+1}^{l}\times I$ is a regular neighborhood of
$c_{l+1}^{l}$ in $E^{l}_{0}$. Obviously, $c_{i}^{l}$ is disjoint
from $c_{l+1}^{l}\times I$ for $i\geq l+2$. Now let
$F^{l+1}=F^{l}-c_{l+1}^{l}\times (-1,1)$. We denote by
$E^{l+1}_{f}$ the disk $E^{l}_{f}$ for $1\leq f\leq l$,
$E^{l+1}_{0}, E^{l+1}_{l+1}$ the two components of
$E^{l}_{0}-c_{l+1}^{l}\times (-1,1)$. Let $c_{i}^{l+1}=c_{i}^{l}$,
$d_{j}^{l+1}=d_{j}^{l}$ for $j\notin m(l)$,
$e^{l+1}_{\gamma}\times I=e^{l}_{\gamma}\times I$ for $\gamma\in
m(l)$. Specially, let $m^{l+1}=\emptyset$,
$e^{l+1}_{m^{l+1}}=\emptyset$. Now if $c^{l}_{i}\subset
E^{l}_{0}$, then $c^{l+1}_{i}$ lies in one of $E^{l+1}_{0}$ and
$E^{l+1}_{l+1}$. Similarly, if $d_{j,0}^{l}\neq \emptyset$, then
$d_{j,0}^{l}$ lies in one of $E^{l+1}_{0}$ and $E^{l+1}_{l+1}$.
Since $m^{l+1}=\emptyset$, $I(w_{\gamma},l+1)=I(w_{\gamma},l)$.
Hence Proposition 4 holds.\qquad Q.E.D.\vskip 0.5mm

{\bf \em The proof of Proposition 5.} \ By the above argument,
$d_{j}^{l+1}-\cup_{\gamma<j}inte^{l+1}_{\gamma}\times
I=d_{j}^{l}-\cup_{\gamma<j}inte^{l}_{\gamma}\times I$,
$c_{i}^{l+1}=c_{i}^{l}$ for $i\geq l+2$. Note that $j\notin
L(c_{l+1}^{l})$. Hence Proposition 5 holds.\qquad Q.E.D.\vskip
1.5mm

{\bf \em The proof of Proposition 6.} \ Since $m^{l+1}=\emptyset$,
$m(l+1)=m(l)$. By Definition 3.2.1, $I(w_{j},l+1)=I(w_{j},l)$. Now
we denote by $W_{j}^{l+1}$ the disk $W_{j}^{l}$ for $j\notin m(l)$
with $s(w_{j})=-$. Then $W_{j}^{l+1}\cap
F^{l+1}=d_{j}^{l+1}\cup_{r\in I(w_{j},l+1)} d_{r}^{l+1}$.

By Definition 2.3.1, $I(v_{i},l+1)=I(v_{i},l)-\bigl\{l+1\bigr\}$.
We denote by $V_{i}^{l+1}$ the disk $V_{i}^{l}$ for $i\geq l+2$
with $s(v_{i})=-$. Since $c_{l+1}^{l}$ is disjoint from $F^{l+1}$
and $d_{j}^{l+1}$. Hence $V_{i}^{l+1}\cap
F^{l+1}=c_{i}^{l+1}\cup_{r\in I(v_{i},l+1)} c_{r}^{l+1}$, and
$V_{i}^{l+1}\cap W_{j}^{l+1}=(c_{i}^{l+1}\cup_{r\in I(v_{i},l+1)}
c_{r}^{l+1})\cap (d_{j}^{l+1}\cup_{r\in I(w_{j},l+1)}
d_{r}^{l+1})=V^{l+1}_{i}\cap W_{j}^{l+1}\cap F^{l+1}$.\qquad
Q.E.D.

\subsection{Constructions}

\ \ \ \ \ By the argument in Section 4.2, in the following
argument, we shall assume that $m^{l+1}\neq\emptyset$. In this
section, we shall construct $c^{l+1}_{i}$ for $i\geq l+2$,
$d_{j}^{l+1}$ for $j\notin m(l+1)=m(l)\cup\bigl\{m^{l+1}\bigr\}$,
$e^{l+1}_{\gamma}\times I$ for $\gamma\in m(l+1)$, and $F^{l+1}$
from $c^{l}_{i}, d_{j}^{l}, e^{l}_{\gamma}\times I, F^{l}$.

Let $\partial W_{j_{0}}^{l}\times I$ be a regular neighborhood of
$\partial W_{j_{0}}^{l}$ in $\partial_{+} \cal V_{-}$  in Lemma
4.1.5. Let $c_{l+1}^{l}\times I$ be a regular neighborhood of
$c_{l+1}^{l}$ in $E^{l}_{0}$.  By Lemma 4.1.3, $\partial
W_{j_{0}}^{l}$ intersects $c_{l+1}^{l}\times I$ in an arc
$a\subset intd_{0,0}$, and $\partial W_{j_{0}}^{l}\times I$
intersects $c_{l+1}^{l}\times I$ in  $a\times I$. (See Figure 13.)

\begin{center}
\includegraphics[totalheight=6cm]{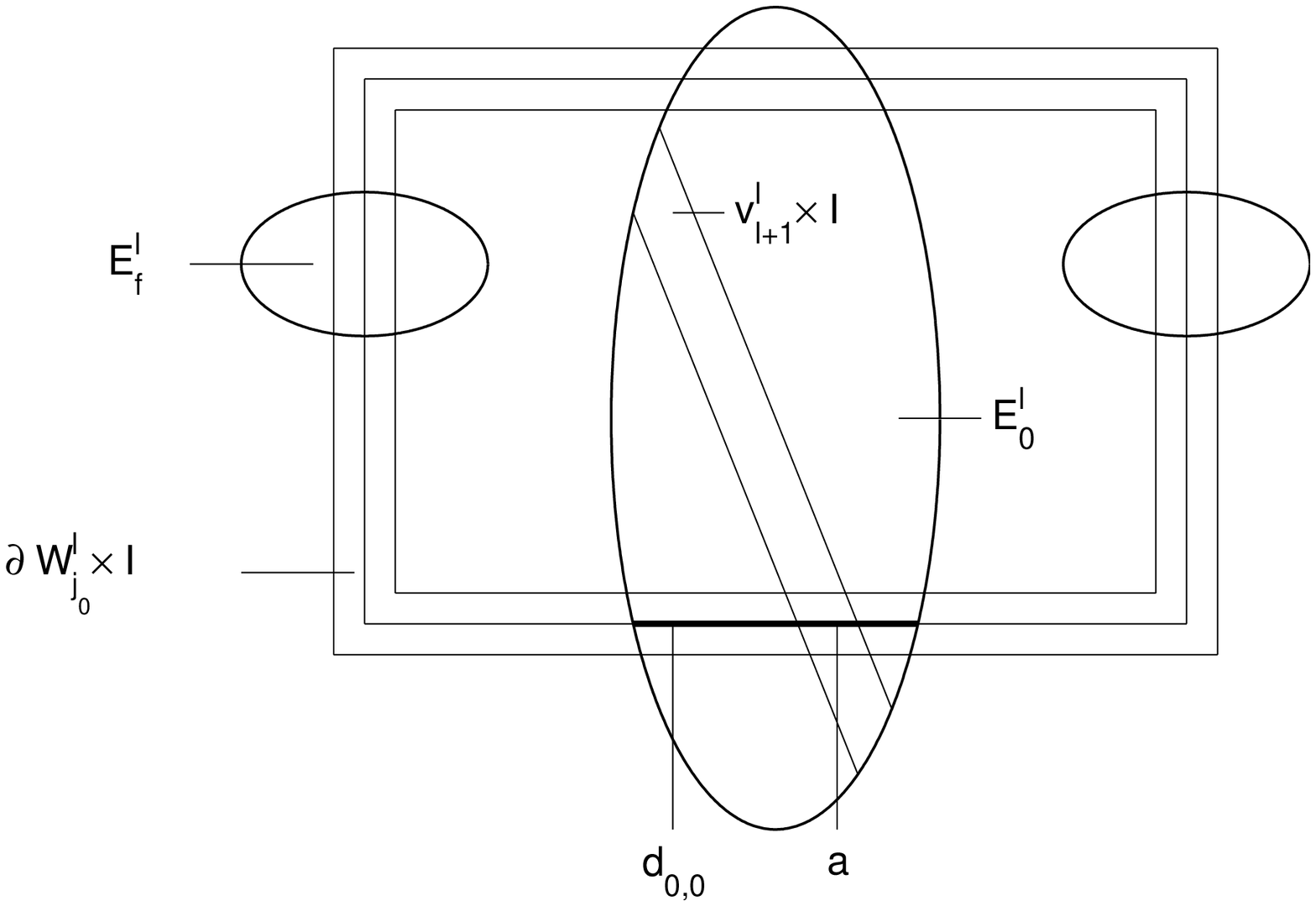}
\begin{center}
Figure 13
\end{center}
\end{center}

{\bf \em Definition 4.3.1.} \ Let $F^{l+1}=(F^{l}
-c_{l+1}^{l}\times (-1,1)) \cup (\partial
W_{j_{0}}^{l}-inta)\times I$.\vskip 0.5mm

{\bf \em Lemma 4.3.2.} \ $F^{l+1}$ is a compact surface in
$\partial_{+} \cal V_{-}$ such that $c_{l+1}^{l}$ is disjoint from
$F^{l+1}$.

{\bf \em Proof.} \  By Lemma 4.1.4, $d_{j_{0}}^{l}$ and
$d_{r}^{l}$ for $r\in I(w_{j_{0}},l)$ are properly embedded in
$F^{l}$. By Lemma 4.1.5(1) and Lemma 4.1.3, the lemma holds.\qquad
Q.E.D.\vskip 0.5mm

{\bf \em Definition 4.3.3.} \ Let $E^{l+1}_{f}=E^{l}_{f}$ for
$1\leq f\leq l$. Let $E^{l+1}_{0}$ be the component of
$E^{l}_{0}-c_{l+1}^{l}\times (-1,1)$ containing $c_{l+1}^{l}\times
\bigl\{-1\bigr\}$, $E_{l+1}^{l+1}$  be the component of
$E_{0}^{l}-c_{l+1}^{l}\times (-1,1)$ containing $c_{l+1}^{l}\times
\bigl\{1\bigr\}$. \vskip 0.5mm

{\bf \em Definition 4.3.4.} \  Let $c_{i}^{l+1}=c_{i}^{l}$ for
$i\geq l+2$.\vskip 0.5mm

%\begin{center}
%\includegraphics[totalheight=5.5cm]{s6.eps}
%\begin{center}
%Figure 14
%\end{center}
%\end{center}

{\bf \em Construction (*)}

By Lemma 3.3.3, in $F^{l}$, $(\cup_{j\notin m(l)}
d_{j}^{l}\cup_{\gamma\in m(l)} e^{l}_{\gamma}\times I)\cap
c_{l+1}^{l}\times I=(\cup_{\alpha} (d_{0,\alpha}\times
I)_{\gamma_{\alpha}})\cap c_{l+1}^{l}\times I$ where
$(d_{0,\alpha}\times I)_{\gamma_{\alpha}}$ is defined in
Definition 3.2.1. By Lemma 4.1.4, $(d_{j_{0}}^{l}\times
I)_{\gamma_{0}}=d_{j_{0}}^{l}$. Note that $j_{0}=m^{l+1}$. Hence,
$(\cup_{j\notin m(l+1)} d_{j}^{l}\cup_{\gamma\in m(l)}
e^{l}_{\gamma}\times I)\cap c_{l+1}^{l}\times I=(\cup_{\alpha\neq
0} (d_{0,\alpha}\times I)_{\gamma_{\alpha}})\cap c_{l+1}^{l}\times
I$. By Definition 3.5.1, $d_{0,\alpha}\subset a_{\alpha}$.

Suppose that $\alpha\neq 0$. If $j_{\alpha}\in
I(w_{\gamma_{\alpha}})$ for $\gamma_{\alpha}\in m(l)$, then
$\gamma_{\alpha}>j_{\alpha}>j_{0}$. If $j_{\alpha}\notin
I(w_{\gamma},l)$ for each $\gamma\in m(l)$, then
$(d_{j_{\alpha}}^{l}\times
I)_{\gamma_{\alpha}}=d_{j_{\alpha}}^{l}$.  By Lemma 4.1.5(2) and
(4), $\partial W_{j_{0}}^{l}\times I$ is disjoint from
$(a_{\alpha}\times I)_{\gamma_{\alpha}}\subset
(d_{j_{\alpha}}^{l}\times I)_{\gamma_{\alpha}}$ for $\alpha\neq
0$. By Lemma 4.1.5(8), $\partial W_{j_{0}}^{l}\times I\cap
D_{0,\alpha}^{*}\subset a^{0}_{\alpha}\times I$. Now
$a_{1}^{0}\cup a_{-1}^{0}$ separates $\cup_{\alpha\neq 0}
D_{0,\alpha}^{*}$ into two disks $\cup_{\alpha< 0}
D_{0,\alpha}^{*}$ and $\cup_{\alpha>0} D_{0,\alpha}^{*}$. We may
assume that

(0) \ $D_{0,\alpha}^{*}\cap a^{0}_{\alpha}\times I\subset
a^{0}_{\alpha}\times [0,1]$ for $\alpha>0$, and
$D_{0,\alpha}^{*}\cap a^{0}_{\alpha}\times I\subset
a^{0}_{\alpha}\times [-1,0]$ for $\alpha<0$.

Now $\partial W_{j_{0}}^{l}\times I\cup_{\alpha} D_{0,\alpha}^{*}$
is as in Figure 14.
\begin{center}
\includegraphics[totalheight=4.5cm]{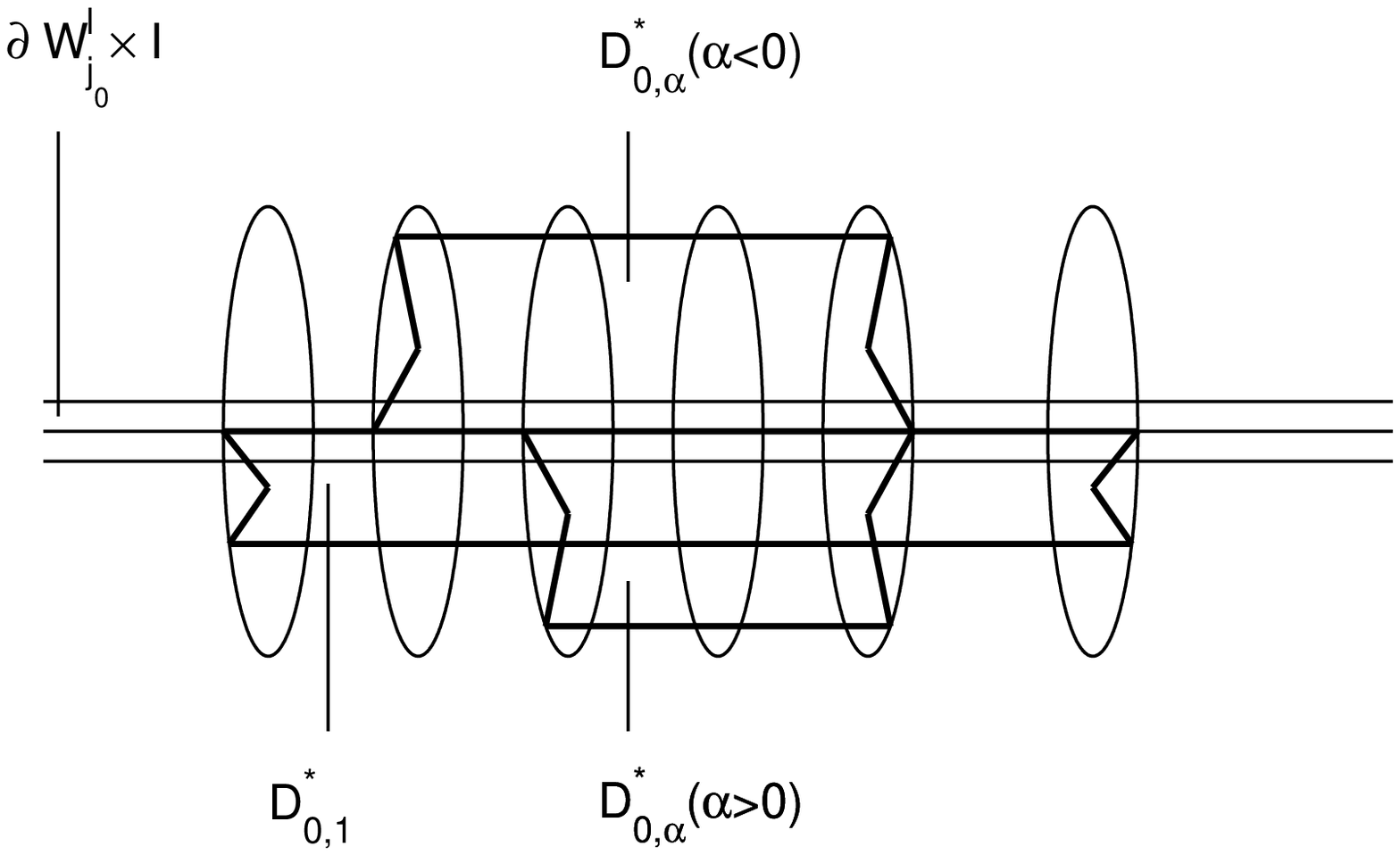}
\begin{center}
Figure 14
\end{center}
\end{center}

\begin{center}
\includegraphics[totalheight=4.5cm]{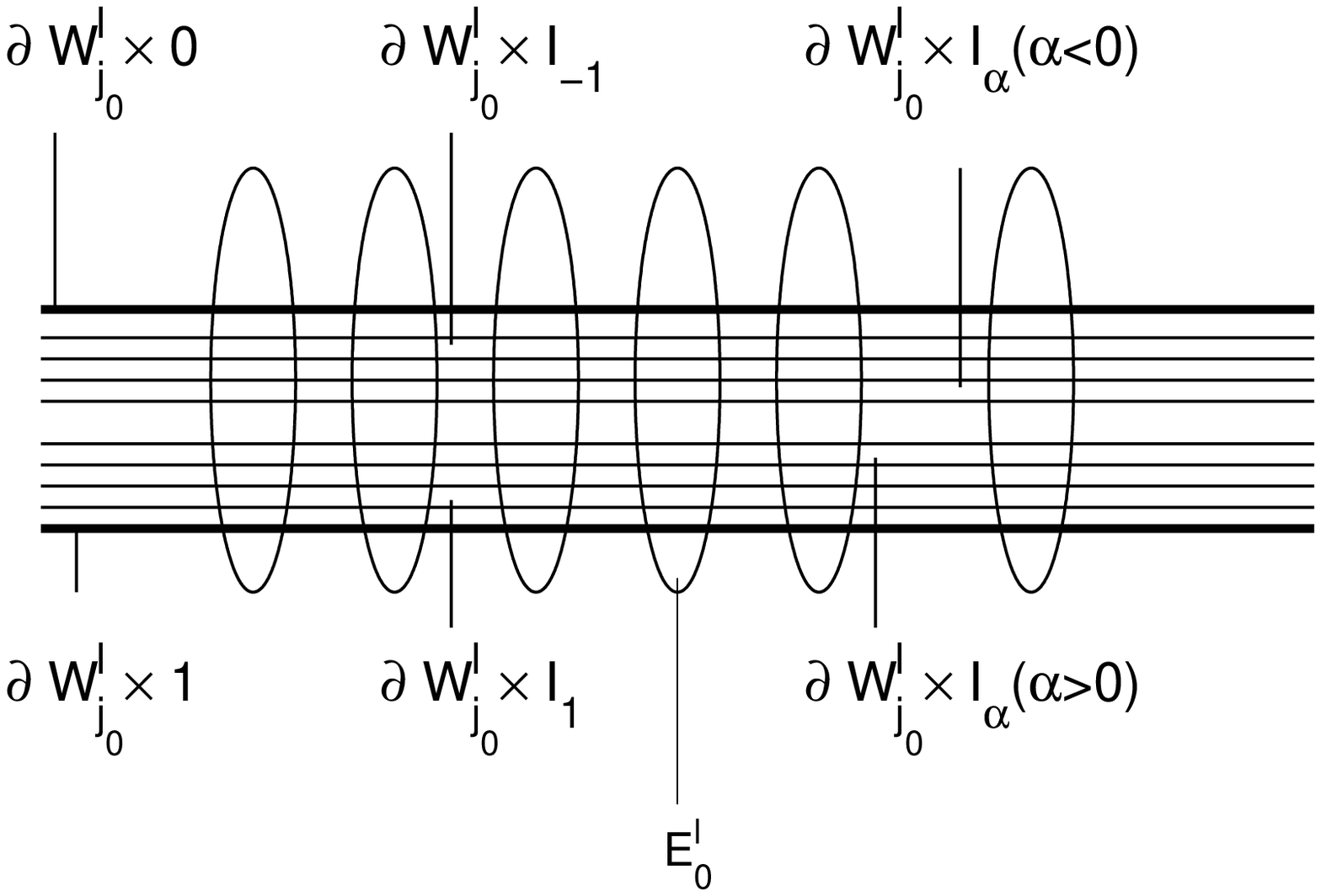}
\begin{center}
Figure 15
\end{center}
\end{center}

Now let $I_{\alpha}=[t_{1}^{\alpha},t_{2}^{\alpha}]$ be a
sub-interval of $[0,1]$ such that

(1) \ If $0<\alpha<\beta$, then
$t_{2}^{\alpha}>t_{1}^{\alpha}>t_{2}^{\beta}>t_{1}^{\beta}>1/2$.

(2) \ If $\alpha<\beta<0$, then
$1/2>t_{2}^{\alpha}>t_{1}^{\alpha}>t_{2}^{\beta}>t_{1}^{\beta}$.

Now $\partial W_{j_{0}}^{l}\times I_{\alpha}$ is as in Figure 15.

Let $t_{\alpha}$ be the center point of $I_{\alpha}$, and
$a_{0,\alpha}^{*}$ be an arc connecting $\partial_{1} a_{\alpha}$
to $\partial_{1} a_{\alpha}^{0}\times t_{\alpha}$ in
$E^{l}_{f_{\delta(0,\alpha)},\alpha}$ such that

(3) $a_{0,\alpha}-\partial W_{j_{0}}^{l}\times
I=a_{0,\alpha}^{*}-\partial W_{j_{0}}^{l}\times I$, where
$a_{0,\alpha}$ is as in Lemma 3.5.5 and
$a_{0,\alpha}^{*}\cap\partial W_{j_{0}}^{l}\times I$ is an arc in
$d_{\delta(0,\alpha),0}\times I\subset \partial
W_{j_{0}}^{l}\times I$ as in Figure 16. By Definition 3.5.1,
$d_{\delta(0,\alpha),0}\subset a^{0}_{\alpha}$.

Similarly, $b_{0,\alpha}^{*}$ be an arc connecting $\partial_{2}
a_{\alpha}$ to $\partial_{2} a_{\alpha}^{0}\times t_{\alpha}$ in
$E^{l}_{f_{\theta(0,\alpha)},\alpha}$ such that

(4) $b_{0,\alpha}-\partial W_{j_{0}}^{l}\times
I=b_{0,\alpha}^{*}-\partial W_{j_{0}}^{l}\times I$ and
$b_{0,\alpha}^{*}\cap\partial W_{j_{0}}^{l}\times I$ is an arc in
$d_{\theta(0,\alpha),0}\times I\subset \partial
W_{j_{0}}^{l}\times I$.

\begin{center}
\includegraphics[totalheight=4cm]{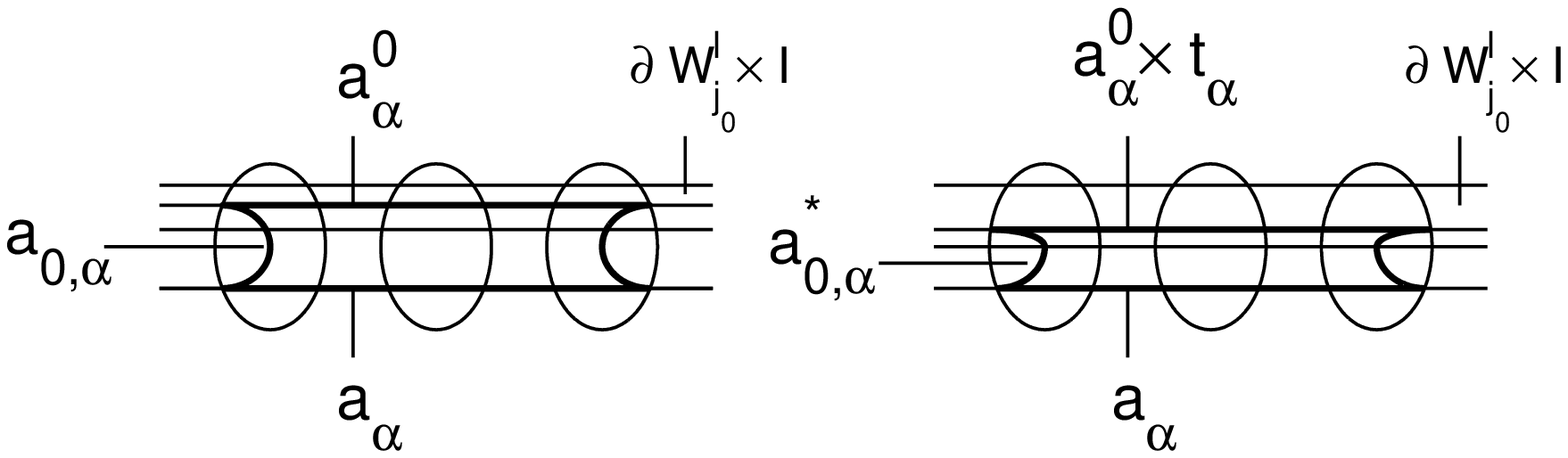}
\begin{center}
Figure 16
\end{center}
\end{center}

\begin{center}
\includegraphics[totalheight=6cm]{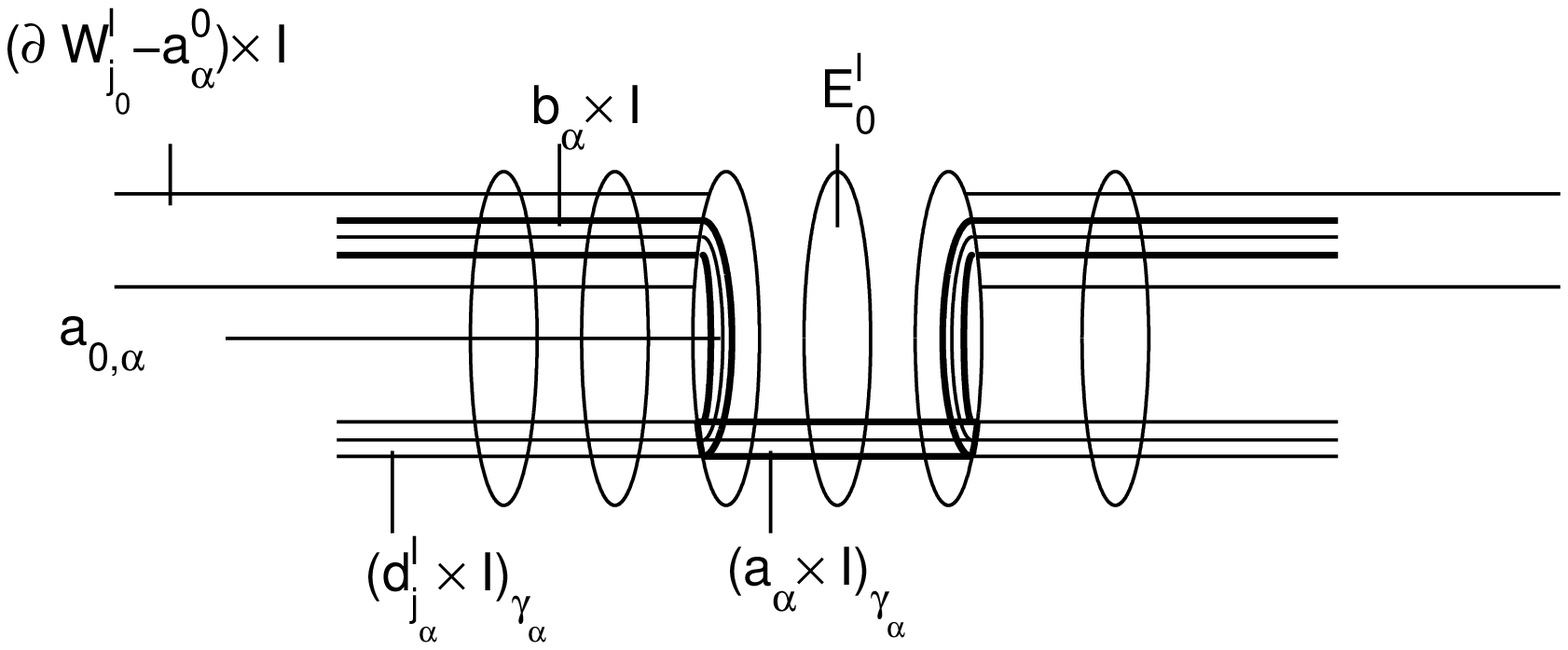}
\begin{center}
Figure 17
\end{center}
\end{center}

By Definition 3.5.1, $d_{\theta(0,\alpha),0}\subset
a^{0}_{\alpha}$. For $\alpha\neq 0$, let $b_{\alpha}=(\partial
W_{j_{0}}-a_{\alpha}^{0})\times \bigl\{t_{\alpha}\bigr\}\cup
a_{0,\alpha}^{*}\cup b_{0,\alpha}^{*}$ as in Figure 17. Then
$b_{\alpha}$ is an arc in $F^{l+1}$ by Lemma 3.5.8. Now let
$b_{\alpha}\times I$ be a  neighborhood of $b_{\alpha}$ in
$F^{l+1}$ satisfying the following conditions:

(5) \ $b_{\alpha}\times I\cap \partial W_{j_{0}}^{l}\times
I=(\partial W_{j_{0}}-inta_{\alpha}^{0})\times I_{\alpha}$.

(6) \ If $s(w_{\gamma_{\alpha}})=-$, then $(\partial
b_{\alpha})\times I=((\partial a_{\alpha})\times
I)_{\gamma_{\alpha}}$; if $s(w_{\gamma_{\alpha}})=+$ or
$\gamma_{\alpha}=\emptyset$, then $\partial b_{\alpha}=((\partial
a_{\alpha})\times I)_{\gamma_{\alpha}}$.

If $j_{\alpha}\in I(w_{\gamma_{\alpha}},l)$ with
$s(w_{\gamma_{\alpha}})=-$, then there is a homeomorphism
$H_{\alpha}$ from $(a_{\alpha}\times I)_{\gamma_{\alpha}}$ to
$b_{\alpha}\times I$ such that $H_{\alpha}$ is an identifying map
on $(\partial b_{\alpha})\times I=(\partial a_{\alpha})\times I$.
If $s(w_{\gamma_{\alpha}})=+$ or $\gamma_{\alpha}=\emptyset$, then
there is a homeomorphism $H_{\alpha}$ from $(a_{\alpha}\times
I)_{\gamma_{\alpha}}$ to $b_{\alpha}=b_{\alpha}\times
\bigl\{0\bigr\}$. In this case, $(a_{\alpha}\times
I)_{\gamma_{\alpha}}=a_{\alpha}$. (See Definition 3.2.1.)\qquad
Q.E.D.(Construction(*))

{\bf \em Lemma 4.3.5.} \ (1) If $\alpha\neq \beta$, then
$b_{\alpha}\times I\cap b_{\beta}\times I=\emptyset$.

(2) \ $H=\cup H_{\alpha}$ is an injective map from $\cup_{\alpha}
(a_{\alpha}\times I)_{\gamma_{\alpha}}$ to $\cup_{\alpha}
b_{\alpha}\times I$.

(3) \ If $j\notin m(l+1)$, and $d_{j}^{l}\cap
(a_{0,\alpha}^{*}\times I\cup b_{0,\alpha}^{*}\times
I)\neq\emptyset$, then $d_{j}^{l}\cap (a_{0,\alpha}^{*}\times
I\cup b_{0,\alpha}^{*}\times I)\subset\cup_{\beta\neq 0}
(a_{\beta}\times I)_{\gamma_{\beta}}$. Furthermore, either
$j=j_{\lambda}$ or $j>\gamma_{\lambda}$ for some $\lambda$.

(4) \ $e^{l}_{\gamma}\times I\cap (a_{0,\alpha}^{*}\times I\cup
b_{0,\alpha}^{*}\times I)\subset\cup_{\beta\neq 0}
(a_{\beta}\times I)_{\gamma_{\beta}}$. If $e^{l}_{\gamma}\times
I\cap (inta_{0,\alpha}^{*}\times I\cup intb_{0,\alpha}^{*}\times
I)\neq\emptyset$, then $e^{l}_{\gamma}\times I\cap
(inta_{0,\alpha}^{*}\times I\cup intb_{0,\alpha}^{*}\times
I)\subset\cup_{\beta\neq 0} (a_{\beta}\times I)_{\gamma_{\beta}}$,
and $\gamma\geq\gamma_{\lambda}$ for some $\lambda$.

(5) \ $b_{\alpha}\times I\subset F^{l+1}$.

{\bf \em Proof.} \ (1) \ By Construction(*) (5) and (6), we only
need to prove that $(a_{0,\alpha}^{*}\cup b_{0,\alpha}^{*})\cap
(a_{0,\beta}^{*}\cup b_{0,\beta}^{*})=\emptyset$. Without loss of
generality, we assume that $0<\alpha<\beta$. By Lemma 3.4.7,
$\delta(0,\alpha)\leq\delta(0,\beta)\leq\theta(0,\beta)\leq\theta(0,\alpha)$.

Suppose that $\delta(0,\alpha)<\delta(0,\beta)\leq
\theta(0,\beta)$. Since $f_{i,\alpha}=f_{i,0}$ for
$\delta(0,\alpha)\leq i\leq \theta(0,\alpha)$
$f_{i,\beta}=f_{i,0}$ for $\delta(0,\beta)\leq i\leq
\theta(0,\beta)$. By Lemma 3.3.4, $f_{\delta(0,\alpha),\alpha}\neq
f_{\delta(0,\beta),\beta},f_{\theta(0,\beta),\beta}$.
$a_{0,\alpha}^{*}\cap (a_{0,\beta}^{*}\cup
b_{0,\beta}^{*})=\emptyset$.

Suppose that $\delta(0,\alpha)=\delta(0,\beta)\leq
\theta(0,\beta)$. Since $a_{\alpha}$ separates
$a^{0}_{\alpha}=a^{0}_{\beta}$ and $a_{\beta}$ in
$D_{0,\beta}^{*}$. Now by the choice of $I_{\alpha}$, (1) holds.

(2) \ By Lemma 3.5.4, $(a_{\alpha}\times I)_{\gamma_{\alpha}}\cap
(a_{\beta}\times I)_{\gamma_{\beta}}=\emptyset$ for
$\alpha\neq\beta$. By (1), $b_{\alpha}\times I\cap b_{\beta}\times
I=\emptyset$. Hence (2) holds.

(3, 4) \ By Construction(*) (3) and (4), $a_{0,\alpha}-\partial
W_{j_{0}}^{l}\times I=a^{*}_{0,\alpha}-\partial
W_{j_{0}}^{l}\times I$ and $a_{0,\alpha}^{*}\cap
\partial W_{j_{0}}^{l}\times I$ is an arc  in
$d_{\delta(0,\alpha)}\times I$.  By Lemma 4.1.5(2), (4) and (6),
$d_{j}^{l}$ and $e^{l}_{\gamma}\times I$ are disjoint from
$intd_{i,0}\times I$ for $j\neq j_{0}=m^{l+1}$ and $\gamma\in
m(l)$. Now by Lemma 3.5.5, (3) and (4) holds.

(5) \ Since $d_{0,0}\subset a_{\alpha}^{0}$, $(\partial
W_{j_{0}}^{l}-inta_{\alpha}^{0})\times I_{\alpha}\subset F^{l+1}$.
By Lemma 3.5.8, $a_{0,\alpha}^{*}\cup b_{0,\alpha}^{*}$ is
disjoint from $c_{l+1}^{l}\times I$. Hence (5) holds. \qquad
Q.E.D.\vskip 3mm

{\bf \em Definition 4.3.6.} \ Let
$d_{j}^{l+1}=(d_{j}^{l}-\cup_{\alpha\neq 0} (a_{\alpha}\times
I)_{\gamma_{\alpha}})\cup H(d_{j}^{l}\cap(\cup_{\alpha\neq 0}
(a_{\alpha}\times I)_{\gamma_{\alpha}}))$  for $j\notin m(l+1)$.
\vskip 0.5mm

{\bf \em Lemma 4.3.7.} \ $\bigl\{c_{i}^{l+1} \ | \ i\geq
l+2\bigr\}$ is a set of pairwise disjoint arcs properly embedded
in $F^{l+1}$ which lies in one of $E^{l+1}_{f}$.

{\bf \em Proof.} \ Since $i\geq l+2$, $c_{i}^{l+1}=c_{i}^{l}$ is
disjoint from $c_{l+1}^{l}\times I$. Hence $c_{i}^{l+1}$ lies in
$E^{l+1}_{f}$ for some $f$ by Proposition 4(4). Furthermore, if
$c_{i}^{l}\subset E^{l}_{0}$, then $c_{i}^{l+1}$ lies in one of
$E^{l+1}_{0}$ and $E^{l+1}_{l+1}$. By Proposition 4(4) and Lemma
4.1.5(5), $c_{i}^{l+1}$ is properly embedded in $F^{l+1}$. \qquad
Q.E.D.\vskip 0.5mm

{\bf \em Lemma 4.3.8.} \ (1) \ $\bigl\{d_{j}^{l+1} \ | \ j\notin
m(l+1)\bigr\}$ is a set of pairwise disjoint arcs in $F^{l+1}$
such that $\partial d_{j}^{l+1}=\partial d_{j}^{l}$.

(2) \ If $j\in \bigl\{1,\ldots,n\bigr\}-m(l+1)$ and $j\neq
j_{\alpha}$ for each $\alpha$, then
$d_{j}^{l+1}=(d_{j}^{l}-\cup_{\gamma_{\alpha}<j} (a_{\alpha}\times
I)_{\gamma_{\alpha}})\cup_{\gamma_{\alpha}<j}
H_{\alpha}(d_{j}^{l}\cap (a_{\alpha}\times I)_{\gamma_{\alpha}})$.

(3) \ If $j=j_{\beta}$ for some $\beta\neq 0$, then
$d_{j}^{l+1}=(d_{j}^{l}-a_{\beta}\cup_{\gamma_{\alpha}<j}
(a_{\alpha}\times I)_{\gamma_{\alpha}})\cup_{\gamma_{\alpha}<j}
H_{\alpha}(d_{j}^{l}\cap(a_{\alpha}\times
I)_{\gamma_{\alpha}})\cup b_{\beta}$.

(4) \ If $j<j_{0}$, then $d_{j}^{l}=d_{j}^{l+1}$.

{\bf \em  Proof.} \ (1) \ By Proposition 4(2) and (3), each
component $d_{j}^{l}\cap e^{l}_{\gamma_{\alpha}}\times I$ is a
core of $e^{l}_{\gamma_{\alpha}}\times (0,1)$ if $j\notin
I(w_{\gamma_{\alpha}},l)$. Hence each component of $d_{j}^{l}\cap
(a_{\alpha}\times I)_{\gamma_{\alpha}}$ is a core of
$(a_{\alpha}\times I)_{\gamma_{\alpha}}$ even if $j=j_{\alpha}$.
By Lemma 4.3.5, $d_{j}^{l}\cap (\cup_{\alpha\neq 0}
a_{0,\alpha}^{*}\cup b_{0,\alpha}^{*})\subset \cup_{\alpha\neq 0}
(a_{\alpha}\times I)_{\gamma_{\alpha}}$ and $b_{\alpha}\times
I\cap b_{\beta}\times I=\emptyset$. By Lemma 3.3.3 and Lemma
4.1.4, $(\cup_{j\notin m(l+1)} d_{j}^{l})\cap c_{l+1}^{l}\times
I=(\cup_{\alpha\neq 0} (d_{0,\alpha}\times
I)_{\gamma_{\alpha}})\cap c_{l+1}^{l}\times I$. By Definition
3.5.1, $d_{0,\alpha}\subset a_{\alpha}$. By Lemma 4.3.5(5),
$d_{j}^{l+1}$ is an arc in $F^{l+1}$. By Proposition 4,
$d_{j}^{l}\cap d_{r}^{l}=\emptyset$ for $j\neq r$. Since $H$ is
injective, $d_{j}^{l+1}\cap d_{r}^{l+1}=\emptyset$.

(2) \ If $j\neq j_{\alpha}$ for each $\alpha\neq 0$ and
$j<\gamma_{\alpha}$, then, by Proposition 4(3), $d_{j}^{l}$ is
disjoint from $(a_{\alpha}\times I)_{\gamma_{\alpha}}\subset
e^{l}_{\gamma_{\alpha}}\times I$. Hence (2) holds.

(3) \ Suppose that $j=j_{\beta}$ for some $\beta\neq 0$. Then
$a_{\beta}=a_{\beta}\times \bigl\{0\bigr\}\subset d_{j}^{l}$.
Furthermore, $H_{\beta}(a_{\beta})=b_{\beta}$. By the argument in
(2), (3) holds.

(4)  If $j<j_{0}$, then $j<j_{\alpha}$ for each $\alpha\neq 0$. By
Proposition 4(3), $d_{j}^{l}\cap (a_{\alpha}\times
I)_{\gamma_{\alpha}}=\emptyset$ for each $\alpha$. By (2), (4)
holds.  \qquad Q.E.D.\vskip 0.5mm

In fact, if $j>j_{0}$, then $d_{j}^{l+1}$ is obtained by doing
band sums with copies of $\partial W_{j_{0}}^{l}$ to $d_{j}^{l}$.
See Lemma 4.6.3.\vskip 0.5mm

{\bf \em Definition 4.3.9.} \ (1) \ Let $e^{l+1}_{\gamma}\times
I=e^{l}_{\gamma}\times I$ for $\gamma\in m(l)$ with
$s(w_{\gamma})=+$.

(2) \ Suppose that $r\in I(w_{\gamma},l)$ for some $\gamma\in
m(l)$ with $s(w_{\gamma})=-$. If $\gamma>j_{0}$, then let
$(d_{r}^{l+1}\times I)_{\gamma}= ((d_{r}^{l}\times
I)_{\gamma}-\cup_{\alpha} (a_{\alpha}\times
I)_{\gamma_{\alpha}})\cup H((d_{r}^{l}\times
I)_{\gamma}\cap(\cup_{\alpha} (a_{\alpha}\times
I)_{\gamma_{\alpha}}))$. If $\gamma\leq j_{0}$, then let
$(d_{r}^{l+1}\times I)_{\gamma}=(d_{r}^{l}\times I)_{\gamma}$.

(3) Suppose that $\gamma\in m(l)$ with $s(w_{\gamma})=-$. Let
$e^{l+1}_{\gamma}=(e^{l}_{\gamma}-\cup_{r\in I(w_{\gamma},l)}
d_{r}^{l})\cup_{r\in I(w_{\gamma},l)} d_{r}^{l+1}$,
$e^{l+1}_{\gamma}\times I=(e^{l}_{\gamma}\times I-\cup_{r\in
I(w_{\gamma},l)} (d_{r}^{l}\times I)_{\gamma})\cup_{r\in
I(w_{\gamma},l+1)} (d_{r}^{l+1}\times I)_{\gamma}$. Specially, let
$e^{l+1}_{m^{l+1}}=\partial W_{j_{0}}^{l}-intd_{j_{0}}^{l}$,
$e^{l+1}_{m^{l+1}}\times I=(\partial
W_{j_{0}}^{l}-intd_{j_{0}}^{l})\times I$. \vskip 0.5mm

{\bf \em Lemma 4.3.10.} \ (1) \ For each $\gamma\in m(l+1)$,
$e^{l+1}_{\gamma}\times I$ is a disk in $F^{l+1}$ such that
$(\partial e^{l}_{\gamma})\times I=(\partial
e^{l+1}_{\gamma})\times I$ for $\gamma\in m(l)$ and $(\partial
e^{l+1}_{m^{l+1}})\times I=(\partial d_{j_{0}}^{l})\times I$.

(2) \ Suppose that $r\in I(w_{\gamma},l)$ for some
$j_{0}<\gamma\in m(l)$ with $s(w_{\gamma})=-$. If $r\neq
j_{\alpha}$ for each $\alpha$, then $(d_{r}^{l+1}\times
I)_{\gamma}= ((d_{r}^{l}\times
I)_{\gamma}-\cup_{\gamma_{\alpha}<r} (a_{\alpha}\times
I)_{\gamma_{\alpha}})\cup_{\gamma_{\alpha}<r}
H_{\alpha}((d_{r}^{l}\times I)_{\gamma}\cap(a_{\alpha}\times
I)_{\gamma_{\alpha}})$.

(3) \ Suppose that $s(w_{\gamma_{\beta}})=-$. Then
$(d_{j_{\beta}}^{l+1}\times
I)_{\gamma_{\beta}}=((d_{j_{\beta}}^{l}\times
I)_{\gamma_{\beta}}-(a_{\beta}\times
I)_{\gamma_{\beta}}\cup_{\gamma_{\alpha}<j_{\beta}}
(a_{\alpha}\times
I)_{\gamma_{\alpha}})\cup_{\gamma_{\alpha}<j_{\beta}}
H_{\alpha}((d_{j_{\beta}}^{l}\times I)_{\gamma_{\beta}}\cap
(a_{\alpha}\times I)_{\gamma_{\alpha}})\cup b_{\beta}\times I$.

{\bf \em Proof.} \ \ Suppose that $s(w_{\gamma})=+$. By
Proposition 4(1) and (5), $e^{l}_{\gamma}\times I$ is disjoint
from $c_{l+1}^{l}\times I\subset F^{l}$. Hence
$e^{l+1}_{\gamma}\times I=e^{l}_{\gamma}\times I\subset F^{l+1}$
is a disk.

Suppose that  $s(w_{\gamma})=-$. There are four cases:

Case 1. \ $\gamma<j_{0}$.

Now if $r\in I(w_{\gamma})$, then $r<\gamma<j_{0}$. By Lemma 4.1.2
and Lemma 4.3.7, $(d_{r}^{l+1}\times I)_{\gamma}=(d_{r}^{l}\times
I)_{\gamma}\subset e^{l}_{\gamma}\times I$ is disjoint from
$c_{l+1}^{l}\times I$. By Definition 4.3.9,
$e^{l+1}_{\gamma}\times I=e^{l}_{\gamma}\times I\subset F^{l+1}$
is a disk.

Case 2. \ By the construction, $e^{l+1}_{m^{l+1}}\times I$ is a
disk in $F^{l+1}$.

Case 3. \ $\gamma>j_{0}$ and $\gamma\neq \gamma_{\alpha}$ for each
$\alpha$.

Now $\gamma\geq Max\bigl\{ \gamma_{i,0} \ | \ \delta(0)\leq i\leq
\theta(0)\bigr\}$. Since $\delta(0,\alpha)\leq i\leq
\theta(0,\alpha)$, $\gamma_{i,\alpha}=\gamma_{i,0}$. By Lemma
3.5.3,
 if $e^{l}_{\gamma}\times
I\cap(inta_{\alpha}\times I)_{\gamma_{\alpha}}\neq\emptyset$,
then, $\gamma>\gamma_{\alpha}$.  By Proposition 4(1),
$e^{l}_{\gamma}\times I\cap e^{l}_{\gamma_{\alpha}}\times
I=(\cup_{r\in I(w_{\gamma},l)} (d_{r}^{l}\times I)_{\gamma})\cap
e^{l}_{\gamma_{\alpha}}\times I$.  By Definition 2.1.4, each
component of $e^{l}_{\gamma}\times I\cap(a_{\alpha}\times
I)_{\gamma_{\alpha}}$ is $(c\times I)_{\gamma}\subset
(a_{\alpha}\times (0,1))_{\gamma_{\alpha}}$ where $c$ is a core of
$(a_{\alpha}\times I)_{\gamma_{\alpha}}$. By Lemma 3.2.6, if $r\in
I(w_{\gamma},l)$ and $r<\gamma_{\alpha}$, then $(d_{r}^{l}\times
I)_{\gamma}$ is disjoint from $(a_{\alpha}\times
I)_{\gamma_{\alpha}}$. Hence (2) holds.  By Lemma 4.1.5(2),
$e^{l}_{\gamma}\times I$ is disjoint from $W_{j_{0}}^{l}\times
I_{\alpha}\subset W_{j_{0}}^{l}\times I$. By Lemma 4.3.5, if
$e^{l}_{\gamma}\times I\cap((inta_{0,\alpha}^{*}\cup
intb_{0,\alpha}^{*})\times I)_{\gamma_{\alpha}}\neq\emptyset$,
then $e^{l}_{\gamma}\times I\cap((inta_{0,\alpha}^{*}\cup
intb_{0,\alpha}^{*})\times I)_{\gamma_{\alpha}}\subset
\cup_{\beta\neq 0}((inta_{\beta}\cup intb_{\beta})\times
I)_{\gamma_{\beta}}$. Since $H$ is injective,
$e^{l+1}_{\gamma}\times I$ is a disk in $F^{l+1}$.

Case 4. $\gamma=\gamma_{\beta}$ for some $\beta$.

Now $(a_{\beta}\times
I)_{\gamma_{\beta}}\subset(d_{j_{\beta}}^{l}\times
I)_{\gamma_{\beta}} \subset e^{l}_{\gamma_{\beta}}\times I$. By
Definition 4.3.9 and the argument in Case 3,
$e^{l+1}_{\gamma}\times I$ is a disk in $F^{l+1}$ and (3) holds.
\qquad Q.E.D.

%(2) \ Since $\gamma>j_{0}$, by Lemma 3.5.3, if
%$e^{l}_{\gamma}\times I\cap(inta_{\alpha}\times
%I)_{\gamma_{\alpha}}\neq\emptyset$, then,
%$\gamma\geq\gamma_{\alpha}$. Note that $a_{\alpha}\subset
%d_{j_{\alpha}}^{l}$. Since $r\in I(w_{\gamma},l)$ and $r\neq
%j_{\alpha}$ for each $\alpha$, by Lemma 3.2.6 and Definition
%4.3.9(2), (2) holds.

%(3) \ Since $a_{\beta}\subset d_{j_{\beta}}^{l}$, $b_{\beta}\times
%I=H_{\beta}((a_{\beta}\times I)_{\gamma_{\beta}})\subset
%(d_{j_{\beta}}^{l+1}\times I)_{\gamma_{\beta}}$. By the same
%argument in (2), (3) holds.\qquad Q.E.D.

\subsection{$F^{l+1}$ is a surface generated by an abstract tree
in $\partial{+} V_{-}$}\vskip 0.5mm

\ \ \ \ \ By the construction, $F^{l+1}$ is a surface in
$\partial{+} V_{-}$. In this section, we shall prove that
$\cup_{f}E^{l+1}_{f}\cup_{\gamma\in m(l+1)} e^{l+1}_{\gamma}$ is
an abstract tree and $F^{l+1}$ is a surface generated by
$\cup_{f}E^{l+1}_{f}\cup_{\gamma\in m(l+1)} e^{l+1}_{\gamma}$,
which satisfies Proposition 4(1).\vskip 0.5mm

{\bf \em Lemma 4.4.1.} \ (1) \ $F^{l+1}=\cup_{f=0}^{l+1}
E^{l+1}_{f}\cup_{\gamma\in m(l+1)} e^{l+1}_{\gamma}\times I$.

(2) \ For $\gamma\neq\lambda\in m(l+1)$, $(\partial
e^{l+1}_{\gamma})\times I\cap (\partial e^{l+1}_{\lambda})\times I
=\emptyset$.

{\bf \em Proof.} \ By the construction, $F^{l+1}=(F^{l}-c\times
(-1,1))\cup (\partial W_{j_{0}}^{l}-inta)\times I$, where
$a\subset d_{0,0}$. Since $d_{j_{0}}^{l}\subset F^{l}$,
$F^{l+1}=(F^{l}-c\times (-1,1))\cup e^{l+1}_{m^{l+1}}\times I$.
Hence (1) follows from the construction.

Since $F^{l}$ is generated by $\cup_{f}E^{l}_{f}\cup_{\gamma\in
m(l)}e^{l}_{\gamma}$. Hence, by Lemma 4.3.10,  $(\partial
e^{l+1}_{\gamma})\times I\cap (\partial e^{l+1}_{\lambda})\times I
=\emptyset$ for $\gamma\neq\lambda\in m(l)$. By Lemma 4.3.10 and
Lemma 4.1.5(1), $(\partial e^{l+1}_{\gamma})\times I\cap (\partial
e^{l+1}_{m^{l+1}})\times I =\emptyset$ for $\gamma\in m(l)$.
\qquad Q.E.D.\vskip 0.5mm

{\bf Lemma 4.4.2.} \ If $\lambda\in m(l+1)$ with
$s(w_{\lambda})=+$, then  $inte^{l+1}_{\gamma}\times I$ is
disjoint from $\cup_{f}
E^{l+1}_{f}\cup_{\gamma<\lambda}e^{l+1}_{\gamma}\times I$.

{\bf \em Proof.} \ Since $s(w_{\lambda})=+$, by Lemma 4.1.3,
$\lambda\neq j_{0}=m^{l+1}$. By Proposition 4(1) and (5),
$e^{l}_{\gamma}\times I$ is disjoint from $c_{l+1}^{l}\times
I\subset E^{l}_{0}$. By Definition 4.3.9, $e^{l+1}_{\gamma}\times
I=e^{l}_{\gamma}\times I$.

Suppose that $\lambda<j_{0}$. By Definition 4.3.9,
$e^{l+1}_{\gamma}\times I=e^{l}_{\gamma}\times I$ for
$\gamma<\lambda$. Note that $E^{l+1}_{0}, E^{l+1}_{l+1}\subset
E^{l}_{0}$ and $E^{l+1}_{f}=E^{l}_{f}$ for $1\leq f\leq l$.  In
this case, $\cup_{f}
E^{l+1}_{f}\cup_{\gamma<\lambda}e^{l+1}_{\gamma}\times I =\cup_{f}
E^{l}_{f}\cup_{\gamma<\lambda}e^{l}_{\gamma}\times
I-c_{l+1}^{l}\times (-1,1)$. Hence $inte^{l+1}_{\lambda}\times I$
is disjoint from $\cup_{f}
E^{l+1}_{f}\cup_{\gamma<\lambda}e^{l+1}_{\gamma}\times I$.

Suppose that $\lambda>j_{0}$. By Proposition 4(1),
$inte^{l}_{\lambda}\times I$ is disjoint from
$\cup_{f}E^{l}_{f}\cup_{\gamma<\lambda} e^{l}_{\gamma}\times I$.
By Lemma 4.1.5(1) and (2), $e^{l}_{\lambda}\times I$ is disjoint
from $e^{l+1}_{m^{l+1}}\times I\subset \partial
W_{j_{0}}^{l}\times I$. By Definition 4.3.9, $\cup_{f}
E^{l+1}_{f}\cup_{\gamma<\lambda}e^{l+1}_{\gamma}\times I=(\cup_{f}
E^{l}_{f}\cup_{\gamma<\lambda}e^{l}_{\gamma}\times
I-c_{l+1}^{l}\times (-1,1))\cup e^{l+1}_{m^{l+1}}\times I$. Hence
$inte^{l+1}_{\lambda}\times I=inte^{l}_{\lambda}\times I$ is
disjoint from $\cup_{f}
E^{l+1}_{f}\cup_{\gamma<\lambda}e^{l+1}_{\gamma}\times I$. \qquad
Q.E.D.\vskip 0.5mm

{\bf Lemma 4.4.3.} \ Suppose that $\lambda\in m(l+1)$ with
$s(w_{\lambda})=-$. Then   $inte^{l+1}_{\lambda}\times
I\cap(\cup_{f}
E^{l+1}_{f}\cup_{\gamma<\lambda}e^{l+1}_{\gamma}\times
I)=\cup_{r\in I(w_{\lambda},l+1)} d_{r}^{l+1}\times I$.

{\bf \em Proof.} \  By assumption, Proposition 4 holds for $k=l$.
Hence $inte^{l}_{\lambda}\times I\cap(\cup_{f}
E^{l}_{f}\cup_{\gamma<\lambda}e^{l}_{\gamma}\times I)=\cup_{r\in
I(w_{\lambda},l)} d_{r}^{l}\times I$ for $\lambda\in m(l)$ with
$s(w_{\lambda})=-$.

Since $s(w_{\lambda})=-$, by Lemma 2.2.5 and Lemma 4.1.3,
$j_{0}=m^{l+1}\notin I(w_{\lambda},l)$. By Definition 2.3.1,
$I(w_{\lambda},l+1)=I(w_{\lambda},l)-\bigl\{m^{l+1}\bigr\}=I(w_{\lambda},l)$.
Now there are four cases:

Case 1. \ $\lambda<j_{0}=m^{l+1}$.

In this case, $\bigl\{\gamma \ | \ \lambda>\gamma\in
m(l+1)\bigr\}=\bigl\{\gamma \ | \ \lambda>\gamma\in m(l)\bigr\}$.

By Definition 4.3.9, $e^{l+1}_{\lambda}\times
I=e^{l}_{\lambda}\times I$ and $e^{l+1}_{\gamma}\times
I=e^{l}_{\gamma}\times I$ for $\gamma\leq\lambda$. By Lemma 4.3.8
and Lemma 2.2.4, $d_{r}^{l+1}=d_{r}^{l}$ for $r\in
I(w_{\lambda},l)$ or $r\in I(w_{\gamma},l)$.  Hence
$inte^{l+1}_{\lambda}\times I\cap (\cup_{f}
E^{l+1}_{f}\cup_{\gamma<\lambda} e^{l+1}_{\gamma}\times
I)=inte^{l}_{\lambda}\times I\cap (\cup_{f}
E^{l}_{f}\cup_{\gamma<\lambda} e^{l}_{\gamma}\times
I-c_{l+1}^{l}\times (-1,1))$. By Lemma 4.1.2,
$e^{l}_{\gamma}\times I$ is disjoint from $c_{l+1}^{l}\times I$.
By Definition 4.3.9, the lemma holds.

Case 2. \ $\lambda=j_{0}=m^{l+1}$.

By Definition 4.3.9, $e^{l+1}_{\gamma}\times
I=e^{l}_{\gamma}\times I$ for $\gamma<\lambda$. and
$e^{l+1}_{m^{l+1}}\times I=(\partial
W_{j_{0}}-intd_{j_{0}}^{l})\times I$. If $r\in I(w_{j_{0}},l)$,
then $r<j_{0}$. By Lemma 4.3.8, $d_{r}^{l+1}=d_{r}^{l}$ is
disjoint from $c_{l+1}^{l}\times I$.  By Lemma 4.1.5(1) and (2),
$d_{r}^{l+1}\times I\subset \cup_{f}
E^{l+1}_{f}\cup_{\gamma<j_{0}} e^{l+1}_{\gamma}\times
I=\cup_{f}E^{l}_{f}\cup_{\gamma<j_{0}} e^{l}_{\gamma}\times
I-c_{l+1}^{l}\times (-1,1)$.  Hence the lemma holds.

Case 3. \ $\lambda>j_{0}=m^{l+1}$ and $\lambda\neq
\gamma_{\alpha}$ for each $\alpha$.

In this case, $\bigl\{\gamma \ | \ \gamma<\lambda\in
m(l+1)\bigr\}=\bigl\{\gamma \ | \ \gamma<\lambda\in
m(l)\bigr\}\cup\bigl\{m^{l+1}\bigr\}$.

By Lemma 4.1.5(1) and (2), $e^{l}_{\lambda}\times I$ is disjoint
from $e^{l+1}_{m^{l+1}}\times I\subset \partial
W_{j_{0}}^{l}\times I$. By Proposition 4(1),
$inte^{l}_{\lambda}\times I-\cup_{r\in I(w_{\lambda},l)}
d_{r}^{l}\times I$ is disjoint from
$\cup_{f}E^{l}_{f}\cup_{\gamma<\lambda\in m(l)
}e^{l}_{\gamma}\times I$. By Definition 4.3.9, for each $\gamma\in
m(l)$, $e^{l+1}_{\gamma}\times I\subset e^{l}_{\gamma}\times I\cup
\partial W_{j_{0}}^{l}\times I$. Hence $inte^{l}_{\lambda}\times
I-\cup_{r\in I(w_{\lambda},l)} d_{r}^{l}\times I$ is disjoint from
$\cup_{f}E^{l+1}_{f}\cup_{\gamma<\lambda\in m(l+1)
}e^{l+1}_{\gamma}\times I$. Since $\lambda>j_{0}$, by Lemma
4.3.10(2), $(d_{r}^{l+1}\times I)_{\lambda}= ((d_{r}^{l}\times
I)_{\lambda}-\cup_{\gamma_{\alpha}<r} (a_{\alpha}\times
I)_{\gamma_{\alpha}})\cup_{\gamma_{\alpha}<r}
H_{\alpha}((d_{r}^{l}\times I)_{\lambda}\cap (a_{\alpha}\times
I)_{\gamma_{\alpha}})\subset\cup_{f}E^{l+1}_{f}\cup_{\gamma<\lambda\in
m(l+1) }e^{l+1}_{\gamma}\times I$.  Hence the lemma holds.

Case 4. $\lambda=\gamma_{\beta}$.

By Lemma 4.3.5 and Construction(*), $b_{\beta}\times I\subset
(\partial W_{j_{0}}^{l}-a_{\beta}^{0})\times I\cup
E^{l+1}_{f_{\delta(0,\beta),0}}\cup
E^{l+1}_{f_{\theta(0,\beta),0}}$. By Lemma 4.1.5 and Definition
4.3.9, $(\partial W_{j_{0}}^{l}-a_{\beta}^{0})\times I\subset
\cup_{f} E^{l+1}_{f}\cup_{\gamma<j_{0}} e^{l+1}_{\gamma}\times I$.
By the argument in Case 3 and Lemma 4.3.10(3), the lemma holds
 \qquad Q.E.D.\vskip 0.5mm

%\begin{center}
%\includegraphics[totalheight=5cm]{r1.eps}
%\begin{center}
%Figure 14
%\end{center}
%\end{center}

{\bf \em Lemma 4.4.4.} \ If $\lambda>\gamma\in m(l+1)$, then each
component   of $e^{l+1}_{\lambda}\cap e^{l+1}_{\gamma}\times I$ is
an arc $c\subset inte^{l+1}_{\lambda}$ which is a core of
$e^{l+1}_{\gamma}\times (0,1)$, and each component of
$e^{l+1}_{\lambda}\times I \cap e^{l+1}_{\gamma}\times I$ is
$(c\times I)_{\lambda}\subset e^{l+1}_{\gamma}\times (0,1)$.

{\bf \em Proof.} \ Recalling the assumption that Proposition 4
holds for $k=l$. There are five cases:

Case 1. $\gamma<\lambda<j_{0}$.

Now $e^{l+1}_{\gamma}\times I=e^{l}_{\gamma}\times I$ and
$e^{l+1}_{\lambda}\times I=e^{l}_{\lambda}\times I$. Since $F^{l}$
is a surface generated by the abstract tree $\cup_{f}
E^{l}_{f}\cup_{\gamma\in m(l)} e^{l}_{\gamma}$ , the lemma holds.
(See Definition 2.1.4.)

Case 2. $\gamma<\lambda=j_{0}$.

Now $e^{l+1}_{\gamma}\times I=e^{l}_{\gamma}\times I\subset
F^{l}-c_{l+1}^{l}\times (-1,1)$. By Lemma 4.4.3,
$e^{l+1}_{m^{l+1}}\times
I\cap(\cup_{f}E^{l+1}_{f}\cup_{\gamma<m^{l+1}}
e^{l+1}_{\gamma}\times I)= (\partial d_{j_{0}}^{l})\times
I\cup_{r\in I(w_{j_{0}},l+1)} d_{r}^{l+1}\times I$. By Lemma
4.3.8, $d_{r}^{l+1}=d_{r}^{l}$. By Lemma 4.1.5(1), (2),(3) and
Lemma 4.4.1(2), the lemma holds.

Case 3. $\gamma<j_{0}=m^{l+1}<\lambda$.

By Definition 4.3.9, $e^{l+1}_{\gamma}\times
I=e^{l}_{\gamma}\times I\subset F^{l}-c_{l+1}^{l}\times [-1,1]$.
By Lemma 4.1.5, each component of $\partial W_{j_{0}}^{l}\times
I\cap e^{l}_{\gamma}\times I$ is $c\times I\subset
e^{l}_{\gamma}\times (0,1)$ where $c$ is a core of
$e^{l}_{\gamma}\times (0,1)$. By  Lemma 4.1.5(6),
$(intd_{i,0}\times I)_{j_{0}}$ is disjoint from
$e^{l}_{\gamma}\times I=e^{l}_{\gamma}\times I$. Note that
$a_{\alpha}^{0}=\cup_{i=\delta(0,\alpha)}^{\theta(0,\alpha)}
d_{i,0}\cup_{i=\delta(0,\alpha)}^{\theta(0,\alpha)} e_{i,0}$.
Hence each component of $\partial W_{j_{0}}^{l}\times I\cap
e^{l}_{\gamma}\times I$ is either in $(\partial
W_{j_{0}}^{l}-inta_{\alpha}^{0})\times I$ or in
$inta_{\alpha}^{0}\times I$.

Since $\gamma<j_{0}<j_{\alpha}$, by Lemma 4.3.5(4),
$inta_{0,\alpha}^{*}\times I\cup intb_{0,\alpha}^{*}\times I$ is
disjoint from $e^{l}_{\gamma}\times I$. Since $j_{\alpha}>j_{0}$,
by Lemma 3.2.6, $(intd_{i,\alpha}\times I)_{\gamma_{\alpha}}$ is
also disjoint from $e^{l}_{\gamma}\times I$. Hence each component
of $e^{l}_{\gamma_{\alpha}}\times I\cap e^{l}_{\gamma}\times I$ is
either in $(e^{l}_{\gamma_{\alpha}}-inta_{\alpha})\times I$ or in
$(inta_{\alpha}\times I)_{\gamma_{\alpha}}$.

Since $\lambda>j_{0}$, $\lambda>Max\bigl\{\gamma_{i,0} \ | \
\delta(0)\leq i\leq \theta(0)\bigr\}$ by Lemma 3.3.1. Hence if
$e^{l}_{\lambda}\times I\cap (a_{\alpha}\times
I)_{\gamma_{\alpha}}\neq\emptyset$, then, by Lemma 3.5.3,
$\lambda\geq\gamma_{\alpha}$. By Definition 2.1.4 and Proposition
4(1), each component of $e^{l}_{\lambda}\times I\cap
(a_{\alpha}\times I)_{\gamma_{\alpha}}$ is $(c\times
I)_{\lambda}\subset (a_{\alpha}\times I)_{\gamma_{\alpha}}$ where
$c$ is a core of $(a_{\alpha}\times I)_{\gamma_{\alpha}}$. Hence
each component of $e^{l}_{\lambda}\times I\cap
e^{l}_{\gamma}\times I$ is either in $(a_{\alpha}\times
I)_{\gamma_{\alpha}}$ or in $e^{l}_{\lambda}\times
I-(inta_{\alpha}\times I)_{\gamma_{\alpha}}$. By Lemma 4.3.10(2),
(3) and Lemma 4.4.3, $e^{l+1}_{\lambda}\times I\cap
e^{l+1}_{\gamma}\times I=S_{1}\cup S_{2}$ where
$S_{1}=(e^{l}_{\lambda}\times
I-\cup_{\gamma_{\alpha}\leq\lambda}(inta_{\alpha}\times
I)_{\gamma_{\alpha}})\cap e^{l}_{\gamma}\times I$ and
$S_{2}=(\cup_{\gamma_{\alpha}\leq\lambda}
H_{\alpha}(e^{l}_{\lambda}\times I\cap (a_{\alpha}\times
I)_{\gamma_{\alpha}}))\cap e^{l}_{\gamma}\times I$.  By the
construction, $H_{\alpha}((c\times
I)_{\lambda})=(H_{\alpha}(c)\times I)_{\lambda}\subset
b_{\alpha}\times I$ where $H_{\alpha}(c)$ is a core of
$b_{\alpha}\times I$. By the above argument,  the lemma holds.

Case 4. $\gamma=j_{0}<\lambda$.

By Lemma 4.1.5, $e^{l}_{\lambda}\times I$ is disjoint from
$e^{l+1}_{m^{l+1}}\times I=(\partial
W_{j_{0}}^{l}-intd_{j_{0}}^{l})\times I$. By Lemma 4.3.10(2), (3)
and Lemma 4.4.3, $e^{l+1}_{\lambda}\times I\cap
e^{l+1}_{m^{l+1}}\times I=(\cup_{\gamma_{\alpha}\leq\lambda}
H_{\alpha}(e^{l}_{\lambda}\times I\cap (a_{\alpha}\times
I)_{\gamma_{\alpha}}))\cap e^{l+1}_{m^{l+1}}\times I$. Note that
each component of $e^{l}_{\lambda}\times I\cap (a_{\alpha}\times
I)_{\gamma_{\alpha}}$ is $(c\times I)_{\lambda}\subset
(a_{\alpha}\times I)_{\gamma_{\alpha}}$. By Definition 4.3.9,
$H_{\alpha}((c\times I)_{\lambda})$ intersects $b_{\alpha}\times
I$ in $(H_{\alpha}(c)\times I)_{\lambda}$. Since
$a_{\alpha}^{0}\subset d_{j_{0}}^{l}$, by Definition 4.3.9,
$e^{l+1}_{m^{l+1}}\times I\subset (\partial
W_{j_{0}}^{l}-inta_{\alpha}^{0})\times I$. By Lemma 4.1.5(8) and
Construction(*)(3),(4), $inta_{0,\alpha}^{*}\times I\cup
b_{0,\alpha}^{*}\times I$ is disjoint from
$e^{l+1}_{m^{l+1}}\times I$. By the construction,
$b_{\alpha}\times I$ intersects $e^{l+1}_{m^{l+1}}\times I$ in
$e^{l+1}_{m^{l+1}}\times I_{\alpha}$. Since $I_{\alpha}\subset
(0,1)$, by Proposition 4 for $k=l$, the lemma holds.

Case 5. $j_{0}<\gamma<\lambda$. By Lemma 4.1.5,
$e^{l}_{\gamma}\times I,e^{l}_{\lambda}\times I$ are disjoint from
$\partial W_{j_{0}}\times I$.  By Lemma 4.3.5, $b_{\alpha}\times
I\cap b_{\beta}\times I=\emptyset$. By  Definition 4.3.10(2),
$e^{l+1}_{\lambda}\times I\cap e^{l+1}_{\gamma}\times I=S_{1}\cup
S_{2}$ where $S_{1}= (e^{l}_{\lambda}\times
I-\cup_{\gamma_{\alpha}\leq\lambda,\gamma}(e^{l}_{\lambda}\times
I\cap (a_{\alpha}\times I)_{\gamma_{\alpha}}))\cap
(e^{l}_{\gamma}\times
I-\cup_{\gamma_{\alpha}\leq\lambda,\gamma}(e^{l}_{\gamma}\times
I\cap (a_{\alpha}\times I)_{\gamma_{\alpha}}))$ and
$S_{2}=\cup_{\gamma_{\alpha}\leq\gamma,\lambda}
H(e^{l}_{\lambda}\times I\cap (a_{\alpha}\times
I)_{\gamma_{\alpha}})\cap H(e^{l}_{\gamma}\times I\cap
(a_{\alpha}\times I)_{\gamma_{\alpha}})$. Since $H$ is injective,
the lemma holds.\qquad Q.E.D.

{\bf \em Lemma 4.4.5.} \ $F^{l+1}$ is a surface generated by an
abstract tree $\cup_{f} E^{l+1}_{f}\cup_{\gamma\in m(l+1)}
e^{l+1}_{\gamma}$.

{\bf \em Proof.} \ By Lemmas 4.4.1-4.4.4, we only need to prove
that $\cup_{f} E^{l+1}\cup_{\gamma\in m(l+1)} e^{l+1}_{\gamma}$ is
an abstract tree.

By Lemma 4.3.10, $\partial e^{l}_{\gamma}=\partial
e^{l+1}_{\gamma}$ for $\gamma\in m(l)$ and $\partial
e_{m^{l+1}}^{l+1}=\partial d_{j_{0}}^{l}$. Let
$\bigl\{e^{*}_{\gamma} \ | \ \gamma\in m(l+1)\bigr\}$ be  a set of
pairwise disjoint  arcs obtained by pushing $inte^{l+1}_{\gamma}$
off $\cup_{f} E^{l}_{f}$ in $\partial \cal V_{-}$$\times I$. By
Proposition 4(1) for $k=l$, each component of
$\cup_{f=0}^{l}E^{l}_{f}\cup_{\gamma\in m(l)} e_{\gamma}^{*}$ is
an tree. By Definition 4.3.3, $E^{l}_{f}=E^{l+1}_{f}$ for $1\leq
f\leq l$, and $E^{l+1}_{0},E^{l+1}_{l}\subset E_{0}^{l}$. Hence
each component of $\cup_{f=0}^{l+1}E^{l+1}_{f}\cup_{\gamma\in
m(l)} e_{\gamma}^{*}$ is also an tree.  Note that $m^{l+1}=j_{0}$.

By Lemma 3.3.1, $d_{j_{0}}^{l}=\cup_{i=\delta(0)}^{\theta(0)}
d_{i,0}\cup_{i=\delta(0)}^{\theta(0)} e_{i,0}$. Since
$d_{j_{0}}^{l}$ intersects $c_{l+1}^{l}$ in one point lying in
$intd_{0,0}\subset E^{l}_{0}$, and $E^{l+1}_{0}\cup
E^{l+1}_{l+1}=E^{l}_{0}-c_{l+1}^{l}\times (-1,1)$,
$\cup_{f=0}^{l+1}E^{l+1}_{f}\cup_{\gamma\in m(l)} e_{\gamma}^{*}$
contains at least two components. Furthermore, $E^{l+1}_{0}$ and
$E^{l+1}_{f_{\delta(0),0}}$ lies in one component of
$\cup_{f=0}^{l+1}E^{l+1}_{f}\cup_{\gamma\in m(l)} e_{\gamma}^{*}$,
$E^{l}_{l+1}$ and $E^{l+1}_{f_{\theta(0),0}}$ lies in another
component of $\cup_{f=0}^{l+1}E^{l+1}_{f}\cup_{\gamma\in m(l)}
e_{\gamma}^{*}$ as in Figure 18. Since $\partial_{1}
e^{l+1}_{m^{l+1}}=\partial_{1} d_{j_{0}}^{l}\subset
E^{l+1}_{f_{\delta(0),0}}$ and $\partial_{2}
e^{l+1}_{m^{l+1}}=\partial_{2} d_{j_{0}}^{l}\subset
E^{l+1}_{f_{\theta(0),0}}$, so each component of
$\cup_{f=0}^{l+1}E^{l+1}_{f}\cup_{\gamma\in m(l+1)}
e_{\gamma}^{*}$ is also  an  tree. By Definition 2.1.2, $\cup_{f}
E^{l+1}_{f}\cup_{\gamma\in m(l+1)} e^{l+1}_{\gamma}$ is an
abstract tree. \qquad Q.E.D.
\begin{center}
\includegraphics[totalheight=6cm]{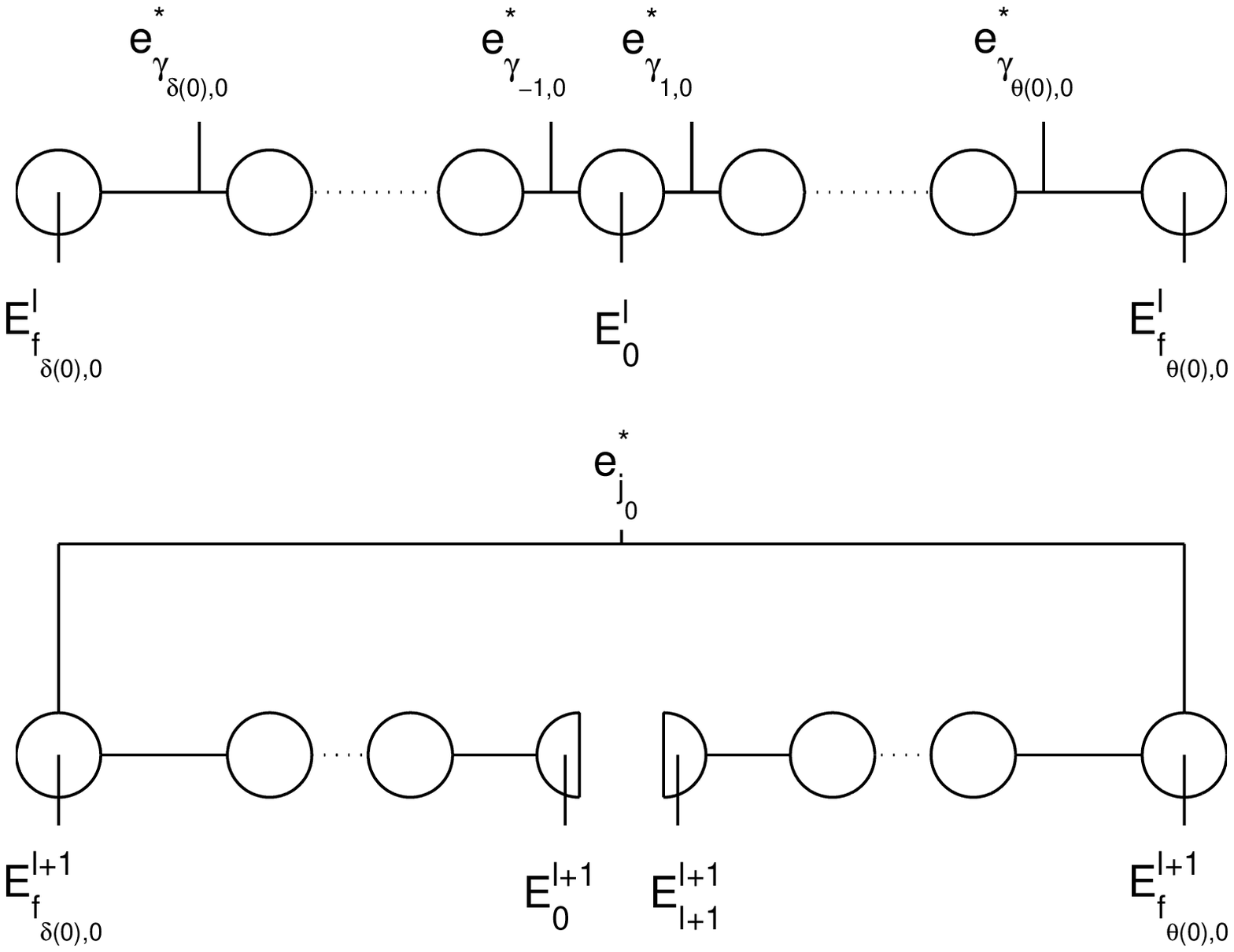}
\begin{center}
Figure 18
\end{center}
\end{center}

\subsection{The properties of $d_{j}^{l+1}$}

\ \ \ \ \ {\bf \em Lemma 4.5.1.} \ Suppose that $j\in
\bigl\{1,\ldots,n\bigr\}-m(l+1)$ and $\lambda\in m(l+1)$. If
$j<\lambda$, then either $d_{j}^{l+1}$ is disjoint from
$e^{l+1}_{\lambda}\times I$ or $j\in I(w_{\lambda},l+1)$ and
$s(w_{\lambda})=-$.

{\bf \em Proof.} \ By assumption, Proposition 4 holds for $k=l$.

Suppose that $j<\lambda$ and $j\notin I(w_{\lambda},l+1)$. Then,
by Proposition 4(3), $d_{j}^{l}$ is disjoint from
$e^{l}_{\lambda}\times I$ for $\lambda\in m(l)$. Now there are two
cases:

Case 1. $j<m^{l+1}=j_{0}$.

Suppose that $\lambda=m^{l+1}$. Since $j\notin
I(w_{m^{l+1}},l+1)$, by Lemma 4.1.5, Lemma 4.3.7,
$d_{j}^{l+1}=d_{j}^{l}$ is disjoint from $e^{l+1}_{m^{l+1}}\times
I\subset\partial W_{j_{0}}^{l}\times I$.

Suppose $\lambda\in m(l)$. Since $j<j_{0}$, by Lemma 4.3.5(3),
$d_{j}^{l+1}=d_{j}^{l}$ is disjoint from $(a_{0,\alpha}^{*}\cup
b_{0,\alpha}^{*})\times I$. Hence $d_{j}^{l+1}$ is disjoint from
$b_{\alpha}\times I$. By Definition 4.3.9,
$e^{l+1}_{\lambda}\times I\subset e^{l}_{\lambda}\times
I\cup_{\alpha} b_{\alpha}\times I$. Hence $d_{j}^{l+1}$ is
disjoint from $e^{l+1}_{\lambda}\times I$.

Case 2. $j>j_{0}=m^{l+1}$.

Since $\lambda>j$, $\lambda>j_{0}$. By Lemma 4.1.5, $d_{j}^{l},
e^{l}_{\lambda}\times I$ are disjoint from $\partial
W_{j_{0}}^{l}\times I$. If $j<\gamma_{\alpha}$, then, by Lemma
3.5.3, $d_{j}^{l}$ is disjoint from $(a_{\alpha}\times
I)_{\gamma_{\alpha}}$ except $j=j_{\alpha}$. Now there are two
sub-cases:

(1) \ $j\neq j_{\alpha}$ for each $j_{\alpha}\in L(c_{l+1}^{l})$.

By Lemma 4.3.7(2),
$d_{j}^{l+1}=(d_{j}^{l}-\cup_{\gamma_{\alpha}<j} (a_{\alpha}\times
I)_{\gamma_{\alpha}})\cup_{\gamma_{\alpha}<j}
H_{\alpha}(d_{j}^{l}\cap(a_{\alpha}\times I)_{\gamma_{\alpha}})$.
By Construction(*)(3), (4) and Lemma 4.3.5(3),
$d_{j}^{l}-\cup_{\gamma_{\alpha}<j} (a_{\alpha}\times
I)_{\gamma_{\alpha}}$ is disjoint from $(a_{0,\alpha}^{*}\cup
b_{0,\alpha}^{*})\times I$. By Lemma 4.3.10 and Lemma 4.4.3,
$e^{l+1}_{\lambda}\times I=(e^{l}_{\lambda}\times
I-\cup_{\gamma_{\alpha}\leq\lambda} (a_{\alpha}\times
I)_{\gamma_{\alpha}})\cup_{\gamma_{\alpha}\leq\lambda}
H_{\alpha}(e^{l}_{\lambda}\times I \cap (a_{\alpha}\times
I)_{\gamma_{\alpha}})$. By Lemma 4.3.5(4), $e^{l}_{\lambda}\times
I-\cup_{\gamma_{\alpha}\leq\lambda} (a_{\alpha}\times
I)_{\gamma_{\alpha}}$ is disjoint from $(a_{0,\alpha}^{*}\cup
b_{0,\alpha}^{*})\times I$. By Lemma 4.3.10(2), $d_{j}^{l+1}\cap
e^{l+1}_{\lambda}\times I=S_{1}\cup S_{2}$ where $S_{1}=
(d_{j}^{l}-\cup_{\gamma_{\alpha}<j} (a_{\alpha}\times
I)_{\gamma_{\alpha}})\cap (e^{l}_{\lambda}\times
I-\cup_{\gamma_{\alpha}\leq\lambda} (a_{\alpha}\times
I)_{\gamma_{\alpha}})$ and $S_{2}= H(e^{l}_{\lambda}\times I
\cap(\cup_{\gamma_{\alpha}\leq\lambda} (a_{\alpha}\times
I)_{\gamma_{\alpha}}))\cap
H(d_{j}^{l}\cap(\cup_{\gamma_{\alpha}<j} (a_{\alpha}\times
I)_{\gamma_{\alpha}}))$. Since $H$ is injective, by Proposition
4(3) for $k=l$, the lemma holds.

(2) \ $j= j_{\beta}$ for some $j_{\beta}\in L(c_{l+1}^{l})$.

Now $(a_{\beta}\times
\bigl\{0\bigr\})_{\gamma_{\beta}}=a_{\beta}\subset d_{j}^{l}$. By
Lemma 4.3.10 and Lemma 4.4.3,
$d_{j}^{l+1}=(d_{j}^{l}-a_{\beta}\cup_{\gamma_{\alpha}<j}
(a_{\alpha}\times I)_{\gamma_{\alpha}})\cup
H(d_{j}^{l}\cap(\cup_{\gamma_{\alpha}<j} (a_{\alpha}\times
I)_{\gamma_{\alpha}}))\cup b_{\beta}$. By assumption, $\lambda\neq
\gamma_{\beta}$. Hence, by the same argument in (1), the lemma
holds.\qquad Q.E.D.\vskip 0.5mm

{\bf \em Lemma 4.5.2.} \ If $j>\lambda$, then each component of
$d_{j}^{l+1}\cap e^{l+1}_{\lambda }\times I$ is a core of
$e^{l+1}_{\lambda}\times (0,1)$.

{\bf \em Proof.}  \ By Proposition 4(2) for $k=l$ and the same
argument as that in the proof of Lemma 4.4.4, the lemma
holds.\qquad Q.E.D.\vskip 0.5mm

{\bf \em Lemma 4.5.3.} \ For each $j\notin L(c_{l+1}^{l})$,
$d_{j}^{l+1}=\cup_{i=1}^{\theta(j)}d_{j,f_{i,j}}^{l}
\cup_{i=1}^{\theta(j)-1} e_{\gamma_{i,j}}^{*}$ satisfying the
following conditions:

(1) \ $\gamma_{i,j}\in m(l)$ is as in Lemma 3.1.2, and
$e_{\gamma_{i,j}}^{*}$ is a core of $e^{l+1}_{\gamma_{i,j}}\times
(0,1)$ for $1\leq i\leq\theta(j)-1$.

(2) \ $f_{i,j}$ is as in Lemma 3.1.2. Furthermore, if $f_{i,j}=0$,
then  $d_{j,0}^{l}$ is either properly embedded in $E^{l+1}_{0}$
or $E^{l+1}_{l+1}$ which is disjoint from $\cup_{\gamma<j}
inte_{\gamma}^{l+1}\times I$, if $f_{i,j}\neq 0$, then
$d_{j,f_{i,j}}^{l}$ is properly embedded in $E^{l+1}_{f_{i,j}}$
which is disjoint from $\cup_{\gamma<j} inte_{\gamma}^{l+1}\times
I$ for $1\leq i\leq \theta(j)$.

(3) \ $d_{j}^{l+1}$ is regular in $\cup_{f}
E^{l+1}_{f}\cup_{\gamma<j} e^{l+1}_{\gamma}\times I$.

{\bf \em Proof.} \ By Lemma 3.1.2, in $F^{l}$,
$d_{j}^{l}=\cup_{i=1}^{\theta(j)}d_{j,f_{i,j}}^{l}
\cup_{i=1}^{\theta(j)-1} e_{\gamma_{i,j}}$ such that
$\gamma_{i,j}<j$ and $d_{j,f_{i,j}}^{l}$ is disjoint from
$\cup_{\gamma<j} inte^{l}_{\gamma}\times I$.  Since $j\notin
L(c_{l+1}^{l})$, $d_{j,f}^{l}$ is disjoint from $c_{l+1}^{l}\times
I$ even if $f=0$. Hence $d_{j,0}^{l}$ lies in one of $E^{l+1}_{0}$
and $E^{l+1}_{l+1}$.

Suppose that $\gamma<j$. There are two cases:

(1) \ $j<j_{0}$.

Now $\gamma_{i,j}<j<j_{0}$. By Lemma 4.3.7 and Definition 4.3.9,
$e^{l+1}_{\gamma_{i,j}}\times I=e^{l}_{\gamma_{i,j}}\times I$ and
$d_{j}^{l+1}=d_{j}^{l}$. Furthermore, $e^{l+1}_{\gamma}\times
I=e^{l}_{\gamma}\times I$ for $\gamma<j$. In this case,
$e_{\gamma_{i,j}}$ is also a core of $e^{l+1}_{\gamma_{i,j}}\times
(0,1)$. By Lemma 3.1.2, the lemma holds.

(2) \ $j>j_{0}$.

By Lemma 4.1.5, $d_{j,f_{i,j}}^{l}\subset d_{j}^{l}$ is disjoint
from $e^{l+1}_{m^{l+1}}\times I\subset\partial W_{j_{0}}^{l}\times
I$.

Now we claim that $intd_{j,f_{i,j}}^{l}$  is disjoint from
$(a_{0,\alpha}^{*}\cup b_{0,\alpha}^{*})\times I$. By
Construction(*)(3) and (4), if $intd_{j,f_{i,j}}^{l}\cap
(a_{0,\alpha}^{*}\cup b_{0,\alpha}^{*})\times I\neq\emptyset$,
then $intd_{j,f_{i,j}}^{l}\cap (a_{0,\alpha}^{*}\cup
b_{0,\alpha}^{*})\times I\subset (a_{\lambda}\times
I)_{\gamma_{\lambda}}$ for some $j_{\lambda}\in L(c_{l+1}^{l})$.
By assumption, $j\neq j_{\lambda}$. By Lemma 4.3.5,
$j>\gamma_{\lambda}$, contradicting Lemma 3.1.2.

If $\gamma<j$, then, by Lemma 3.1.2, $d_{j,f_{i,j}}^{l}$ is
disjoint from $inte^{l}_{\gamma}\times I$. By Definition 4.3.9,
$inte^{l+1}_{\gamma}\times I\subset \cup_{\alpha} b_{\alpha}\times
I\cup inte^{l}_{\gamma}\times I$. Hence $d_{j,f_{i,j}}^{l}$ is
disjoint from $inte^{l+1}_{\gamma}\times I$. Now
$d_{j}^{l+1}=\cup_{i=1}^{\theta(j)}d_{j,f_{i,j}}^{l}
\cup_{i=1}^{\theta(j)-1} e_{\gamma_{i,j}}^{*}$  where
$e^{*}_{\gamma_{i,j}}=(e_{\gamma_{i,j}}-e_{\gamma_{i,j}}\cap
(\cup_{\alpha}(a_{\alpha}\times I)_{\gamma_{\alpha}}))\cup
\cup_{\alpha} H_{\alpha}(e_{\gamma_{i,j}}\cap (a_{\alpha}\times
I)_{\gamma_{\alpha}})$ is a core of $e^{l+1}_{\gamma}\times
(0,1)$. Now each component of $d_{j}^{l+1}\cap
(\cup_{\gamma<j}e^{l+1}_{\gamma}\times I)$ is
$e_{\gamma_{i,j}}^{*}$ for $1\leq i\leq\theta(j)-1$. By Lemma
4.5.1 and 4.5.2, (1) and (2) holds.

By Lemma 4.3.7, $\partial d_{j}^{l+1}=\partial d_{j}^{l}$. By
Lemma 4.3.10, $(\partial e^{l+1}_{\gamma})\times I=(\partial
e^{l}_{\gamma})\times I$ for $\gamma\in m(l)$. Hence $\partial
d_{j}^{l+1}$ is disjoint from $e^{l+1}_{\gamma}\times I$. Since
$j>j_{0}$, $j\notin I(w_{j_{0}},l)$, by Lemma 4.1.5(4), $\partial
d_{j}^{l+1}=\partial d_{j}^{l}$ is disjoint from
$e^{l+1}_{m^{l+1}}\times I\subset
\partial W_{j_{0}}^{l}\times I$. By Lemma 4.3.7(2) and Definition
4.3.9, $d_{j}^{l}$ is properly embedded in $\cup_{f}
E^{l+1}_{f}\cup_{\gamma<j} e^{l+1}_{\gamma}\times I$. By Lemma
3.1.2, $f_{i,j}\neq f_{i,j}$ and $\gamma_{i,j}=\gamma_{r,j}$ for
$i\neq r$. By Lemma 4.5.2, $d_{j}^{l+1}$ is regular in $\cup_{f}
E^{l+1}_{f}\cup_{\gamma<j} e^{l+1}_{\gamma}\times I$.\qquad
Q.E.D.\vskip 3mm

{\bf \em Lemma 4.5.4.} \ If $\gamma<j_{\alpha}$, then
$inta_{0,\alpha}^{*}\cup intb_{0,\alpha}^{*}$ is disjoint from
$e^{l+1}_{\gamma}\times I$.

{\bf \em Proof.} \ Suppose that $\gamma<j_{\alpha}$. There are
three cases:

Case 1. $\gamma<j_{0}=m^{l+1}$.

By Definition 4.3.9, $e^{l+1}_{\gamma}\times
I=e^{l}_{\gamma}\times I$. By Lemma 4.3.5(4), the lemma holds.

Case 2. $\gamma=j_{0}=m^{l+1}$.

By Lemma 3.5.2 and Definition 4.3.9, $a^{0}_{\alpha}\subset
d_{j_{0}}^{l}$. Hence $e^{l+1}_{m^{l+1}}\times I\subset (\partial
W_{j_{0}}^{l}-inta_{\alpha}^{0})\times I$. By Lemma 4.1.5(8) and
Construction(*)(3) and (4), the lemma holds.

Case 3. $\gamma>j_{0}=m^{l+1}$. By Lemma 4.3.5(4),
$e^{l}_{\gamma}\times I\cap (inta_{0,\alpha}^{*}\cup
intb_{0,\alpha}^{*})\subset \cup_{\beta\leq\gamma}
(a_{\beta}\times I)_{\gamma_{\beta}}$. Note that
$\gamma_{\alpha}>j_{\alpha}$ if $\gamma_{\alpha}\neq\emptyset$.
Since $\gamma<j_{\alpha}$, $\gamma<\gamma_{\alpha}$ if
$\gamma_{\alpha}\neq \emptyset$.

By Lemma 4.3.10(2) and (3) and Lemma 4.4.3,
$e^{l+1}_{\gamma}\times I=(e^{l}_{\gamma}\times
I-(e^{l}_{\gamma}\times I\cap(\cup_{\gamma_{\beta}\leq\gamma}
(a_{\beta}\times
I)_{\gamma_{\beta}}))\cup_{\gamma_{\beta}\leq\gamma}
H_{\beta}(e^{l}_{\gamma}\times I\cap (a_{\beta}\times
I)_{\gamma_{\beta}})$. Note that $H_{\beta}(e^{l}_{\gamma}\times
I\cap (a_{\beta}\times I)_{\gamma_{\beta}})\subset b_{\beta}\times
I$. By Lemma 4.3.5, $b_{\alpha}\times I\cap b_{\beta}\times
I=\emptyset$. Hence the lemma holds.\qquad Q.E.D.\vskip 0.5mm

{\bf \em Lemma 4.5.5.} \ Suppose that $j_{\alpha}\in
L(c_{l+1}^{l})$and $\alpha\neq 0$. Then

(1) \ Each component of $d_{j_{\alpha}}^{l+1}\cap
(\cup_{\gamma<j_{\alpha}} e^{l+1}_{\gamma}\times I)$ is  a core of
$e^{l+1}_{\gamma}\times I$ for $\gamma\in\bigl\{m^{l+1}\bigr\}\cup
\bigl\{\gamma_{i,\alpha} \ | \ i<\delta(0,\alpha) \ or \ i>
\theta(0,\alpha)\bigr\}\cup \bigl\{\gamma_{i,0} \ | \
i<\delta(0,\alpha) \ or \ i> \theta(0,\alpha)\bigr\}$.

(2) \ Each  component of $d_{j_{\alpha}}^{l+1}-
\cup_{\gamma<j_{\alpha}} \textrm{int}e^{l+1}_{\gamma}\times I$
lies in $E^{l+1}_{f}$ for $f\in\bigl\{0,l+1\bigr\}\cup
\bigl\{f_{i,\alpha} \ | \ i\leq\delta(0,\alpha) \ or \
i\geq\theta(0,\alpha)\bigr\}\cup \bigl\{f_{i,0} \ | \
i<\delta(0,\alpha) \ or \ i>\theta(0,\alpha)\bigr\}$.

(3) \ $d_{j_{\alpha}}^{l+1}$ is regular in $\cup_{f}
E^{l+1}_{f}\cup_{\gamma<j_{\alpha}} e^{l+1}_{\gamma}\times I$.

{\bf \em Proof.} \ Since $a_{\alpha}=a_{\alpha}\times
\bigl\{0\bigr\}\subset (a_{\alpha}\times I)_{\gamma_{\alpha}}$. By
Lemma 4.3.7(3), $d_{j_{\alpha}}^{l+1}=S\cup b_{\alpha}$ where $S=
(d_{j_{\alpha}}^{l}-a_{\alpha}\cup_{\gamma_{\beta}<j_{\alpha}}
(a_{\beta}\times
I)_{\gamma_{\beta}})\cup_{\gamma_{\beta}<j_{\alpha}}
H_{\beta}(d_{j_{\alpha}}^{l}\cap (a_{\beta}\times
I)_{\gamma_{\beta}})$.

By Definition 3.5.1,
$a_{\alpha}=\cup_{i=\delta(0,\alpha)}^{\theta(0,\alpha)}d_{i,\alpha}\cup_{i=\delta(0,\alpha)}^{\theta(0,\alpha)}e_{i,\alpha}$.
By Lemma 3.3.1,
$d_{j_{\alpha}}^{l}-inta_{\alpha}=\cup_{i=\delta(\alpha)}^{\delta(0,\alpha)-1}(d_{i,\alpha}\cup
e_{i,\alpha})
\cup_{i=\theta(0,\alpha)+1}^{\theta(\alpha)}(d_{i,\alpha}\cup
e_{i,\alpha})$. By the argument in Lemma 4.5.3,
$S_{1}=\cup_{i=\delta(\alpha)}^{\delta(0,\alpha)-1}(d_{i,\alpha}\cup
e_{i,\alpha}^{*})
\cup_{i=\theta(0,\alpha)+1}^{\theta(\alpha)}(d_{i,\alpha}\cup
e_{i,\alpha}^{*})$ satisfying the following conditions:

(1) \ $d_{i,\alpha}$ is an arc in $E^{l+1}_{f_{i,\alpha}}$ such
that $intd_{i,\alpha}$ is disjoint from $e^{l+1}_{\gamma}\times I$
for $\gamma<j_{\alpha}$.

(2) \ $e^{*}_{i,\alpha}=(e_{i,\alpha}-e_{i,\alpha}\cap
(\cup_{\alpha}(a_{\alpha}\times I)_{\gamma_{\alpha}}))\cup
\cup_{\alpha} H_{\alpha}(e_{i,\alpha}\cap (a_{\alpha}\times
I)_{\gamma_{\alpha}})$ is a core of
$e^{l+1}_{\gamma_{i,\alpha}}\times I$.

By the construction, $b_{\alpha}=b_{\alpha}\times
\bigl\{0\bigr\}=a_{0,\alpha}^{*}\cup (\partial
W_{j_{0}}^{l}-a_{\alpha}^{0})\times \bigl\{t_{\alpha}\bigr\}\cup
b_{0,\alpha}^{*}$. Note that $a_{0,\alpha}^{*}\subset
E^{l+1}_{f_{\delta(0,\alpha),\alpha}}$ and
$b_{0,\alpha}^{*}\subset E^{l+1}_{f_{\theta(0,\alpha),\alpha}}$.
Let $S_{2}=(\partial W_{j_{0}}^{l}-a_{\alpha}^{0})\times
\bigl\{t_{\alpha}\bigr\}$. Note that $\gamma_{i,0}<j_{0}$ for
$\delta(0)\leq i\leq \theta(0)$. By Definition 4.3.9,
$e^{l+1}_{\gamma_{i,0}}\times I=e^{l}_{\gamma_{i,0}}\times I$. By
Lemma 3.3.1, $d_{j_{0}}^{l}=\cup_{\delta(0)}^{\theta(0)}
d_{i,0}\cup_{\delta(0)}^{\theta(0)} e_{i,0}$. By Lemma 3.5.1,
$a_{\alpha}^{0}=\cup_{\delta(0,\alpha)}^{\theta(0,\alpha)}
d_{i,0}\cup_{\delta(0)}^{\theta(0)} e_{i,0}$. Hence $S_{2}$
intersects $e^{l+1}_{\gamma}\times I$ in a core of
$e^{l+1}_{\gamma}\times I$ for $\gamma\in \bigl\{\gamma_{i,0} \ |
\ i<\delta(0,\alpha) \ or \ i> \theta(0,\alpha)\bigr\}$, and
$E^{l+1}_{f_{i,\alpha}}$ in an arc
$d_{i,0}\times\bigl\{t_{\alpha}\bigr\}$ for $f\in \bigl\{f_{i,0} \
| \ i<\delta(0,\alpha) \ or \ i>\theta(0,\alpha)\bigr\}$. In
particular, by Definition 4.3.9, $S_{2}$ intersects
$e^{l+1}_{m^{l+1}}\times I$ in a core of $e^{l+1}_{m^{l+1}}\times
(0,1)$. By Lemma 4.5.4, $a_{0,\alpha}^{*}\cup b_{0,\alpha}^{*}$ is
disjoint from $inte^{l+1}_{\gamma}\times I$ for
$\gamma<j_{\alpha}$.

Now we claim that $d_{i,0}\times\bigl\{t_{\alpha}\bigr\}$ is
disjoint from $inte^{l+1}_{\gamma}\times I$ for
$\gamma<j_{\alpha}$. There are two cases:

Case 1. $\gamma\leq j_{0}=m^{l+1}$.

Suppose that $\gamma<j_{0}$. By Definition 4.3.9,
$e^{l+1}_{\gamma}\times I=e^{l}_{\gamma}\times I$. By Lemma
4.1.5(6), $d_{i,0}\times\bigl\{t_{\alpha}\bigr\}\subset
d_{i,0}\times I$ is disjoint from $e^{l+1}_{\gamma}\times I$.
Since $d_{i,0}\subset d_{j_{0}}^{l}$, $inte^{l+1}_{m^{l+1}}\times
I=(\partial W_{j_{0}}^{l}-d_{j_{0}}^{l})\times I$ is disjoint from
$d_{i,0}\times\bigl\{t_{\alpha}\bigr\}$.

Case 2. $\gamma>j_{0}$.

Since $\gamma<j_{\alpha}$, $\gamma<\gamma_{\alpha}$. Since
$\gamma>j_{0}>Max\bigl\{\gamma_{i,0} \ | \ \delta(0)\leq
i\leq\theta(0)\bigr\}$, by Lemma 3.5.3, $e^{l}_{\gamma}\times I$
is disjoint from $(inta_{\alpha}\times I)_{\gamma_{\alpha}}$.

By Lemma 4.3.10, $e^{l+1}_{\gamma}\times I=(e^{l}_{\gamma}\times
I-\cup_{\gamma_{\beta}\leq\gamma} (a_{\beta}\times
I)_{\gamma_{\beta}})\cup_{\gamma_{\beta}\leq\gamma}
H_{\beta}(e^{l}_{\gamma}\times I\cap (a_{\beta}\times
I)_{\gamma_{\beta}})$. Note that $H_{\beta}(e^{l}_{\gamma}\times
I\cap (a_{\beta}\times I)_{\gamma_{\beta}})\subset b_{\beta}\times
I$. By Lemma 4.3.5, $b_{\alpha}\times I\cap b_{\beta}\times
I=\emptyset$. Hence $e^{l+1}_{\gamma}\times I$ is disjoint from
$H_{\alpha}((inta_{\alpha}\times
I)_{\gamma_{\alpha}})=intb_{\alpha}\times I$. Thus
$d_{i,0}\times\bigl\{t_{\alpha}\bigr\}$ is disjoint from
$inte^{l+1}_{\gamma}\times I$.

By Lemma 3.3.4(3), $\bigl\{\gamma_{i,\alpha} \ | \
i>\theta(0,\alpha) \ or \ i<\delta(0,\alpha)\bigr\}\cap
\bigl\{\gamma_{i,0} \ | \ i>\theta(0,\alpha) \ or \
i<\delta(0,\alpha)\bigr\}=\emptyset$. By Lemma s 4.5.1 and 4.5.2,
$d_{j}^{l+1}$ is regular in $\cup_{f} E^{l+1}_{f}\cup_{\gamma<j}
e^{l+1}_{\gamma}\times I$.\qquad Q.E.D.\vskip 0.5mm

\subsection{The proofs of Propositions 4-6 for the case: $k=l+1$ and $s(v_{l+1})=+$}

\ \ \ \ \ In this section, we shall finish the proofs of
Propositions 4-6 for the case: $k=l+1$ and $s(v_{l+1})=+$.

{\bf \em Lemma 4.6.1.} \ $\partial c_{i}^{l+1}$ is disjoint from
$d_{j}^{l+1}$ and $e^{l+1}_{\gamma}$.

{\bf \em Proof.} \ By Definition 4.3.3, $c_{i}^{l+1}=c_{i}^{l}$.
By Lemma 4.1.5(5), $\partial c_{i}^{l+1}$ is disjoint from
$\partial W_{j_{0}}^{l}\times I$. By Lemma 3.5.8, $\partial
c_{i}^{l+1}$ is disjoint from $b_{\alpha}\times I$. By Proposition
4(5), $\partial c_{i}^{l+1}$ is disjoint from $d_{j}^{l}$ and
$e^{l}_{\gamma}\times I$. By Definition 4.3.9, $d_{j}^{l+1}\subset
d_{j}^{l}\cup \partial W_{j_{0}}^{l}\times I\cup_{\alpha}
b_{\alpha}\times I$ and $e^{l+1}_{\gamma}\times I\subset
e^{l}_{\gamma}\times I\cup
\partial W_{j_{0}}^{l}\times I\cup_{\alpha} b_{\alpha}\times I$. Hence the
lemma holds.\qquad Q.E.D.

{\bf \em The proof of Proposition 4.} \ By Lemmas 4.4.1-4.4.5,
Proposition 4(1) holds. By Lemma 4.5.3(3) and Lemma 4.5.5(3),
$d_{j}^{l+1}$ is regular in $\cup_{f}E^{l+1}_{f}\cup_{\gamma<j}
e^{l+1}_{\gamma}\times I$. Hence Proposition 4(2) holds.
Proposition 4(3) follows from Lemma 4.5.1.

If $j\notin L(c_{l+1}^{l})$, then, by Proposition 4(4),
$c_{i}^{l+1}=c_{i}^{l}$ intersects $d_{j,f_{i,j}}^{l}$ in at most
one point. By Lemma 4.5.3 and Lemma 4.3.7,
$d_{j}^{l+1}-\cup_{\gamma<j} inte^{l+1}_{\gamma}\times I$
intersects $c_{i}^{l+1}$ in at most one point.

Suppose that $j=j_{\alpha}$. By Lemma 3.3.4(4),
$\bigl\{f_{i,\alpha} \ | \ i>\theta(\alpha,\beta) \ or \
i<\delta(\alpha,\beta)\bigr\}\cap \bigl\{f_{i,\beta} \ | \
i>\theta(\alpha,\beta) \ or \
i<\delta(\alpha,\beta)\bigr\}=\emptyset$.

Suppose fist that $\delta(0,\alpha)=\theta(0,\alpha)$, then, by
Lemma 3.5.8, $a_{0,\alpha}^{*}\subset E^{l+1}_{0}$
$b_{0,\alpha}^{*}\subset E^{l+1}_{l+1}$. Note that $E^{l+1}_{0},
E^{l+1}_{l+1}\subset E^{l}_{0}$. By Lemma 3.3.1,
$\bigl\{0,l+1\bigr\}\cap(\bigl\{f_{i,\alpha} \ | \
i>\theta(\alpha,\beta) \ or \ i<\delta(\alpha,\beta)\bigr\}\cup
\bigl\{f_{i,\beta} \ | \ i>\theta(\alpha,\beta) \ or \
i<\delta(\alpha,\beta)\bigr\})=\emptyset$.

Suppose now that $\delta(0,\alpha)\neq\theta(0,\alpha)$. Then, by
Lemma 3.3.1,
$\bigl\{f_{\delta(0,\alpha),\alpha},f_{\theta(0,\alpha),\alpha}\bigr\}\cap(\bigl\{f_{i,\alpha}
\ | \ i>\theta(\alpha,\beta) \ or \
i<\delta(\alpha,\beta)\bigr\}\cup \bigl\{f_{i,\beta} \ | \
i>\theta(\alpha,\beta) \ or \
i<\delta(\alpha,\beta)\bigr\})=\emptyset$.

By Lemma 4.5.5,
$d_{j}^{l+1}-\cup_{\gamma<j}inte^{l+1}_{\gamma}\times I
=\cup_{i=\delta(\alpha)}^{\delta(0,\alpha)-1}
d_{i,\alpha}\cup_{i=\theta(0,\alpha)+1}^{\theta(\alpha)}d_{i,\alpha}
\cup_{i=\delta(0)}^{\delta(0,\alpha)-1}
d_{i,0}\times\bigl\{t_{\alpha}\bigr\}\cup_{i=\theta(0,\alpha)+1}^{\theta(0)}d_{i,0}\times\bigl\{t_{\alpha}\bigr\}\cup
a_{0,\alpha}^{*}\cup b_{0,\alpha}^{*}$. By Lemma 4.3.7, there is
at most one component $b^{i}$ of
$d_{j}^{l+1}-\cup_{\gamma<j}inte^{l+1}_{\gamma}\times I$ such that
$b^{i}\cap c^{l+1}\neq \emptyset$ for each $i\geq l+2$.

By Lemma 4.1.5(6), Construction(*)(3), (4) and Lemma 3.5.8,
$c_{i}^{l+1}=c_{i}^{l}$ intersects $a_{0,\alpha}^{*}$ in at most
one point. Furthermore, $c_{i}^{l+1}=c_{i}^{l}$ intersects
$a_{0,\alpha}^{*}$in one point if and only if $c_{i}^{l+1}\cap
a_{0,\alpha}$ in one point, where $a_{0,\alpha}$ is as in Lemma
3.5.5. Similarly, $c_{i}^{l+1}=c_{i}^{l}$ intersects
$b_{0,\alpha}^{*}$in one point if and only if $c_{i}^{l+1}\cap
b_{0,\alpha}$ in one point. By Lemma 4.1.5(6), $c_{i}^{l+1}$
intersects $d_{i,0}\times\bigl\{t_{\alpha}\bigr\}$ in one point if
and only if $c_{i}^{l}$ intersects $d_{i,0}$ in one point. By
Proposition 4(4) and Lemma 3.3.1, $c_{i}^{l+1}=c_{i}^{l}$
intersects $d_{i,\alpha}$ in at most one point. Hence Proposition
4(4) holds.

Proposition 4(5) follows from Lemma 4.6.1.\qquad Q.E.D.\vskip
1.5mm

{\bf \em The proof of Proposition 5.} \ Suppose that $j\notin
L(c_{l+1}^{l})$. By the argument in the proof of Lemma 4.5.3 and
Lemma 3.1.4, $L(d_{j}^{l+1})=L(d_{j}^{l})$.

Suppose that $j=j_{\alpha}\in L(c_{l+1}^{l})$. By the argument in
the proofs of Lemma 4.5.5 and Proposition 4 for $k=l+1$,
$L(d_{j}^{l+1}) =\cup_{i=\delta(\alpha)}^{\delta(0,\alpha)-1}
L(d_{i,\alpha})\cup_{i=\theta(0,\alpha)+1}^{\theta(\alpha)}L(d_{i,\alpha})
\cup_{i=\delta(0)}^{\delta(0,\alpha)-1}
L(d_{i,0})\cup_{i=\theta(0,\alpha)+1}^{\theta(0)}L(d_{i,0})\cup
L(a_{0,\alpha}^{*})\cup L(b_{0,\alpha}^{*})$.

Now there are two cases:

Case 1. $\delta(0,\alpha)=\theta(0,\alpha)=0$.

By Lemma 3.5.8(1), $L(a_{0,\alpha}^{*})\cup
L(b_{0,\alpha}^{*})=L(d_{0,0})\cup L(d_{0,\alpha})- L(d_{0,0})\cap
L(d_{0,\alpha})$. By Lemma 3.3.1, $f_{i,0}\neq f_{r,\alpha}$. By
Proposition 4(4), $L(d_{i,0})\cap L(d_{i,\alpha})=\emptyset$ for
$i\neq 0$. Hence, $L(d_{j}^{l+1})=L(d_{j}^{l})\cup
L(d_{m^{l+1}}^{l})- L(d_{j}^{l})\cap L(d_{m^{l+1}}^{l})$.

Case 2. $\delta(0,\alpha)\neq\theta(0,\alpha)$.

By Lemma 3.3.4, $L(d_{i,0})=L(d_{i,\alpha})$ for
$\delta(0,\alpha)<i<\theta(0,\alpha)$. By Lemma 3.5.8(2),
$L(a_{0,\alpha}^{*})=L(d_{\delta(0,\alpha),0})\cup
L(d_{\delta(0,\alpha),\alpha})-L(d_{\delta(0,\alpha),0})\cap
L(d_{\delta(0,\alpha),\alpha})$. Similarly,
$L(b_{0,\alpha}^{*})=L(d_{\theta(0,\alpha),0})\cup
L(d_{\theta(0,\alpha),\alpha})-L(d_{\theta(0,\alpha),0})\cap
L(d_{\theta(0,\alpha),\alpha})$. Hence
$L(d_{j}^{l+1})=L(d_{j}^{l})\cup L(d_{m^{l+1}}^{l})-
L(d_{j}^{l})\cap L(d_{m^{l+1}}^{l})$. \\ \qquad Q.E.D.\vskip 0.5mm

{\bf \em Lemma 4.6.3.} \ $b_{j}^{l+1}$ is  obtained by doing band
sums with some copies of $\partial W_{j_{0}}^{l}$ to $d_{j}^{l}$.

{\bf \em Proof.} \ Since each component of $d_{j}^{l}\cap
(a_{\alpha}\times I)_{\gamma_{\alpha}}$ is a core of
$(a_{\alpha}\times I)_{\gamma_{\alpha}}$, say
$a_{\alpha}\times\bigl\{t_{j}\bigr\}$. By Definition 4.3.6,
$H_{\alpha}(a_{\alpha}\times\bigl\{t_{j}\bigr\})$ intersects
$a_{0,\alpha}^{*}\times I$ in $a_{0,\alpha}^{*}\times
\bigl\{t_{j}\bigr\}$, $b_{0,\alpha}^{*}\times I$ in
$b_{0,\alpha}^{*}\times \bigl\{t_{j}\bigr\}$, $(\partial
W_{j_{0}}^{l}-inta_{\alpha}^{0})\times I$ in $(\partial
W_{j_{0}}^{l}-inta^{0}_{\alpha})\times \bigl\{t_{j}\bigr\}$.  Now
by Lemma 3.5.5, $a_{0,\alpha}^{*}\times \bigl\{t_{j}\bigr\}\cup
b_{0,\alpha}^{*}\times \bigl\{t_{j}\bigr\}\cup
a_{\alpha}\times\bigl\{t_{j}\bigr\}\cup
a_{\alpha}^{0}\times\bigl\{t_{j}\bigr\}$ bounds a disk in
$F^{l+1}$. Hence,  the lemma holds.\qquad Q.E.D.

{\bf \em The Proof of Proposition 6.} \ Suppose that $i\geq l+2$
and $s(v_{i})=-$. Then, by Proposition 6, there is a properly
embedded disk $V_{i}^{l}$ in $\cal V_{-}$ such that $\partial
V_{i}^{l}\cap F^{l}=c_{i}^{l}\cup_{r\in I(v_{i},l)} c_{r}^{l}$.
Now we denote by $V_{i}^{l+1}$ the disk $V_{i}^{l}$. By Lemma
4.1.5(7),  $\partial W_{j_{0}}^{l}\times I$ is disjoint from
$\partial V_{i}^{l}-c_{i}^{l}\cup_{r\in I(v_{i},l)} v_{r}^{l}$ for
each $i\geq l+2$ with $s(v_{i})=-$. Since
$F^{l+1}=(F^{l}-c_{l+1}^{l}\times (-1,1))\cup
e^{l+1}_{m^{l+1}}\times I$,  $\partial V_{i}^{l+1}\cap
F^{l+1}=c_{i}^{l+1}\cup_{r\in I(v_{i},l+1)} c_{r}^{l+1}$. In this
case, if $l+1\in I(v_{i},l)$, then
$I(v_{i},l+1)=I(v_{i},l)-\bigl\{l+1\bigr\}$, if not, then
$I(v_{i},l+1)=I(v_{i},l)$.

Now suppose that $j\notin m(l+1)$ and $s(w_{j})=-$. By Lemma
4.1.5(7), $\partial W_{j_{0}}^{l}\times I$ is disjoint from
$\partial W_{j}^{l}$.  By Proposition 6,  $\partial W_{j}^{l}\cap
F^{l}=d_{j}^{l}\cup_{r\in I(w_{j},l)}d_{r}^{l}$.  Now let
$C_{j}=(\partial W_{j}^{l}-d_{j}^{l}\cup_{r\in
I(w_{j},l)}d_{r}^{l})\cup d_{j}^{l+1}\cup_{i\in I(w_{j},l)}
d_{r}^{l+1}$. Then $C_{j}\cap F^{l+1}=d_{j}^{l+1}\cup_{r\in
I(w_{j},l)} d_{r}^{l+1}$. Since $s(w_{j_{0}})=-$,  by Lemma 2.2.4
and Lemma 4.1.3, $j_{0}=m^{l+1}\notin I(w_{j},l)$. Hence
$I(w_{j},l+1)=I(w_{j},l)$. By Lemma 4.6.3, $C_{j}$ is obtained by
doing band sums with copies of $\partial W_{j_{0}}^{l}$  to
$\partial W_{j}^{l}$. Hence $C_{j}$ bounds a disk in $\cal W_{-}$,
denoted by $W_{j}^{l+1}$.

Now by Proposition 6, $(\partial V_{i}-c_{i}^{l}\cup_{r\in
I(v_{i},l)} c_{r}^{l})\cap (\partial W_{j}^{l}-d_{j}^{l}\cup_{r\in
I(w_{j},l)}d_{r}^{l})=\emptyset$ for $i\geq l+2$ with $s(v_{i})=-$
and $j\notin m(l)$ with $s(w_{j})=-$. Now there are two case:

(1) \ $l+1\notin I(v_{i},l)$.

Now $I(v_{i},l+1)=I(v_{i},l)$. Hence $V_{i}^{l+1}\cap
W_{j}^{l+1}=(c_{i}^{l+1}\cup_{r\in I(v_{i},l+1)} c_{r}^{l+1})\cap
(d_{j}^{l+1}\cup_{r\in I(w_{j},l+1)} d_{r}^{l+1})$.

(2) \  $l+1\in I(v_{i},l)$.

Now $I(v_{i},l+1)=I(v_{i},l)-\bigl\{l+1\bigr\}$. By Proposition
6(2), $c_{l+1}^{l}$ is disjoint from
$W_{j}^{l}-d_{j}^{l}\cup_{r\in I(w_{j},l)} d_{r}^{l}$. Obviously,
$c_{l+1}^{l}$ is disjoint from $d_{j}^{l+1}$ and $d_{r}^{l+1}$.
Note that $(\partial V_{i}^{l}-c_{i}^{l}\cup_{r\in I(v_{i},l)}
c_{r}^{l})\cup c_{l+1}^{l}=\partial
V_{i}^{l+1}-c_{i}^{l+1}\cup_{r\in I(v_{i},l+1)} c_{r}^{l+1}$.
Hence $V_{i}^{l+1}\cap W_{j}^{l+1}=(c_{i}^{l+1}\cup_{r\in
I(v_{i},l+1)} c_{r}^{l+1})\cap d_{j}^{l+1}\cup_{r\in I(w_{j},l+1)}
d_{r}^{l+1})$. \qquad Q.E.D \vskip 0.5mm

\section{The Proofs of Propositions 4-6 for the case: $k=l+1$ and $s(v_{l+1})=-$}

\ \ \ \ \  In this section,  we shall Propositions 4-6 for the
case: $k=l+1$ and $s(v_{l+1})=-$. In the following argument, we
assume $s(v_{l+1})=-$.

\subsection{The element of constructions}

\ \ \ \ \ Recalling $m^{l+1}=MinL(c_{l+1}^{l})$, and
$L(c_{l+1}^{l})=\bigl\{\ldots,j_{-1},j_{0}=m^{l+1},j_{1},\ldots\bigr\}$.

{\bf \em Lemma 5.1.1.} \ Suppose that $m^{l+1}\neq \emptyset$.
Then

(1) \ $V_{l+1}^{l}\subset \cal V_{-}$ such that $\partial
V_{l+1}^{l}-c_{l+1}^{l}$ is disjoint from $F^{l}$.

(2) \ $\partial V_{l+1}^{l}$ intersects $d_{j_{0}}^{l}$ in one
point $p=intc_{l+1}^{l}\cap intd_{0,0}$.

(3) \ If $j<j_{0}=m^{l+1}$, then $d_{j}^{l}$ is disjoint from
$\partial V_{l+1}^{l}$; if $\gamma<j_{0}=m^{l+1}$, then
$e^{l}_{\gamma}\times I$ is disjoint from $\partial V_{l+1}^{l}$.

(4) \ $s(w_{m^{l+1}})=+$.

{\bf \em Proof.} \ (1) \ Since $s(v_{l+1})=-$,  by Proposition 6,
$V_{l+1}^{l}\subset \cal V_{-}$ and $\partial V_{l+1}^{l}\cap
F^{l}=c_{l+1}^{l}\cup_{r\in I(v_{l+1},l)} c_{r}^{l}$. If $r\in
I(v_{l+1},l)$, then, by Lemma 2.2.4, $r<l+1$. By Definition 2.3.1,
$I(v_{l+1},l)=\emptyset$. Thus (1) holds.

(2) \ By (1),  $\partial V_{l+1}^{l}\cap
d_{j_{0}}^{l}=c_{l+1}^{l}\cap d_{j_{0}}^{l}$. By the proof of
Lemma 4.1.2, $c_{l+1}^{l}$ intersects $d_{j_{0}}^{l}$ in one point
lying in $intd_{0,0}\cap intc_{l+1}^{l}$.

(3) \ Since $d_{j}^{l},e^{l}_{\gamma}\times I\subset F^{l}$. Hence
$d_{j}^{l}\cap \partial V_{l+1}^{l}=d_{j}^{l}\cap c_{l+1}^{l}$ and
$e^{l}_{\gamma}\times I\cap\partial
V_{l+1}^{l}=e^{l}_{\gamma}\times I\cap c_{l+1}^{l}$. By the
minimality of $m^{l+1}$ in $L(c_{l+1}^{l})$, (3) follows from the
proof of Lemma 4.1.2.

(4) \ Suppose that $s(w_{m^{l+1}})=-$. Then, by Proposition 6,
$W_{m^{l+1}}^{l}\subset \cal W_{-}$. By (1), $W_{m^{l+1}}^{l}\cap
V_{l+1}^{l}=c_{l+1}^{l}\cap (d_{m^{l+1}}^{l}\cup_{r\in
I(w_{m^{l+1}},l)} d_{r}^{l})$. Since $r\in I(w_{m^{l+1}},l)$, by
Lemma 2.2.4, $r<j_{0}=m^{l+1}$. By (2) and (3), $V_{l+1}^{l}$
intersects $W_{m^{l+1}}^{l}$ in one point. Hence $\cal V_{-}\cup
W_{-}$ is stabilized, a contradiction. \qquad Q.E.D.\vskip 0.5mm

{\bf \em Lemma 5.1.2.} \ Suppose that
$j_{0}=m^{l+1}\neq\emptyset$. If $j_{0}\in I(w_{\gamma_{0}},l)$
for some $\gamma_{0}\in m(l)$ with $s(w_{\gamma_{0}})=-$, then
there is a neighborhood of $d_{j_{0}}^{l}$ in $F^{l}$, say
$d_{j_{0}}^{l}\times [-2.5,2.5]$, which satisfies the following
conditions:

(1) \ $(d_{j_{0}}^{l}\times I)_{\gamma_{0}}=d_{j_{0}}^{l}\times
[-1,1]$.

(2) \ If $\gamma\leq Max\bigl\{ \gamma_{i,0} \ | \ \delta(0)\leq
i\leq \theta(0)\bigr\}$, then each component of
$d_{j_{0}}^{l}\times [-2.5,2.5]\cap e^{l}_{\gamma}\times I$ is
$c\times [-2.5,2.5]\subset e^{l}_{\gamma}\times (0,1)$ where
$c\subset intd_{j_{0}}^{l}$ is a core of $e^{l}_{\gamma}\times
(0,1)$.

(3) \ If $\gamma_{0}>\gamma >Max\bigl\{ \gamma_{i,0} \ | \
\delta(0)\leq i\leq \theta(0)\bigr\}$, then $d_{j_{0}}^{l}\times
[-2.5,2.5]$ is disjoint $e^{l}_{\gamma}\times I$.

(4) \ If $\gamma> \gamma_{0}$, then $e^{l}_{\gamma}\times I$ is
disjoint from $d_{j_{0}}^{l}\times ([1,2.5]\cup [-2.5,-1])$.

(5) \ If $j\neq j_{0}<\gamma_{0}$, then $d_{j}^{l}$ is disjoint
from $d_{j_{0}}^{l}\times [-2.5,2.5]$, if $j>\gamma_{0}$, then
$d_{j}^{l}$ is disjoint from $d_{j_{0}}^{l}\times ([1,2.5]\cup
[-2.5,-1])$.

(6) \ $intd_{i,0}\times [-2.5,-1.5]\cup [1.5,2.5]$ is disjoint
from $e^{l}_{\gamma}\times I$ for each $\gamma\in m(l)$ and
$\delta(0)\leq i\leq \theta(0)$, $intd_{i,0}\times [-2.5,2.5]$ is
disjoint from $e^{l}_{\gamma}\times I$ for $\gamma<j_{0}$.

(7) \ $d_{j_{0}}^{l}\times [-2.5,2.5]$ is disjoint from
$(a_{\alpha}\times I)_{\gamma_{\alpha}}$ for $\alpha\neq 0$, and
$d_{j}^{l}\times [-2.5,2.5]\cap
D^{*}_{0,\alpha}=a^{0}_{\alpha}\times [-2.5,2.5]\cap
D^{*}_{0,\alpha}$.

(8) \ For each $i\geq l+1$, $\partial c_{i}^{l}$ is disjoint from
$d_{j_{0}}^{l}\times [-2.5,2.5]$.

{\bf \em Proof.} \ (1) is trivial.

(2) \ Since $\gamma<j_{0}$, by  Proposition 4(2) and Definition
2.1.5, each component of $d_{j_{0}}^{l}\cap e^{l}_{\gamma}\times
I$ is a core of $e^{l}_{\gamma}\times I$, say $c$, which is a core
of $e^{l}_{\gamma}\times (0,1)$. Note that $d_{j_{0}}^{l}\subset
e^{l}_{\gamma}$. By Proposition 4(1) and Definition 2.1.4, each
component of $(d_{j_{0}}^{l}\times I)_{\gamma_{0}} \cap
e^{l}_{\gamma}\times I$ is $(c\times I)_{\gamma_{0}}\subset
e^{l}_{\gamma}\times (0,1)$.

(3) By Lemma 3.2.3 and Lemma 3.2.6, $(d_{j_{0}}^{l}\times
I)_{\gamma_{0}}$ is disjoint from $e^{l}_{\gamma}\times I$ for
$\gamma_{0}>\gamma >Max\bigl\{ \gamma_{i,0} \ | \ \delta(0)\leq
i\leq \theta(0)\bigr\}$.

(4) \ If $\gamma> \gamma_{0}\in m(l)$, then, by Proposition 4(1),
each component of $e^{l}_{\gamma}\times I\cap (d_{j_{0}}^{l}\times
I)_{\gamma_{0}}=(c\times I)_{\gamma}\subset (d_{j_{0}}^{l}\times
(0,1))_{\gamma_{0}}$.

(5) \ If $j\neq j_{0}<\gamma_{0}$, then, by Proposition 4(3),
$d_{j}^{l}$ is disjoint from $(d_{j_{0}}^{l}\times
I)_{\gamma_{0}}$. By Proposition 4(2) and Definition 2.1.5, if
$j>\gamma_{0}$, then each component of $d_{j}^{l}\cap
(d_{j_{0}}^{l}\times I)_{\gamma_{0}}$ is $c\subset
(d_{j_{0}}^{l}\times (0,1))_{\gamma_{0}}$.

(6) \ If $j_{0}<\gamma$, then either $\gamma=\gamma_{0}$ or
$(d_{j_{0}}^{l}\times I)_{\gamma_{0}}$ is disjoint from
$e^{l}_{\gamma}\times I$ by Lemma 3.2.6. If $j_{0}>\gamma$, then,
by Lemma 3.3.1 and Definition 2.1.4, $(intd_{i,0}\times
I)_{\gamma_{0}}$ is disjoint from $e^{l}_{\gamma}\times I$. If
$\gamma>\gamma_{0}$, then each component of $e^{l}_{\gamma}\times
I\cap (d_{j_{0}}^{l}\times I)_{\gamma_{0}}=(c\times
I)_{\gamma}\subset (d_{j_{0}}^{l}\times (0,1))_{\gamma_{0}}$.

(7)  By Lemma 3.2.6 and Lemma 3.5.7, $(d_{j_{0}}^{l}\times
I)_{\gamma_{0}}$ is disjoint from $(a_{\alpha}\times
I)_{\gamma_{\alpha}}$ for $\alpha\neq 0$.

(8) \ By Proposition 4(5), $\partial c_{i}^{l}$ is disjoint from
$(d_{j_{0}}^{l}\times I)_{\gamma_{0}}\subset
e^{l}_{\gamma_{0}}\times I$.

Now $d_{j_{0}}^{l}\times [-2,2]$ can be obtained by making
$(d_{j_{0}}^{l}\times I)_{\gamma_{0}}$ wide slightly. Hence Lemma
5.1.2 holds. \qquad Q.E.D. \vskip 0.5mm

{\bf \em Lemma 5.1.3.} \ Suppose $j_{0}\neq\emptyset$. If
$j_{0}\notin I(w_{\gamma},l)$ for each $\gamma\in m(l)$ with
$s(w_{\gamma})=-$, there is a neighborhood of $d_{j_{0}}^{l}$ in
$F^{l}$, say $d_{j_{0}}^{l}\times [-2.5,2.5]$, which satisfies the
following conditions:

1) \ If $j\neq j_{0}$, then $d_{j_{0}}^{l}\times [-2.5,2.5]$ is
disjoint from $d_{j}^{l}$.

2) \ If $\gamma >Max\bigl\{ \gamma_{i,0} \ | \ \delta(0)\leq i\leq
\theta(0)\bigr\}$, then $d_{j_{0}}^{l}\times [-2.5,2.5]$ is
disjoint from $e^{l}_{\gamma}\times I$.

3) \ If $\gamma \leq Max\bigl\{ \gamma_{i,0} \ | \ \delta(0)\leq
i\leq \theta(0)\bigr\}$, then each component of
$d_{j_{0}}^{l}\times [-2.5,2.5]\cap e^{l}_{\gamma}\times I$ is
$c\times [-2.5,2.5]\subset e^{l}_{\gamma}\times (0,1)$ where
$c\subset intd_{j_{0}}^{l}$ is a core of $e^{l}_{\gamma}\times
(0,1)$.

4) \ $d_{i,0}\times [-2.5,2.5]$ is disjoint from
$\textrm{int}e^{l}_{\gamma}\times I$ for each $\gamma\in m(l)$.

(5) \ $d_{i,0}\times [-2.5,2.5]$ is disjoint from $\partial
c_{i}^{l}$ for each $i\geq l+1$.

(6) \ $d_{j_{0}}^{l}\times [-2.5,2.5]$ is disjoint from
$(a_{\alpha}\times I)_{\gamma_{\alpha}}$ for $\alpha$, and
$d_{j}^{l}\times [-2.5,2.5]\cap
D^{*}_{0,\alpha}=a^{0}_{\alpha}\times [-2.5,2.5]\cap
D^{*}_{0,\alpha}$.

{\bf \em Proof.} \ By Lemma 5.1.1(4), $s(w_{m^{l+1}})=+$. Hence if
$j_{0}=m^{l+1}\in I(w_{\gamma},l)$ then $s(w_{\gamma})=-$. By
assumption, $j_{0}\notin I(w_{\gamma},l)$ for each $\gamma\in
m(l)$. Now (1) follows from Proposition 4. (2) follows from Lemma
3.2.3. (3) follows Proposition 4(2) and Definition 2.1.5. (4)
follows from (2) and Lemma 3.3.1. (5) follows from Proposition
4(5). (6) follows from Lemma 3.5.7.\qquad Q.E.D.\vskip 0.5 mm

{\bf \em Lemma 5.1.4.} \ Suppose that $j_{0}\neq\emptyset$. Then
there is a neighborhood of $\partial V_{l+1}^{l}$ in
$\partial_{+}\cal V_{-}$, say $\partial V_{l+1}^{l}\times I$,
satisfying the following conditions:

(1) \ $\partial V_{l+1}^{l}\times I\cap F^{l}=c_{l+1}^{l}\times
I$.

(2) \ If $j, \gamma<j_{0}$, then $d_{j}^{l}$ and
$e^{l}_{\gamma}\times I$ are disjoint from $c_{l+1}^{l}\times I$.

(3) \ $\partial d_{j}^{l}, (\partial e^{l}_{\gamma})\times I$ are
disjoint from $c_{l+1}^{l}\times I$ for each $j\notin m(l)$ and
$\gamma\in m(l)$.

(4) \ $d_{j_{0}}^{l}$ intersects $c_{l+1}^{l}\times I$ in an arc
$a\subset intd_{0,0}$, and $d_{j_{0}}^{l}\times [-2.5,2.5]$
intersects $c_{l+1}^{l}\times I$ in $a\times [-2.5,2.5]$. See
Figure 19.

(5) \ For $i\geq l+2$, $c_{i}^{l}$ is disjoint from
$c_{l+1}^{l}\times I$, if $s(v_{i})=-$, then $V_{i}^{l}$ is
disjoint from $\partial V_{l+1}^{l}\times I$.

(6) \ For each $j\notin m(l)$ with $s(w_{j})=-$,
$W_{j}^{l}-d_{j}^{l}\cup_{r\in I(w_{j},l)} d_{r}^{l}$ is disjoint
from $\partial V_{l+1}^{l}\times I$.

{\bf \em Proof.} \ (1) follows from Lemma 5.1.1(1). (2) follows
from Lemma 5.1.1(3) and (4). (3) follows from Proposition 4(5).
(4) follows from Lemma 5.1.1(2). (5) and (6) follow from
Proposition 6.\qquad Q.E.D.
\begin{center}
\includegraphics[totalheight=4cm]{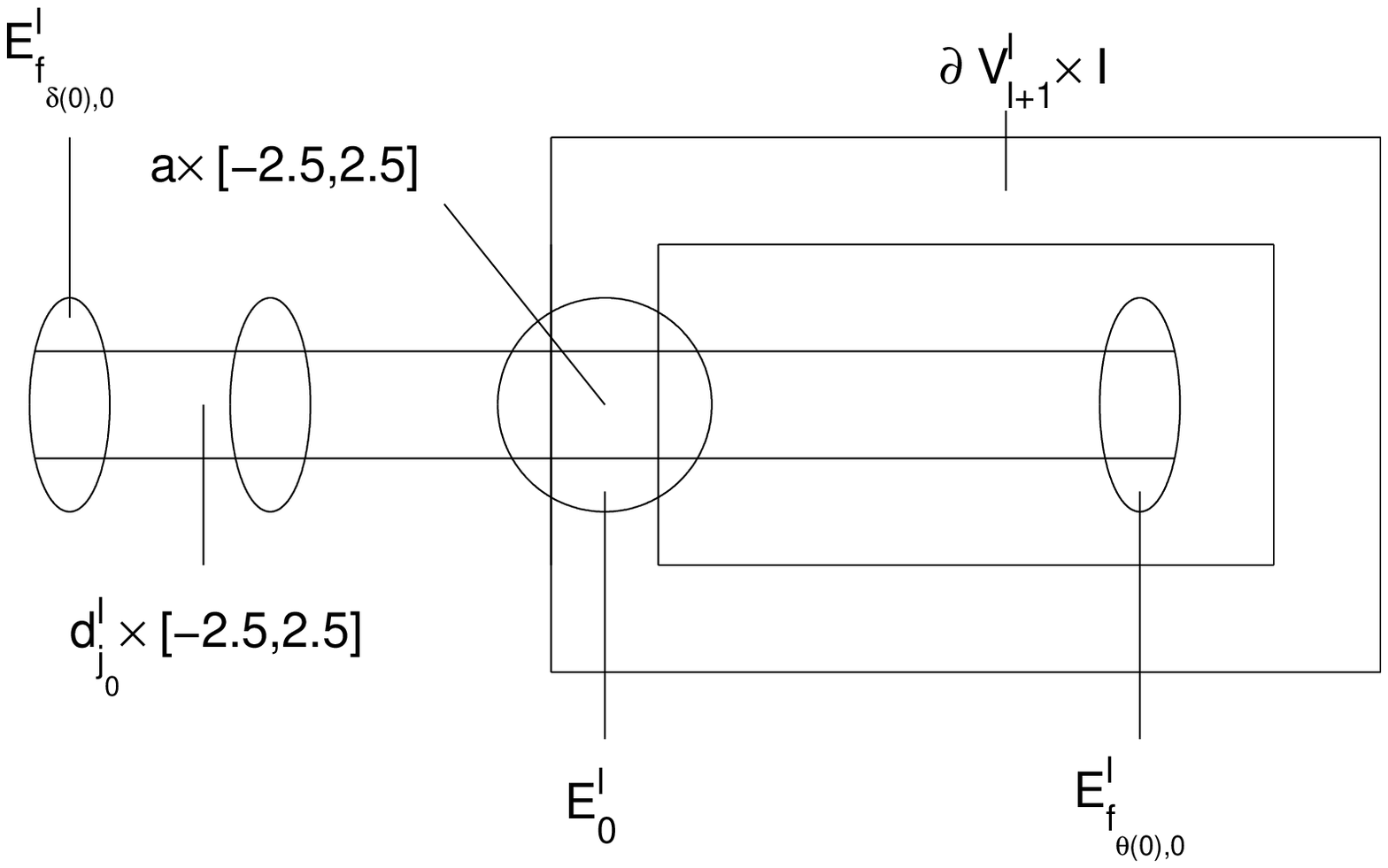}
\begin{center}
Figure 19
\end{center}
\end{center}

{\bf \em Lemma 5.1.5.} \ If $L(c_{l+1}^{l})=\emptyset$, then
Propositions 4-6 hold for the case: $k=l+1$ and $s(v_{l+1})=-$.

{\bf \em Proof.} \ By Corollary 3.2.5, $d_{j}^{l},
e^{l}_{\gamma}\times I$ are disjoint from $c_{l+1}^{l}\times I$
where $c_{l+1}^{l}\times I$ is a regular neighborhood of
$c_{l+1}^{l}$ in $E^{l}_{0}$. Let $F^{l+1}=F^{l}-c_{l+1}^{l}\times
(-1,1)$. We denote by $d_{j}^{l+1}$ the arc $d_{j}^{l}$ for $j\in
\bigl\{1,\ldots,n\bigr\}-m(l)$, $e^{l+1}_{\gamma}\times I$ the
disk $e^{l}_{\gamma}\times I$ for $\gamma\in m(l)$, $E^{l+1}_{f}$
the disk $E^{l}_{f}$ for $1\leq f\leq l$. In particular, we denote
by $E^{l+1}_{0}, E^{l+1}_{l+1}$ the two components of
$E^{l}_{0}-c_{l+1}^{l}\times (-1,1)$. By Proposition 4,
$c_{i}^{l}\cap c_{l+1}^{l}=\emptyset$ for $i\geq l+2$.  Now we
denote by $c^{l+1}_{i}$ the arc $c^{l}_{i}$. Let $m(l+1)=m(l)$.
For $i\geq l+2$ with $s(v_{i})=-$ and $j\in
\bigl\{1,\ldots,n\bigr\}-m(l)$, we denote by $V_{i}^{l+1}$ the
disk $V^{l}_{i}$, $W^{l+1}_{j}$ the disk $W^{l}_{j}$. By the
argument in Section 4.2, the lemma holds.\qquad Q.E.D.

\subsection{The constructions of $c^{l+1}_{i}, d_{j}^{'},
e^{'}_{\gamma}\times I, F^{l+1}$ (I).}

\ \ \ \ \ By Lemma 5.1.5, in the following argument, we may assume
that $L(c_{l+1}^{l})\neq\emptyset$. Hence
$j_{0}=m^{l+1}\neq\emptyset$. In Sections 5.2 and 5.3, we shall
construct $c^{l+1}_{i}, d_{j}^{'}, e^{'}_{\gamma}\times I,
F^{l+1}$ from $d_{j}^{l},c_{i}^{l}, e^{l}_{\gamma}\times I,
F^{l}$, where $c_{i}^{l+1}, F^{l+1}$ are just $c_{i}^{l+1},
F^{l+1}$ in Propositions 4-6, but $d_{j}^{'}, e^{'}_{\gamma}\times
I$ are not $d_{j}^{l+1},e^{l+1}_{\gamma}\times I$ in Propositions
4-6. In this section, we first assume that $j_{0}\notin
I(w_{\gamma},l)$ for each $\gamma\in m(l)$ with $s(w_{\gamma})=-$.

Let $\partial V_{l+1}^{l}\times I$ be a neighborhood of $\partial
V_{l+1}^{l}$ in $\partial_{+} \cal V_{-}$ satisfying  Lemma 5.1.4,
and $d_{j_{0}}^{l}\times [-2.5,2.5]$ be a neighborhood of
$d_{j_{0}}^{l}$ in $F^{l}$ satisfying Lemma 5.1.3.\vskip 0.5mm

{\bf \em Definition 5.2.1.}  \ We denote by $E_{0}, E_{l+1}$ the
two components of $E^{l}_{0}-c_{l+1}^{l}\times (-1,1)$.
\begin{center}
\includegraphics[totalheight=4cm]{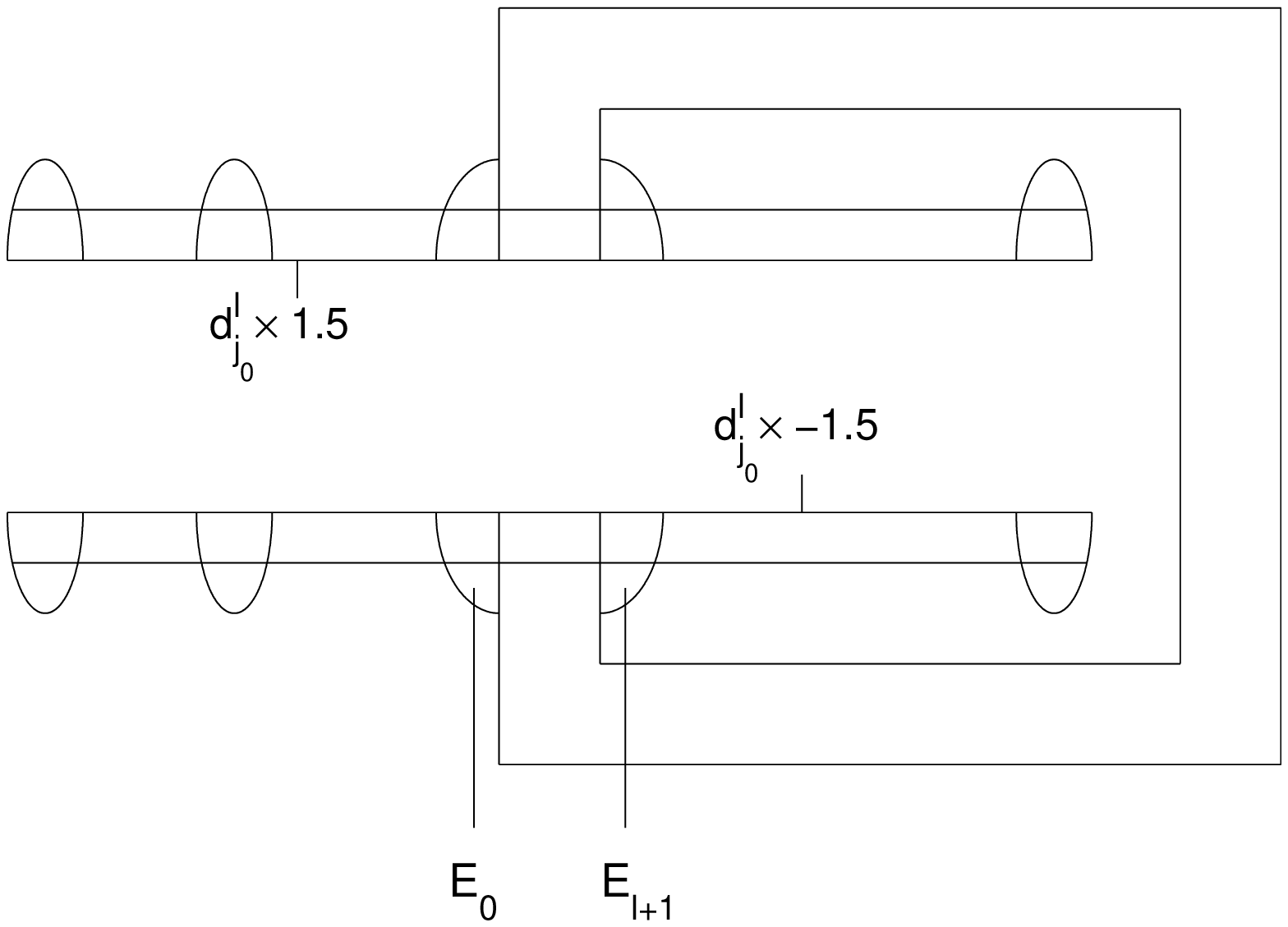}
\begin{center}
Figure 20
\end{center}
\end{center}
\begin{center}
\includegraphics[totalheight=6cm]{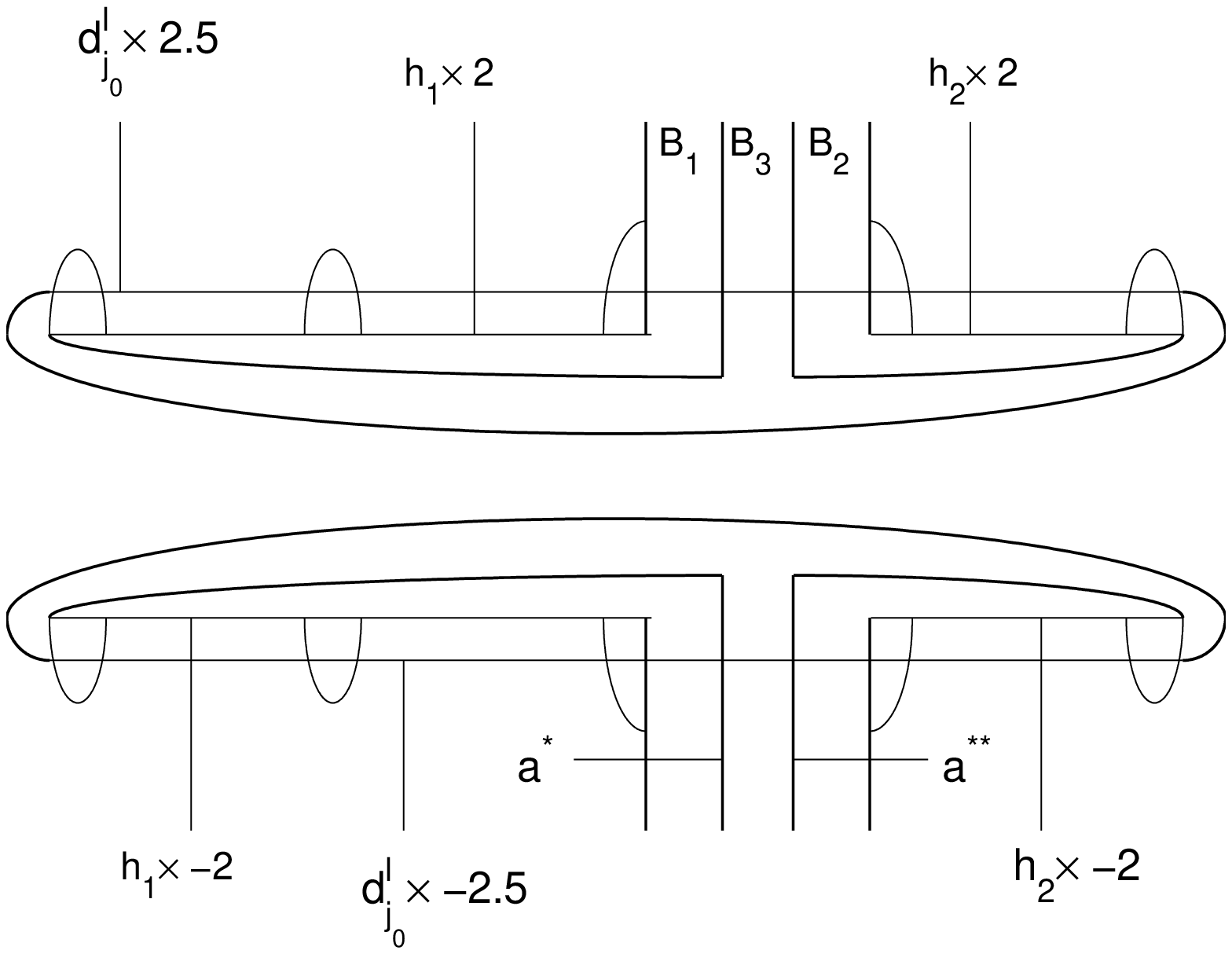}
\begin{center}
Figure 21
\end{center}
\end{center}

{\bf \em Definition 5.2.2.}  \ (1) \ Let $F^{*}=F^{l}\cup
\partial V_{l+1}^{l}\times I-d_{j_{0}}^{l}\times (-1.5,1.5)$. See
Figure 20.

(2) \ Let $F^{l+1}$ be the surface obtained by identifying
$\partial_{1} d_{j_{0}}^{l}\times [1.5,2]$ to $\partial_{1}
d_{j_{0}}^{l}\times [2,2.5]$ and $\partial_{1} d_{j_{0}}^{l}\times
[-2,-1.5]$ to $\partial_{1} d_{j_{0}}^{l}\times [-2.5,-2]$ in
$\partial_{+} \cal V_{-}$$-F^{l+1}$. See Figure 21.

Let $b=\partial V_{l+1}^{l}-a\times (-2,2)$. Then $b\times
I=(\partial V_{l+1}^{l}-a\times (-2,2))\times I$. Now by Lemma
5.1.4(4), $B=(d_{j_{0}}^{l}\times [-2,-1.5]\cup [1.5,2])\cup
b\times I$ is a disk.
%We may assume that $\partial_{1} b\subset
%d_{j_{0}}^{l}\times\bigl\{-2\bigr\}$ and $\partial_{2} b\subset
%d_{j_{0}}^{l}\times\bigl\{2\bigr\}$.
We denote by $h_{1},h_{2}$ the  two components of
$d_{j_{0}}^{l}-inta$ where $a$ is as in Lemma 5.1.4. Now there is
an arc $a^{*}$ in $B=(d_{j_{0}}^{l}\times [-2,-1.8]\cup
[1.8,2])\cup b\times I$ connecting $\partial_{1}
d_{j_{0}}^{l}\times \bigl\{2\bigr\}$ to $\partial_{1}
d_{j_{0}}^{l}\times \bigl\{-2\bigr\}$. Furthermore, $a^{*}$
intersects $b\times I$ in $b\times \bigl\{-1/2\bigr\}$. Similarly,
there is an arc $a^{**}$ in $(d_{j_{0}}^{l}\times [-2,-1.8]\cup
[1.8,2])\cup b\times I$ connecting $\partial_{2}
d_{j_{0}}^{l}\times \bigl\{2\bigr\}$ to $\partial_{2}
d_{j_{0}}^{l}\times \bigl\{-2\bigr\}$. Furthermore, $a^{**}$
intersects $b\times I$ in $b\times \bigl\{1/2\bigr\}$. Now
$a^{*}\cup a^{**}$ separates $B=(d_{j_{0}}^{l}\times [-2,-1.5]\cup
[1.5,2])\cup b\times I$ into three disks $B_{1},B_{2},B_{3}$ as in
Figure 21. Without loss of generality, we may assume that
$h_{1}\times\bigl\{-2,2\bigr\}\subset B_{1}$ and
$h_{2}\times\bigl\{-2,2\bigr\}\subset B_{2}$ as in Figure 21. Note
that $B_{1}\cup B_{2}\subset (d_{j_{0}}^{l}\times [-2,-1.8]\cup
[1.8,2])\cup b\times I$.

Now there is a homeomorphism $H_{i}$ form $h_{i}\times
[-2,2]\subset d_{j_{0}}^{l}\times [-2,2]$ to $B_{i}$ such that
$H_{i}$ is an identifying map on $h_{i}\times\bigl\{-2,2\bigr\}$
for $i=1,2$. Let $H=H_{1}\cup H_{2}$.
\begin{center}
\includegraphics[totalheight=6cm]{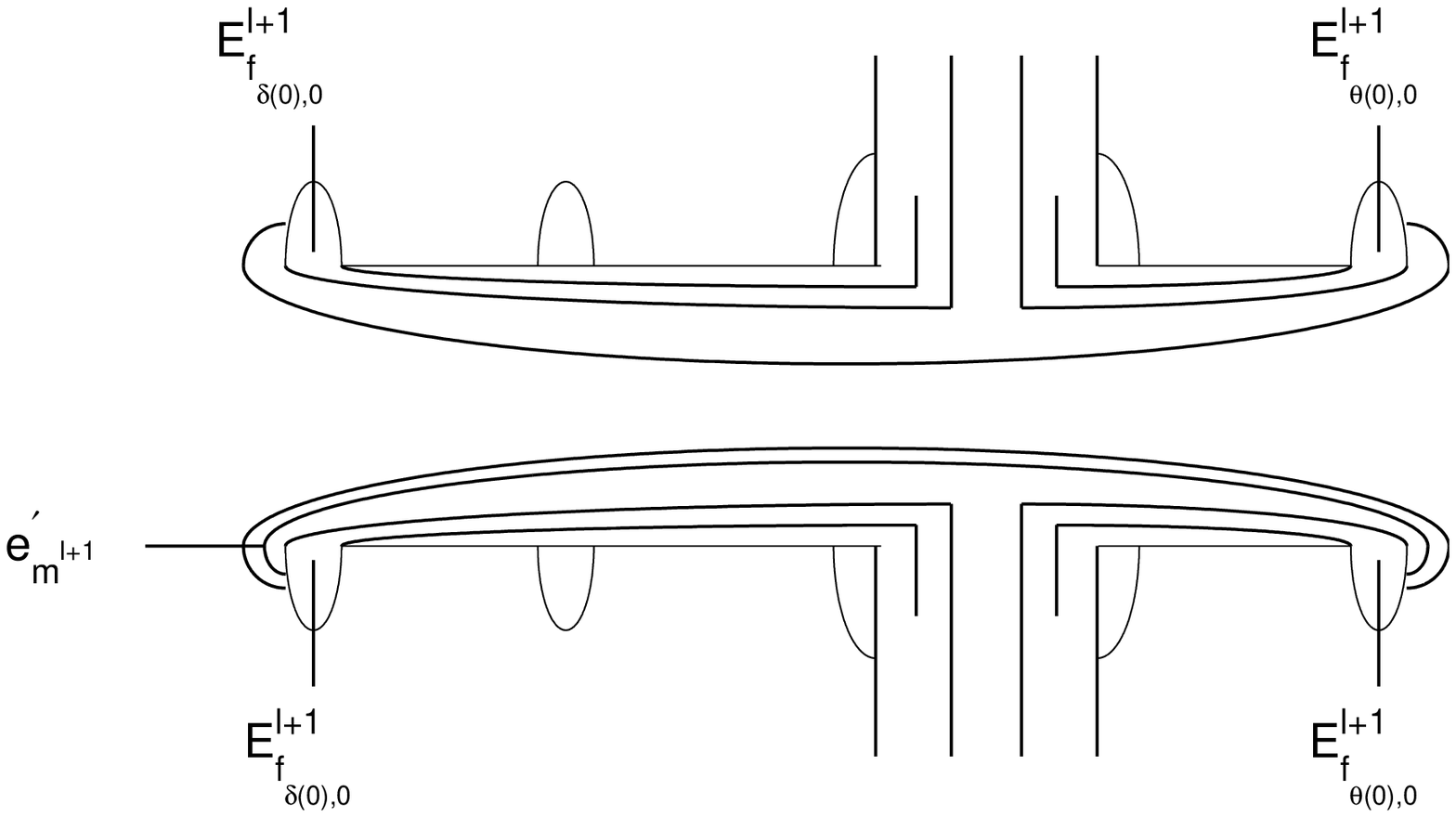}
\begin{center}
Figure 22
\end{center}
\end{center}

By Lemma 3.3.1 and Remark 3.3.2,
$d_{j_{0}}^{l}=\cup_{i=\delta(0)}^{\theta(0)}
d_{i,0}\cup_{i=\delta(0)}^{\theta(0)}e_{i,0}$.

{\bf \em Definition 5.2.3.} \ (1) \ For each $\gamma\in m(l)$, let
$e^{'}_{\gamma}=e^{l}_{\gamma}$, $e^{'}_{\gamma}\times
I=(e^{l}_{\gamma}\times I-e^{l}_{\gamma}\times I\cap (h_{1}\cup
h_{2})\times [-2,2])\cup H(e^{l}_{\gamma}\times I\cap (h_{1}\cup
h_{2})\times [-2,2])$.

(2) \ For $j\notin m(l+1)$, let $d_{j}^{'}=d_{j}^{l}$.

(3) \ For $i\geq l+2$, let $c^{l+1}_{i}=(c_{i}^{l}-c_{i}^{l}\cap
(h_{1}\cup h_{2})\times [-2,2])\cup H(c_{i}^{l}\cap (h_{1}\cup
h_{2})\times [-2,2])$.

(4) \ For $1\leq f\leq l$, let
$E^{l+1}_{f}=(E^{l}_{f}-E^{l}_{f}\cap (h_{1}\cup h_{2})\times
[-2,2])\cup H(E^{l}_{f}\cap (h_{1}\cup h_{2})\times [-2,2])$, let
$E^{l+1}_{0}=(E_{0}-E_{0}\cap (h_{1}\cup h_{2})\times [-2,2])\cup
H(E_{0}\cap (h_{1}\cup h_{2})\times [-2,2])$, let
$E^{l+1}_{l+1}=(E_{l+1}-E_{l+1}\cap (h_{1}\cup h_{2})\times
[-2,2])\cup H(E_{l+1}\cap (h_{1}\cup h_{2})\times [-2,2])$.

(5) \ Let
$e^{l+1}_{m^{l+1}}=d_{j_{0}}^{l}\times\bigl\{-1.6\bigr\}$, and
$e^{l+1}_{m^{l+1}}\times I=B_{3}$ with $e^{l+1}_{m^{l+1}}\times
[-1,0]=d_{j_{0}}^{l}\times [-1.5,-1.6]$. See Figure 22.

\subsection{The constructions of $c^{l+1}_{i}, d_{j}^{'},
e^{'}_{\gamma}\times I, F^{l+1}$ (II).}

\ \ \ \ \ In this section, we shall construct $c^{l+1}_{i},
d_{j}^{'}, e^{'}_{\gamma}\times I, F^{l+1}$ from
$d_{j}^{l},c_{i}^{l}, e^{l}_{\gamma}\times I, F^{l}$ for the case:
 $j_{0}\in I(w_{\gamma_{0}},l)$ for
some $\gamma_{0}\in m(l)$ with $s(w_{\gamma_{0}})=-$.

Let $\partial V_{l+1}^{l}\times I$ be a neighborhood of $\partial
V_{l+1}^{l}$ in $\partial_{+} \cal V_{-}$ satisfying  Lemma 5.1.4,
and $d_{j_{0}}^{l}\times [-2,2]$ be a neighborhood of
$d_{j_{0}}^{l}$ in $F^{l}$ satisfying Lemma 5.1.2.

Recalling the disks $E_{0}, E_{l+1}$, $B_{1}, B_{2}, B_{3}$ and
the arcs $h_{1}, h_{2}, a^{*}, a^{**}$ in Section 5.2.
\begin{center}
\includegraphics[totalheight=6cm]{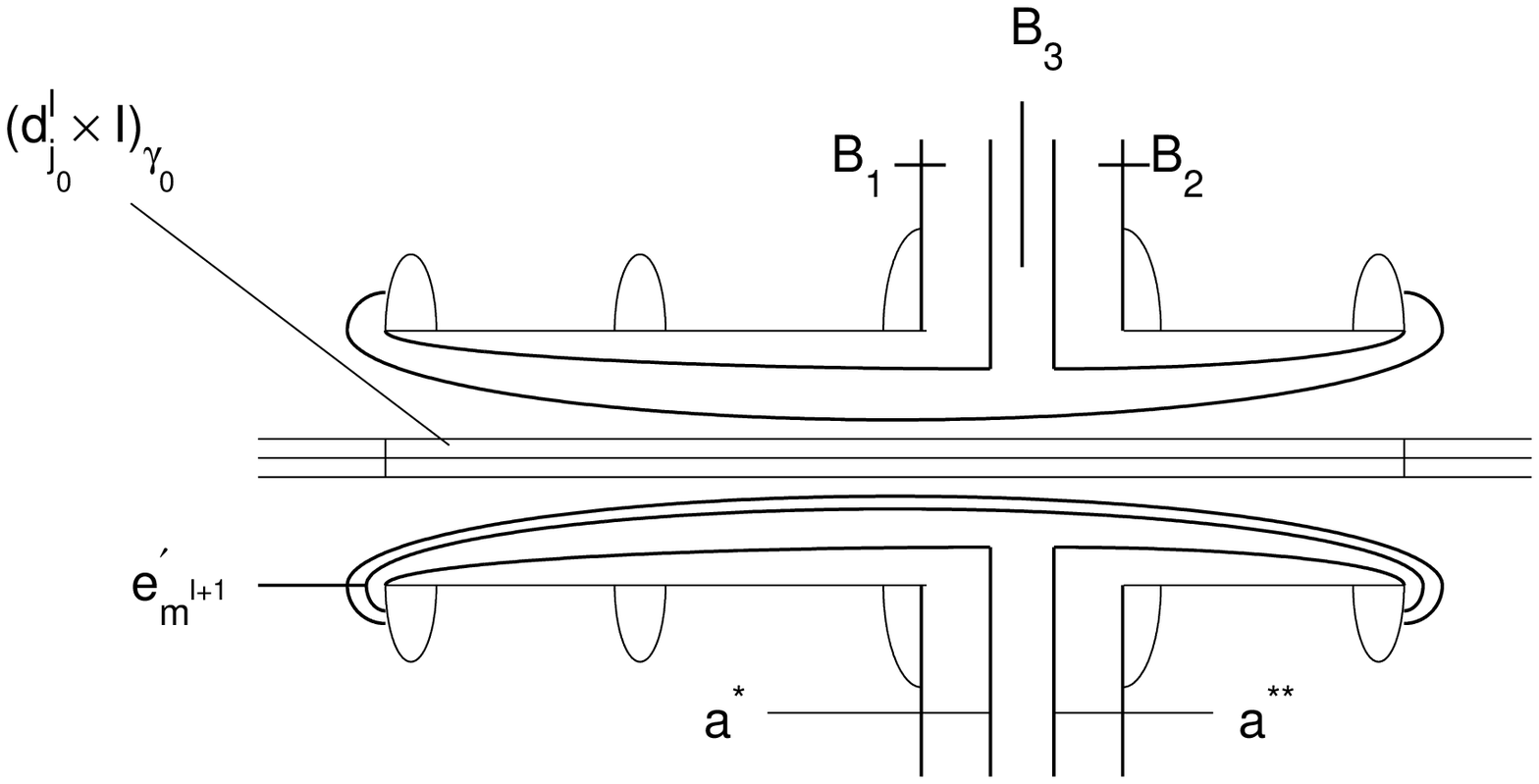}
\begin{center}
Figure 23
\end{center}
\end{center}

{\bf \em Definition 5.3.1.}  \ (1) \ Let $F^{*}=(F^{l}\cup
\partial V_{l+1}^{l}\times I-d_{j_{0}}^{l}\times (-1.5,1.5))\cup (d_{j_{0}}^{l}\times I)_{\gamma_{0}}$.

(2) \ Let $F^{l+1}$ be the surface obtained by identifying
$\partial_{1} d_{j_{0}}^{l}\times [1.5,2]$ to $\partial_{1}
d_{j_{0}}^{l}\times [2,2.5]$ and $\partial_{1} d_{j_{0}}^{l}\times
[-2,-1.5]$ to $\partial_{1} d_{j_{0}}^{l}\times [-2.5,-2]$ in
$\partial_{+} \cal V_{-}$$-F^{l+1}$. See Figure 23.

{\bf \em Definition 5.3.2.} \ (1) \ For each $\gamma\in m(l)$ with
$\gamma<m^{l+1}=j_{0}$, let $e^{'}_{\gamma}=e^{l}_{\gamma}$,
$e^{'}_{\gamma}\times I=(e^{l}_{\gamma}\times
I-e^{l}_{\gamma}\times I\cap (h_{1}\cup h_{2})\times [-2,2])\cup
H(e^{l}_{\gamma}\times I\cap (h_{1}\cup h_{2})\times [-2,2])$.

(2) \ For each $\gamma\in m(l)$ with $\gamma>m^{l+1}=j_{0}$, let
$e^{'}_{\gamma}=e^{l}_{\gamma}$, $e^{'}_{\gamma}\times
I=e^{l}_{\gamma}\times I$.

(3) \ Let
$e^{l+1}_{m^{l+1}}=d_{j_{0}}^{l}\times\bigl\{-1.6\bigr\}$, and
$e^{l+1}_{m^{l+1}}\times I=B_{3}$ with $e^{l+1}_{m^{l+1}}\times
[-1,0]=d_{j_{0}}^{l}\times [-1.5,-1.6]$. See Figure 24.

{\bf \em Definition 5.3.3.} \ (1) \ For $j\notin m(l+1)$, let
$d_{j}^{'}=d_{j}^{l}$.

(2) \ For $i\geq l+2$, let $c^{l+1}_{i}=(c_{i}^{l}-c_{i}^{l}\cap
(h_{1}\cup h_{2})\times [-2,2])\cup H(c_{i}^{l}\cap (h_{1}\cup
h_{2})\times [-2,2])$.

(3) \ For $1\leq f\leq l$, let
$E^{l+1}_{f}=(E^{l}_{f}-E^{l}_{f}\cap (h_{1}\cup h_{2})\times
[-2,2])\cup H(E^{l}_{f}\cap (h_{1}\cup h_{2})\times [-2,2])$, let
$E^{l+1}_{0}=(E_{0}-E_{0}\cap (h_{1}\cup h_{2})\times [-2,2])\cup
H(E_{0}\cap (h_{1}\cup h_{2})\times [-2,2])$, let
$E^{l+1}_{l+1}=(E_{l+1}-E_{l+1}\cap (h_{1}\cup h_{2})\times
[-2,2])\cup H(E_{l+1}\cap (h_{1}\cup h_{2})\times [-2,2])$. See
Figures 23 and 24.

\begin{center}
\includegraphics[totalheight=6cm]{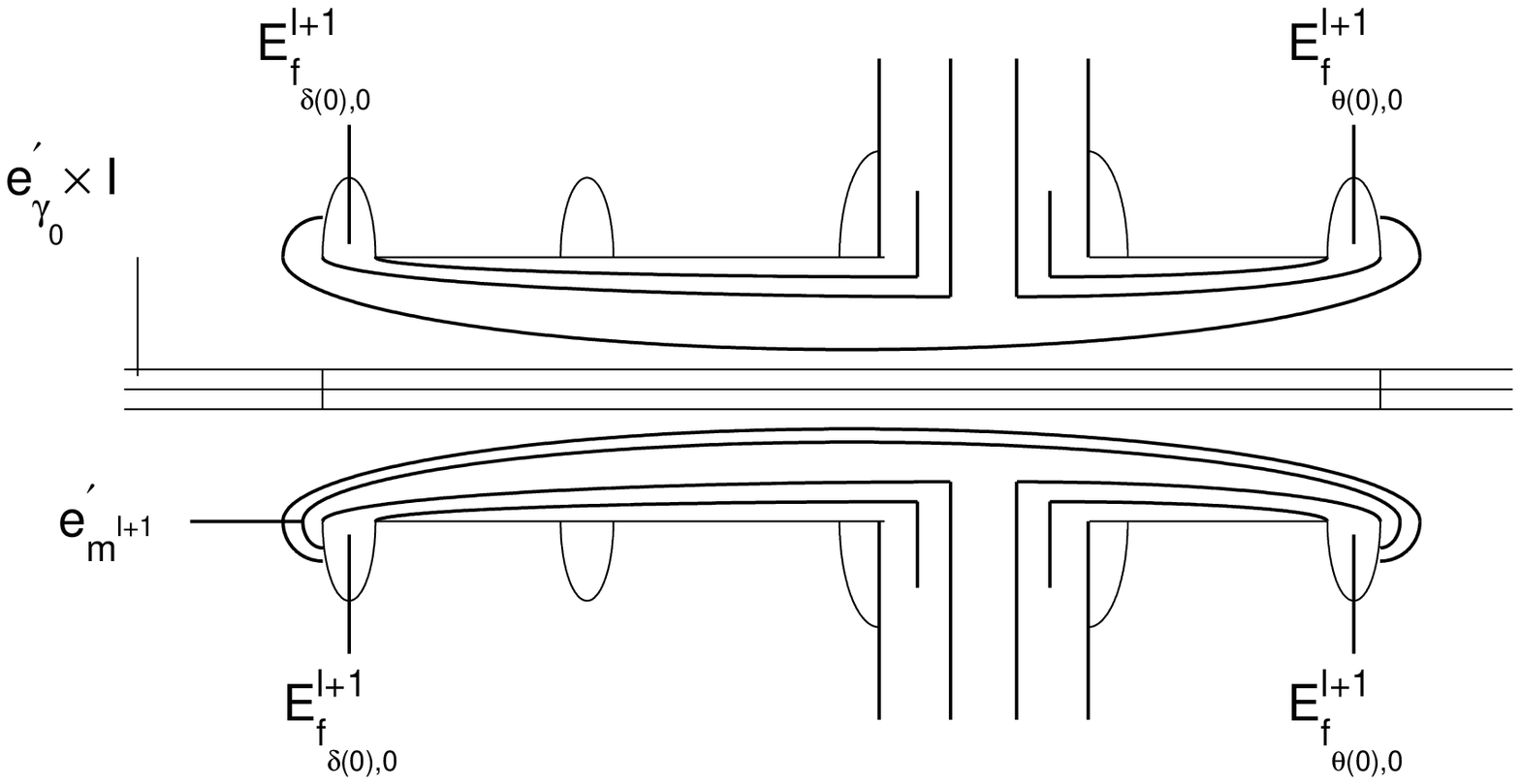}
\begin{center}
Figure 24
\end{center}
\end{center}

\subsection{Properties of $c^{l+1}_{i}, d_{j}^{'},
e^{'}_{\gamma}\times I, F^{l+1}$}

\ \ \ \ \  In this section, we shall introduce simple properties
of $c^{l+1}_{i}, d_{j}^{'}, e^{'}_{\gamma}\times I, F^{l+1}$.

{\bf \em Lemma 5.4.1.} \ $(d_{j_{0}}^{l}\times [-2,-1.5]\cup
[1.5,2])\cup b\times I\subset E^{l+1}_{0}\cup E^{l+1}_{l+1}
\cup_{i=\delta(0)}^{\theta(0)}
E^{l+1}_{f_{i,0}}\cup_{i=\delta(0)}^{\theta(0)}
e^{'}_{\gamma_{i,0}}\times I\cup e^{'}_{m^{l+1}}\times I$.

{\bf \em Proof.} \ By Lemma 3.3.1,
$d_{j_{0}}^{l}=\cup_{i=\delta(0)}^{\theta(0)}
d_{i,0}\cup_{i=\delta(0)}^{\theta(0)}e_{i,0}$, where $d_{i,0}$ is
an properly embedded arc in $E^{l}_{f_{i,0}}$, $e_{i,0}$ is a core
of $e^{l}_{\gamma_{i,0}}\times (0,1)$ for $i\neq 0$, and $d_{0,0}$
is an arc in $E^{l}_{0}$. By Lemma 5.1.4, $a\subset intd_{0,0}$.
We denote by $h^{*}, h^{**}$ the two components of $d_{0,0}-inta$.
We may assume that $h^{*}\subset E_{0}$ and $h^{**}\subset
E_{l+1}$. Then $h_{1}\times [-2,2]=(\cup_{i=\delta(0)}^{-1}
d_{i,0}\cup h^{*}\cup_{i=\delta(0)}^{-1} e_{i,0})\times [-2,2]$
and $h_{2}\times [-2,2]=(h^{**}\cup_{i=1}^{\theta(0)}
d_{i,0}\cup_{i=1}^{\theta(0)} e_{i,0})\times [-2,2]$. Note that
$(d_{j_{0}}^{l}\times [-2,-1.5]\cup [1.5,2])\cup b\times
I=B_{1}\cup B_{2}\cup B_{3}$. Now by Definitions 5.2.3, 5.3.2 and
5.3.3, the lemma holds.\qquad Q.E.D. \vskip 2mm

{\bf \em Lemma 5.4.2.} \ $\bigl\{d_{j}^{'} \ | \
j\in\bigl\{1,2,\ldots, n\bigr\}-m(l+1)\bigr\}$ is a set of
pairwise disjoint arcs in $F^{l+1}$. Furthermore, $d_{j}^{'}$ is
properly embedded in $\cup_{f} E^{l+1}_{f}\cup_{\gamma<j}
e^{'}_{\gamma}\times I$.

{\bf \em Proof.} \ Now there are two cases:

Case 1. \ $j_{0}\in I(w_{\gamma_{0}},l)$ for some $\gamma_{0}\in
m(l)$ with $s(w_{\gamma_{0}})=-$.

By Lemma 5.1.2, $d_{j}^{'}=d_{j}^{l}$ is disjoint from
$d_{j_{0}}^{l}\times [-1.5,1.5]$. By Proposition 4(2), $d_{j}^{l}$
is properly embedded in $\cup_{f} E^{l}_{f}\cup_{\gamma<j}
e^{l}_{\gamma}\times I$. If $j<j_{0}$, then,  by Lemma 5.1.4,
$d_{j}^{l}$ is disjoint from $\partial V_{l+1}^{l}\times I$. Hence
$d_{j}^{l}$ is also properly embedded in $\cup_{f}
E^{l+1}_{f}\cup_{\gamma<j} e^{'}_{\gamma}\times I$. If
$j>j_{0}=m^{l+1}$, then $j>\gamma_{i,0}$ for each $\delta(0)\leq
i\leq\theta(0)$. By Definitions 5.3.2, 5.3.3 and Lemma 5.4.1,
$d_{j}^{l}\cap
\partial V_{l+1}^{l}\times I=d_{j}^{l}\cap c_{l+1}^{l}\times
I\subset E^{l+1}_{0}\cup E^{l+1}_{l+1}
\cup_{i=\delta(0)}^{\theta(0)}
E^{l+1}_{f_{i,0}}\cup_{i=\delta(0)}^{\theta(0)}
e^{'}_{\gamma_{i,0}}\times I\cup e^{'}_{m^{l+1}}\times I\cup
(d_{j_{0}}^{l}\times I)_{\gamma_{0}}$. By Lemma 5.1.2, if
$d_{j}^{l}\cap (d_{j_{0}}^{l}\times I)_{\gamma_{0}}\neq\emptyset$,
then $j>\gamma_{0}$. Hence the lemma holds.

Case 2. \ $s(w_{\gamma_{0}})=+$ or $\gamma_{0}=\emptyset$.

By Lemma 5.1.3 and the argument in Case 1, the lemma holds. \qquad
Q.E.D.\vskip 0.5mm

{\bf \em Lemma 5.4.3.} \ (1) \ $E^{l+1}_{f}$ is a disk in
$F^{l+1}$.

(2) \ $\bigl\{c_{i}^{l+1} \ | \ i\geq l+2\bigr\}$ is a set of
pairwise disjoint arcs properly embedded in $F^{l+1}$.
Furthermore, $c_{i}^{l+1}$ lies in one of $E^{l+1}_{f}$ for some
$f$.

(3) \ $\partial c_{i}^{l+1}\cap e^{'}_{\gamma}\times I=\emptyset$
for each $\gamma\in m(l+1)$, and $\partial c_{i}^{l+1}\cap
d_{j}^{'}=\emptyset$ for each $j\in\bigl\{1,2,\ldots,
n\bigr\}-m(l+1)$.

(4) \ $c_{i}^{l+1}$ is obtained by doing band sums with copies of
$\partial V_{l+1}^{l}$ to $c_{i}^{l}$.

{\bf \em Proof.} \ (1) \ For each $1\leq f\leq l$, by Lemma 3.1.6,
each component of $d_{j_{0}}^{l}\cap E^{l}_{f}$ is a properly
embedded arc $c$ in $E^{l}_{f}$, and each component of
$d_{j_{0}}^{l}\times [-2,2] \cap E^{l}_{f}$ is  $c\times [-2,2]$
in $E^{l}_{f}$. Since $E^{l}_{f}\cap E^{l}_{0}=\emptyset$,
$E^{l+1}_{f}$ is a disk in $F^{l+1}$.

Since each component of $d_{j_{0}}^{l}\cap E^{l}_{0}$ is a
properly embedded arc $c$ in $E^{l}_{0}$. By Lemma 5.1.4(4), each
component of $d_{j_{0}}^{l}\times [-2,2] \cap E_{f}$ is  $(c\cap
E_{f})\times [-2,2]$ in $E^{l}_{f}$, where $c\cap E_{f}$ is a
properly embedded arc in $E_{f}$ for $f=0,l+1$. Hence
$E^{l+1}_{0}, E^{l+1}_{l+1}$ are two disks in $F^{l+1}$.

(2) \ Let $i\geq l+2$. By Proposition 4(4) and Lemma 5.1.4(5),
either $c_{i}^{l}$ lies in $E^{l}_{f}$ for some $1\leq f\leq l$ or
$c_{i}^{l}$ lies in $E_{f}$ for $f=0,l+1$. We may assume that
$c_{i}^{l}\subset E^{l}_{f}$. By Lemma 3.1.6, each component of
$c_{i}^{l}\cap d_{j_{0}}^{l}$ is a point $p$. Furthermore,
$p=c_{i}^{l}\cap c$ where $c$ is one component of
$d_{j_{0}}^{l}\cap E^{l}_{f}$. By Lemma 5.1.2(8) and Lemma
5.1.3(5), each component of $d_{j_{0}}^{l}\times [-2,2]\cap
c_{i}^{l}=p\times [-2,2]\subset intc_{i}^{l}$. By the
construction, $c_{i}^{l+1}$ is a properly embedded arc in
$F^{l+1}$ which lies in $E^{l+1}_{f}$.

(3) \ By Lemma 5.1.2(8) and Lemma 5.1.3(5), $\partial c^{l}_{i}$
is disjoint from $d_{j_{0}}^{l}\times [-2.5,2.5]$.  By Proposition
4(5) and the construction, (3) holds.

(4) \ Since $\partial V_{l+1}^{l}\times I\cup d_{j_{0}}^{l}\times
[-2,2]-d_{j_{0}}^{l}\times (-1.5,1.5)$ is a disk. By (2),  each
component of $c^{l}_{i}\cap d_{j_{0}}^{l} \times [-2,2]$ is an
arc. Hence the homeomorphism $H$ means band sums.\qquad
Q.E.D.\vskip 0.5mm

{\bf \em Lemma 5.4.4.} \ For each $\gamma\in m(l+1)$,
$e^{'}_{\gamma}\times I$ is a disk in $F^{l+1}$ such that

(1) \ $(\partial e^{'}_{\gamma})\times I\subset \cup_{f}\partial
E^{l+1}_{f}$ for $\gamma\in m(l)$,

(2) \ $(\partial e^{'}_{\gamma})\times I\cap (\partial
e^{'}_{\lambda})\times I=\emptyset$ for $\gamma\neq \lambda$.

(3) \ $\partial_{1} e^{'}_{m^{l+1}}\subset
E^{l+1}_{f_{\delta(0),0}}$ and $\partial_{2}
e^{'}_{m^{l+1}}\subset E^{l+1}_{f_{\theta(0),0}}$.

(4) \ $e^{'}_{\gamma}\times I=e^{l}_{\gamma}\times I$ for
$\gamma\neq m^{l+1}> Max\bigl\{\gamma_{i,0} \ | \ \delta(0)\leq
i\leq \theta(0)\bigr\}$.

{\bf \em Proof.} \ There are two cases:

Case 1. \ $s(w_{\gamma_{0}})=-$.

Suppose that $\gamma\neq m^{l+1}>Max\bigl\{\gamma_{i,0} \ | \
\delta(0)\leq i\leq \theta(0)\bigr\}$. If $\gamma\geq \gamma_{0}$,
then $\gamma>j_{0}$.  By Definition 5.3.2, $e^{'}_{\gamma}\times
I=e^{l}_{\gamma}\times I$. By Lemma 5.1.2(4),
$e^{l}_{\gamma}\times I$ is disjoint from $d_{j_{0}}^{l}\times
[-2.5,-1.5]\cup [1.5,2.5]$. Hence $e^{'}_{\gamma}\times I\subset
F^{l+1}$. If $\gamma<\gamma_{0}$, then, by Lemma 5.1.2(3),
$e^{l}_{\gamma}\times I$ is disjoint from $d_{j_{0}}^{l}\times
[-2,2]$. By Definition 5.3.2, $e^{'}_{\gamma}\times
I=e^{l}_{\gamma}\times I\subset F^{l+1}$.

Suppose that $\gamma\leq Max\bigl\{\gamma_{i,0} \ | \
\delta(0)\leq i\leq \theta(0)\bigr\}$. By Lemma 5.1.2, each
component of $d_{j_{0}}^{l}\times [-2,2]\cap e^{l}_{\gamma}\times
I$ is $c\times [-2,2]\subset e^{l}_{\gamma}\times (0,1)$. By Lemma
5.1.4, $e^{l}_{\gamma}\times I$ is disjoint from $\partial
V_{l+1}^{l}\times I$. By Definition 5.3.2, $e^{'}_{\gamma}\times
I$ is a disk in $F^{l+1}$.

It is easy to see that $e^{'}_{m^{l+1}}\times I=B_{3}$ is a disk
in $F^{l+1}$.

By Proposition 4(1) and Definition 2.1.4, $(\partial
e^{l}_{\gamma})\times I\subset \cup_{f}\partial E^{l}_{f}$,
$(\partial e^{l}_{\gamma})\times I\cap (\partial
e^{l}_{\lambda})\times I=\emptyset$ for $\gamma,\lambda\in m(l)$.
Hence (1) follows from the construction, and $(\partial
e^{'}_{\gamma})\times I\cap (\partial e^{'}_{\lambda})\times
I=\emptyset$ for $\gamma\neq \lambda\in m(l)$. By Proposition 4(1)
and Lemma 5.1.2, $\partial d_{j_{0}}^{l}\times [-2,2]$ is disjoint
from $e^{l}_{\gamma}\times I$ for $\gamma\in m(l)$ even if
$\gamma=\gamma_{0}$. By Lemma 5.3.2, $(\partial
e^{'}_{m^{l+1}})\times I=H(\partial d_{j_{0}}^{l}\times [-2,2])$.
Hence $(\partial e^{'}_{m^{l+1}})\times I\cap (\partial
e^{'}_{\gamma})\times I=\emptyset$ for $\gamma\in m(l)$. Thus (2)
holds.

By Lemma 3.3.1, $d_{j_{0}}^{l}=\cup_{i=\delta(0)}^{\theta(0)}
d_{i,0}\cup_{i=\delta(0)}^{\theta(0)}e_{i,0}$. Hence $\partial_{1}
d_{j_{0}}^{l}\subset E^{l}_{f_{\delta(0),0}}$ and $\partial_{2}
d_{j_{0}}^{l}\subset E^{l}_{f_{\theta(0),0}}$. By Definitions
5.3.1, 5.3.2 and 5.3.3, (3) holds. See Figure 24.

Case 2.  $s(w_{\gamma_{0}})=+$ or $\gamma_{0}=\emptyset$.

By Lemma 5.1.3 and the argument in Case 1, the lemma holds for
this case.\qquad Q.E.D.\vskip 0.5mm

\subsection{$F^{l+1}$ is a surface generated by
$\cup_{f} E^{l+1}_{f}\cup_{\gamma\in m(l+1)} e^{'}_{\gamma}$}

\ \ \ \ \  {\bf \em Lemma 5.5.1.} \ $F^{l+1}=\cup_{f}
E^{l+1}_{f}\cup_{\gamma\in m(l+1)} e^{'}_{\gamma}\times I$.

{\bf \em Proof.}  \ By Lemma 5.4.1, $(d_{j_{0}}^{l}\times
([-2,-1.5]\cup[1.5,2]))\cup b\times I\subset \cup_{f}
E^{l+1}_{f}\cup_{i=\delta(0)}^{\theta(0)}
e^{'}_{\gamma_{i,0}}\times I\cup e^{'}_{m^{l+1}}\times I$. By
Proposition 4(1), the lemma holds. \qquad Q.E.D.\vskip 0.5mm

{\bf \em Lemma 5.5.2.}  \ If $s(w_{\lambda})=+$, then
$inte^{'}_{\lambda}\times I$ is disjoint from $\cup_{f}
E^{l+1}_{f}\cup_{\gamma<\lambda} e^{'}_{\gamma}\times I$.

{\bf \em Proof.} \ By Proposition 4(1), if $\lambda\in m(l)$ with
$s(w_{\lambda})=+$, then $inte^{l}_{\lambda}\times I$ is disjoint
from $\cup_{f} E^{l}_{f}\cup_{\gamma<\lambda} e^{l}_{\gamma}\times
I$. There are two cases:

Case 1. $s(w_{\gamma_{0}})=-$.

There are three sub-cases:

Case 1.1. $\lambda\neq m^{l+1}>Max\bigl\{\gamma_{i,0} \ | \
\delta(0)\leq i\leq \theta(0)\bigr\}$.

By Lemma 5.4.3, $e^{'}_{\lambda}\times I=e^{l}_{\lambda}\times I$.
By Proposition 4(1), $inte^{l}_{\lambda}\times I$ is disjoint from
$\cup_{f} E^{l}_{f}\cup_{\gamma<\lambda}e^{l}_{\gamma}\times I$.
By Lemma 5.1.4(1), Proposition 4(5), $inte^{l}_{\lambda}\times I$
is disjoint from $\partial V_{l+1}^{l}\times I$. By Lemma 5.1.2(3)
and (4), $d_{j_{0}}^{l}\times [-2-1.5]\cup[1.5,2]$ is disjoint
from $e^{l}_{\lambda}\times I$. By Definitions 5.3.2 and 5.3.3,
$\cup_{f}E^{l+1}_{f}\cup_{\gamma<\lambda} e^{'}_{\gamma}\times
I\subset
\cup_{f}E^{l}_{f}\cup_{\gamma<\lambda}e^{l}_{\gamma}\times I\cup
\partial V_{l+1}^{l}\times I\cup d_{j_{0}}^{l}\times [-2-1.5]\cup[1.5,2]$.
Hence $inte^{'}_{\lambda}\times I$ is
disjoint from $\cup_{f} E^{l+1}_{f}\cup_{\gamma<\lambda}
e^{'}_{\gamma}\times I$.

Case 1.2. $\lambda=m^{l+1}$. Now if $\gamma<\lambda$, then,  by
Lemma 5.1.4(2), $e^{l}_{\gamma}\times I$ is disjoint from
$\partial V_{l+1}^{l}\times I$. By Definitions 5.3.2 and 5.3.3,
$\cup_{f}E^{l+1}_{f}\cup_{\gamma<\lambda} e^{'}_{\gamma}\times
I\subset
(\cup_{f}E^{l}_{f}\cup_{\gamma<\lambda}e^{l}_{\gamma}\times
I-d_{j_{0}}^{l}\times [-2,2])\cup B_{1}\cup B_{2}$,
$e^{'}_{m^{l+1}}\times I=B_{3}$ where $B_{i}$ is as in Section
5.2. Hence the lemma holds.

Case 1.3. $\lambda\leq Max\bigl\{\gamma_{i,0} \ | \ \delta(0)\leq
i\leq \theta(0)\bigr\}$.  Now if $\gamma\leq\lambda$, then
$\gamma\leq\lambda<j_{0}=m^{l+1}$. By Lemma 5.1.4,
$e^{l}_{\gamma}\times I$ is disjoint from $\partial
V_{l+1}^{l}\times I$. By Proposition 4(1),
$inte^{l}_{\lambda}\times I$ is disjoint from $\cup_{f}
E^{l}_{f}\cup_{\gamma<\lambda}e^{l}_{\gamma}\times I$. By Lemma
5.1.2, each component of $d_{j_{0}}^{l}\times [-2,2]\cap
e^{l}_{\lambda}\times I$ is $c\times [-2,2]\subset
e^{l}_{\lambda}\times (0,1)$ where $c$ is a core of $
e^{l}_{\lambda}\times (0,1)$. Hence $intc\times [-2,2]$ is
disjoint from $\cup_{f}
E^{l}_{f}\cup_{\gamma<\lambda}e^{l}_{\gamma}\times I$. By
Definitions 5.3.2 and 5.3.3, the lemma  holds.

Case 2. $s(w_{\gamma_{0}})=+$ or $\gamma_{0}=\emptyset$.

By Lemma 5.1.3, Definition 5.2.3 and the argument in Case 1, the
lemma holds.\qquad Q.E.D.\vskip 0.5mm

{\bf \em Lemma 5.5.3.}  \ If $s(w_{\lambda})=-$, then
$inte^{'}_{\lambda}\times I\cup(\cup_{f}
E^{l+1}_{f}\cup_{\gamma<\lambda} e^{'}_{\gamma}\times
I)=\cup_{r\in I(w_{\lambda},l+1)} (d_{r}^{'}\times I)_{\lambda}$.

{\bf \em Proof.} \ Suppose that $\lambda\in m(l+1)$ and
$s(w_{\lambda})=-$. By Lemma 5.1.4(4), $\lambda\neq m^{l+1}$. By
Proposition 4(1), $inte^{l}_{\lambda}\times I\cup(\cup_{f}
E^{l}_{f}\cup_{\gamma<\lambda} e^{l}_{\gamma}\times I)=\cup_{r\in
I(w_{\lambda},l)} (d_{r}^{l}\times I)_{\lambda}$.  Now there are
two cases:

Case 1. $s(w_{\gamma_{0}})=-$. There are four sub-cases:

Case 1.1. $\lambda\neq \gamma_{0}>m^{l+1}$.

In this case, $\bigl\{\gamma<\lambda\in
m(l+1)\bigr\}=\bigl\{\gamma<\lambda\in
m(l)\bigr\}\cup\bigl\{m^{l+1}\bigr\}$.

By Definition 2.3.1,
$I(w_{\lambda},l+1)=I(w_{\lambda},l)-\bigl\{m^{l+1}\bigr\}$. Since
$m^{l+1}=j_{0}\in I(w_{\gamma_{0}},l)$, by Lemma 2.2.5,
$I(w_{\lambda},l+1)=I(w_{\lambda},l)$. By Lemma 5.4.4,
$e^{'}_{\gamma}\times I=e^{l}_{\gamma}\times I$. By Lemmas 5.1.2,
$d_{j_{0}}^{l}\times [-2,-1]\cup [1,2]$ is disjoint from
$e^{l}_{\lambda}\times I$. By Proposition 4(1),
$e^{l}_{\lambda}\times I\cap (\cup_{f}
E^{l}_{f}\cup_{\gamma<\lambda} e^{l}_{\gamma}\times
I-d_{j_{0}}^{l}\times [-2,-1)\cup (1,2])=\cup_{r\in
I(w_{\lambda},l)} (d_{r}^{l}\times I)_{\lambda}$. By Proposition
4(5) and  Lemma 5.1.4(1), $e^{l}_{\lambda}\times I$ is disjoint
from $(\partial V_{l+1}^{l}-c_{l+1}^{l})\times I$.  By Lemma
5.4.1, the lemma holds.

Case 1.2. $\lambda=\gamma_{0}$.

In this case, $\bigl\{\gamma<\lambda\in
m(l+1)\bigr\}=\bigl\{\gamma<\lambda\in
m(l)\bigr\}\cup\bigl\{m^{l+1}\bigr\}$.

By Definition 2.3.1,
$I(w_{\lambda},l+1)=I(w_{\lambda},l)-\bigl\{m^{l+1}\bigr\}$.

Now by Lemma 5.1.2(2), (3) and Definitions 5.3.2, 5.3.3,
$(d_{j_{0}}^{l}\times I)_{\gamma_{0}}$ is disjoint from $\cup_{f}
E^{l+1}_{f}\cup_{\gamma<\gamma_{0}} e^{'}_{\gamma}\times I$. By
the argument in Case 1.1, the lemma holds.

Case 1.3. $m^{l+1}>\lambda>Max\bigl\{\gamma_{i,0} \ | \
\delta(0)\leq i\leq\theta(0)\bigr\}$.

In this case, $\bigl\{\gamma<\lambda\in
m(l+1)\bigr\}=\bigl\{\gamma<\lambda\in m(l)\bigr\}$.

By Lemma 5.1.4 and Lemma 5.4.4, $e^{'}_{\lambda}\times
I=e^{l}_{\lambda}\times I$ is disjoint from $\partial
V_{l+1}^{l}\times I$. By Lemma 5.1.2(3), $e^{'}_{\lambda}\times I$
is disjoint from $d_{j_{0}}^{l}\times [-2,2]$. By Definition 5.3.2
and 5.3.3, $e^{'}_{\lambda}\times
I\cap(\cup_{f}E^{l+1}_{f}\cup_{\gamma<\lambda}e^{'}_{\gamma}\times
I )=e^{l}_{\lambda}\times
I\cap(\cup_{f}E^{l}_{f}\cup_{\gamma<\lambda}e^{l}_{\gamma}\times I
)$. Hence the lemma holds.

Case 1.4. $\lambda\leq Max\bigl\{\gamma_{i,0} \ | \ \delta(0)\leq
i\leq\theta(0)\bigr\}$.

In this case, $\bigl\{\gamma<\lambda\in
m(l+1)\bigr\}=\bigl\{\gamma<\lambda\in m(l)\bigr\}$.

Now we denote by $P_{\lambda}$ the surface
$\cup_{f}E^{l}_{f}\cup_{\gamma<\lambda}e^{l}_{\gamma}\times I$.

Now $e^{'}_{\lambda}\times I
\cap(\cup_{f}E^{l+1}_{f}\cup_{\gamma<\lambda}e^{'}_{\gamma}\times
I )=((e^{l}_{\lambda}\times I-e^{l}_{\lambda}\times I\cap
d_{j_{0}}^{l}\times [-2,2])\cap (P_{\lambda}-P_{\lambda}\cap
d_{j_{0}}^{l}\times [-2,2]))\cup (H(e^{l}_{\lambda}\times I\cap
d_{j_{0}}^{l}\times [-2,2])\cap H(P_{\lambda}\cap
d_{j_{0}}^{l}\times [-2,2]))$. Hence the lemma holds.

Case 2. $s(w_{\gamma_{0}})=+$ or $\gamma_{0}=\emptyset$.

By Lemma 5.1.3 and Definition 5.2.3, and the argument in Case 1,
the lemma holds.\qquad Q.E.D.\vskip 0.5mm

{\bf \em Lemma 5.5.4.} \ If $\gamma<\lambda\in m(l+1)$, then each
component of $e^{'}_{\lambda}\times I\cap e^{'}_{\gamma}\times I$
is $(c\times I)_{\lambda}\subset e^{'}_{\gamma}\times (0,1)$ where
$c\subset inte^{'}_{\lambda}$ is a core of $e^{'}_{\gamma}\times
(0,1)$.

{\bf \em Proof.} \ By Proposition 4(1) and Definition 2.1.4, if
$\gamma<\lambda\in m(l)$, then each component of
$e^{l}_{\lambda}\times I\cap e^{l}_{\gamma}\times I$ is $(c\times
I)_{\lambda}\subset e^{l}_{\gamma}\times (0,1)$ where $c\subset
inte^{l}_{\lambda}$ is a core of $e^{l}_{\gamma}\times (0,1)$.

Without loss of generality, we may assume that
$s(w_{\gamma_{0}})=-$. Now there are six cases:

Case 1. $\lambda, \gamma\in m(l)$ and
$\lambda>\gamma>Max\bigl\{\gamma_{i,0} \ | \ \delta(0)\leq
i\leq\theta(0)\bigr\}$.

By Lemma 5.4.3, $e^{'}_{\lambda}\times I=e^{l}_{\lambda}\times I$
and $e^{'}_{\gamma}\times I=e^{l}_{\gamma}\times I$. Hence the
lemma holds.

Case 2. $\lambda=m^{l+1}>\gamma$.

Since $s(w_{m^{l+1}})=+$, by Lemma 5.4.4(2) and Lemma 5.5.2,
$e^{'}_{m^{l+1}}\times I$ is disjoint from $e^{'}_{\gamma}\times
I$.

Case 3. $\lambda>\gamma=m^{l+1}$.

By Lemma 3.2.4, and Lemma 5.1.4(1), each component of
$e^{l}_{\lambda}\times I\cap c_{l+1}^{l}\times
I=e^{l}_{\lambda}\times I\cap \partial V_{l+1}^{l}\times I$ is
$(c\times I)_{\lambda}$ where $c\subset inte^{l}_{\lambda}$.
Furthermore, by Lemma 5.1.2(3) and (4), $(c\times I)_{\lambda}$
lies either in $F^{l}-d_{j_{0}}^{l}\times [-2,2]$ or in
$(d_{j_{0}}^{l}\times I)_{\gamma_{0}}$. If $(c\times I)_{\lambda}$
lies in $(d_{j_{0}}^{l}\times I)_{\gamma_{0}}$. Then $(c\times
I)_{\lambda}$ is disjoint from $e^{'}_{m^{l+1}}\times I$. If
$(c\times I)_{\lambda}$   lies in $F^{l}-d_{j_{0}}^{l}\times
[-2,2]$, then, by Definitions 5.3.3, $(c\times I)_{\lambda}$
intersects $e^{'}_{m^{l+1}}\times I$ in $(c_{*}\times
I)_{\lambda}\subset e^{'}_{m^{l+1}}\times (0,1)$, where
$c_{*}\subset inte^{'}_{\lambda}$ is a core of
$e^{'}_{m^{l+1}}\times (0,1)$. Hence the lemma holds.

Case 4. $\gamma_{0}>\lambda>Max\bigl\{\gamma_{i,0} \ | \
\delta(0)\leq i\leq\theta(0)\bigr\}\geq\gamma$.

By Lemma 5.4.4(4), $e^{l}_{\lambda}=e^{'}_{\lambda}\times I$. By
Lemma 5.1.2(3), $e^{l}_{\lambda}\times I$ is disjoint from
$d_{j_{0}}^{l}\times [-2,2]$. By Lemma 5.1.4,
$e^{l}_{\lambda}\times I\cap \partial V_{l+1}^{l}\times
I=e^{l}_{\lambda}\times I\cap c_{l+1}^{l}\times I$. By Lemma
3.2.4, each component of  $e^{l}_{\lambda}\times I\cap \partial
V_{l+1}^{l}\times I$ is $(c\times I)_{\lambda}$ where $c\subset
inte^{l}_{\lambda}$. By Lemma 5.1.4, $e^{l}_{\gamma}\times I$ is
disjoint from $\partial V_{l+1}^{l}\times I$. By Definition 5.3.2,
$e^{'}_{\lambda}\times I\cap e^{'}_{\gamma}\times I=S_{1}\cup
S_{2}$, where  $S_{1}=e^{l}_{\lambda}\times I\cap
e^{l}_{\gamma}\times I$ and $S_{2}= (e^{l}_{\lambda}\times I\cap
\partial V_{l+1}^{l}\times I)\cap H(e^{l}_{\gamma}\cap (h_{1}\cup h_{2})\times [-2,2])$.
Hence the lemma holds.

Case 5. $\lambda\geq\gamma_{0}>Max\bigl\{\gamma_{i,0} \ | \
\delta(0)\leq i\leq\theta(0)\bigr\}\geq\gamma$.

Now by Lemma 5.1.2 and Lemma 5.4.4, $e^{'}_{\lambda}\times
I=e^{l}_{\lambda}\times I$ is disjoint from $d_{j_{0}}^{l}\times
[-2,-1)\cup (1,2]$ even if $\lambda=\gamma_{0}$. By Definition
5.3.2, $e^{'}_{\lambda}\times I\cap e^{'}_{\gamma}\times
I=S_{1}\cup S_{2}$, where $S_{1}= (e^{l}_{\lambda}\times
I-e^{l}_{\lambda}\times I\cap d_{j_{0}}^{l}\times [-1,1])\cap
e^{l}_{\gamma}\times I$ and $S_{2}=(e^{l}_{\lambda}\times I\cap
\partial V_{l+1}\times I)\cap H(e^{l}_{\gamma}\times I\cap
(h_{1}\cup h_{2})\times I)$.

By Lemma 3.3.1,
$d_{j_{0}}^{l}=\cup_{i=\delta(0)}^{\theta(0)}d_{i,0}\cup_{i=\delta(0)}^{\theta(0)}e_{i,0}$.
Since $\gamma<j_{0}$,  by Lemma 5.1.2(6),
$intd_{\delta(0),0}\times [-1,1]$ and $intd_{\theta(0),0}\times
[-1,1]$ are disjoint from $e^{l}_{\gamma}\times I$. Hence each
component of $e^{l}_{\lambda}\times I\cap e^{l}_{\gamma}\times I$
is either in $e^{l}_{\lambda}\times I-e^{l}_{\lambda}\times I\cap
d_{j_{0}}^{l}\times [-1,1]$ or in $d_{j_{0}}^{l}\times [-1,1]$. By
the argument in Case 4, the lemma holds.

Case 6. $Max\bigl\{\gamma_{i,0} \ | \ \delta(0)\leq
i\leq\theta(0)\bigr\}\geq\lambda>\gamma$.

By Definition 5.3.2, $e^{'}_{\lambda}\times I\cap
e^{'}_{\gamma}\times I=S_{1}\cup S_{2}$ where $S_{1}=
(e^{l}_{\lambda}\times I-e^{l}_{\lambda}\times I\cap (h_{1}\cup
h_{2})\times I)\cap (e^{l}_{\gamma}\times I-e^{l}_{\gamma}\times
I\cap (h_{1}\cup h_{2})\times I)$ and
$S_{2}=H(e^{l}_{\lambda}\times I\cap (h_{1}\cup h_{2})\times
I)\cap H(e^{l}_{\gamma}\times I\cap (h_{1}\cup h_{2})\times I)$.
Since $H$ is a homeomorphism, the lemma holds. \qquad Q.E.D.\vskip
1.5mm

{\bf \em Lemma 5.5.5.} \ $\cup_{f} E^{l+1}_{f}\cup_{\gamma\in
m(l+1)} e^{'}_{\gamma}$ is a abstract tree, and $F^{l+1}$ is
generated by $\cup_{f} E^{l+1}_{f}\cup_{\gamma\in m(l+1)}
e^{'}_{\gamma}$.

{\bf \em Proof.} \ By Lemma 5.4.4 and Lemmas 5.5.1-5.5.4, we only
need to prove that $\cup_{f} E^{l+1}_{f}\cup_{\gamma}
e^{'}_{\gamma}$ is a abstract tree. By Lemma 5.4.4, $\partial_{1}
e^{'}_{m^{l+1}}\subset \partial E^{l+1}_{f_{\delta(0),0}}$ and
$\partial_{2} e^{'}_{m^{l+1}}\subset \partial
E^{l+1}_{f_{\theta(0),0}}$. By the argument in Lemma 4.4.5, the
lemma holds.\qquad Q.E.D.

\subsection{Properties of $d_{j}^{'}$}

 \ \ \ \ \ \ {\bf \em Lemma 5.6.1.} \ Suppose that
$j\in\bigl\{1,\ldots,n\bigr\}-m(l+1)$, $\gamma\in m(l+1)$ and
$j\notin I(w_{\gamma},l+1)$. If $j<\gamma$, then $d_{j}^{'}$ is
disjoint from $e^{'}_{\gamma}\times I$.

{\bf \em Proof.} \ Suppose that $j<\gamma$ and $j\notin
I(w_{\gamma},l+1)$. By Proposition 4(3), $d_{j}^{l}$ is disjoint
from $e^{l}_{\gamma}\times I$ for $\gamma\in m(l)$. Without loss
of generality, we may assume that $s(w_{\gamma_{0}})=-$. Now there
are three cases:

Case 1. \ $\gamma\neq m^{l+1}$, $j>Max\bigl\{\gamma_{i,0} \ | \
\delta(0)\leq i\leq\theta(0)\bigr\}$.

Now $\gamma>j>Max\bigl\{\gamma_{i,0} \ | \ \delta(0)\leq
i\leq\theta(0)\bigr\}$. By Definition 5.3.2 and Lemma 5.4.4,
$e^{'}_{\gamma}\times I=e^{l}_{\gamma}\times I$ and
$d_{j}^{'}=d_{j}^{l}$. Hence $d_{j}^{'}$ is disjoint from
$e^{l}_{\gamma}\times I$.

Case 2. $\gamma= m^{l+1}$.

Since $j<m^{l+1}$, $j<\gamma_{0}$.  By Lemma 5.1.2(5),
$d_{j}^{'}=d_{j}^{l}$ is disjoint from $d_{j_{0}}^{l}\times
[-2,2]$. By Lemma 5.1.4(2), $d_{j}^{'}$ is disjoint from $\partial
V_{l+1}^{l}$. Hence $d_{j}^{'}$ is disjoin from
$e^{'}_{m^{l+1}}\times I\subset d_{j_{0}}^{l}\times [-2,2]\cup
\partial V_{l+1}^{l}\times I$.

Case 3. \ $j<Max\bigl\{\gamma_{i,0} \ | \ \delta(0)\leq
i\leq\theta(0)\bigr\}$.

Now by Lemma 5.1.2(5) and Lemma 5.1.4(2), $d_{j}^{l}$ is disjoint
from $\partial V_{l+1}^{l}\times I$ and $d_{j_{0}}^{l}\times
[-2,2]$. By Definition 5.3.2, $e^{'}_{\gamma}\times I\subset
e^{l}_{\gamma}\times I\cup d_{j_{0}}^{l}\times [-2,2]\cup \partial
V_{l+1}^{l}\times I$. Hence $d_{j}^{'}=d_{j}^{l}$ is disjoint from
$e^{'}_{\gamma}\times I$. \qquad Q.E.D.\vskip 0.5mm

{\bf \em Lemma 5.6.2.} \ Suppose that
$j\in\bigl\{1,\ldots,n\bigr\}-m(l+1)$, $\gamma\in m(l+1)$. If
$j>\gamma$, then each component of $d_{j}^{'}\cap
e^{'}_{\gamma}\times I$ is a core $c\subset intd_{j}^{'}$ of
$e^{'}_{\gamma}\times (0,1)$.

{\bf \em Proof.} \ Suppose that $j>\gamma$. By Proposition 4(2),
each component of $d_{j}^{l}\cap e^{l}_{\gamma}\times I$ is a core
$c\subset intd_{j}^{l}$ of $e^{l}_{\gamma}\times (0,1)$ for
$\gamma\in m(l)$. Without loss of generality, we may assume that
$s(w_{\gamma_{0}})=-$. Now there are four cases:

Case 1. \ $\gamma\neq m^{l+1}>Max\bigl\{\gamma_{i,0} \ | \
\delta(0)\leq i\leq \theta(0)\bigr\}$.

By Definition 5.3.2 and Lemma 5.4.4, $d_{j}^{'}=d_{j}^{l}$ and
$e^{'}_{\gamma}\times I=e^{l}_{\gamma}\times I$. Hence the lemma
holds.

Case 2. \ $\gamma=m^{l+1}$.

By Definition 5.3.2, $e^{'}_{m^{l+1}}\times I=B_{3}$. By Lemma
5.1.4, $d_{j}^{l}\cap
\partial V_{l+1}^{l}\times I=d_{j}^{l}\cap c_{l+1}^{l}\times I$.
By Lemma 5.1.2 and Lemma 3.1.6, each component of $d_{j}^{l}\cap
\partial V_{l+1}^{l}\times I$ is an arc $c\subset intd_{j}^{l}$
which lies either in $d_{j_{0}}^{l}\times [-1,1]$ or in
$c_{l+1}^{l}\times I-d_{j_{0}}^{l}\times [-2,2]$. If $c$ lies in
$d_{j_{0}}^{l}\times [-1,1]$, then $c$ is disjoint from
$e^{'}_{m^{l+1}}\times I$. If $c$ lies in $c_{l+1}^{l}\times
I-d_{j_{0}}^{l}\times [-2,2]$, then $c$ intersects
$e^{'}_{m^{l+1}}\times I$ in a core of $e^{'}_{m^{l+1}}\times
(0,1)$. Hence the lemma holds.

Case 3. \ $\gamma<j<Max\bigl\{\gamma_{i,0} \ | \ \delta(0)\leq
i\leq \theta(0)\bigr\}$.

Now $j<m^{l+1}$. By Lemma 5.1.4(2), $d_{j}^{'}=d_{j}^{l}$ is
disjoint from $\partial V_{l+1}^{l}\times I$. By Lemma 5.1.2(5),
$d_{j}^{'}$ is disjoint from $d_{j_{0}}^{l}\times [-2,2]$. By
Definition 5.3.2, $e^{'}_{\gamma}\times I\subset
e^{l}_{\gamma}\times I\cup d_{j_{0}}^{l}\times [-2,2]\cup\partial
V_{l+1}^{l}\times I$. Hence $d_{j}^{'}\cap e^{'}_{\gamma}\times
I=d_{j}^{l}\cap e^{l}_{\gamma}\times I$, and the lemma holds.

Case 4. \ $j>Max\bigl\{\gamma_{i,0} \ | \ \delta(0)\leq i\leq
\theta(0)\bigr\}\geq\gamma$.

Now by Lemma 5.1.2(5), $d_{j}^{'}=d_{j}^{l}$ is disjoint from
$d_{j_{0}}^{l}\times [-2,-1]\cup [1,2]$.  By Lemma 5.1.2(6),
$intd_{\delta(0),0}\times [-2,2]$ and $intd_{\theta(0),0}\times
[-2,2]$ are disjoint from $e^{l}_{\gamma}\times I$. Hence each
component of  $d_{j}^{l}\cap e^{l}_{\gamma}\times I$ is either in
$d_{j_{0}}^{l}\times [-1,1]$ or disjoint from $d_{j_{0}}^{l}\times
[-2,2]$. By Lemma 3.1.6, each component of $d_{j}^{l}\cap
c_{l+1}^{l}\times I$ is $c\subset intd_{j}^{l}$ with $\partial_{1}
c\subset c_{l+1}^{l}\times\bigl\{1\bigr\}$ and $\partial_{2}
c\subset c_{l+1}^{l}\times\bigl\{-1\bigr\}$.  By Lemma 5.1.4,
$e^{l}_{\gamma}\times I$ is disjoint from $\partial
V_{l+1}^{l}\times I$. Now $d_{j}^{'}\cap e^{'}_{\gamma}\times
I=S_{1}\cup S_{2}$, where $S_{1}=(d_{j}^{l}-d_{j_{0}}^{l}\times
[-1,1])\cap e^{l}_{\gamma}\times I$ and $S_{2}=
d_{j}^{l}\cap\partial V_{l+1}^{l}\times I\cap
H(e^{l}_{\gamma}\times I\cap (h_{1}\cup h_{2})\times [-2,2])$.
Hence the lemma holds. \qquad Q.E.D.\vskip 0.5mm

{\bf \em Lemma 5.6.3.} \ Suppose that $j\notin L(c_{l+1}^{l})$.
Then

(1) \ $d_{j}^{'}$ is regular in $\cup_{f}
E^{l+1}_{f}\cup_{\gamma<j} e^{'}_{\gamma}\times I$.

(2) \ For each $i\geq l+2$, $d_{j}^{'}-\cup_{\gamma<j}
inte^{'}_{\gamma}\times I$ intersects $c_{i}^{l+1}$ in at most one
point. Furthermore, $d_{j}^{'}-\cup_{\gamma<j}
inte^{'}_{\gamma}\times I$ intersects $c_{i}^{l+1}$ in  one point
if and only if $d_{j}^{l}-\cup_{\gamma<j} inte^{l}_{\gamma}\times
I$ intersects $c_{i}^{l}$ in  one point.

{\bf \em Proof.} \ Without loss of generality, we may assume that
$s(w_{\gamma_{0}})=-$.

By Lemma 3.1.2, in $F^{l}$,
$d_{j}^{l}=\cup_{i=1}^{\theta(j)}d_{j,f_{i,j}}^{l}\cup_{i=1}^{\theta(j)-1}
e_{\gamma_{i,j}}$, where $d_{j,f_{i,j}}^{l}$ is a properly
embedded arc in $E^{l}_{f_{i,j}}$ which is disjoint from
$\cup_{\gamma<j} inte^{l}_{\gamma}\times I$ and $e_{\gamma_{i,j}}$
is a core of $e^{l}_{\gamma_{i,j}}\times (0,1)$ for some
$j>\gamma_{i,j}\in m(l)$. Furthermore, $f_{i,j}\neq f_{r,j}$ and
$\gamma_{i,j}\neq \gamma_{r,j}$ for $i\neq r$.

By Lemma 5.1.2, $d_{j}^{l}$ is disjoint from $d_{j_{0}}^{l}\times
[-2-1]\cup [1,2]$.  Hence $e_{\gamma_{i,j}}$ is also a core of
$e^{'}_{\gamma_{i,j}}\times (0,1)$. Since $j\notin
L(c_{l+1}^{l})$, $d_{j,f_{i,j}}^{l}$ is disjoint from
$c_{l+1}^{l}\times I$ and $\partial V_{l+1}^{l}\times I$ even if
$f_{i,j}=0$. Furthermore, if $f_{i,j}=0$, then $d_{j,0}^{l}$ lies
in in one of $E_{0}$ and $E_{l+1}$. See Definition 5.2.1. By
Definition 5.3.2, $e^{'}_{\gamma}\times I\subset
e^{l}_{\gamma}\times I\cup d_{j_{0}}^{l}\times
[-2,-1]\cup[1,2]\cup \partial V_{l+1}^{l}\times I$. Hence if
$\gamma<j$, then $d_{j,f_{i,j}}^{l}$ is disjoint from
$inte^{'}_{\gamma}\times I$. Now each component of
$d_{j}^{'}\cap(\cup_{\gamma<j} e^{'}_{\gamma}\times I)$ is
$e_{\gamma_{i,j}}$ for $1\leq i\leq\theta(j)-1$. By Lemmas 5.6.1
and 5.6.2, (1) holds.

Now $d_{j}^{'}-\cup_{\gamma<j} inte^{'}_{\gamma}\times
I=\cup_{i=1}^{\theta(j)}d_{j,f_{i,j}}^{l}$. Since
$c_{i}^{l+1}=(c_{i}^{l}-d_{j_{0}}^{l}\times [-2,2])\cup
H(c_{i}^{l}\cap d_{j_{0}}^{l}\times [-2,2])$. Since $j\neq j_{0}$,
$d_{j,f_{i,j}}^{l}$ is disjoint from $d_{j_{0}}^{l}\times [-2,2]$.
Hence (2) holds. \qquad Q.E.D.\vskip 0.5mm

{\bf \em Lemma 5.6.4.} \ Suppose that $j=j_{\alpha}\in
L(c_{l+1}^{l})$ and $\alpha\neq 0$. Then

(1) \
$d_{j_{\alpha}}^{'}-inta_{\alpha}=\cup_{\delta(\alpha)}^{\delta(0,\alpha)-1}
d_{i,\alpha}\cup_{\delta(\alpha)}^{\delta(0,\alpha)-1}
e_{i,\alpha}\cup_{\theta(0,\alpha)+1}^{\delta(\alpha)}
d_{i,\alpha}\cup_{\theta(0,\alpha)+1}^{\theta(\alpha)}
e_{i,\alpha}$.

(2) \ $d_{i,\alpha}$ is disjoint from $\cup_{\gamma<j_{\alpha}}
inte^{'}_{\gamma}\times I$.

(3) \ $d_{i,\alpha}$ intersects $c^{l+1}_{i}$ in one point if and
only if $d_{i,\alpha}$ intersects $c^{l}_{i}$ in one point for
$\delta(\alpha)\leq i\leq\delta(0,\alpha)-1$ or
$\theta(0,\alpha)+1\leq i\leq\theta(\alpha)$.

(4) \ For $\gamma<j_{\alpha}$, each component of
$d_{j_{\alpha}}^{'}\cap e^{'}_{\gamma}\times I$ is contained
either in $d_{j_{\alpha}}^{'}-inta_{\alpha}$ or in
$inta_{\alpha}$.

{\bf \em Proof.} \ By Definition 3.5.1,
$a_{\alpha}=\cup_{i=\delta(0,\alpha)}^{\theta(0,\alpha)}
d_{i,\alpha} \cup_{i=\delta(0,\alpha)}^{\theta(0,\alpha)}
e_{i,\alpha}$. Hence  $d_{0,\alpha}\subset a_{\alpha}$,
$d_{i,\alpha}$ is disjoint from $c_{l+1}^{l}\times I$ and
$\partial V_{l+1}^{l}\times I$ for $\delta(0)\leq
i\leq\delta(0,\alpha)-1$ or $\theta(0,\alpha)+1\leq i\leq
\theta(\alpha)$. Now by the argument in Lemma 5.6.3, (1), (2) and
(3) hold.

Suppose that $\delta(0,\alpha)\leq -1$. Then
$f_{\delta(0,\alpha),\alpha}\neq 0$. Hence
$d_{\delta(0,\alpha),\alpha}\subset
E^{l+1}_{f_{\delta(0,\alpha),\alpha}}$ is disjoint from
$c_{l+1}^{l}\times I$ and $\partial V_{l+1}^{l}\times I$. By the
proof of Lemma 5.6.3, $intd_{\delta(0,\alpha),\alpha}$ is disjoint
from $e^{'}_{\gamma}\times I$ for $\gamma<j_{\alpha}$.

Now we assume that $\delta(0,\alpha)=\theta(0,\alpha)=0$. Since
$\alpha\neq 0$, By Proposition 4(3) and Lemma 3.3.1,
$d_{i,\alpha}$ is disjoint from $d_{j_{0}}^{l}\times [-2.5,2.5]$.
Now $a_{\alpha}=a^{1}\cup a^{2}\cup a^{3}$ where $a^{1}$ is an arc
in $E_{0}\subset E^{l+1}_{0}$, $a^{2}$ is an arc in
$c_{l+1}^{l}\times I-d_{j_{0}}^{l}\times [-2.5,2.5]$ (by Lemma
5.1.2(7)), $a^{3}$ is an arc in $E_{l+1}\subset E^{l+1}_{l+1}$. By
Definition 5.2.3 and Definitions 5.3.2 and 5.3.3,
$e^{'}_{\gamma}\times I\subset e^{l}_{\gamma}\times I\cup\partial
V_{l+1}^{l}\times I\cup d_{j_{0}}^{l}\times [-2,2]$. Hence
$inta^{1}, inta^{3}$ are disjoint from $e^{'}_{\gamma}\times I$
for $\gamma<j_{\alpha}$. Thus (4) holds. \qquad Q.E.D.

\subsection{Properly embedded disks in $\cal V_{-}$ and $\cal
W_{-}$}

 \ \ \ \ \ In this section, we shall prove the following Lemma:

{\bf \em Lemma 5.7.1.} \ There are two sets of pairwise disjoint
disks $\bigl\{V^{l+1}_{i} \ | \ i\geq l+2 \ with \
s(v_{i})=-\bigr\}$ properly embedded in $\cal V_{-}$ and
$\bigl\{W^{'}_{j} \ | \ j\in\bigl\{1,\ldots,n\bigr\}-m(l+1) \ with
\ s(w_{j})=-\bigr\}$ properly embedded in $\cal W_{-}$ such that

(1) \ $\partial V_{i}^{l+1}\cap F^{l+1}=c^{l+1}_{i}\cup_{r\in
I(v_{i},l+1)} c^{l+1}_{i}$, $\partial W_{j}^{'}\cap
F^{l+1}=d_{j}^{'}\cup_{r\in I(w_{j},l+1)} d_{r}^{'}$;

(2) \ $V_{i}^{l+1}\cap W_{j}^{'}=V_{i}^{l+1}\cap W_{j}^{'}\cap
F^{l+1}$.

{\bf \em Proof.} \ Suppose that $i\geq l+2$ and $s(v_{i})=-$.
Then, by Proposition 6, $V_{i}^{l}$ is a properly embedded disk in
$\cal V_{-}$ such that $V_{i}^{l}\cap F^{l}=c^{l}_{i}\cup_{i\in
I(v_{i},l)} c_{r}^{l}$. By Lemma 5.1.4, $V_{i}^{l}$ is disjoint
from $\partial V_{l+1}^{l}\times I$. By assumption,
$s(v_{l+1})=-$. By Lemma 2.2.4, $l+1\notin I(v_{i},l)$. By
Definition 2.3.1, $I(v_{i},l+1)=I(v_{i},l)$. Now let
$C_{i}=(\partial V_{i}^{l}-c^{l}_{i}\cup_{i\in I(v_{i},l)}
c_{r}^{l})\cup c^{l+1}_{i}\cup_{i\in I(v_{i},l+1)} c_{r}^{l+1}$.
By Lemma 5.4.3(4), $c^{l+1}_{i}$ and $c_{r}^{l+1}$ are obtained by
doing band sums with copies $\partial V_{l+1}^{l}$ to $c_{i}^{l}$
and $c_{r}^{l}$. Hence $C_{i}$ bounds a disk in $\cal V_{-}$,
denoted by $V_{i}^{l+1}$. Since $F^{l+1}\subset F^{l}\cup\partial
V_{l+1}^{l}\times I$. Hence $V_{i}^{l+1}\cap
F^{l+1}=c^{l+1}_{i}\cup_{i\in I(v_{i},l+1)} c_{r}^{l+1}$.

Suppose now that $j\in\bigl\{1,\ldots,n\bigr\}-m(l+1)$ and
$s(w_{j})=-$. Then, by Proposition 6, $W_{j}^{l}$ is a properly
embedded disk in $F^{l}$ such that $W_{j}^{l}\cap
F^{l}=d_{j}^{l}\cup_{r\in I(w_{j},l)} d_{r}^{l}$. By Lemma 5.1.4,
$W_{j}^{l}-d_{j}^{l}\cup_{r\in I(w_{j},l)} d_{r}^{l}$ is disjoint
from $\partial V_{l+1}^{l}\times I$. We denote by $W_{j}^{'}$ the
disk $W_{j}^{l}$. There are two cases:

Case 1. $\gamma_{0}\neq\emptyset$. See Definition 3.2.1.

Since $j\notin m(l+1)$ and $\gamma_{0}\in m(l)$, by Lemma 2.2.5,
$j_{0}\notin I(w_{j},l)$. By Definition 2.3.1,
$I(w_{j},l+1)=I(w_{j},l)$. Since $d_{j}^{'}=d_{j}^{l}$ and
$d_{r}^{'}=d_{r}^{l}$, $W_{j}^{'}\cap F^{l+1}=d_{j}^{'}\cup_{r\in
I(w_{j},l+1)} d_{r}^{'}$. By Proposition 6, $\partial V_{l+1}^{l}-
c^{l}_{i}\cup_{i\in I(v_{i},l)} c_{r}^{l}$ is disjoint from
$W_{j}^{l}-d_{j}^{l}\cup_{r\in I(w_{j},l)} d_{r}^{l}$. Hence
$V_{i}^{l+1}\cap W_{j}^{'}=(c^{l+1}_{i}\cup_{i\in I(v_{i},l+1)}
c_{r}^{l+1})\cap(d_{j}^{'}\cup_{r\in I(w_{j},l+1)} d_{r}^{'})$.
Hence (2) holds.

Case 2. $\gamma_{0}=\emptyset$.

In this case, $I(w_{j},l+1)=I(w_{j},l)-\bigl\{m^{l+1}\bigr\}$. Now
by Proposition 6, $d_{j_{0}}^{l}$ is disjoint from $\partial
V_{i}^{l}-c^{l}_{i}\cup_{i\in I(v_{i},l)} c_{r}^{l}$. By
Definition 5.2.1, $d_{j_{0}}^{l}$ is disjoint from $F^{l+1}$. By
Definition 5.2.3, $d_{j_{0}}^{l}$ is disjoint from $c_{i}^{l+1}$.
Hence $W_{j}^{'}\cap F^{l+1}=d_{j}^{'}\cup_{r\in I(w_{j},l+1)}
d_{r}^{'}$ and $V_{i}^{l+1}\cap W_{j}^{'}=(c^{l+1}_{i}\cup_{i\in
I(v_{i},l+1)} c_{r}^{l+1})\cap(d_{j}^{'}\cup_{r\in I(w_{j},l+1)}
d_{r}^{'})$.\qquad Q.E.D. \vskip 0.5mm

\subsection{The proofs of Propositions 4-6}

\ \ \ \ \  In this section, we shall first construct $d_{j}^{l+1},
e^{l+1}_{\gamma}\times I$ from $d_{j}^{'}$ and
$e^{'}_{\gamma}\times I$ for $j\in\bigl\{1,\ldots,n\bigr\}-m(l+1)$
and $\gamma\in m(l+1)$. Then we shall prove Propositions 4-6 for
the case: $k=l+1$ and $s(v_{l+1})=-$.

{\bf \em Construction(**).}

Since $j_{\alpha}>j_{0}=m^{l+1}$ for $\alpha\neq 0$,
$\gamma_{\alpha}>m^{l+1}$ if $\gamma_{\alpha}\neq\emptyset$. By
Lemma 5.1.2(7) and Lemma 5.1.3(6), $d_{j_{0}}^{l}\times
[-2.5,2.5]$ is disjoint from $(a_{\alpha}\times
I)_{\gamma_{\alpha}}$. By Definitions 3.2.1 and Lemma 5.4.4,
$a_{\alpha}\subset d_{j}^{'}=d_{j}^{l}$, $(a_{\alpha}\times
I)_{\gamma_{\alpha}}\subset (d_{j_{\alpha}}^{'}\times
I)_{\gamma_{\alpha}}\subset e^{'}_{\gamma_{\alpha}}\times
I=e^{l}_{\gamma_{\alpha}}\times I$ if $s(w_{\gamma_{\alpha}})=-$,
and $(a_{\alpha}\times I)_{\gamma_{\alpha}}=a_{\alpha}\subset
d_{j_{\alpha}}^{'}$ if $s(w_{\gamma_{\alpha}})=+$ or
$\gamma_{\alpha}=\emptyset$.

Without loss of generality, we may assume that, in $F^{l}$,
$d_{j_{0}}^{l}\times [-2.5,2.5]\cap D^{*}_{0,\alpha}\subset
d_{j_{0}}^{l}\times [0,2.5]$ for $\alpha>0$, and
$d_{j_{0}}^{l}\times [-2.5,2.5]\cap D^{*}_{0,\alpha}\subset
d_{j_{0}}^{l}\times [-2.5,0]$ for $\alpha<0$. See Lemma 3.5.5.

By Definition 5.2.2(2) and Definition 5.3.1(2), in $F^{l+1}$,
$d_{j_{0}}^{l}\times [-2.5,-1.5]$ is a disk such that $\partial
d_{j_{0}}^{l}\times [-2.5,-2]=\partial d_{j_{0}}^{l}\times
[-2,-1.5]$, and $d_{j_{0}}^{l}\times [1.5,2.5]$ is a disk such
that $\partial d_{j_{0}}^{l}\times [1.5,2]=\partial
d_{j_{0}}^{l}\times [2,2.5]$.

For $\alpha>0$, let $C_{\alpha}$ be a simple closed curve in
$d_{j_{0}}^{l}\times [1.5,2.5]$ satisfying the following
conditions:

(1) \ $C_{\alpha}$ intersects $d_{j_{0}}^{l}\times [2,2.5]$ in a
core $d_{j_{0}}^{l}\times \bigl\{t_{\alpha}\bigr\}$ of
$d_{j_{0}}^{l}\times [2,2.5]$ where $t_{\alpha}\in (2,2.5)$,
$C_{\alpha}$ intersects $e^{'}_{m^{l+1}}\times I$ in a core of
$e^{'}_{m^{l+1}}\times (0,1)$ lying in  $d_{j_{0}}^{l}\times
(1.6,1.7)$, say $e_{m^{l+1},\alpha}$.

(2) \ If $0<\alpha<\beta$, then $t_{\alpha}>t_{\beta}$.

(3) \ If $\alpha\neq \beta$, then $C_{\alpha}\cap
C_{\beta}=\emptyset$. See Figure 25.

For $\alpha<0$, let $C_{\alpha}$ be a simple closed curve in
$d_{j_{0}}^{l}\times [-2.5,-2]$ satisfying the following
conditions:

(4) \ $C_{\alpha}$ intersects $d_{j_{0}}^{l}\times [-2.5,-2]$ in a
core $d_{j_{0}}^{l}\times \bigl\{t_{\alpha}\bigr\}$ of
$d_{j_{0}}^{l}\times [-2.5,-2]$ where $t_{\alpha}\in (-2.5,-2)$,
$C_{\alpha}$ intersects $e^{'}_{m^{l+1}}\times I$ in a core of
$e^{'}_{m^{l+1}}\times (0,1)$ lying in  $d_{j_{0}}^{l}\times
(-1.7,-1.6)$, say $e_{m^{l+1},\alpha}$.

(5) \ If $\beta<\alpha<0$, then $t_{\beta}>t_{\alpha}$.

(6) \ If $\alpha\neq \beta$, then $C_{\alpha}\cap
C_{\beta}=\emptyset$. See Figure 25.

\begin{center}
\includegraphics[totalheight=6cm]{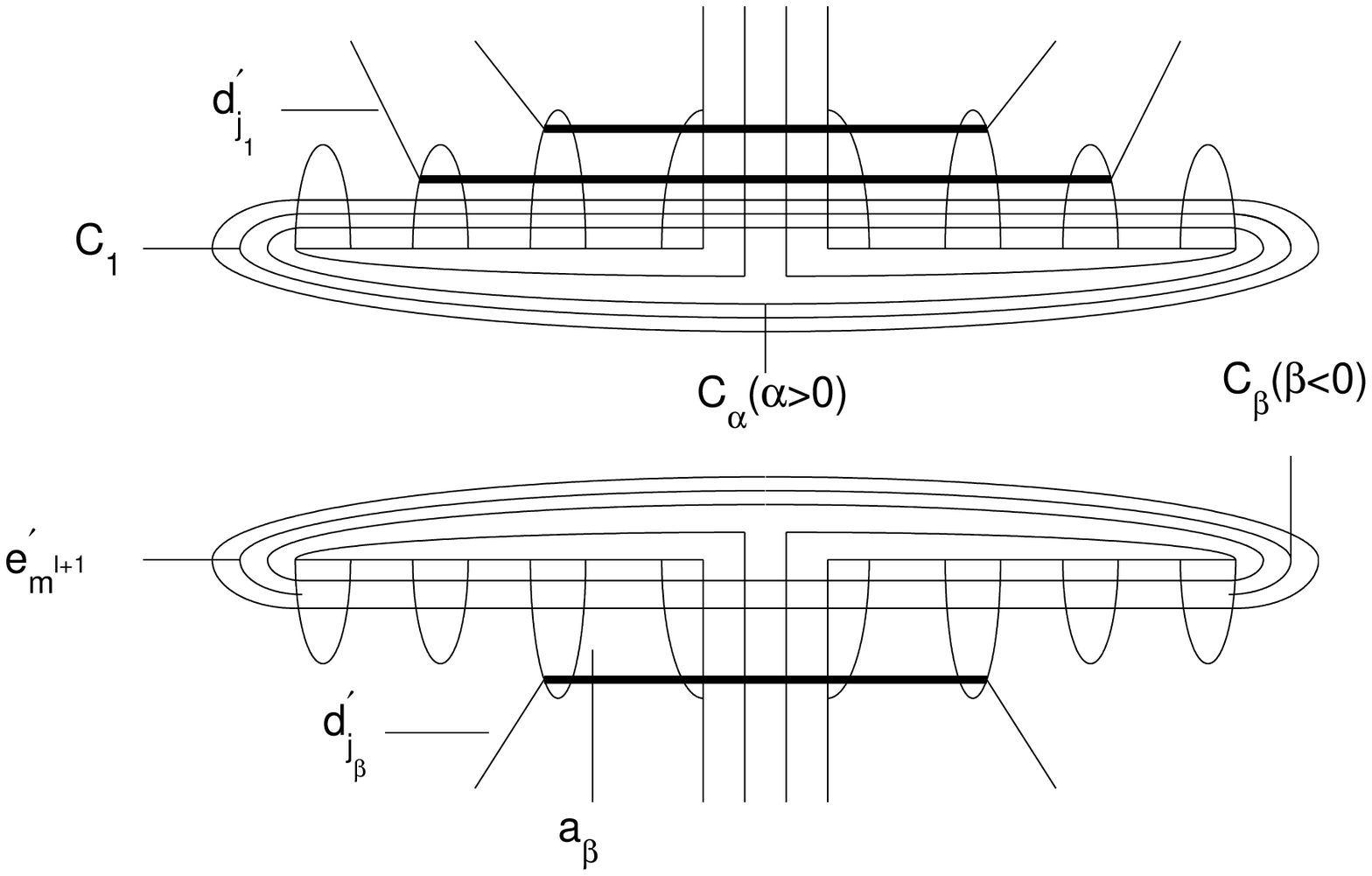}
\begin{center}
Figure 25
\end{center}
\end{center}

Now $C_{\alpha}$ bounds a disk in $d_{j_{0}}^{l}\times
[-2.5,-1.5]\subset F^{l+1}$ for $\alpha<0$, and $C_{\alpha}$
bounds a disk in $d_{j_{0}}^{l}\times [1.5,2.5]\subset F^{l+1}$
for $\alpha>0$.

Note that
$a_{\alpha}=\cup_{i=\delta(0,\alpha)}^{\theta(0,\alpha)}d_{i,\alpha}
\cup_{i=\delta(0,\alpha)}^{\theta(0,\alpha)}e_{i,\alpha}$, and
$a_{\alpha}^{0}\times
\bigl\{t_{\alpha}\bigr\}=\cup_{i=\delta(0,\alpha)}^{\theta(0,\alpha)}d_{i,0}\times\bigl\{t_{\alpha}\bigr\}
\cup_{i=\delta(0,\alpha)}^{\theta(0,\alpha)}e_{i,0}\times\bigl\{t_{\alpha}\bigr\}$.
Since $a_{\alpha}$ and $a_{\alpha}^{0}\times
\bigl\{t_{\alpha}\bigr\}$ is disjoint from $d_{j_{0}}^{l}\times
[-2,2]$. By Definition 5.2.3, Definitions 5.3.2 and 5.3.3,
$\partial_{1} a_{\alpha},
\partial_{1} a^{0}_{\alpha}\times \bigl\{t_{\alpha}\bigr\}\subset
E^{l+1}_{f_{\delta(0,\alpha),i}}$ and $\partial_{2} a_{\alpha},
\partial_{2} a^{0}_{\alpha}\times \bigl\{t_{\alpha}\bigr\}\subset
E^{l+1}_{f_{\theta(0,\alpha),i}}$. In particular, if
$\delta(0,\alpha)=\theta(0,\alpha)=0$. By Lemma 3.5.8, Definition
5.2.3 and Definition 5.3.3, we may assume that
$E^{l+1}_{f_{\delta(0,\alpha),i}}=E^{l+1}_{0}$ and
$E^{l+1}_{f_{\theta(0,\alpha),i}}=E^{l+1}_{l+1}$.
\begin{center}
\includegraphics[totalheight=5cm]{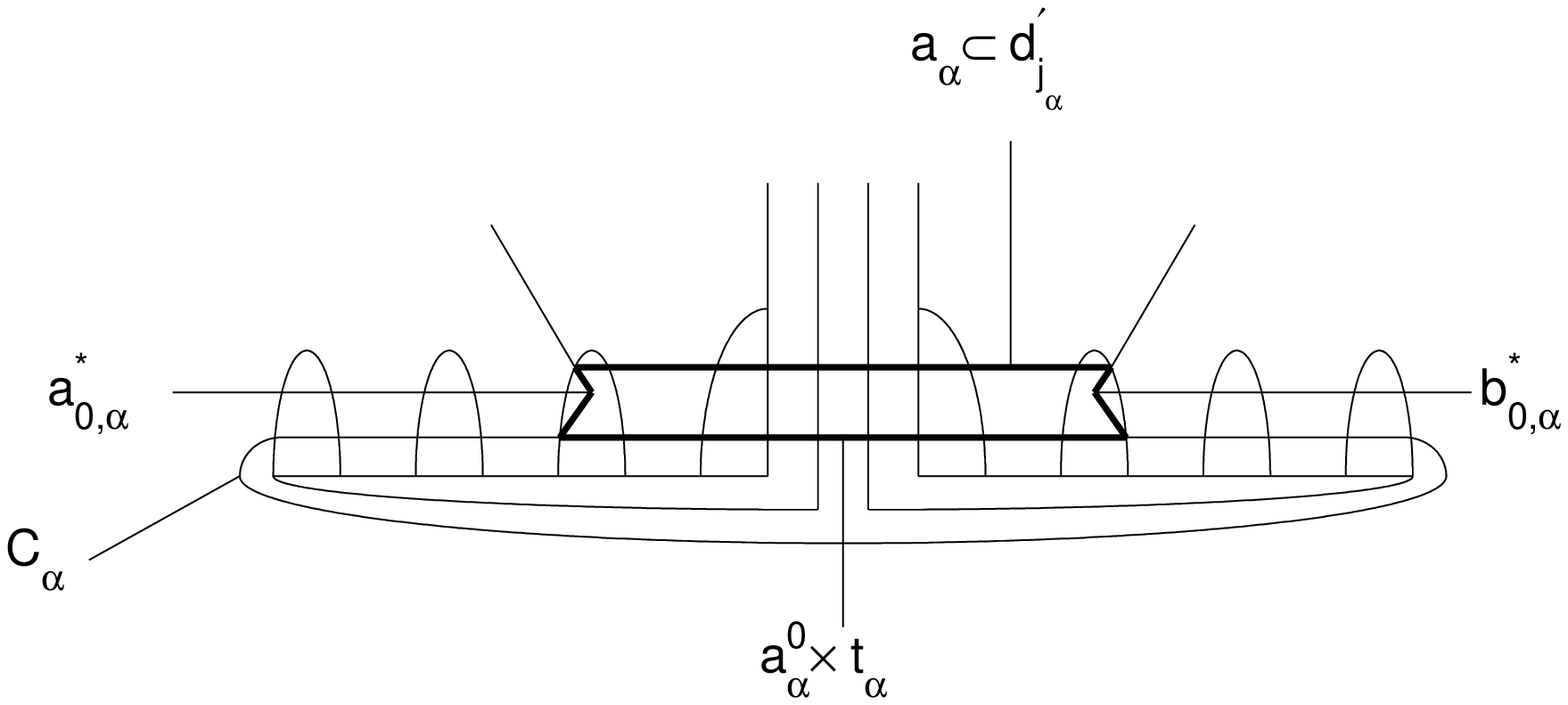}
\begin{center}
Figure 26
\end{center}
\end{center}

Now let $a_{0,\alpha}^{*}$ be an arc in
$E^{l+1}_{f_{\delta(0,\alpha),i}}$ connecting $\partial_{1}
a_{\alpha}$ to, $\partial_{1} a^{0}_{\alpha}\times
\bigl\{t_{\alpha}\bigr\}$ such that

(7) \ $a_{0,\alpha}^{*}-d_{j_{0}}^{l}\times
[-2.5,2.5]=a_{0,\alpha}-d_{j_{0}}^{l}\times [-2.5,2.5]$ and
$a_{0,\alpha}^{*}\cap d_{j_{0}}^{l}\times [-2.5,2.5]$ is an arc
lying in $d_{\delta(0,\alpha),0}\times [-2.5,-2]\cup [2,2.5]$.

Now let $b_{0,\alpha}^{*}$ be an arc in
$E^{l+1}_{f_{\theta(0,\alpha),i}}$ connecting $\partial_{2}
a_{\alpha}$ to, $\partial_{2} a^{0}_{\alpha}\times
\bigl\{t_{\alpha}\bigr\}$ such that

(8) \ $b_{0,\alpha}^{*}-d_{j_{0}}^{l}\times
[-2.5,2.5]=b_{0,\alpha}-d_{j_{0}}^{l}\times [-2.5,2.5]$ and
$b_{0,\alpha}^{*}\cap d_{j_{0}}^{l}\times [-2.5,2.5]$ is an arc
lying in $d_{\theta(0,\alpha),0}\times [-2.5,-2]\cup [2,2.5]$.\\
Where $a_{0,\alpha}, b_{0,\alpha}$ are as in Lemma 3.5.5. See
Figure 26.

By Lemma 5.1.4(1) and the proof of Lemma 4.3.5,
$a_{0,\alpha}^{*}\cup b_{0,\alpha}^{*}$ is disjoint from
$c_{l+1}^{l}\times I$ and $\partial V_{l+1}^{l}\times I$.
Furthermore, $(a_{0,\alpha}^{*}\cup b_{0,\alpha}^{*})\cap
(a_{0,\beta}^{*}\cup b_{0,\beta}^{*})=\emptyset$ for $\alpha,
\beta\neq 0$, $\alpha\neq\beta$.

Now let
$b_{\alpha}=a_{0,\alpha}^{*}\cup(C_{\alpha}-a_{\alpha}^{0}\times\bigl\{t_{\alpha}\bigr\})\cup
b_{0,\alpha}^{*}$, and $b_{\alpha}\times I$ be a neighborhood of
$b_{\alpha}$ in $F^{l+1}$ satisfying the following conditions:

(9) \ $(\partial b_{\alpha})\times I=((\partial a_{\alpha})\times
I)_{\gamma_{\alpha}}$ if $s(w_{\gamma_{\alpha}})=-$.

(10) \ $\partial b_{\alpha}=((\partial a_{\alpha})\times
I)_{\gamma_{\alpha}}$ if $s(w_{\gamma_{\alpha}})=+$ or
$\gamma_{\alpha}=\emptyset$.

(11) \ $b_{\alpha}\times I\subset d_{j_{0}}^{l}\times
(-2.5,-2)\cup (2,2.5)\cup (-1.7,-1.6)\cup (1.6,1.7)$.

(12) \ For $\alpha\neq \beta$, $b_{\alpha}\times I\cap
b_{\beta}\times I=\emptyset$. \qquad \qquad
Q.E.D.(Construction(**))\vskip 1mm

{\bf \em Lemma 5.8.1.} \ (1) \ If
$d_{j}^{'}\cap(a_{0,\beta}^{*}\times I\cup b_{0,\beta}^{*}\times
I)\neq\emptyset$, then  $d_{j}^{'}\cap(a_{0,\beta}^{*}\times I\cup
b_{0,\beta}^{*}\times I)\subset \cup_{\alpha\neq 0}
(a_{\alpha}\times I)_{\gamma_{\alpha}}$. Furthermore, either
$j=j_{\lambda}$, or $j>\gamma_{\lambda}$ for some $\lambda\neq 0$.

(2) \ $e^{'}_{\gamma}\times I\cap(a_{0,\beta}^{*}\times I\cup
b_{0,\beta}^{*}\times I)\subset \cup_{\alpha\neq 0}
(a_{\alpha}\times I)_{\gamma_{\alpha}}$. If $e^{'}_{\gamma}\times
I\cap(inta_{0,\beta}^{*}\times I\cup intb_{0,\beta}^{*}\times
I)\neq\emptyset$, then $e^{'}_{\gamma}\times
I\cap(inta_{0,\beta}^{*}\times I\cup intb_{0,\beta}^{*}\times
I)\subset \cup_{\alpha\neq 0} (a_{\alpha}\times
I)_{\gamma_{\alpha}}$. Furthermore,  $\gamma \geq\gamma_{\lambda}$
for some $\lambda\neq 0$.

{\bf \em Proof.} \ By the proof of Lemma 4.3.5 and Lemma 5.1.4(1),
$a_{0,\alpha}^{*}\times I\cup b_{0,\alpha}^{*}\times I$ is
disjoint from $\partial V_{l+1}^{l}\times I$. By
Construction(**)(7) and (8), $inta_{0,\alpha}^{*}\times I\cup
intb_{0,\alpha}^{*}\times I$ is disjoint from $d_{j_{0}}^{l}\times
[-2,2]$. Hence $inta_{0,\alpha}^{*}\times I\cup
intb_{0,\alpha}^{*}\times I$ is disjoint from
$e^{'}_{m^{l+1}}\times I\subset d_{j_{0}}^{l}\times
[-2,2]\cup\partial V_{l+1}^{l}\times I$. Suppose $j\notin m(l+1)$
and $\gamma\in m(l)$.   By Definitions 5.2.3, 5.3.2 and 5.3.3,
$d_{j}^{'}\cap (a_{0,\alpha}^{*}\times I\cup
b_{0,\alpha}^{*}\times I)=d_{j}^{l}\cap (a_{0,\alpha}^{*}\times
I\cup b_{0,\alpha}^{*}\times I)$ and $e^{'}_{\gamma}\times
I\cap(a_{0,\alpha}^{*}\times I\cup b_{0,\alpha}^{*}\times
I)=e^{l}_{\gamma}\times I\cap (a_{0,\alpha}^{*}\times I\cup
b_{0,\alpha}^{*}\times I)$. By Lemma 4.3.5, the lemma holds.\qquad
\qquad Q.E.D.\vskip 0.5mm

{\bf \em Lemma 5.8.2.} \
$b_{\alpha}=a_{0,\alpha}^{*}\cup_{i=\delta(0)}^{\delta(0,\alpha)-1}
d_{i,0}\times
\bigl\{t_{\alpha}\bigr\}\cup_{i=\delta(0)}^{\delta(0,\alpha)-1}
e_{i,0}\times \bigl\{t_{\alpha}\bigr\}\cup
e_{m^{l+1},\alpha}\cup_{i=\theta(0,\alpha)+1}^{\theta(0)}
d_{i,0}\times
\bigl\{t_{\alpha}\bigr\}\cup_{i=\theta(0,\alpha)+1}^{\theta(0)}
e_{i,0}\times \bigl\{t_{\alpha}\bigr\}\cup b_{0,\alpha}^{*}$
satisfying the following conditions:

(1) \ $intd_{i,0}\times\bigl\{t_{\alpha}\bigr\}$ is disjoint from
$e^{'}_{\gamma}\times I$ for $\gamma\in m(l+1)$.

(2) \ $e_{i,0}\times\bigl\{t_{\alpha}\bigr\}$ is a core of
$e^{'}_{\gamma}\times (0,1)$.

(3) \ $e_{m^{l+1},\alpha}$ is a core of $e^{'}_{m^{l+1}}\times
(0,1)$. \\ Where $e_{m^{l+1},\alpha}$  is as in
Construction(**)(1).

{\bf \em Proof.} \ (1) \ By Definition 3.5.1, $d_{0,0}\subset
a_{\alpha}^{0}$ for each $\alpha\neq 0$. Hence $f_{i,0}\neq 0$ for
$\delta(0)\leq i\leq\delta(0,\alpha)-1$ or $\theta(0,\alpha)+1\leq
i\leq \theta(0)$. Hence $d_{i,0}\times [-2.5,2.5]$ is disjoint
from $c_{l+1}^{l}\subset E^{l}_{0}$. Since $t_{\alpha}\in
[-2.5,-2]\cup [2,2.5]$, by Definitions 5.2.3, 5.3.2 and Lemmas
5.1.2(8), 5.1.3(4), $intd_{i,0}\times \bigl\{t_{\alpha}\bigr\}$ is
disjoint from $e^{'}_{\gamma}\times I$ for each $\gamma\in
m(l+1)$.

(2) By Definitions 5.2.3, 5.3.2,
$e_{i,0}\times\bigl\{t_{\alpha}\bigr\}$ is a core of
$e^{'}_{\gamma}\times (0,1)$.

(3) follows from Construction(**), Definitions 5.2.3 and 5.3.2.
\qquad Q.E.D.\vskip 0.5mm

{\bf \em Lemma 5.8.3.} \ (1) \ If $\gamma\neq m^{l+1}
>Max\bigl\{\gamma_{i,0} \ | \ \delta(0)\leq
i\leq\theta(0)\bigr\}$, then
$C_{\alpha}-inta_{\alpha}^{0}\times\bigl\{t_{\alpha}\bigr\}$ is
disjoint from $e^{'}_{\gamma}\times I$.

(2) \
$(C_{\alpha}-inta_{\alpha}^{0}\times\bigl\{t_{\alpha}\bigr\})\cap
e^{'}_{m^{l+1}}\times I=e_{m^{l+1},\alpha}$.

(3) \ If $\gamma\leq Max\bigl\{\gamma_{i,0} \ | \ \delta(0)\leq
i\leq\theta(0)\bigr\}$, then each component of
$(C_{\alpha}-inta_{\alpha}^{0}\times\bigl\{t_{\alpha}\bigr\})\cap
e^{'}_{\gamma}\times I$ is a core of $e^{'}_{\gamma}\times (0,1)$.

{\bf \em Proof.} \ (1) \ Suppose that $\gamma\neq m^{l+1}
>Max\bigl\{\gamma_{i,0} \ | \ \delta(0)\leq
i\leq\theta(0)\bigr\}$. Since
$(C_{\alpha}-inta_{\alpha}^{0}\times\bigl\{t_{\alpha}\bigr\})\subset
d_{j_{0}}^{l}\times (-2.5,-1.5)\cup (1.5,2.5)$, by Lemma 5.1.2(2)
and (4), Lemma 5.1.3(2), $e^{l}_{\gamma}\times I$ is disjoint from
$C_{\alpha}-inta_{\alpha}^{0}\times\bigl\{t_{\alpha}\bigr\}$. By
Lemma 5.4.4, $e^{'}_{\gamma}\times I=e^{l}_{\gamma}\times I$.
Hence (1) holds.

(2) \ By Construction(**)(1),
$C_{\alpha}-inta_{\alpha}^{0}\times\bigl\{t_{\alpha}\bigr\}=\cup_{i=\delta(0)}^{\delta(0,\alpha)-1}
d_{i,0}\times
\bigl\{t_{\alpha}\bigr\}\cup_{i=\delta(0)}^{\delta(0,\alpha)-1}
e_{i,0}\times \bigl\{t_{\alpha}\bigr\}\cup
e_{m^{l+1},\alpha}\cup_{i=\theta(0,\alpha)+1}^{\theta(0)}
d_{i,0}\times
\bigl\{t_{\alpha}\bigr\}\cup_{i=\theta(0,\alpha)+1}^{\theta(0)}
e_{i,0}\times \bigl\{t_{\alpha}\bigr\}$.  Since $d_{0,0}\subset
a_{\alpha}^{0}$, by Lemma 5.1.4(4),
$C_{\alpha}-inta_{\alpha}^{0}\times\bigl\{t_{\alpha}\bigr\}-e_{m^{l+1},\alpha}$
is disjoint from $c_{l+1}^{l}\times I$ and $\partial
V_{l+1}^{l}\times I$. Furthermore,
$C_{\alpha}-inta_{\alpha}^{0}\times\bigl\{t_{\alpha}\bigr\}-e_{m^{l+1},\alpha}\subset
d_{j_{0}}^{l}\times (-2.5,-2)\cup (2,2.5)$, and
$e^{'}_{m^{l+1}}\times I\subset d_{j_{0}}^{l}\times [-2,2]\cup
\partial V_{l+1}^{l}\times I$. Hence (2) holds.

(3) \ Suppose that $\gamma\leq Max\bigl\{\gamma_{i,0} \ | \
\delta(0)\leq i\leq\theta(0)\bigr\}$. Then $\gamma<j_{0}=m^{l+1}$.
By (2),
$C_{\alpha}-inta_{\alpha}^{0}\times\bigl\{t_{\alpha}\bigr\}-e_{m^{l+1},\alpha}$
is disjoint from $c_{l+1}^{l}\times I$. By Lemma 5.5.2,
$inte_{m^{l+1},\alpha}$ is disjoint from $e^{'}_{\gamma}\times I$.
By Definition 5.2.3 and Definition 5.3.3,
$(C_{\alpha}-inta_{\alpha}^{0}\times\bigl\{t_{\alpha}\bigr\}-inte_{m^{l+1},\alpha})\cap
e^{'}_{\gamma}\times I=
(C_{\alpha}-inta_{\alpha}^{0}\times\bigl\{t_{\alpha}\bigr\}-inte_{m^{l+1},\alpha})\cap
e^{l}_{\gamma}\times I$. By Lemma 5.1.2 and Lemma 5.1.3, (3)
holds.\qquad Q.E.D.\vskip 0.5mm

{\bf \em Definition 5.8.4.} \ (1) If $s(w_{\gamma_{\alpha}})=-$,
then let $H_{\alpha}$ be a homeomorphism from $(a_{\alpha}\times
I)_{\gamma_{\alpha}}$ to $b_{\alpha}\times I$ such that
$H_{\alpha}$ is an identifying map on $((\partial
a_{\alpha})\times I)_{\gamma_{\alpha}}$.

(2) \ If $s(w_{\gamma})=+$ or $\gamma_{\alpha}=\emptyset$, let
$H_{\alpha}$ be a homeomorphism from $a_{\alpha}=(a_{\alpha}\times
I)_{\gamma_{\alpha}}$ to $b_{\alpha}$ such that $H_{\alpha}$ is an
identifying map on $\partial a_{\alpha}$.\vskip 0.5mm

{\bf \em Definition 5.8.5.} \ (1) \ For $j\notin m(l+1)$, let
$d_{j}^{l+1}=(d_{j}^{'}-\cup_{\alpha\neq 0} (a_{\alpha}\times
I)_{\gamma_{\alpha}})\cup_{\alpha\neq 0} H_{\alpha}(d_{j}^{'}\cap
(a_{\alpha}\times I)_{\gamma_{\alpha}})$.

(2) \ For $\gamma\in m(l+1)$, if $\gamma\leq m^{l+1}$, let
$e^{l+1}_{\gamma}\times I=e^{'}_{\gamma}\times I$; if
$\gamma>m^{l+1}$, let $e^{l+1}_{\gamma}\times
I=(e^{'}_{\gamma}\times I-\cup_{\alpha\neq 0} (a_{\alpha}\times
I)_{\gamma_{\alpha}})\cup_{\alpha\neq 0}
H_{\alpha}(e^{'}_{\gamma}\times I\cap (a_{\alpha}\times
I)_{\gamma_{\alpha}})$.\vskip 0.5mm

{\bf \em Lemma 5.8.6.} \ $d_{j}^{l+1}$ is isotopic to $d_{j}^{'}$
in $F^{l+1}\subset \partial_{+} \cal V_{-}$.

{\bf \em Proof.} \ By Lemma 3.5.5, $a_{0,\alpha}^{*}\cup
a_{\alpha}\cup a_{\alpha}^{0}\times\bigl\{t_{\alpha}\bigr\}\cup
b_{0,\alpha}^{*}$ bounds a disk in $F^{l+1}$. See Figure 26. By
Construction(**), $C_{\alpha}$ bounds a disk in $F^{l+1}$. Hence
$a_{\alpha}$ is isotopic to $b_{\alpha}$. By Lemma 5.6.1 and Lemma
5.6.2, each component of $d_{j}^{'}\cap (a_{\alpha}\times
I)_{\gamma_{\alpha}}$ is a core of $(a_{\alpha}\times
I)_{\gamma_{\alpha}}$ even if $j=j_{\alpha}$, say $c_{j}$. By
Definitions 5.8.4 and 5.8.5, $H_{\alpha}(c_{j})$ is a core of
$b_{\alpha}\times I$. Hence $d_{j}^{l+1}$ is isotopic to
$d_{j}^{'}$. \qquad Q.E.D.\vskip 2mm

Now we prove Propositions 4-6 for the case: $k=l+1$ and
$s(v_{l+1})=-$.

{\bf \em The proofs of Propositions 4-6.} \ By Lemma 5.4.4, Lemmas
5.5.1-5.5.5, $F^{l+1}$ is generated by the abstract tree $\cup_{f}
E^{l+1}_{f}\cup_{\gamma\in m(l+1)} e^{'}_{\gamma}$ satisfying the
following conditions:

(1) \ If $s(w_{\lambda})=+$, then $inte^{'}_{\lambda}\times I$ is
disjoint from $\cup_{f} E^{l+1}_{f}\cup_{\gamma<\lambda}
e^{'}_{\gamma}\times I$.

(2) \ If $s(w_{\lambda})=-$, then $inte^{'}_{\lambda}\times
I\cup(\cup_{f} E^{l+1}_{f}\cup_{\gamma<\lambda}
e^{'}_{\gamma}\times I)=\cup_{r\in I(w_{\lambda},l+1)}
(d_{r}^{'}\times I)_{\lambda}$.

By Lemma 5.6.1, if $j<\gamma$ and $j\notin I(w_{\gamma},l+1)$,
then $d_{j}^{'}$ is disjoint from $e^{'}_{\gamma}\times I$.

By Lemma 5.6.2, each component of  $d_{j}^{'}\cap
e^{'}_{\gamma}\times I$ is a core of $e^{'}_{\gamma}\times (0,1)$
for $j\notin m(l+1), \gamma\in m(l+1)$ and $j>\gamma$.

By Lemma 5.6.3, $d_{j}^{'}$ is regular in $\cup_{f}
E^{l+1}_{f}\cup_{\gamma<j} e^{'}_{\gamma}\times I$ for $j\notin
L(c_{l+1}^{l})$.

By Lemma 5.4.3, $c_{i}^{l+1}$ is an  arc properly embedded in
$F^{l+1}$ which lies in one of $E^{l+1}_{f}$ for some $f$;
$\partial c_{i}^{l+1}\cap e^{'}_{\gamma}\times I=\emptyset$ for
each $\gamma\in m(l+1)$, and $\partial c_{i}^{l+1}\cap
d_{j}^{'}=\emptyset$ for each $j\in\bigl\{1,2,\ldots,
n\bigr\}-m(l+1)$.

Now if we take place of $c^{l+1}_{i}, d_{j}^{l+1}, E^{l+1}_{f},
e^{l+1}_{\gamma}\times I, F^{l+1}$ with $c^{l+1}_{i}, d_{j}^{'},
E^{l+1}_{f}, e^{'}_{\gamma}\times I, F^{l+1}$, then Proposition 4
holds except that

(i) \ $d_{j_{\alpha}}^{'}$ is regular in $\cup_{f}
E^{l+1}_{f}\cup_{\gamma<j_{\alpha}}e^{'}_{\gamma}\times I$;

(ii)
$d_{j_{\alpha}}^{'}-\cup_{\gamma<j_{\alpha}}inte^{'}_{\gamma}\times
I$ intersects $c^{l+1}_{i}$ in at most one point. \\ Where
$j_{\alpha}\in L(c_{l+1}^{l})$.

Now by Lemmas 5.6.4, 5.8.1, 5.8.2, 5.8.3 and the argument in
Chapter 4, Proposition 4 holds for the case: $k=l+1$ and
$s(v_{l+1})=-$.

By Lemmas 5.6.3, 5.6.4, 5.8.2 and the proof of Proposition 5 in
Section 4.6, Proposition 5 holds for the case: $k=l+1$ and
$s(v_{l+1})=-$.

Proposition 6 follows from Lemmas 5.8.6, 5.7.1 and the proof of
Proposition 6 in Section 4.6.\qquad Q.E.D.\vskip 0.5mm

\section{The proofs of Propositions 1-3}

\ \ \ \ \ Now we  prove Propositions 1-3 for $k=l+1$ under the
assumptions that Propositions 1-6 hold for $k\leq l$. Then we
finish the proofs of Propositions 1-6.

{\bf \em The Proofs of Propositions 1-3.} \ By Lemma 3.1.5,
$L(v_{i}^{l})=L(c^{l}_{i})$ and $L(d_{j}^{l})=L(w_{j}^{l})$. Hence
$m^{l+1}=MinL(c_{l+1}^{l})=MinL(v_{l+1}^{l})$. Now there are two
cases:

Case 1. \ $s(v_{l+1})=-$.

Now by Lemma 5.1.1(4), $s(w_{m^{l+1}})=+$. Now by  Remark 3.5.9
and the argument in Chapter 4, Propositions 1-3 hold.

Case 2. \ $s(v_{l+1})=+$.

By Remark 3.5.9 and the argument in Chapter 5, Propositions 1-3
hold.\qquad Q.E.D.\vskip 0.5mm

Now Propositions 1-6 hold for $0\leq k\leq m$. Hence Theorem 1 is
true.\vskip 0.5mm

{\bf Acknowledgement:} \ The author is grateful to Professor
Shicheng Wang for many helpful discussions with the author on this
problem and some comments on this paper. The author is  grateful
to Professors C. Gordon, T. Kobayashi, Fengchun Lei and Ying-qing
Wu for some helpful discussions on this problem. The author thank
Doctor Jiming Ma for careful reading of this paper. \vskip 4mm

{\bf References.} \vskip 2mm

[B] \ David Bachman, Connected sums of unstabilized Heegaard
splittings are unstabilized, Math.GT/0404058.

[CG1] \ A. Casson and C. Gordon, 3-Manifolds with irreducible
Heegaard splittings of arbitrarily high genus, unpublished.

[CG2] \ A. Casson and C. Gordon, Reducing Heegaard splittings,
Topology and Its Applications, 27(1987), 275-283.

[H] \ W. Haken, \ Some results on surfaces in 3-manifolds, Studies
in Modern Topology (Math. Assoc. Amer., distributed by:
prentice-Hall, 1968), 34-98.

[Ki] \ R. Kirby, Problems in Low-Dimensional Topology, Geometric
Topology, Edited by H. Kazez, AMS/IP Vol. 2, International Press,
1997.

[Kn] \ H. Kneser, \ Geschlossene Fl$\ddot{a}$chen in
dreidimensionalen Mannig-flatigkeien, Jahresbericht der Deut.
Math. Verein. 38(1929), 248-260.

[Ko] \ T. Kobayashi, A construction of 3-manifolds whose
homeomorphism classes of Heegaard splittings have polynomial
growth, Osaka Journal of Mathematics, 29(1992), 653-674.

[M] \ J. Milnor,  A unique factorization theorem for 3-manifolds,
Amer. J. Math. 84(1962), 1-7.

[QM] \ Ruifeng Qiu and Jiming Ma, Heegaard splittings of annular
3-manifolds and $\partial$-reducible 3-manifolds, Preprint.

[RS] \ H. Rubinstein and M. Scharlemann, Comparing Heegaard
splittings of non-Haken 3-manifolds, Topology, 35(1996),
1005-1026.

[S] \ J. Singer, Three-dimensional manifolds and their
Heegaard-diagrams, Trans. Amer. Math. Soc., 35(1933),88-111.

[ST] \ M. Scharlemann and A. Thompson, Heegaard splittings of
${\rm Surfaces}\times I$ are standard, Math. Ann., 295(1993),
549-564.

[W] \   F.Waldhausen, \ Heegaard-Zerlegungen der 3-sphere,
Topology, 7(1968), 195-203. \vskip 3mm

Ruifeng Qiu

Department of Mathematics

Dalian University of Technology

Dalian 116024, China

Email: qiurf@dlut.edu.cn

\enddocument
\bye

\end{document}